\documentclass[12pt]{article}

\usepackage{amsfonts,amsmath,amssymb,epsf}
\usepackage{mathrsfs}
\usepackage{graphics} 

\usepackage{latexsym,amsmath,amssymb,amsfonts}
\usepackage{rotating}
\usepackage{mathrsfs}
\usepackage{xypic} \xyoption{all}
\usepackage{amscd}
\usepackage{hyperref}
\usepackage{euscript}
\usepackage{hhline}
\usepackage{graphicx,epstopdf}
\usepackage{xy}
\usepackage{epsfig}

\allowdisplaybreaks

\newlength{\fighskip} \fighskip=2pt
\newlength{\figvskip} \figvskip=3pt

\newcommand*{\figbox}[2]{{
  \def\figscale{#1}
  \def\arraystretch{0.8}
  \arraycolsep=0pt
  \begin{array}{c}
    \vbox{\vskip\figscale\figvskip
      \hbox{\hskip\figscale\fighskip
        \includegraphics[scale=\figscale]{#2}}}
  \end{array}}}

\advance\textheight by \topskip
\newcounter{thm}

\newcommand\sect[1]{\section{#1}\setcounter{equation}0\setcounter{thm}0} 

\newtheorem{thm}{Theorem}[section]
\newtheorem{defn}[thm]{Definition}
\newtheorem{prop}[thm]{Proposition}
\newtheorem{cor}[thm]{Corollary}
\newtheorem{rema}[thm]{Remark}
\newtheorem{lemma}[thm]{Lemma}

\newcommand\void[1]       {}

\newcommand\be            {\begin{equation}}
\newcommand\ee            {\end{equation}}
\newcommand\bea           {\begin{equation}\begin{array}l\displaystyle}
\newcommand\bearll        {\begin{array}{ll}\displaystyle}
\newcommand\eear          {\end{array}}
\newcommand\bnu          {\begin{enumerate}}
\newcommand\enu          {\end{enumerate}}

\newcommand\enl           {\\[1em]\displaystyle}

\newcommand\etb           {& \hspace*{-.5em} \displaystyle}
\newcommand\erf[1]        {(\ref{#1})}
\newcommand\id            {{\mathrm{id}}}
\newcommand\labl[1]       {\label{#1}\ee}

\newcommand\nxt{\noindent\raisebox{.08em}{\rule{.44em}{.44em}}%
                          \hspace{.4em}}

\newcommand\arxiv[2]      {\href{http://arXiv.org/abs/#1}{#2}}
\newcommand\doi[2]        {\href{http://dx.doi.org/#1}{#2}}
\newcommand\httpurl[2]    {\href{http://#1}{#2}}

\renewcommand\cir          {\,{\circ}\,}
\newcommand\eq            {\,{=}\,}
\newcommand\In            {\,{\in}\,}
\newcommand\oti           {\,{\otimes}\,}
\newcommand\ti            {\,{\times}\,}
\newcommand\To            {\,{\rightarrow}\,}
\newcommand\nn          {\nonumber \\}

\newcommand\Aop           {A_\mathrm{op}}
\newcommand\Acl           {A_\mathrm{cl}}
\newcommand\Bop           {B_\mathrm{op}}
\newcommand\Bcl           {B_\mathrm{cl}}
\newcommand\eop           {e_\mathrm{op}}
\newcommand\ecl           {e_\mathrm{cl}}
\newcommand\rop           {r_\mathrm{op}}
\newcommand\rcl           {r_\mathrm{cl}}
\newcommand\ord          {\mathrm{ord}}

\newcommand\Bl            {B\hspace*{-.5pt}\ell}
\newcommand\BlxBr         {B_l{\times}B_r}
\newcommand\Cor           {C\hspace*{-.9pt}o\hspace*{-.6pt}r}
\newcommand\CxC           {{\mathcal{C}_\pm^2}}
\newcommand\Dim           {\mathrm{Dim}}

\newcommand\eps           {\varepsilon}
\newcommand\Fun           {\mathcal{F}\hspace*{-.8pt}\mbox{\sl un}}
\newcommand\Hom           {\mathrm{Hom}}
\newcommand\ico           {\iota_{\text{{\rm cl-op}}}}
\newcommand\Ob            {\mathrm{Ob}}
\newcommand\one           {{\bf1}}
\newcommand\One           {O\hspace*{-.6pt}n\hspace*{-.6pt}e}
\newcommand\sew           {{\mathscr S}}
\newcommand\tftC          {t\hspace*{-1.2pt}f\hspace*{-1.2pt}t_\Cc}
\newcommand\Vect          {{\mathcal V}\hspace*{-.8pt}{\sl ect}}
\newcommand\WSh           {{\mathcal{WS}}\hspace*{-.5pt}{\sl h}}
\newcommand\Xhat          {\widehat{\Xr}}
\newcommand\Xtil          {\widetilde{\Xr}}

\newcommand{\halmos}{\rule{1ex}{1.4ex}}
\newcommand{\pf}{\noindent{\it Proof.}\hspace{2ex}}
\newcommand{\epf}{\hspace*{\fill}\mbox{$\halmos$}}

\def\Cb            {\mathbb{C}}
\def\Hb            {\mathbb{H}}

\def\Rb            {\mathbb{R}}
\def\Zb            {\mathbb{Z}}

\def\Cc            {\mathcal{C}}

\def\Ic            {\mathcal{I}}

\def\Sc            {\mathcal{S}}

\def\Fr            {\mathrm{F}}
\def\Mr            {\mathrm{M}}
\def\Rr            {\mathrm{R}}
\def\Tr            {\mathrm{T}}
\def\Xr            {\mathrm{X}}
\def\Yr            {\mathrm{Y}}
\def\Zr            {\mathrm{Z}}

\begin{document}
\thispagestyle{empty}
\def\thefootnote{\fnsymbol{footnote}}
\begin{flushright}
\end{flushright}
\vskip 1.0em
\begin{center}\LARGE
Cardy algebras and sewing constraints, II
\end{center}\vskip 1.5em
\begin{center}\large
  Liang Kong$\,{}^{a,b}$\footnote{Email: {\tt kong.fan.liang@gmail.com }}
  Qin Li$\,{}^{c,d}$\footnote{Email: {\tt qli@math.cuhk.edu.hk}}
  and
  Ingo Runkel$\,{}^{e}$\footnote{Email: {\tt ingo.runkel@uni-hamburg.de}}%
\end{center}
\begin{center}\it$^a$
Institute for Advanced Study\\
Tsinghua University, Beijing, 100084, China
\end{center}
\begin{center}\it$^b$
Department of Mathematics and Statistics\\
University of New Hampshire\\
33 Academic Way, Durham, 03824, USA
\end{center}
\begin{center}\it$^c$
School of Mathematical Sciences\\
Wu Wen-Tsun Key Laboratory of Mathematics\\
 University of Science and Technology of China\\
Hefei, Anhui, 230026, China
\end{center}
\begin{center}\it$^d$
  Department of Mathematics\\
  The Chinese University of Hong Kong\\
  Shatin, Hong Kong, China
\end{center}
\begin{center}\it$^e$
  Fachbereich Mathematik,
Universit$\ddot{a}$t Hamburg \\
  Bundesstrasse 55,
D - 20146 Hamburg, Germany
\end{center}

\vskip 1em

\begin{abstract}
This is the part II of a two-part work started in \cite{part1}. In part I, Cardy algebras were studied, a notion which arises from the classification of genus-0,1 open-closed rational conformal field theories. In this part, we prove that a Cardy algebra also satisfies the higher genus factorisation and modular-invariance properties formulated in \cite{unique} in terms of the notion of a solution to the sewing constraints. We present the proof by showing that the latter notion, which is defined as a monoidal natural transformation, can be expressed in terms of generators and relations, which correspond exactly to the defining data and axioms of a Cardy algebra. 
\end{abstract}

\newpage

\tableofcontents

\setcounter{footnote}{0}
\def\thefootnote{\arabic{footnote}}

\sect{Introduction and Summary}

This is the part II of the two-part work started in \cite{part1}. In part I, the first and third authors have studied so-called Cardy algebras,  which were
first introduced in \cite{cardy} by the first author. Cardy algebras arises in a careful analysis of the generators and relations of genus 0 and 1 open-closed
rational conformal field theory (RCFT). The RCFT is required to satisfy a so-called $V$-invariant boundary condition (see \cite{cardy} for details), where $V$
is a rational vertex operator algebra such that the category $\Cc_V$ of $V$-modules is a modular tensor category according to Huang's Theorem \cite{huang}. In
this case, the chiral data, which is complex-analytic and infinite-dimensional in nature, is included completely in $V$, and can be proved to be separated
entirely from the rest of the data defining an open-closed RCFT. The latter 
is finite-dimensional, topological and combinatoric in nature and is formulated entirely in
the modular tensor category $\Cc_V$. As a consequence, a Cardy algebra can be defined as a purely categorical notion, namely a triple $(\Aop|\Acl,
 \iota_{\mathrm{cl}-\mathrm{op}})$ based on a modular tensor category $\Cc$, without referring to the chiral data $V$. In particular, the open-string (resp.\ closed-string) sector of a RCFT is determined by the data $\Aop$ (resp.\ $\Acl$),  which can be viewed as an algebra object over a ``stringy ground field" given by a vertex operator algebra $V$ (resp.\ $V\otimes_{\Cb} V$). The precise definition of a Cardy algebra will be recalled in Definition \ref{def:cardy-alg}. More details on the above statements can be found in the the introductory sections of \cite{cardy} and \cite{part1}.
 
\medskip
 
Let us now summarise the pertinent results in \cite{unique} and \cite{part1} and explain how the present paper complements these results.
\begin{itemize}\setlength{\leftskip}{-1em}
\item \cite{unique}:
The separation of the complex-analytic data of an open-closed RCFT from the finite-topological-combinatoric (or just combinatoric for short) data has already been emphasised in \cite{fffs,tft1}. The key tool used in \cite{fffs,tft1,tft5,unique} to efficiently formulate the combinatoric part of the all-genus factorisation and modular invariance properties is the 3-d topological field theory defined by the modular tensor category $\Cc_V$ \cite{turaev-bk}. It is proved in \cite[Thms.\,2.1,\,2.6,\,2.9]{tft5} that for each special symmetric Frobenius algebra in $\Cc_V$ one obtains a solution to the combinatoric all-genus factorisation and modular invariance conditions.

The converse question is addressed in \cite{unique}. To do so, the combinatoric all-genus factorisation and modular invariance properties for $\Cc_V$ are reformulated in terms of a so-called {\em solution to the sewing constraints}, see \cite[Def.\,3.14]{unique} and Definition \ref{def:sol-sew} below. 
This notion can be viewed as an instance of a relative quantum field theory as recently discussed in \cite{relative-qft} (see Remark\,\ref{rema:relative-qft} below). It is proved in \cite[Thm.\,4.26]{unique} that every solution to the sewing constraints for which (i) there is a unique closed-string vacuum, (ii) the disc two-point amplitude is non-degenerate, (iii) the sphere two-point amplitude is non-degenerate, and (iv) the quantum dimension of the open-string sector is non-zero\footnote{
        In \cite{unique} there is an additional condition on the quantum dimensions of the simple objects of $\Cc_V$, but as shown in \cite[Thm.\,3.21]{part1} this condition is not needed.
},
comes from a simple special symmetric Frobenius algebra via the construction in \cite{tft1,tft5}.

We stress here conditions (i) and (iv), which turn out to imply that the algebra in $\Cc_V$ describing the open sector of the RCFT is semi-simple.
 
\item Part I:
The relation between the works \cite{tft1,tft5,unique} and the present two-part work is analogous to the relation between the state sum construction of open-closed 2-d {\em topological} field theory \cite{tft1,Lauda:2006mn} and the classification of all such theories \cite{Lz,Moore:2000nn} (see also \cite{AN,LP,MS}). Roughly speaking, the state sum construction relies on a semi-simplicity assumption which need not be satisfied by a general open-closed 2-d TFT. Therefore, not all open-closed 2-d TFTs arise from a state sum construction.

In the present setting, if we do not require conditions (i) and (iv) from above, it is no longer true that all solutions to the sewing constraints are obtained by applying the construction in \cite{tft1,tft5} to a simple special symmetric Frobenius algebra. The point of the present two-part work is to show that in this more general situation, one has to work with Cardy algebras. 

An important result in part I is that if a Cardy algebra $(A|B, \iota)$ satisfies $\dim(A) \neq 0$ and $B$ is simple, then it is isomorphic to $(A,Z(A),e)$, where $A$ is simple special symmetric Frobenius and $Z(A)$ is the full centre of $A$ \cite[Prop.\,3.16\,\&\,Thm.\,3.21]{part1}. This the Cardy-algebra version of \cite[Thm.\,4.26]{unique}.

\item Part II: The main result of the present paper is a one-to-one correspondence between solutions to the sewing constraints and Cardy algebras; the precise statement is given in Theorem \ref{thm:cardy-sew}. Our approach is as follows:

A solution to the sewing constraints is by definition a natural monoidal transformation from a trivial topological modular functor to a (variant of a) topological modular functor obtained from the 3-d TFT associated to a modular tensor category $\Cc$. The key technical result, maybe even of some interest of its own, is the generators and relations description of such natural monoidal transformations given in Theorem \ref{thm:gen-rel}. This can be understood as one way of making precise the arguments with which the sewing constraints in open-closed RCFT where first analysed in \cite{Le}.

Once the generators and relations presentation of solutions to the sewing constraints in the above sense has been established, one only has to match them to the generators and relations of a Cardy algebra. This is done in Theorem \ref{thm:cardy-sew}.
\end{itemize}

The difference between the construction in \cite{unique} and the present work becomes most striking if we look at the special case of choosing the modular tensor category $\Cc$ to be $\Vect$, the category of finite dimensional complex vector spaces. This amounts to considering 2-d open-closed topological field theories instead of conformal ones.
Conditions (i) and (iv) in \cite[Thm.\,4.26]{unique} now imply that the open and closed sector of the 2-d TFT are just given by $\Cb$ itself. On the other hand, the definition of a Cardy algebra reduces, for $\Cc=\Vect$, to the algebraic object classifying open-closed 2-d TFTs as given in \cite{Lz,Moore:2000nn} and \cite{AN,LP,MS}.

It is worthwhile to point out that the generators and relations description we find (in agreement with \cite{Le}),  differs from the one used in \cite{LP,MS}.
This is not surprising since in 2-d CFTs one assigns to a surface a vector in the space of conformal blocks which is equipped with an action of the mapping
class group, while in 2-d TFT there are no non-trivial actions of mapping class groups. In particular, the condition of modular invariance of the 1-point
genus-one correlation function is trivial in a 2-d TFT. 

\bigskip
\noindent {\bf Acknowledgements}: 
We want to thank Yun Gao, Sen Hu and Yi-Zhi Huang for valuable comments and Jens Fjelstad for helpful discussions. LK wants to thank
IQI at Caltech for its support during his post-doctor year 2008-2009 at the early stage of this work, and the support from the
Gordon and Betty Moore Foundation through Caltech's Center for the Physics of Information, and NSF under Grant No.PHY-0803371. LK
is
also supported by the Basic Research Young Scholars Program and the Initiative Scientific Research Program of Tsinghua University,
and NSFC under Grant No.11071134. QL would like to thank his Ph.D. thesis advisor Peter Teichner for his guidance and
encouragement. QL is supported by Chinese Universities Scientific Fund WK0010000030.  IR is supported in part by the German Science
Foundation (DFG) within the Collaborative Research Center 676 ``Particles, Strings and the Early Universe''.

\sect{Open-closed world sheets}

\subsection{The category of open-closed world sheets}

Up to some minor modifications, the content of this section follows Definitions 3.1, 3.3 and 3.4 of \cite{unique}.

\medskip

We denote by $S^1$ the unit circle $\{|z|\eq 1\}$
in the complex plane, with
counter-clockwise orientation. The map that assigns to a complex number its
complex conjugate is denoted by $\text{Conj} : z \mapsto \bar{z}$.
\\[-2em]

\begin{defn} \label{def:ws} {\rm
An {\em oriented open/closed topological world sheet\/}, or
{\em world sheet\/} for short, is a tuple
  \be
  \Xr \,\equiv \big(\, \Xtil, \imath_\Xr, \delta_\Xr,
  b^\text{in}_\Xr, b^\text{out}_\Xr, \text{or}_\Xr, \ord \,\big)
  \ee
consisting of:
\\[.3em]
\nxt An oriented compact two-dimensional topological manifold $\Xtil$.
The (possibly empty) boundary $\partial\Xtil$ of
$\Xtil$ is oriented by the inward-pointing normal.
\\[.3em]
\nxt A continuous orientation-reversing involution
  \be
  \imath_\Xr{:}\quad \Xtil \to \Xtil \,,
  \ee
such that fixed point set of $\imath_\Xr$ is a submanifold and the quotient surface $\dot\Xr \,{:=}\, \Xtil / \langle \imath_\Xr \rangle$ is a manifold with boundaries and corners. We denote by $\tilde\pi_\Xr {:}\ \Xtil \To \dot\Xr$ the canonical projection.
\\[.3em]
\nxt A partition of the set $\pi_0(\partial \Xtil)$ of connected components
of $\partial\Xtil$ into two subsets $b^\text{in}_\Xr$ and $b^\text{out}_\Xr$
(i.e.\ $b^\text{in}_\Xr \,{\cap}\, b^\text{out}_\Xr \eq \emptyset$ and
$b^\text{in}_\Xr \,{\cup}\, b^\text{out}_\Xr \eq \pi_0(\partial \Xtil)$).
The subsets $b^\text{in}_\Xr$ and $b^\text{out}_\Xr$ are required to be fixed
(as sets, not necessarily element-wise) under the involution ${\imath_\Xr}_*$
on $\pi_0(\partial \Xtil)$ that is induced by $\imath_\Xr$.
We denote the set of fix points of the induced map $:{\imath_\Xr}_*: \pi_0(\partial \Xtil) \to \pi_0(\partial \Xtil)$ by
$\text{op}_{\partial\Xr}$ and its complement in $\pi_0(\partial \Xtil)$ by $\text{cl}_{\partial\Xr}$. We also denote
$$
(b^\text{in}_\Xr)_\text{op}:=\text{op}_{\partial\Xr} \cap b^\text{in}_\Xr, \,\,
(b^\text{out}_\Xr)_\text{op}:=\text{op}_{\partial\Xr} \cap b^\text{out}_\Xr, \,\,
(b^\text{in}_\Xr)_\text{cl}:=\text{cl}_{\partial\Xr} \cap b^\text{in}_\Xr, \,\,
(b^\text{out}_\Xr)_\text{cl}:=\text{cl}_{\partial\Xr} \cap b^\text{out}_\Xr.
$$
\nxt $\text{or}_\Xr$ is a global section of
the bundle $\Xtil\,{\overset{\tilde\pi_\Xr}\longrightarrow}\,\dot\Xr$,
i.e.\ $\text{or}_\Xr {:}\ \dot\Xr \To \Xtil$ is a continuous
map such that $\tilde\pi_\Xr \cir \text{or}_\Xr \eq \id_{\dot\Xr}$.
In particular, a global section exists. We also denote the image of ${\text{or}_\Xr}_*: \pi_0(\partial \dot\Xr) \to \pi_0(\partial \Xtil)$ in $(b^\text{in}_\Xr)_\text{cl}$ by $(b^\text{in}_\Xr)_\text{cl}^+$ and its complement in $(b^\text{in}_\Xr)_\text{cl}$ by $(b^\text{in}_\Xr)_\text{cl}^-$.
\\[.3em]
\nxt A map
  \be
  \delta_\Xr{:}\quad \partial \Xtil \to S^1
  \ee
with the following properties. For $a \In \pi_0(\partial \Xtil)$ denote by
  \be
  \delta_a := \delta\big|_{a}
  \ee
the restriction of $\delta$ to the connected
component $a$ of $\partial\Xtil$. Then
\\
-- $\delta_a$ is an orientation preserving homeomorphism for each $a \In \pi_0(\partial \Xtil)$.
\\
-- $\delta_\Xr\cir\imath_\Xr \eq \text{Conj}\cir\delta_\Xr$.
\\
-- if $a \In \pi_0(\partial \Xtil)$ is fixed under ${\imath_\Xr}_*$, the
image of $S^1 \cap \{\text{Im}(z) \ge 0\}$ under $\delta_a^{-1}$ must
coincide with the image of $\text{or}_\Xr {:}\ \dot\Xr \To \Xtil$
restricted to $a$.
\\[.3em]
\nxt An order map: $\ord: \pi_0(\partial \Xtil) \to \{ 1, 2, 3, \dots, |\pi_0(\partial \Xtil)| \}$ is a bijection such that
\\
-- $\ord(a) < \ord(b)$ if $a\in \text{op}_{\partial\Xr}$ and $b \in \text{cl}_{\partial\Xr}$;
\\
-- Let $x$ be a connected component of the boundary of $\dot\Xr$ which has non-zero intersection with the fixed point set of $\imath_X$. Let further $U \subset \text{op}_{\partial\Xr}$ be the subset of boundary components of $\Xtil$ which get projected to $x$. We require the ordering restricted to $U$ to be $\{n,n+1,\dots, n+|U|-1\}$ for some $n$, increasing cyclically in the direction of the orientation of $\partial\dot\Xr$.
}
\end{defn}

\medskip
For a given world sheet $\Xr$, we will use either $\Xtil$ or $(\Xr)\widetilde{~}$ to denote the underlying surface.

\begin{rema}   \label{rem:unor} {\rm
(i) A stretch of boundary of the quotient surface $\dot\Xr$ can be of five different types, which is shown in the
following picture:

\begin{center}
\includegraphics[width=100mm]{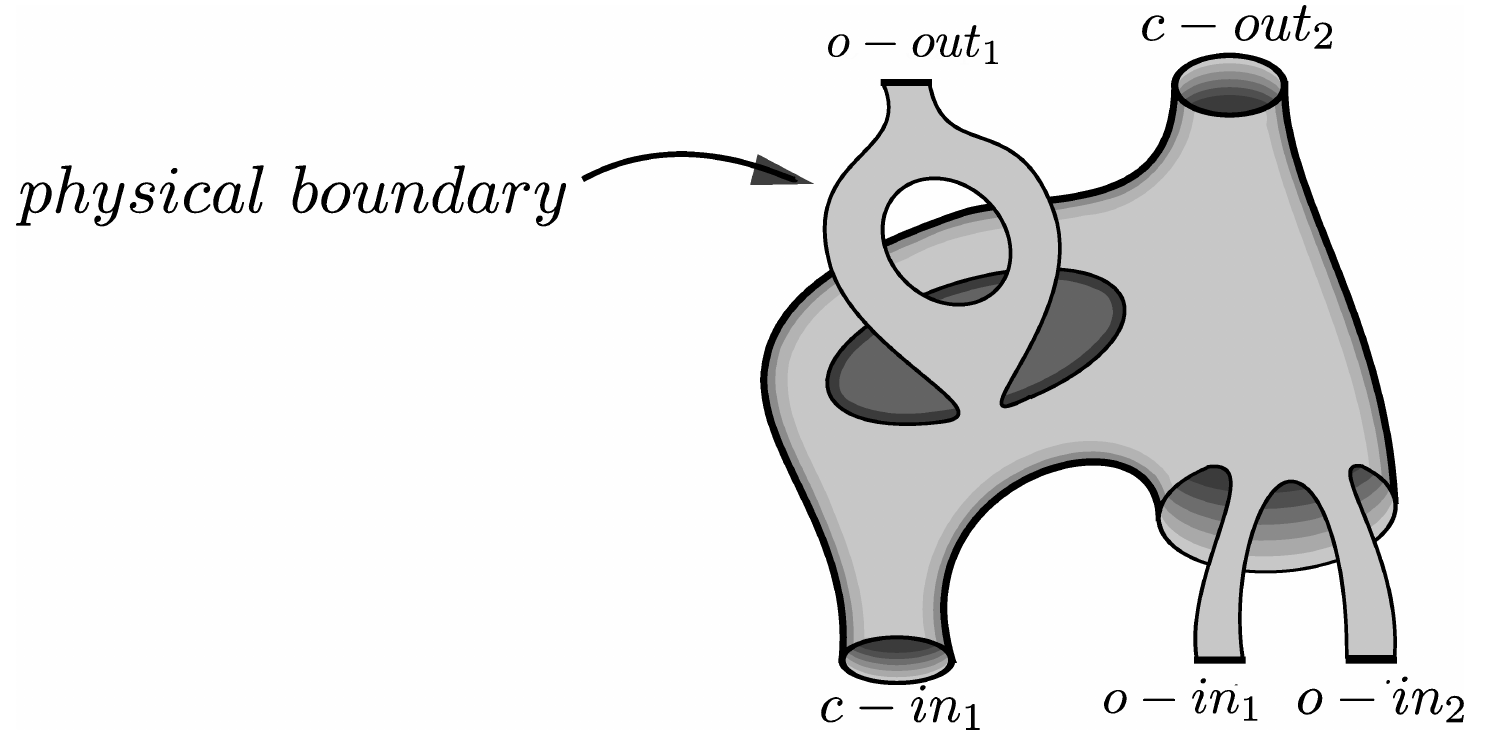}
\end{center}

\begin{list}{-}{\topsep .4em \leftmargin 2em \itemsep 0em}
\item {\em A physical boundary} is a connecting component of the set of fixed points of $\imath$.
\item {\em An in/out-going open state boundary} is the image of a boundary component $a\in (b^\text{in/out}_\Xr)_\text{op}$ of $\Xtil$ under the map $\tilde\pi_\Xr {:}\ \Xtil \To \dot\Xr$.
\item {\em An in/out-going closed state boundary} is the image of a boundary component in $(b^\text{in/out}_\Xr)_\text{cl}$ of $\Xtil$ under the map $\tilde\pi_\Xr {:}\ \Xtil \To \dot\Xr$.
\end{list}
(ii) The first condition on $\delta$ in Definition \ref{def:ws} means that $\delta_a^{-1}$ provides an orientation preserving
 parametrisation of the boundary component $a$ of $\partial\Xtil$ by the unit circle $S^1$. The second condition says that
$\delta$ intertwines the involutions on $\Xtil$ and $\Cb$, while the third condition ensures that in the next definition the gluing of
 open state boundaries is compatible with the orientation.
\\[.3em]
(iii) If one wants to treat unoriented open/closed topological world sheets instead, the section $\text{or}_\Xr$ would be
dropped from the defining data, as well as the third conditions on $\delta_\Xr$.
}
\end{rema}

\begin{figure}[bt]
\vspace{2em}
\begin{center}
\includegraphics[width=100mm]{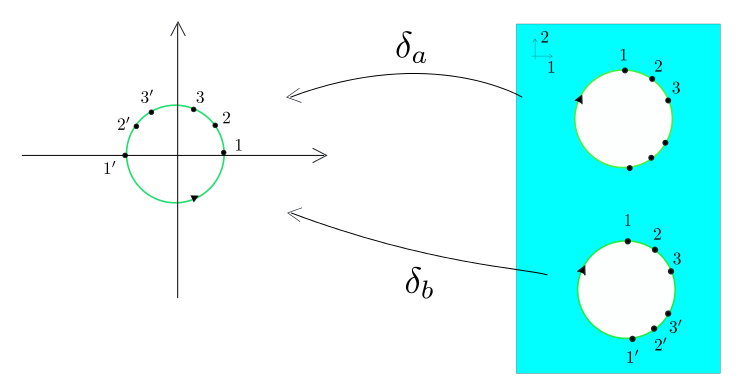}
\end{center}
\vspace{2em}
\caption{Identification of points on $\partial \Xr$ in the sewing procedure. The figure shows part of a world sheet with bulk and boundary orientations as indicated. The boundary parametrisations $\delta_a$ and $\delta_b$ are such that they map the point 1 on $\partial \Xr$ to 1 on the standard $S^1$, as well as 2 to 2 and 3 to 3. When sewing boundaries $a$ and $b$, the points 1,2 and 3 on $a$ get identified with $1'$, $2'$, and $3'$ on $b$, respectively. Note that this is compatible with the world sheet orientation.}
\label{fig:glue-example}
\end{figure}

\begin{defn}   \label{def:oc-sewing} 
{\rm
Let $\Xr \eq \big( \Xtil, \imath, \delta, b^\text{in},
b^\text{out}, \text{or}, \ord \big)$ and $\Yr$ be two world sheets.
\\[.3em]
(i)~\,{\em Sewing data for\/} $\Xr$, or a {\em sewing of\/} $\Xr$,
is a (possibly empty) subset $\sew$
   of $b^\text{out} \ti b^\text{in}$ such that if $(a,b) \In \sew$ then\\
-- $\sew$ does not contain any other pair of the form $(a,\cdot)$ or
   $(\cdot,b)$,\\
-- also $(\imath_*(a),\imath_*(b)) \In \sew$,\\
-- the pair $(a,b)$ either lies in $(b^\text{out})_\text{op} \times (b^\text{in})_\text{op}$ or
$(b^\text{out})_\text{cl}^+ \times (b^\text{in})_\text{cl}^+$ or
$(b^\text{out})_\text{cl}^- \times (b^\text{in})_\text{cl}^-$.
\\[.3em]
(ii)~For a sewing $\sew$ of $\Xr$, the {\em sewn world sheet\/} $\sew(\Xr)$ is
the tuple
$$
\sew(\Xr) \,{\equiv}\, \big( \Xtil',\imath',\delta',{b^\text{in}}',
{b^\text{out}}',\text{or}', \ord' \big),
$$
where
$\Xtil' \,{:=}\, \Xtil / {\sim}$ with $\delta_a^{-1}(z) \,{\sim}\,
\delta_b^{-1}\cir \text{Conj}(-z)$ for all $(a,b) \In \sew$ and $z \In S^1$ (see Figure \ref{fig:glue-example}). Next,
denote by $\pi_{\sew,\Xr}$ the projection from $\Xtil$ to $\Xtil'$ that takes a
point of $\Xtil$ to its equivalence class in $\Xtil'$. Then $\imath'{:}\ \Xtil'
\To \Xtil'$ is the unique involution
such that $\imath'\cir\pi_{\sew,\Xr} \eq \pi_{\sew,\Xr}\cir\imath$. Further,
$\delta'$ is the restriction of $\delta$ to $\partial\Xtil'$,
${b^\text{out}}'\eq \{ a \In b^\text{out}| (a,\cdot) \,{\notin}\, \sew \}$,
${b^\text{in}}' \eq \{ b \In b^\text{in} | (\cdot,b) \,{\notin}\, \sew \}$,
and $\text{or}'$ is the unique continuous section of
$\Xtil'\,{\overset{\tilde\pi_{\Xr'}}{\longrightarrow}}\,{\dot\Xr}'$ such that
the image of $\,\text{or}'$ coincides with the image of
$\pi_{\sew,\Xr} \cir \text{or}$. The order map $\ord'$ restricting on $({b^\text{in/out}}')_\text{cl}$ for $\sew(\Xr)$ is defined by shifting the elements in $({b^\text{in/out}}')_\text{cl}$ as a subset of $({b^\text{in/out}})_\text{cl}$ to fill the empty spots left by those involved in sewing data $\sew$ so that the original order in $(b^\text{in/out})_\text{cl}$ is preserved. The order $\ord'$ of elements in $({b^\text{in/out}}')_\text{op}$ is always smaller than that of the elements in $({b^\text{in/out}}')_\text{cl}$ and is defined the following way:
\begin{list}{-}{\topsep .4em \leftmargin 2em \itemsep 0em}
\item We first locate the element $a_0\in ({b^\text{in/out}}')_\text{op}$ such that
$$
\ord(a_0)=\text{min}(\{ \ord(a) | a\in ({b^\text{in/out}}')_\text{op}\}).
$$
We define $\ord'(a_0):=1$; then define the order $\ord'$ for those lying in the same boundary component of $\partial \dot\Xr'$ according to the induced orientation of $\partial \dot\Xr'$.
\item We repeat the procedure for the rest of $a\in ({b^\text{in/out}}')_\text{op}$ until we exhaust all elements in $({b^\text{in/out}}')_\text{op}$.
\end{list}
This completes the definition of the sewn world sheet $\sew(\Xr)$.
\\[.3em]
(iii) A {\em homeomorphism of world sheets\/} is a homeomorphism $f{:}\ \Xtil
\To \widetilde\Yr$ that is compatible with all chosen structures on $\Xtil$,
i.e.\ with orientation, involution, boundary parametrisation and the section $\text{or}:\dot\Xr\rightarrow\Xtil$ and $\ord$. That is, $f$ satisfies
  \be
  f \circ \imath_\Xr = \imath_\Yr \circ f ~,~~~
  \delta_\Yr \circ f = \delta_\Xr ~,~~~
  f_* b^\text{in/out}_\Xr = b^\text{in/out}_\Yr ~,~~~
  f \circ \text{im}(\text{or}_\Xr)
  = \text{im}(\text{or}_\Yr) ~,
  \ee
where $f_*{:}\ \pi_0(\partial\Xtil) \To \pi_0(\partial\widetilde\Yr)$ is  an order preserving bijection induced by $f$.
}
\end{defn}

\begin{rema}  \label{rem:sewing-data}  {\rm
The conditions in part (i) of Definition \ref{def:oc-sewing} are necessary for the sewing procedure in part (ii) to result in a valid world sheet. The first condition ensures that an in-coming state boundary is glued to precisely one out-going state-boundary, the second condition allows $\imath$ to be lifted to the sewn world sheet, and the third condition ensures that the sewing is compatible with the orientation. The order of the sewn world sheet is defined so that the sewing operation is associative.
}
\end{rema}

The category of open-closed topological world sheets $\WSh$ is the
symmetric monoidal category defined as follows.
\\[.2em]
\nxt The objects are world sheets.
\\[.2em]
\nxt Let $\Xr$ and $\Yr$ be two world sheets.
A {\em morphism\/} $\varpi {:}\ \Xr \To \Yr$ is a pair $\varpi \eq (\sew,f)$
where $\sew$ is a sewing of $\Xr$ and
$f{:}\ \widetilde{\sew(\Xr)} \To \widetilde\Yr$
is a homeomorphism of world sheets.
The set of all morphisms from $\Xr$ to $\Yr$ is denoted by $\Hom(\Xr,\Yr)$.
\\[.2em]
\nxt
Given two morphisms $\varpi \eq (\sew,f){:}\ \Xr \To \Yr$ and $\varpi'\eq (\sew',g){:}
\ \Yr \To \Zr$, the composition $\varpi'\cir\varpi$ is defined as follows. The union
$\sew'' \eq \sew \,{\cup}\, (f \cir\pi_{\sew,\Xr})_*^{\,-1}(\sew')$
is again a sewing of $\Xr$. Furthermore there exists a unique isomorphism
$h{:}\ \widetilde{\sew''(\Xr)} \To \widetilde{\Zr}$ such that the diagram
  \be
  \xymatrix@C=2em{
  \Xtil \ar[d]^{\pi_{\sew'',\Xr}} \ar[r]^{\pi_{\sew,\Xr}} &
  \widetilde{\sew(\Xr)} \ar[r]^{f} &
  \widetilde{\Yr} \ar[r]^{\pi_{\sew',\Yr}} & \widetilde{\sew'(\Yr)} \ar[r]^{g}
  & \widetilde{\Zr} \\
  \widetilde{\sew''(\Xr)} \ar[urrrr]_{h} }
  \ee
commutes. We define the composition
$\circ\,{:}\ \Hom(\Yr,\Zr) \ti \Hom(\Xr,\Yr) \To \Hom(\Xr,\Zr)$ as
  \be
  (\sew',g) \circ (\sew,f) = (\sew'',h) \,.
  \ee
One verifies that the composition is associative. The identity morphism on
$\Xr$ is the pair $\id_\Xr \eq (\emptyset,\id_{\Xtil})$.
\\[.2em]
\nxt The tensor product $\Xr \oti \Yr$ is given by the disjoint union. This means that the surface underlying $\Xr \oti \Yr$ is
\be
 \big(\Xr \oti \Yr\big)\widetilde{~} = \Xtil \sqcup \widetilde\Yr =
\{ (x,1) | x \in \Xtil \} \cup \{ (y,2) | y \in \widetilde\Yr \}~,
\ee
and the other data is transported accordingly. (The notation $(~)\widetilde{~}$ for the underlying surface was introduced below Definition \ref{def:ws}.)
Note that the indexing implied by the disjoint
union makes the tensor product non-strict: We have
\bea
\big(\Xr \oti (\Yr \oti \Zr)\big)\widetilde{~} =
\big\{ (x,1) \,\big|\, x \In \Xtil \big\} \cup
\big\{ ((y,1),2) \,\big|\, y \In \widetilde\Yr \big\} \cup
\big\{ ((z,2),2) \,\big|\, z \In \Xtil \big\} \quad \text{and}
\enl
\big((\Xr \oti \Yr) \oti \Zr\big)\widetilde{~} =
\big\{ ((x,1),1) \,\big|\, x \In \Xtil \} \cup
\big\{ ((y,2),1) \,\big|\, y \In \widetilde\Yr \} \cup
\big\{ (z,2) \,\big|\, z \In \Xtil \}~.
\eear\ee
The associator is then the homeomorphism which takes $(x,1)$ to $((x,1),1)$, etc. The unit object is the empty set (and the unit isomorphisms consist of
forgetting the indexing of the disjoint union).
\\[.2em]
\nxt The symmetric braiding isomorphism
$c_{\Xr,\Yr}{:}\ \Xr \oti \Yr \To \Yr \oti \Xr$
is the homeomorphism that exchanges the two
factors of the disjoint union of the underlying surfaces,
$c_{\Xr,\Yr}\big((x,1)\big) = (x,2)$ and
$c_{\Xr,\Yr}\big((y,2)\big) = (y,1)$
where $x \In \Xtil$ and $y \In \widetilde\Yr$. It is obvious that $c_{\Xr,\Yr}\circ c_{\Yr,\Xr}=id_{\Xr\otimes\Yr}$.

\begin{defn}    \label{def:homotop-morph}  {\rm
Let $\Xr,\Yr\in\WSh$ and let $\xi_0,\xi_1 : \Xr \rightarrow \Yr$
be morphisms of the form $\xi_0 = (\sew_0,f_0)$ and $\xi_1=(\sew_1,f_1)$. We call $\xi_0$ and $\xi_1$ {\em homotopic}
 if $\sew_0=\sew_1$ and if there is a homotopy of maps $\{ f_t | t \in [0,1] \}$ between $f_0$ and $f_1$ such that for
each $t$, $f_t$ is a homeomorphism of world sheets.
}
\end{defn}

\subsection{Symmetric monoidal functors from $\WSh$ to $\Vect$}

Denote by $\Vect$ the symmetric monoidal category of finite-dimensional complex vector spaces and let $\Fun_\otimes^H(\WSh,\Vect)$
be the following category. The objects of $\Fun_\otimes^H(\WSh,\Vect)$ are symmetric monoidal functors $F : \WSh \rightarrow \Vect$
with the property that for any two homotopic morphisms $\varpi, \varpi' : \Xr \rightarrow \Yr$ we have $F(\varpi) = F(\varpi')$ (hence the superscript `$H$' in $\Fun_\otimes^H(\WSh,\Vect)$). A
morphism $\mu : F \rightarrow G$ is a monoidal natural transformation from $F$ to $G$.
The set of morphisms will be denoted by
$\text{Nat}_\otimes(F,G)$.

An example of an object in $\Fun_\otimes^H(\WSh,\Vect)$ is the functor $\One : \WSh \rightarrow \Vect$. It is defined on objects
and morphisms as
\be
  \One(\Xr) = \Cb
  \quad , \quad \One( \Xr\overset{\varpi}{\longrightarrow}\Yr ) = \id_\Cb
  ~~,
\ee
together with the isomorphism $\phi_0 : \Cb \rightarrow \One(\emptyset)$, $\phi_0 = \id_\Cb$ and the natural isomorphism
 $\phi_2 : \One(\Xr) \otimes_\Cb \One(\Yr) \rightarrow \One(\Xr \oti \Yr)$, $\phi_2(z \oti w) = zw$.
One verifies that this is a symmetric monoidal functor.

In fact, $\Fun_\otimes^H(\WSh,\Vect)$ inherits a tensor product from $\Vect$, and $\One$ is the tensor unit with respect to this tensor product. But we will not need this extra structure here.

\subsection{A generating set of world sheets}

By a generating set of world sheets we mean a set
$\{ \Xr_\alpha | \alpha \In \Sc \}$ of world sheets such that for each
world sheet $\Xr$ there exists
a (typically non-unique) list $(\alpha_1,\dots,\alpha_n)$ together with
a (typically non-unique) morphism
$\varpi : \Xr_{\alpha_1} \oti \cdots \oti \Xr_{\alpha_n} \rightarrow \Xr$.
We will use the following generating set. The index set consists of the 12 elements
\be
  \Sc = \{ mo, \Delta o, \eta o, \eps o, mc, \Delta c, \eta c, \eps c,
  po, pc, \iota, \iota^* \}~,
\labl{eq:standard-F}
and the corresponding world sheets are
(for a world sheet $\Xr$ we draw the quotient surface $\dot\Xr$)
\be\begin{array}{lllll}
  \Xr_{mo} =
  \raisebox{-23pt}{\begin{picture}(60,46)
  \put(0,0){\scalebox{0.15}{\includegraphics{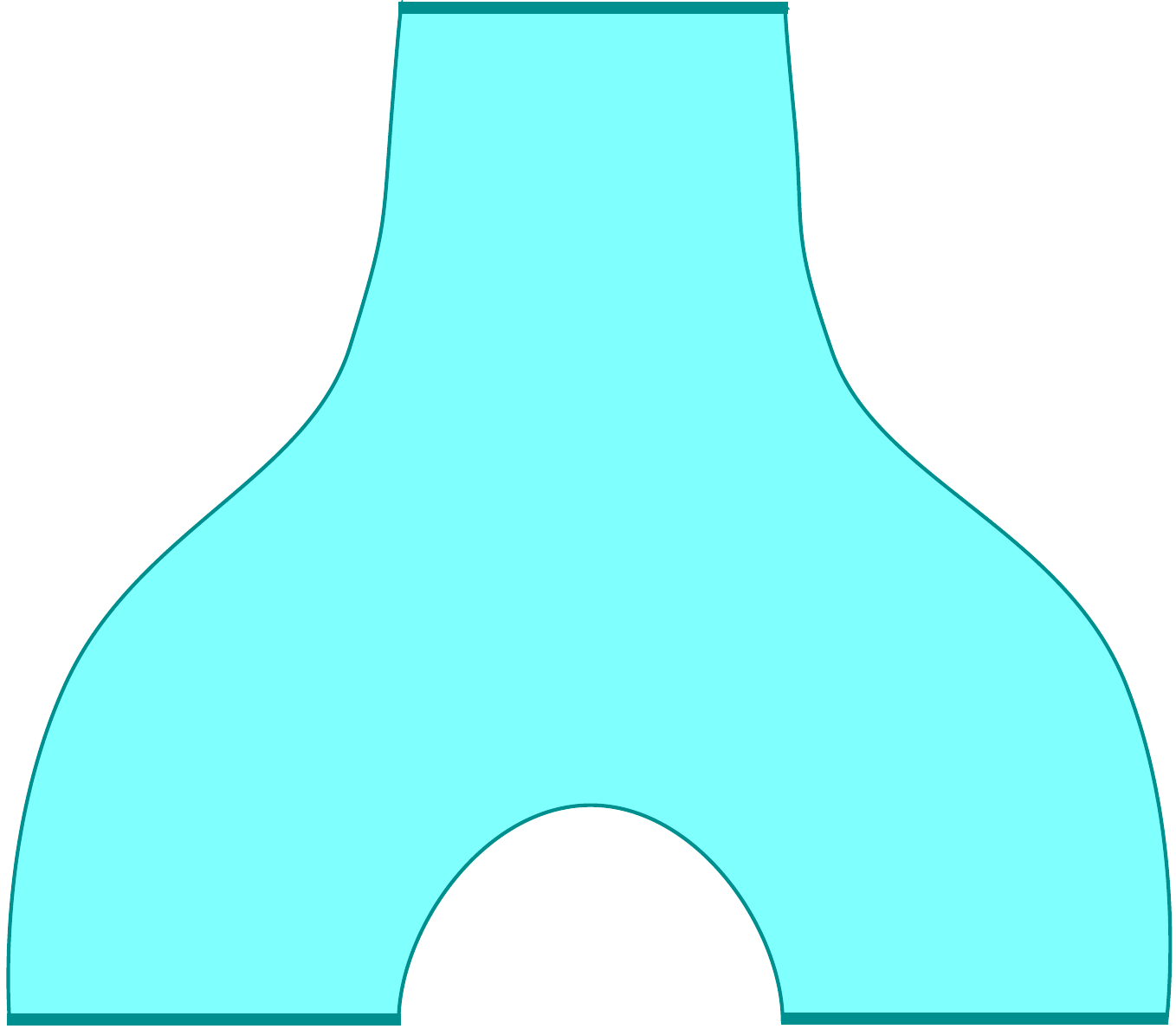}}}\end{picture}}
\etb
  \Xr_{\Delta o} =
  \raisebox{-23pt}{\begin{picture}(60,46)
  \put(0,0){\scalebox{0.15}{\includegraphics{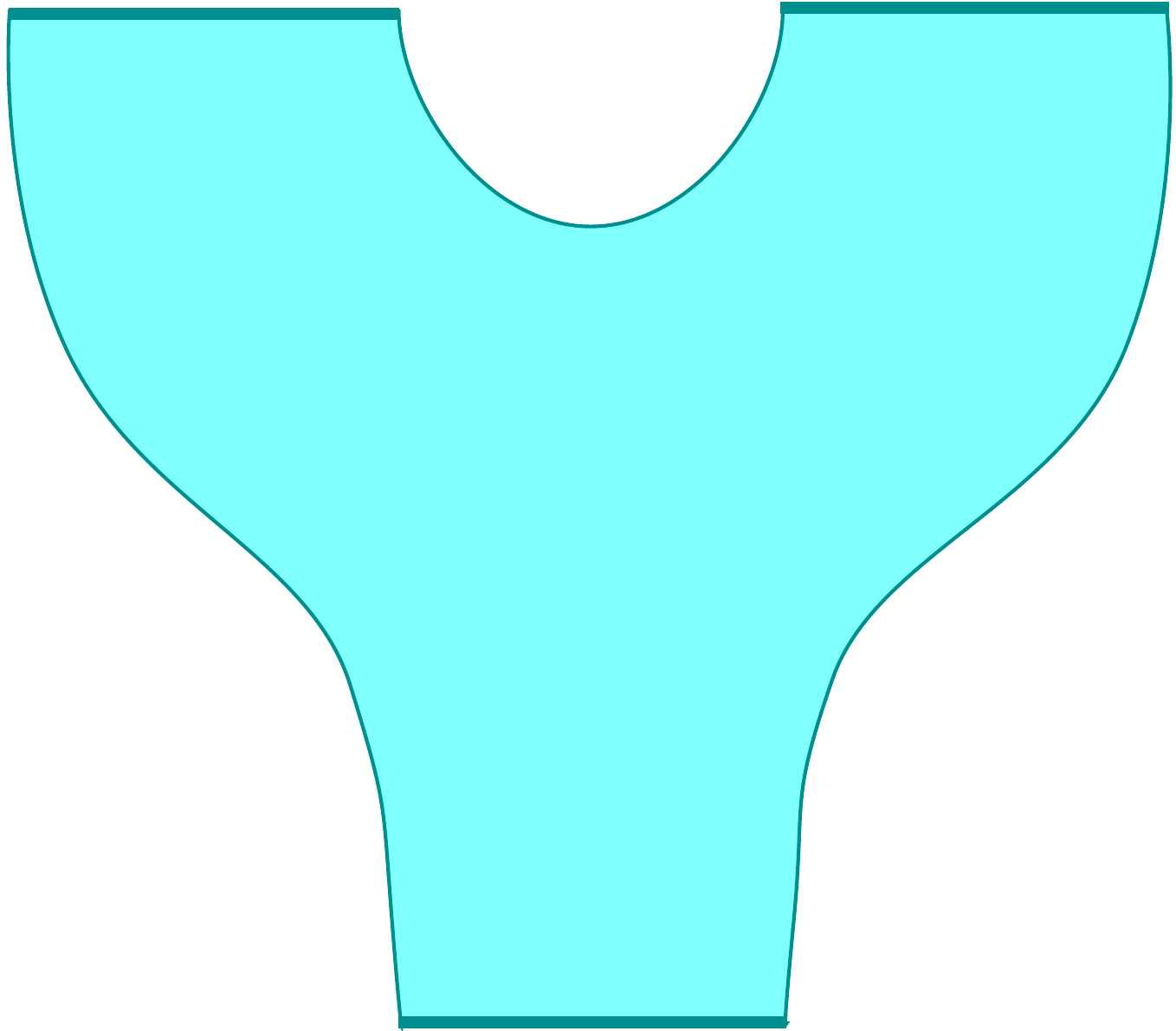}}}\end{picture}}
\etb
  \Xr_{\eta o} =
  \raisebox{-5pt}{\begin{picture}(23,15)
  \put(0,2){\scalebox{0.15}{\includegraphics{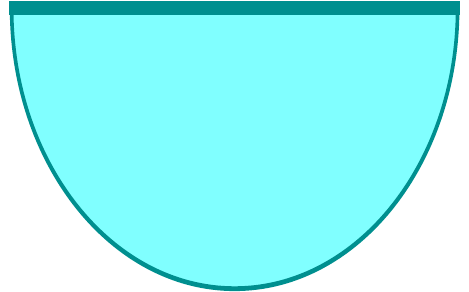}}}\end{picture}}
\etb
  \Xr_{\eps o} =
  \raisebox{-5pt}{\begin{picture}(23,15)
  \put(0,2){\scalebox{0.15}{\includegraphics{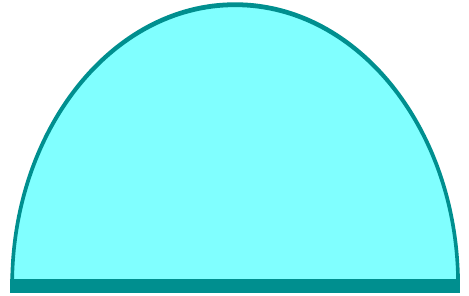}}}\end{picture}}
\\[3em] \displaystyle
  \Xr_{mc} =
  \raisebox{-23pt}{\begin{picture}(60,46)
  \put(0,0){\scalebox{0.15}{\includegraphics{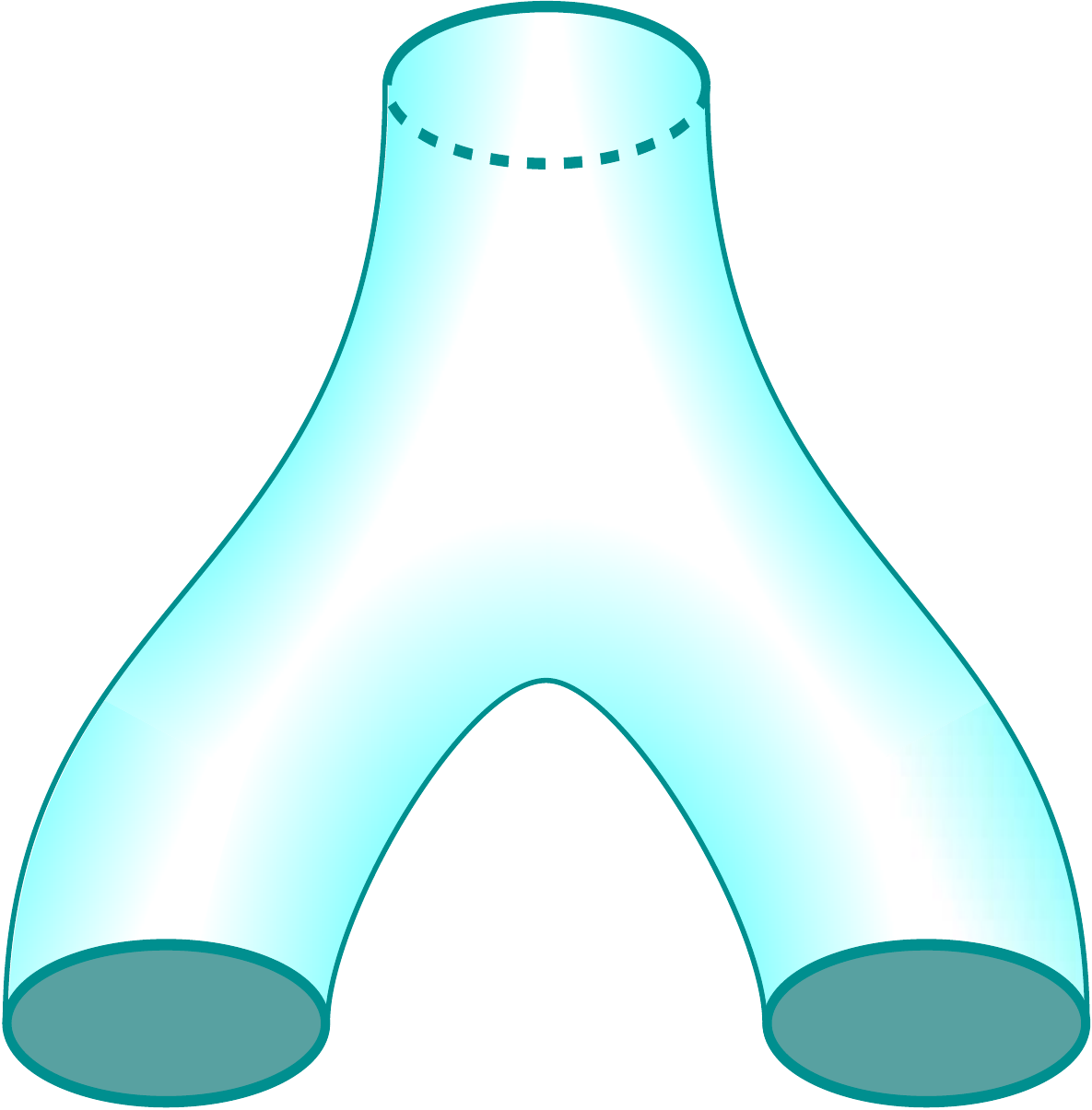}}}\end{picture}}
\etb
  \Xr_{\Delta c} =
  \raisebox{-23pt}{\begin{picture}(60,46)
  \put(0,0){\scalebox{0.15}{\includegraphics{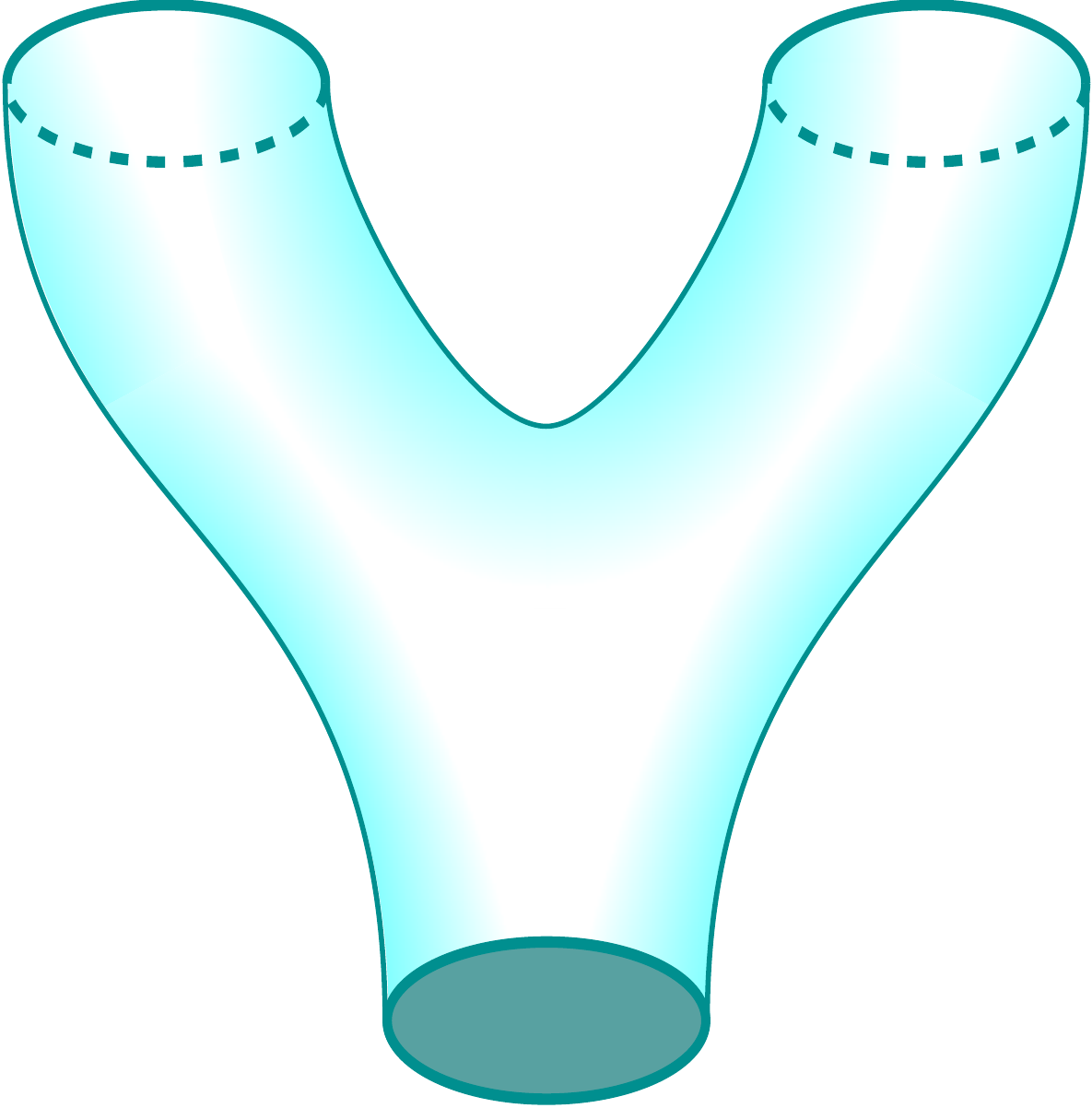}}}\end{picture}}
\etb
  \Xr_{\eta c} =
  \raisebox{-13pt}{\begin{picture}(37,30)
  \put(0,0){\scalebox{0.15}{\includegraphics{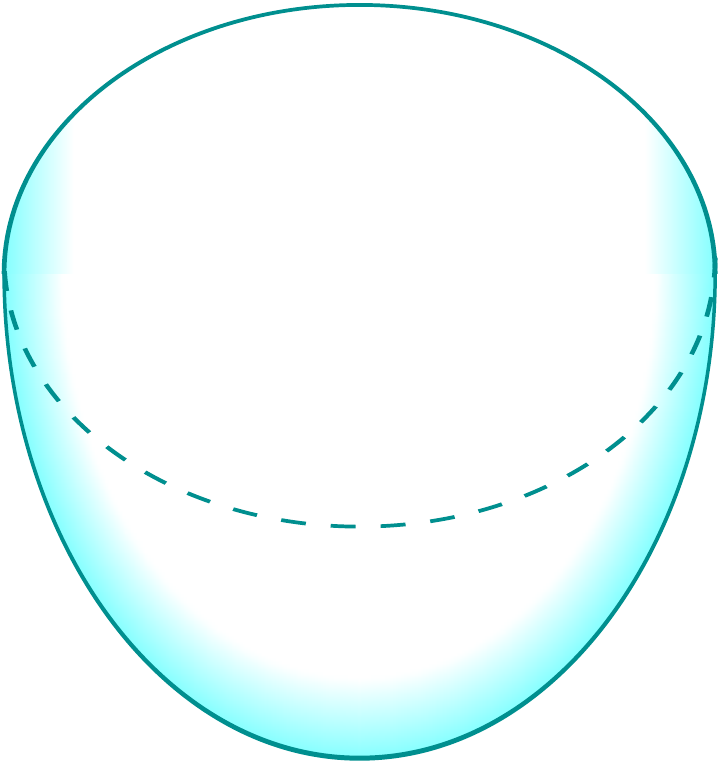}}}\end{picture}}
\etb
  \Xr_{\eps c} =
  \raisebox{-13pt}{\begin{picture}(37,30)
  \put(0,0){\scalebox{0.15}{\includegraphics{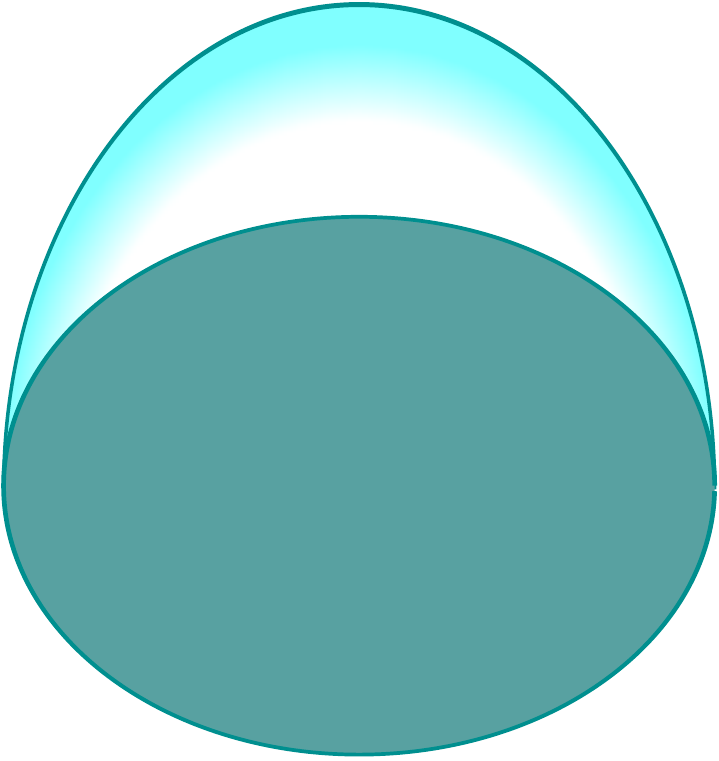}}}\end{picture}}
\\[3em] \displaystyle
  \Xr_{po} =
  \raisebox{-23pt}{\begin{picture}(26,46)
  \put(0,0){\scalebox{0.15}{\includegraphics{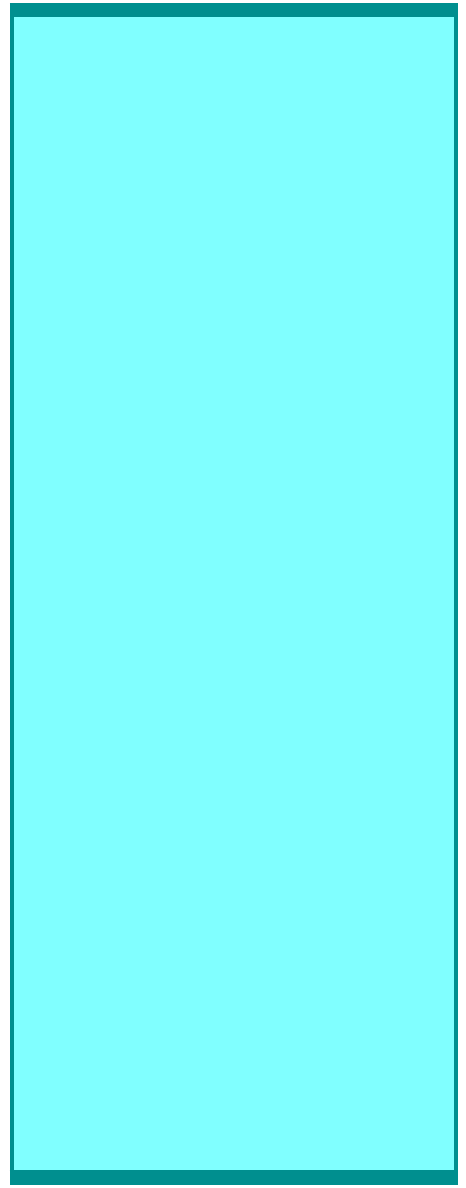}}}\end{picture}}
\etb
  \Xr_{pc} =
  \raisebox{-29pt}{\begin{picture}(26,58)
  \put(0,0){\scalebox{0.15}{\includegraphics{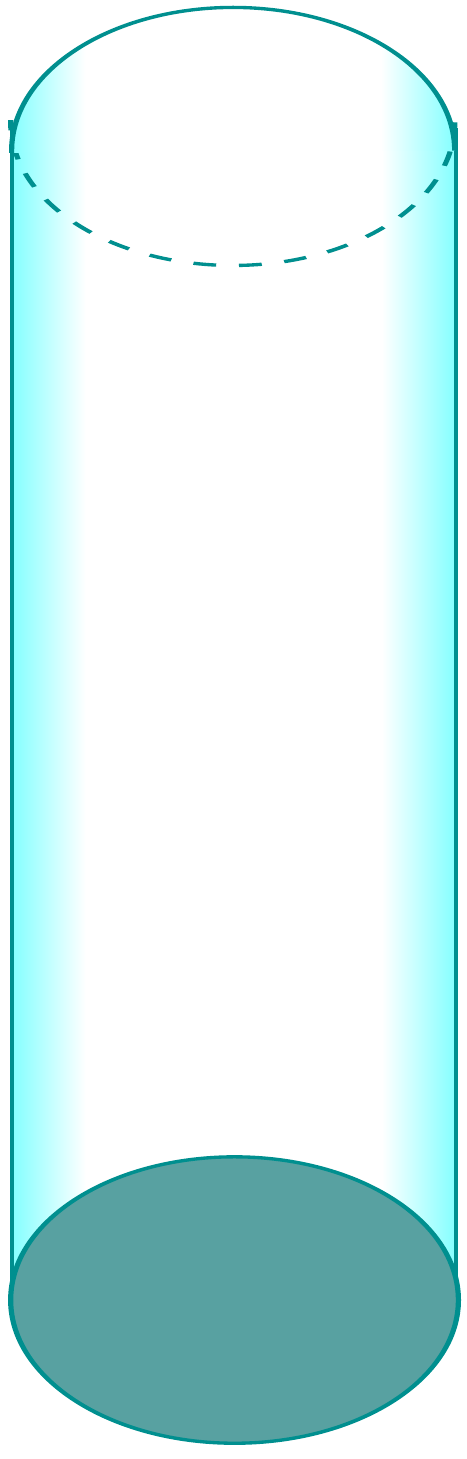}}}\end{picture}}
\etb
  \Xr_{\iota} =
  \raisebox{-26pt}{\begin{picture}(26,52)
  \put(0,0){\scalebox{0.15}{\includegraphics{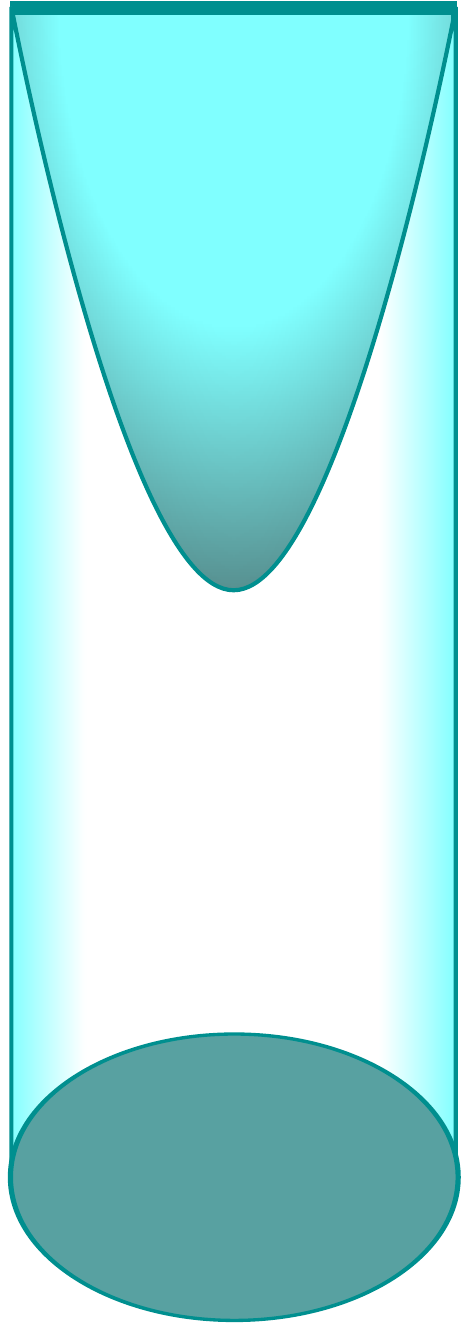}}}\end{picture}}
\etb
  \Xr_{\iota^*} =
  \raisebox{-26pt}{\begin{picture}(26,52)
  \put(0,0){\scalebox{0.15}{\includegraphics{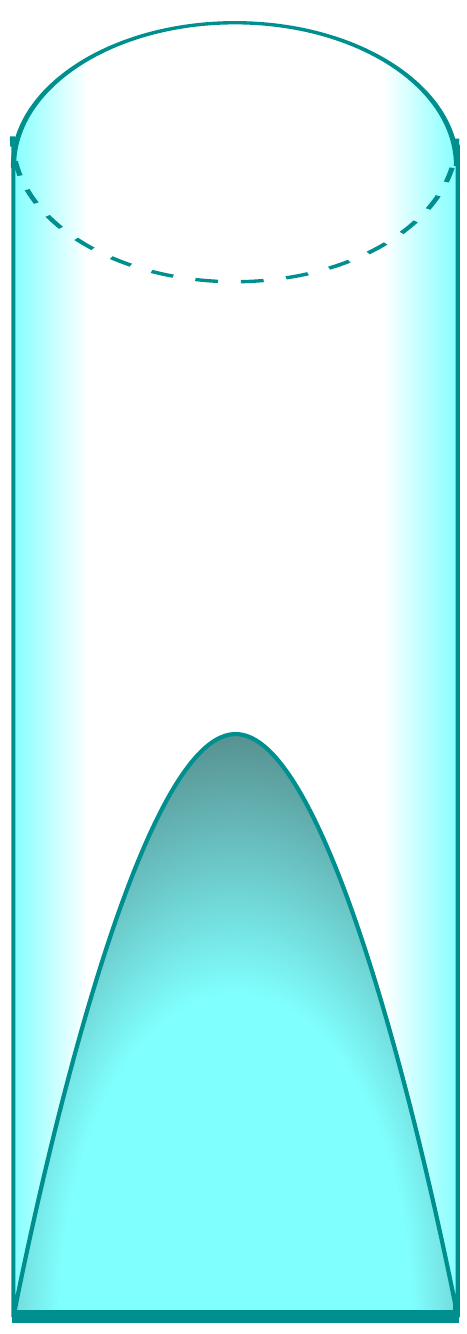}}}\end{picture}}
\end{array}
\labl{eq:standard-gen-ws}
In these pictures, state boundaries at the bottom are in-coming and state-boundaries at the top are out-going.
For each of the world sheets $\Xr_\alpha$, $\alpha \in \Sc$, it is understood (but not indicated in the pictures) that a particular choice for the parametrisation of the state boundaries has been made once and for all. A concrete choice for $\Xr_{\eta o}$, $\Xr_{po}$, and $\Xr_{mo}$ will be described in Section \ref{sec:R1-R13} below. It is enough to pick one parametrisation for the state boundaries to obtain a generating set.\footnote{
   For example, if $\Xr_{po}'$ coincides with $\Xr_{po}$ up to the parametrisation $\delta$, one can always find a morphism $(\emptyset,f) : \Xr_{po} \rightarrow \Xr_{po}'$ by first choosing the unique homeomorphism $f_\partial : \partial \Xtil_{po} \rightarrow\partial \Xtil_{po}$ such that $\delta_{\Xr_{po}} : \partial \Xtil_{po} \rightarrow S^1$ is equal to $\delta_{\Xr_{po}'} \circ f_\partial$ and such that $f_\partial$ maps a given connected component of $\partial \Xtil_{po}$ back to itself. The map $f_\partial$ is homotopic to the identity, and this homotopy can be extended to a neighbourhood of $\Xtil_{po}$, giving an extension of $f_\partial$ to a homeomorphism from $\Xtil_{po}$ to $\Xtil_{po}$ with the desired properties.}

Note that the generating set $\Sc$ is in fact overcomplete, for example one could omit
$\Xr_{\eta o}$, $\Xr_{po}$, $\Xr_{pc}$ and $\Xr_{\iota^*}$. It is, however, convenient to include them because this makes the formulation of the relations in Section \ref{sec:gen-rel} easier.

	The following statement can be found in \cite[App.\,A]{MS} or \cite[Prop.\,3.9]{LP}. 

\begin{prop} \label{prop:gen-set} {\rm
The world sheets in \erf{eq:standard-gen-ws} are a generating set.
}
\end{prop}

\subsection{Natural transformations via
generators and relations}\label{sec:gen-rel}

Let $\{\Xr_\alpha|\alpha\in\mathcal{F}\}$ be a fixed generating set of world sheets.
The category of {\em decomposed open-closed topological world sheets} $\WSh^d$ for the given generating set of world sheets
is defined as follows. As objects it has triples $\Xr^d = (\Xr , (\alpha_1,\dots,\alpha_m), \varpi)$, where $\Xr$ is an object
of $\WSh$, $(\alpha_1,\dots,\alpha_m)$ is an ordered list of elements of $\Sc$, and
\be   \label{eq:varpi}
  \varpi : \big( \cdots \big( (\Xr_{\alpha_1} \oti \Xr_{\alpha_2})
  \oti
  \Xr_{\alpha_3}\big) \oti \cdots \big) \oti \Xr_{\alpha_m}
  \longrightarrow \Xr ~.
\ee
Let $\Xr^d = (\Xr , (\alpha_1,\dots,\alpha_m), \varpi_\Xr)$ and $\Yr^d = (\Yr , (\beta_1,\dots,\beta_n), \varpi_\Yr)$
be two decomposed world sheets.
The set of morphisms $\Hom^d(\Xr^d,\Yr^d)$ is defined to be empty unless $(\alpha_1,\dots,\alpha_m) = (\beta_1,\dots,\beta_n)$. If the two lists are equal the set of morphisms consists of all morphisms $f : \Xr \rightarrow \Yr$ in $\WSh$ such that $f \cir \varpi_\Xr = \varpi_\Yr$. Note that if $f$ exists it is unique.

The tensor product and symmetry operation on $\WSh^d$ are again given by disjoint union and exchange of factors.
As we did for $\WSh$ we will consider the category of symmetric monoidal functors $\Fun_\otimes^H(\WSh^d,\Vect)$
and monoidal natural transformations between
them.

\medskip

The forgetful functor $\text{For}^d : \WSh^d \rightarrow \WSh$ takes an object $\Xr^d$ to the underlying object $\Xr$ of $\WSh$
 and acts as the identity on morphisms. It is surjective, faithful, and symmetric monoidal.
It also gives rise to a functor $(\,\cdot\,)^\delta : \Fun_\otimes^H(\WSh,\Vect) \rightarrow \Fun_\otimes^H(\WSh^d,\Vect)$.
 On objects it acts as $F^\delta = F \cir \text{For}^d$, and to a morphism $\mu \In \text{Nat}_\otimes(F,G)$ it assigns the
 morphism $\mu^\delta \In \text{Nat}_\otimes(F^\delta,G^\delta)$ given by the monoidal natural transformation
$\mu^\delta_{\Xr^d} = \mu_\Xr : F^\delta(X^d) \rightarrow G^\delta(X^d)$.

\medskip

Let $U,F \In \Ob\Fun_\otimes^H(\WSh,\Vect)$ and suppose that $U(f)$ is an invertible linear map for all morphisms $f$ in $\WSh$.
For example, this is the case for $\One$, but will in general not be the case for the functor $\Bl$ constructed in Section
 \ref{sec:Bl-def} below. Let
\be
  \{ c_\alpha : U(\Xr_\alpha) \rightarrow F(\Xr_\alpha) | \alpha \In \Sc \}
\ee
be a collection of linear maps, one for each world sheet in the generating set $\Sc$. Given a decomposed world sheet $\Xr^d = (\Xr , (\alpha_1,\dots,\alpha_m), \varpi_\Xr)$, define the linear map $C_{\Xr^d}( \{ c_\alpha | \alpha \In \Sc \} ) : U(\Xr) \rightarrow F(\Xr)$ as the composition
\bea
U(\Xr)
\xrightarrow{U(\varpi)^{-1}}
U\big(\big( \cdots \big( (\Xr_{\alpha_1} \oti \Xr_{\alpha_2})
  \oti
  \Xr_{\alpha_3}\big) \oti \cdots \big) \oti \Xr_{\alpha_m}\big)
\enl
\xrightarrow{(\phi_2^U)^{-1}}
U\big( \cdots ( (\Xr_{\alpha_1} \oti \Xr_{\alpha_2})
  \oti \cdots ) \oti \Xr_{\alpha_{m-1}} \big) \oti U(\Xr_{\alpha_m})
\enl
\xrightarrow{(\phi_2^U)^{-1} \oti \id} \cdots
\longrightarrow
\big( \cdots \big( (U(\Xr_{\alpha_1}) \otimes_\Cb U(\Xr_{\alpha_2}))
  \otimes_\Cb
  U(\Xr_{\alpha_3})\big) \otimes_\Cb \cdots \big) \otimes_\Cb U(\Xr_{\alpha_m})
\enl
\xrightarrow{(\cdots (c_{\alpha_1} \otimes_\Cb c_{\alpha_2})
\otimes_\Cb \cdots ) \otimes_\Cb c_{\alpha_m}}
\big( \cdots \big( (F(\Xr_{\alpha_1}) \otimes_\Cb F(\Xr_{\alpha_2}))
  \otimes_\Cb
  F(\Xr_{\alpha_3})\big) \otimes_\Cb \cdots \big) \otimes_\Cb F(\Xr_{\alpha_m})
\enl
\xrightarrow{\phi_2^F \otimes_\Cb \id}
\big( \cdots \big( F(\Xr_{\alpha_1} \oti \Xr_{\alpha_2})
  \otimes_\Cb
  F(\Xr_{\alpha_3})\big) \otimes_\Cb \cdots \big) \otimes_\Cb F(\Xr_{\alpha_m})
\enl
\longrightarrow \cdots \xrightarrow{\phi_2^F}
F\big(\big( \cdots \big( (\Xr_{\alpha_1} \oti \Xr_{\alpha_2})
  \oti
  \Xr_{\alpha_3}\big) \oti \cdots \big) \oti \Xr_{\alpha_m}\big)
\xrightarrow{F(\varpi)}
F(\Xr) \ .
\eear\ee
Because the expressions get somewhat cumbersome, we will from here on no longer spell out the isomorphisms $\phi_2$. We will also omit any $\phi_0$'s and associators. The above formula then becomes
\be
  C_{\Xr^d}( \{ c_\alpha | \alpha \In \Sc \} )
  = F(\varpi) \circ ( c_{\alpha_1} \otimes_\Cb \cdots \otimes_\Cb c_{\alpha_m} )
  \circ U(\varpi)^{-1}
  ~.
\labl{eq:nat-C-def}
Note that we need $U(f)$ to be invertible for morphisms $f$ of $\WSh$ for this definition to make sense. As shown in the following lemma, the
linear maps \erf{eq:nat-C-def} provide all monoidal natural
transformations from $U^\delta$ to $F^\delta$.

\begin{lemma}   \label{lem:all-nat-C}  {\rm
$C(\{ c_\alpha | \alpha \In \Sc \} ) \In \text{Nat}_\otimes(U^\delta,F^\delta)$, and all elements in $\text{Nat}_\otimes(U^\delta,F^\delta)$ are of this form for a suitable collection of linear maps $\{ c_\alpha | \alpha \In \Sc \}$.
}
\end{lemma}
\noindent \pf
Let $\Xr^d = (\Xr , (\alpha_1,\dots,\alpha_m), \varpi_\Xr)$
and $\Yr^d = (\Yr , (\alpha_1,\dots,\alpha_m), \varpi_\Yr)$ be two decomposed world sheets with the same collection of generators $\alpha_i \in \Sc$ so
that we can have a morphism from $\Xr^d$ to $\Yr^d$.
To show that $C \equiv C( \{ c_\alpha | \alpha \In \Sc \} )$ is natural
we need to check that for a given $f : \Xr^d \rightarrow \Yr^d$ the diagram
  \be
  \xymatrix{
  U^\delta(\Xr^d) \ar[d]^{C_{\Xr^d}} \ar[r]^{U^\delta(f)} & U^\delta(\Yr^d) \ar[d]^{C_{\Yr^d}} \\
  F^\delta(\Xr^d) \ar[r]^{F^\delta(f)} & F^\delta(\Yr^d)}
  \labl{eq:C-is-natural}
commutes.
This is equivalent to the statement that
\begin{equation}
\begin{aligned}
&F(f) \cir F(\varpi_{\Xr}) \cir
  ( c_{\alpha_1} \otimes_\Cb \cdots \otimes_\Cb c_{\alpha_m} ) \cir U(\varpi_\Xr)^{-1} \\
  = \ &F(\varpi_{\Yr}) \cir
  ( c_{\alpha_1} \otimes_\Cb \cdots \otimes_\Cb c_{\alpha_m} ) \cir U(\varpi_\Yr)^{-1}
  \cir U(f) ~.
\end{aligned}
\end{equation}

By functoriality of $F$ we have $F(f) \cir F(\varpi_{\Xr}) = F(f \cir \varpi_{\Xr})$, which by definition of a morphism in $\WSh^d$ is equal to $F(\varpi_{\Yr})$. Similarly, $U(f) \cir U(\varpi_{\Xr}) = U(\varpi_{\Yr})$.
By assumption $U(\varpi_{\Xr})$ and $U(\varpi_{\Yr})$ are invertible, so that
this implies $U(\varpi_\Yr)^{-1} \cir U(f) = U(\varpi_\Xr)^{-1}$, as required.

To check that $C$ is monoidal one has to write out explicitly the associators and the isomorphisms $\phi^F_2$, $\phi^U_2$ in \erf{eq:nat-C-def}. One also needs to set $C_\emptyset = \phi^F_0 \circ (\phi^U_0)^{-1} : U(\emptyset) \rightarrow F(\emptyset)$. That $C$ is monoidal then follows from the
coherence properties of the $\phi^F$ and $\phi^U$. We do not give the details.

Conversely, let $\mu \in \text{Nat}_\otimes(U^\delta,F^\delta)$. For $\alpha \in \Sc$, let $\Xr_\alpha^d = (\Xr_\alpha,(\alpha),id)$ and set
\be
  c_\alpha = \mu_{\Xr_\alpha^d} : U(\Xr_\alpha) \longrightarrow F(\Xr_\alpha) ~.
\ee
Then, with $C \equiv C( \{ c_\alpha | \alpha \In \Sc \} )$ and
world sheet $\Xr^d = (\Xr , (\alpha_1,\dots,\alpha_m), \varpi_\Xr)$,
\begin{equation}
 \begin{aligned}
  C_{\Xr^d}
  =& \ F(\varpi) \circ ( c_{\alpha_1} \otimes_\Cb \cdots \otimes_\Cb c_{\alpha_m} )
  \circ U(\varpi)^{-1}\\
  \overset{(1)}{=}&\ F(\varpi) \circ ( \mu_{\Xr_{\alpha_1}^d} \otimes_\Cb \cdots
  \otimes_\Cb \mu_{\Xr_{\alpha_m}^d} ) \circ U(\varpi)^{-1}\\
  \overset{(2)}{=}&\ F(\varpi) \circ \mu_{\Xr_{\alpha_1}^d \oti \cdots
  \oti \Xr_{\alpha_m}^d} \circ U(\varpi)^{-1}\\
  \overset{(3)}{=}&\ \mu_{\Xr^d} \circ U(\varpi) \circ U(\varpi)^{-1} = \mu_{\Xr^d}
\end{aligned}
\end{equation}
where (1) is the definition of $c_\alpha$, (2) is monoidality of $\mu$ and (3) follows from the commutativity
of the following diagram,
\begin{equation}
\xymatrix{
U^\delta(\Xr_{\alpha_1}^d\otimes\cdots\otimes \Xr_{\alpha_m}^d) \ar[d]_{U^\delta(\varpi)} \ar[r]
&F^\delta(\Xr_{\alpha_1}^d\otimes\cdots\otimes \Xr_{\alpha_m}^d)\ar[d]_{F^\delta(\varpi)} \\
U^\delta(\Xr^d) \ar[r]^{\mu_{\Xr}} &F^\delta(\Xr^d)}
\end{equation}
 in which the upper horizontal arrow is $\mu_{\Xr_{\alpha_1}^d\otimes\cdots\otimes \Xr_{\alpha_m}^d}$.
\epf

\medskip
The above lemma allows us to express all natural transformations
$U^\delta \rightarrow F^\delta$ by their effect on the generating
set of world sheets. The next question is
which elements in $\text{Nat}_\otimes(U^\delta,F^\delta)$
are of the form $\mu^\delta$ for some
$\mu \in \text{Nat}_\otimes(U,F)$, where $\mu^\delta$ the pull-back
of the monoidal natural transformation by the forgetful functor as
described above.
To have $C(\{ c_\alpha | \alpha \In \Sc \} ) = \mu^\delta$, the
linear maps $\{ c_\alpha | \alpha \In \Sc\}$ have to
obey a number of relations. The relations will be expressed using the decomposed world sheets
 $\Yr^d_{\Rr n,l}$ and $\Yr^d_{\Rr n,r}$, $n=1,\dots,32$, listed in
the following pictures.

\begin{align*}
&\Yr_{\Rr1,l}^d=\hspace{2mm}\figbox{0.36}{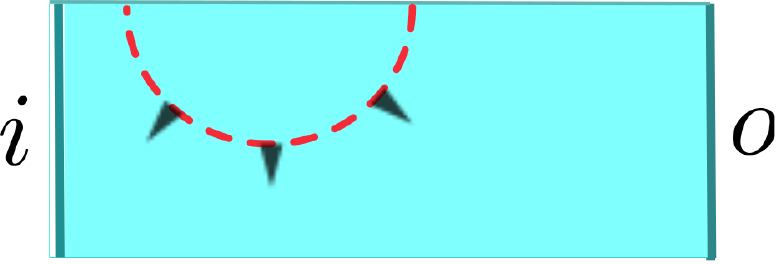}\hspace{30mm}\Yr_{\Rr1,r}^d=\hspace{2mm}\figbox{0.36}{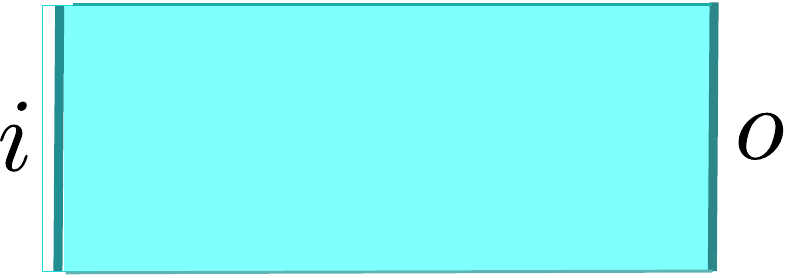}\\
&\Yr_{\Rr2,l}^d=\hspace{2mm}\figbox{0.36}{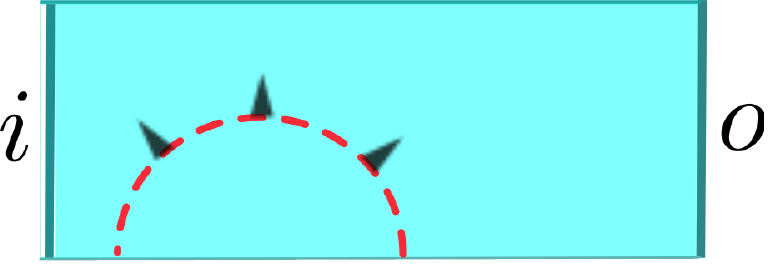}\hspace{30mm}\Yr_{\Rr2,r}^d=\hspace{2mm}\figbox{0.36}{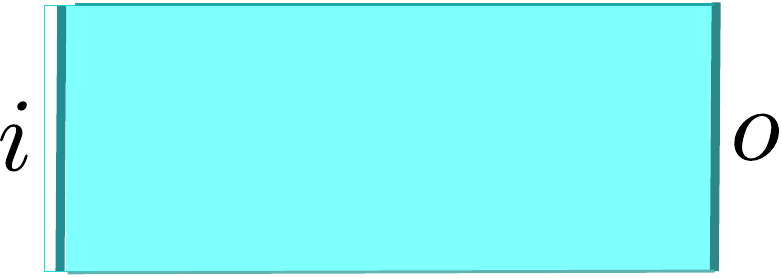}\\
&\Yr_{\Rr3,l}^d=\hspace{2mm}\figbox{0.36}{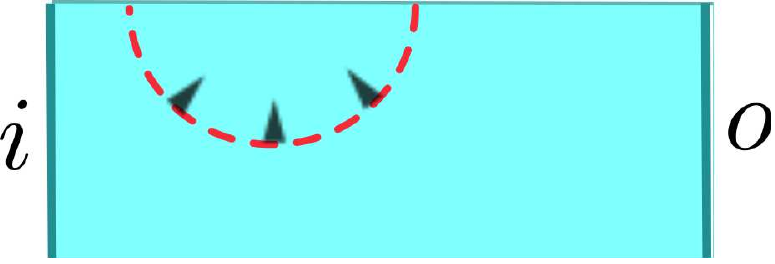}\hspace{30mm}\Yr_{\Rr3,r}^d=\hspace{2mm}\figbox{0.36}{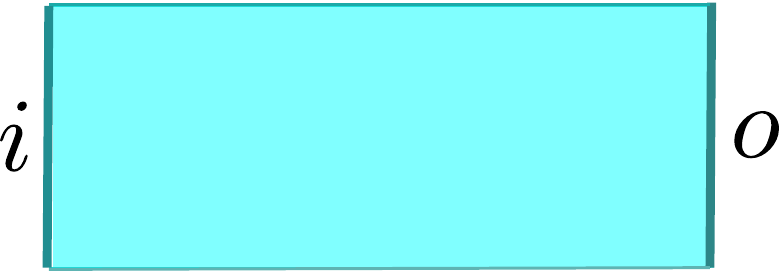}\\
&\Yr_{\Rr4,l}^d=\hspace{2mm}\figbox{0.36}{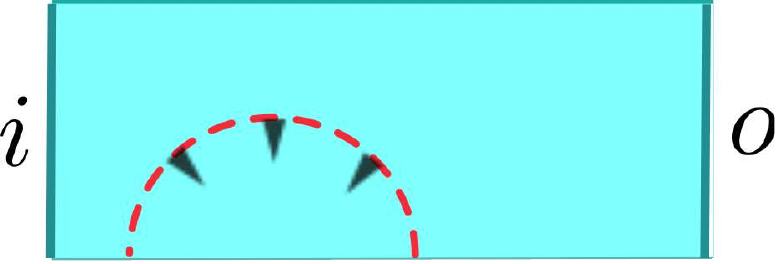}\hspace{30mm}\Yr_{\Rr4,r}^d=\hspace{2mm}\figbox{0.36}{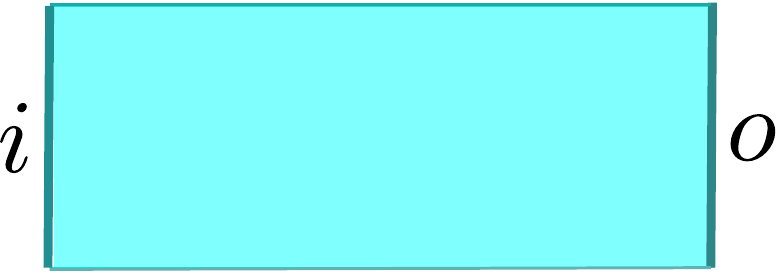}\\
&\Yr_{\Rr5,l}^d=\hspace{2mm}\figbox{0.36}{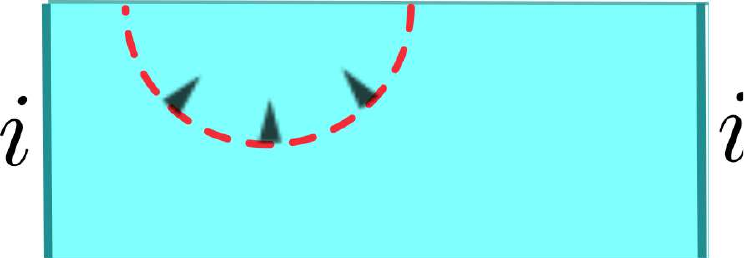}\hspace{31mm}\Yr_{\Rr5,r}^d=\hspace{2.5mm}\figbox{0.36}{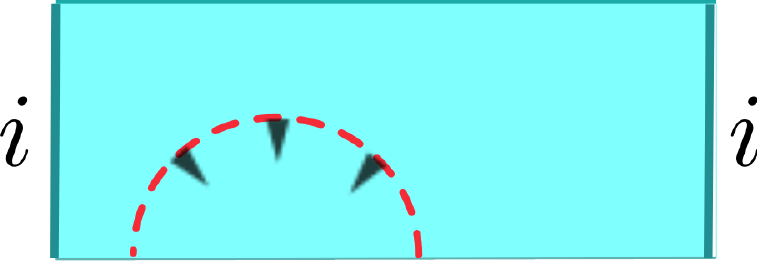}\\
&\Yr_{\Rr6,l}^d=\figbox{0.28}{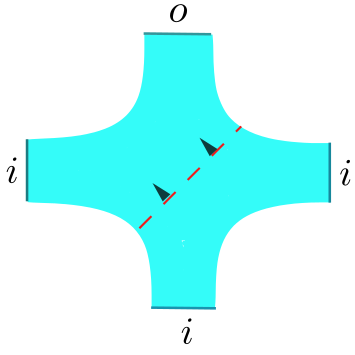}\hspace{26.5mm}\Yr_{\Rr6,r}^d=\figbox{0.28}{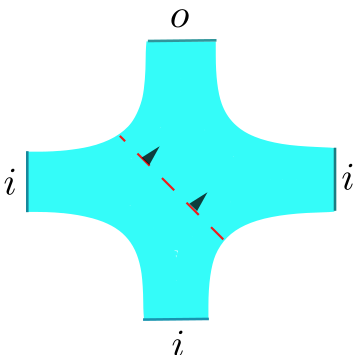}\\
&\Yr_{\Rr7,l}^d=\figbox{0.28}{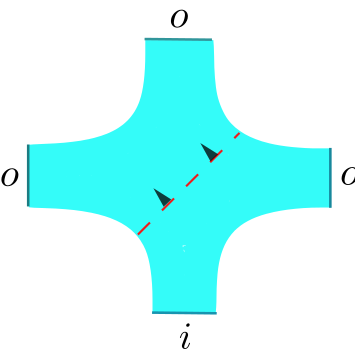}\hspace{25mm}\Yr_{\Rr7,r}^d=\figbox{0.28}{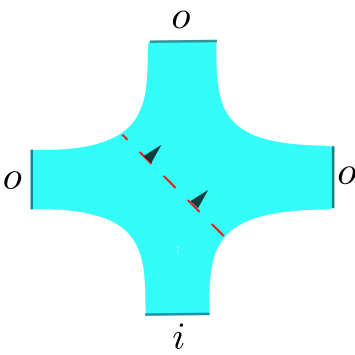}\\
&\Yr_{\Rr8,l}^d=\figbox{0.28}{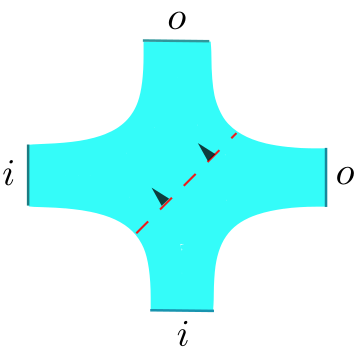}\hspace{26mm}\Yr_{\Rr8,r}^d=\figbox{0.28}{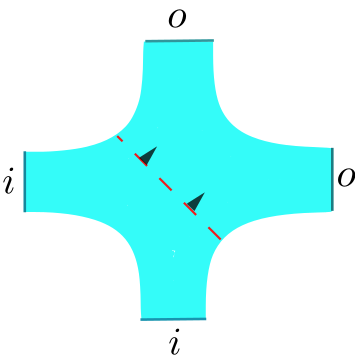}\\
&\Yr_{\Rr9,l}^d=\figbox{0.28}{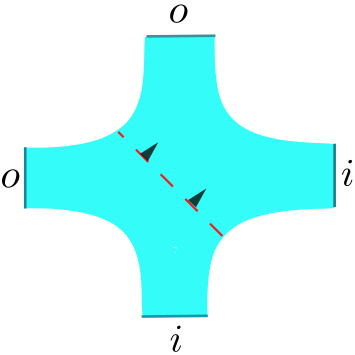}\hspace{26mm}\Yr_{\Rr9,r}^d=\figbox{0.28}{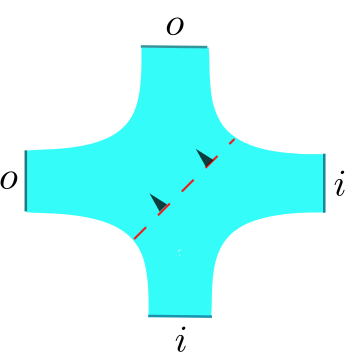}\\
&\Yr_{\Rr10,l}^d=\hspace{3mm}\figbox{0.26}{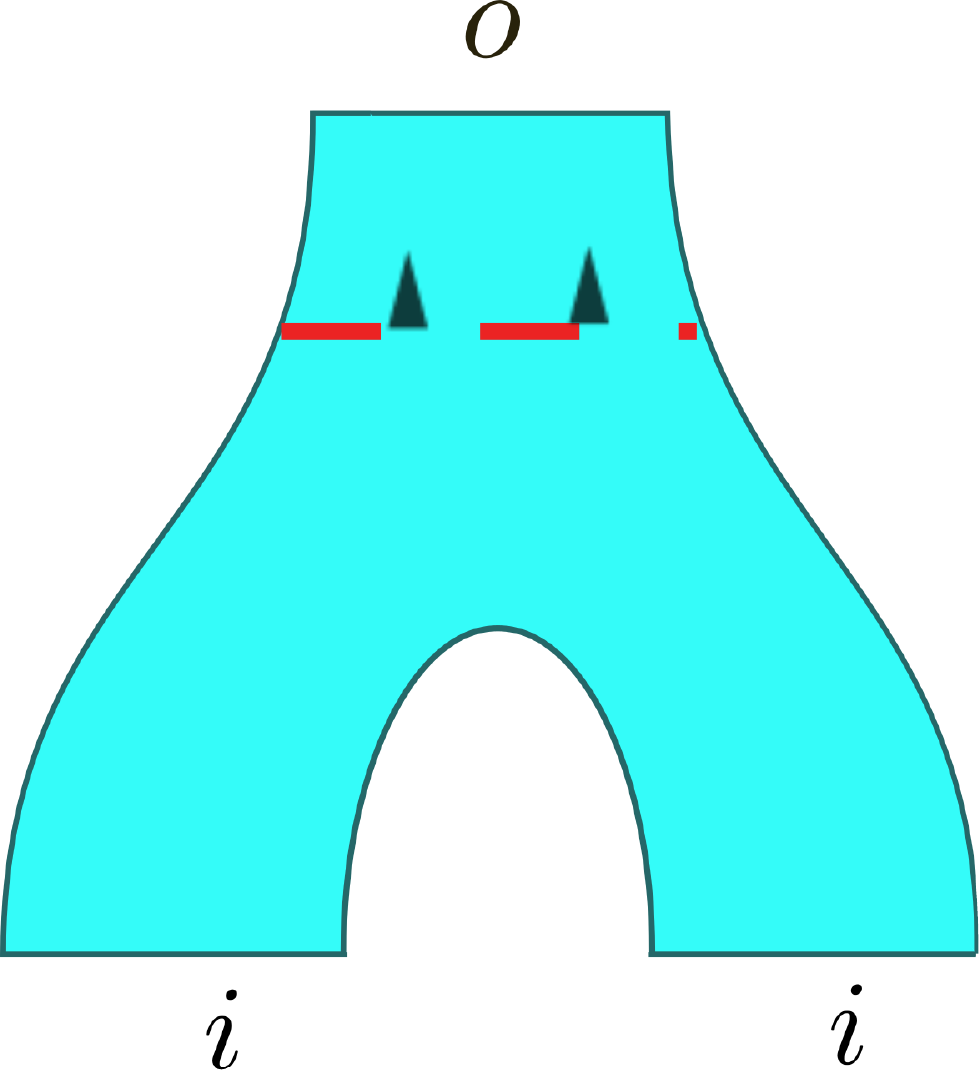}\hspace{31mm}\Yr_{\Rr10,r}^d=\hspace{3mm}\figbox{0.26}{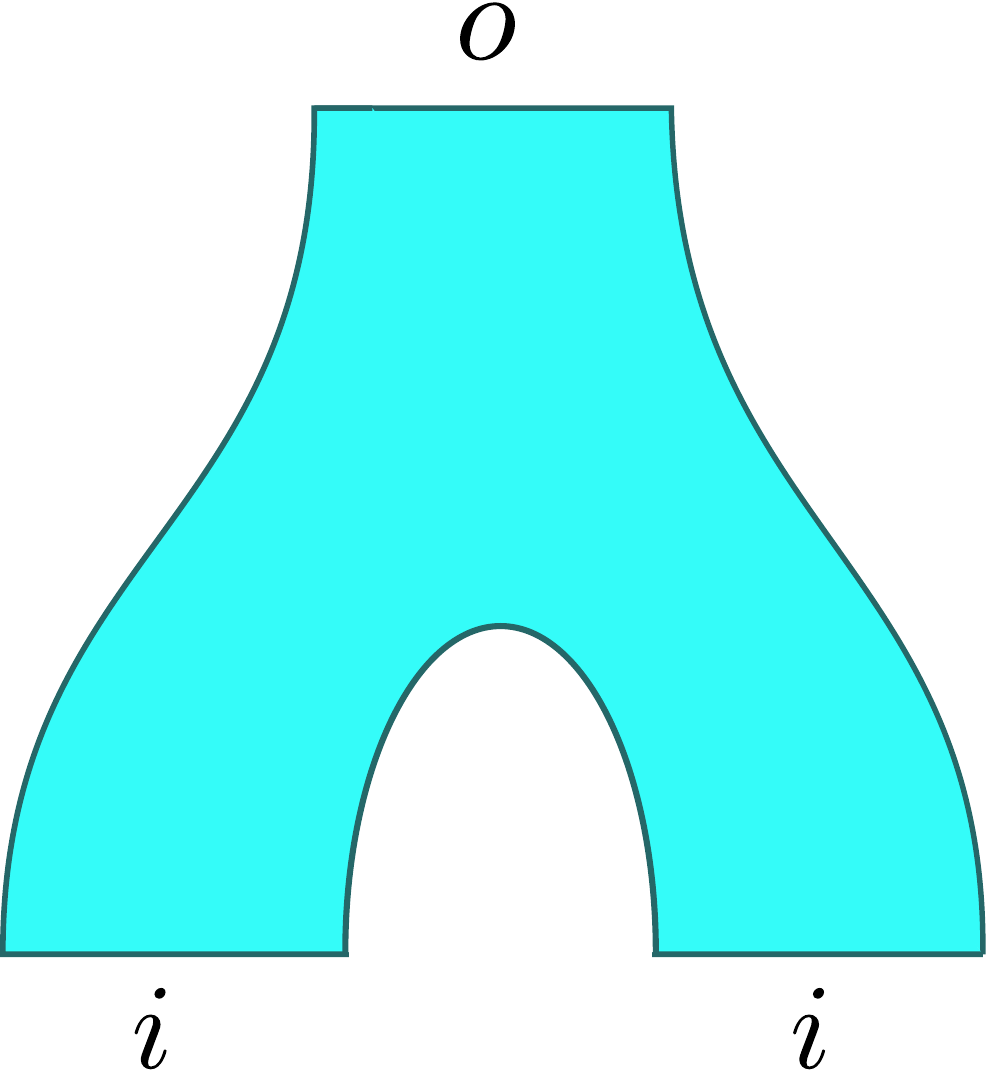}\\
&\Yr_{\Rr11,l}^d=\hspace{12mm}\figbox{0.07}{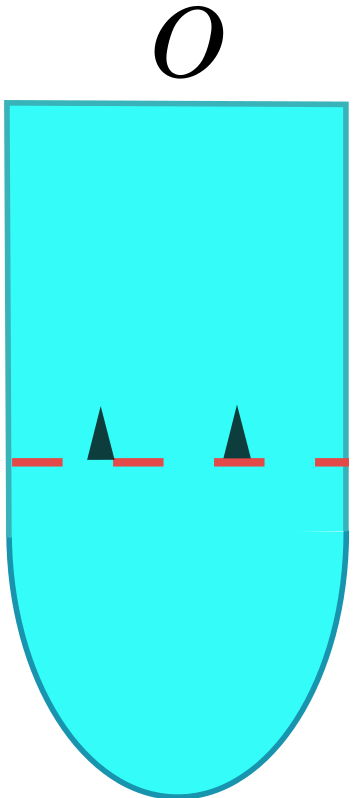}\hspace{38mm}\Yr_{\Rr11,r}^d=\hspace{13mm}\figbox{0.09}{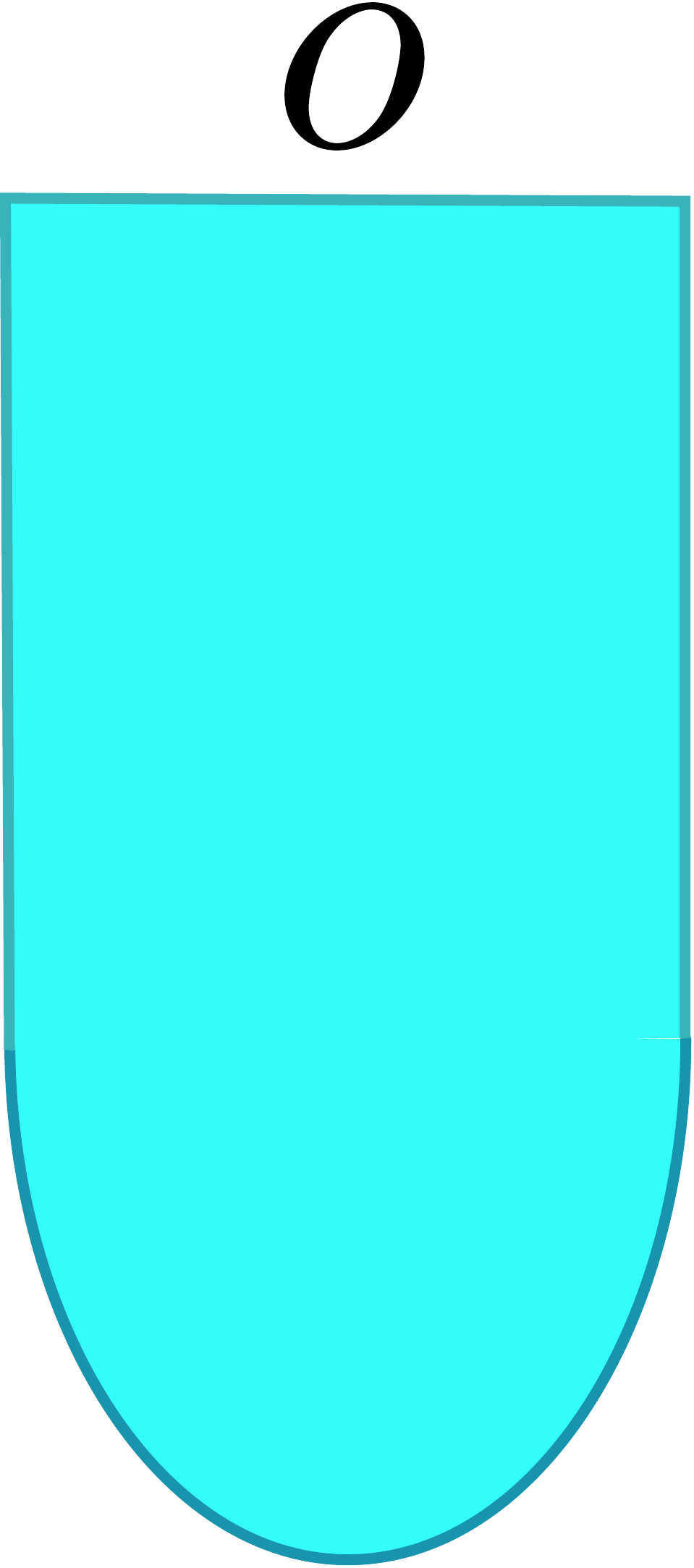}\\
&\Yr_{\Rr12,l}^d=\hspace{3.5mm}\figbox{0.25}{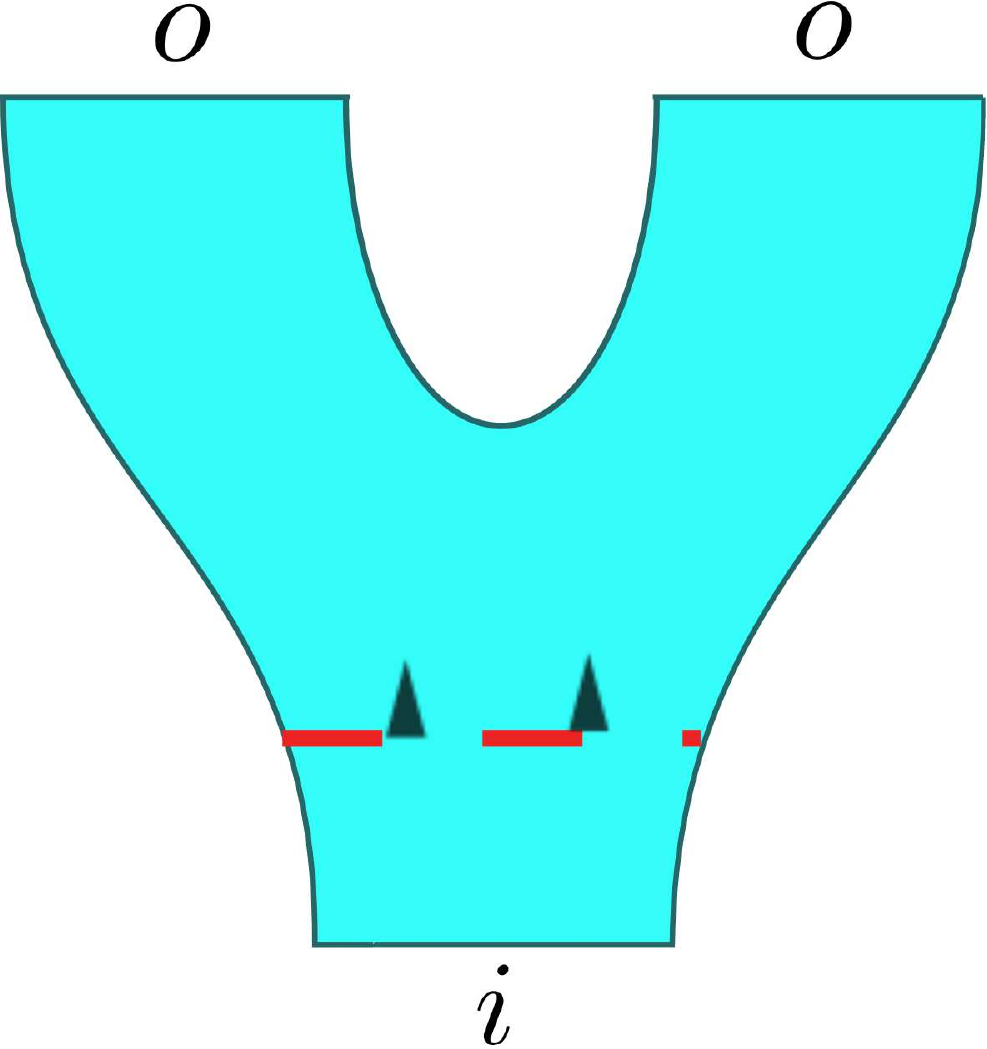}\hspace{31mm}\Yr_{\Rr12,r}^d=\hspace{3.5mm}\figbox{0.26}{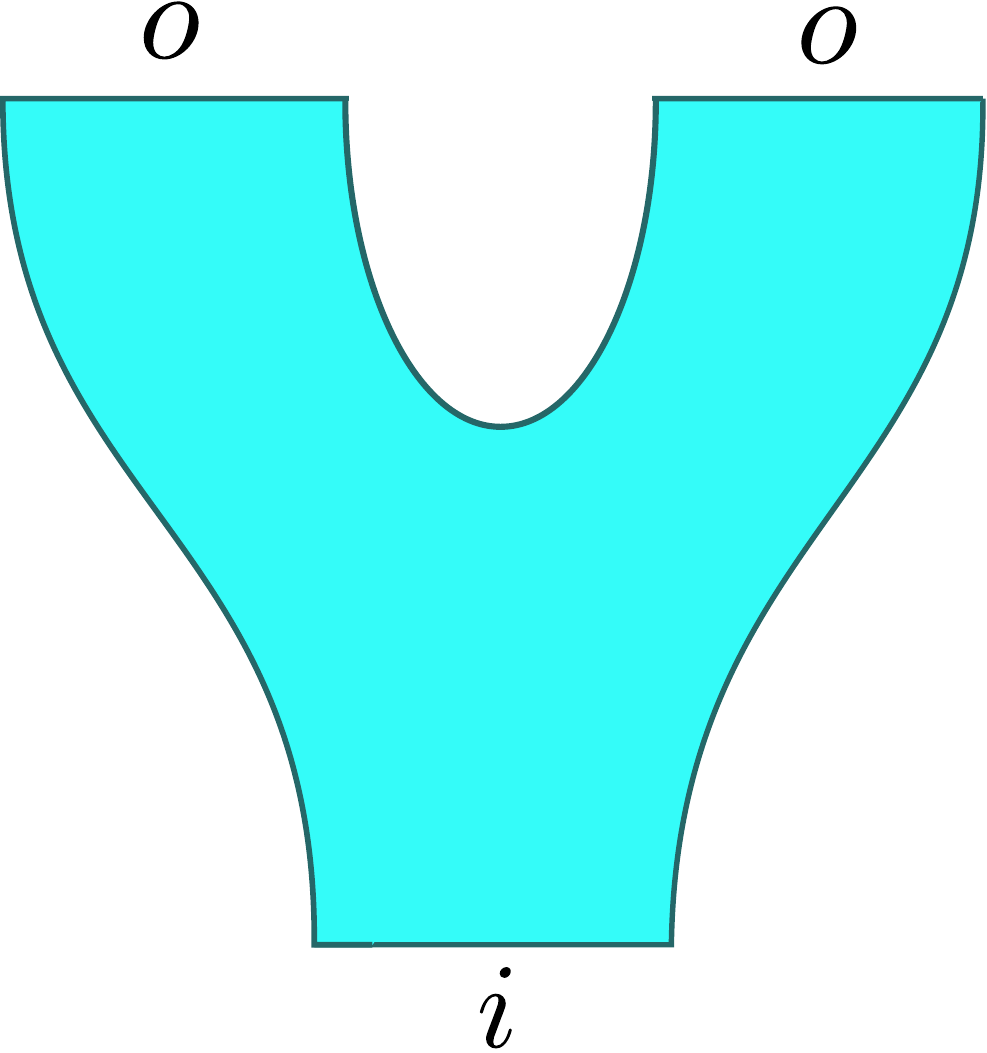}\\
&\Yr_{\Rr13,l}^d=\hspace{12mm}\figbox{0.07}{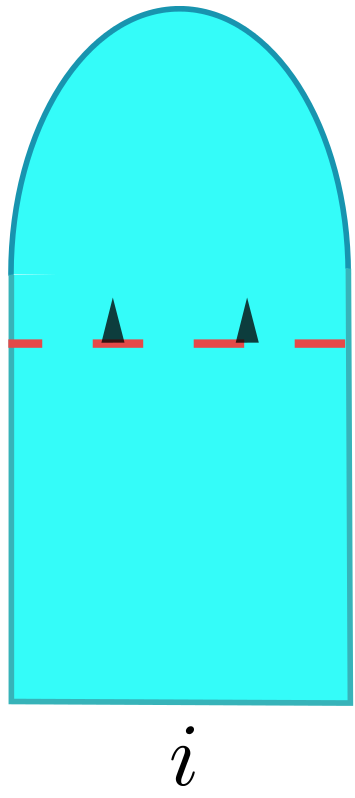}\hspace{39mm}\Yr_{\Rr13,r}^d=\hspace{12mm}\figbox{0.09}{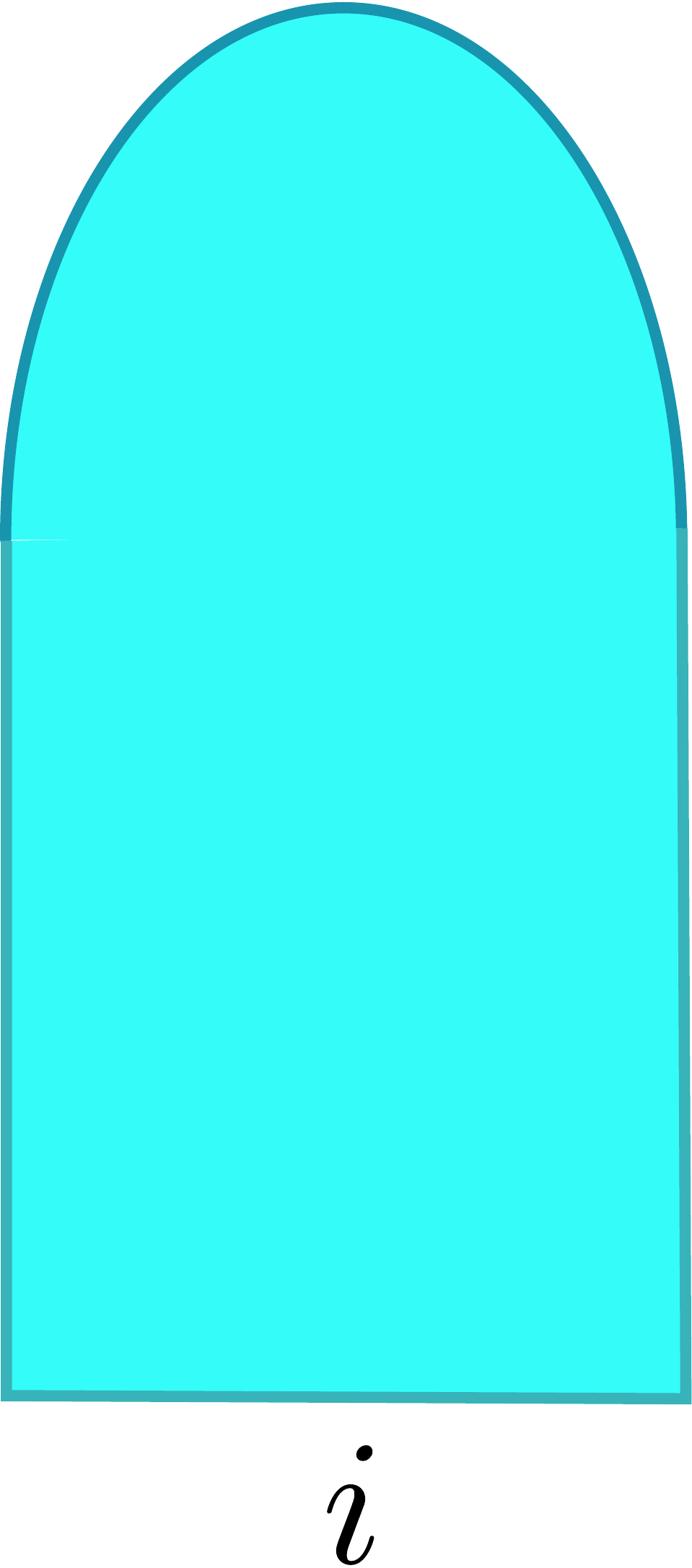}\\
&\Yr_{\Rr14,l}^d=\ \figbox{0.3}{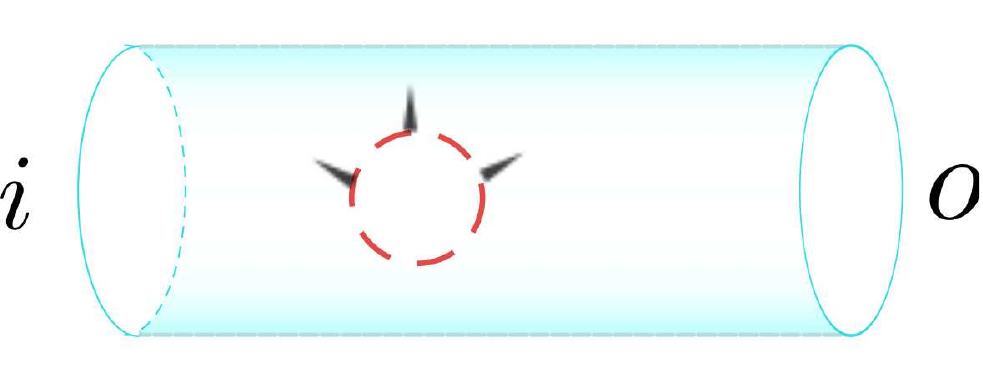}\hspace{28.5mm}\Yr_{\Rr14,r}^d=\ \figbox{0.31}{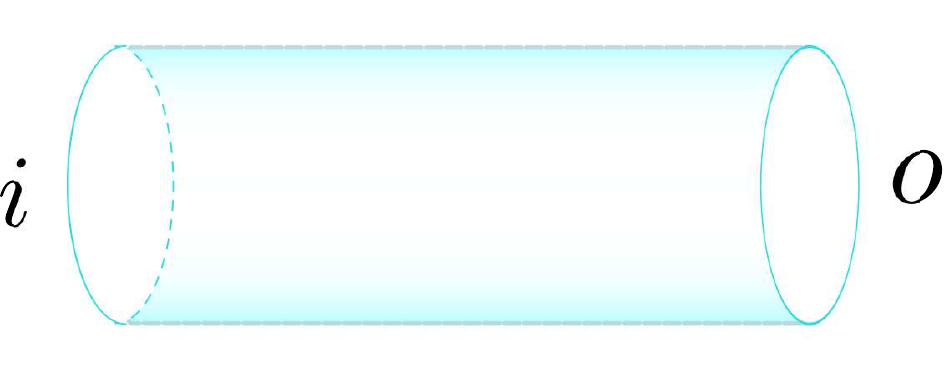}\\
&\Yr_{\Rr15,l}^d=\ \figbox{0.31}{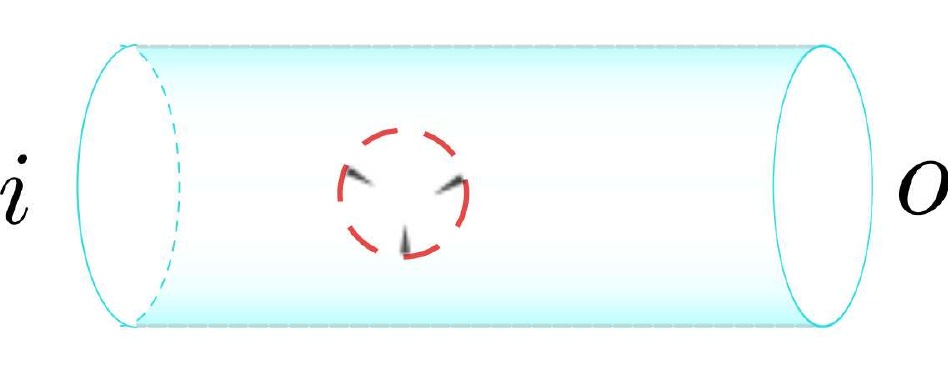}\hspace{28.5mm}\Yr_{\Rr15,r}^d=\ \figbox{0.3}{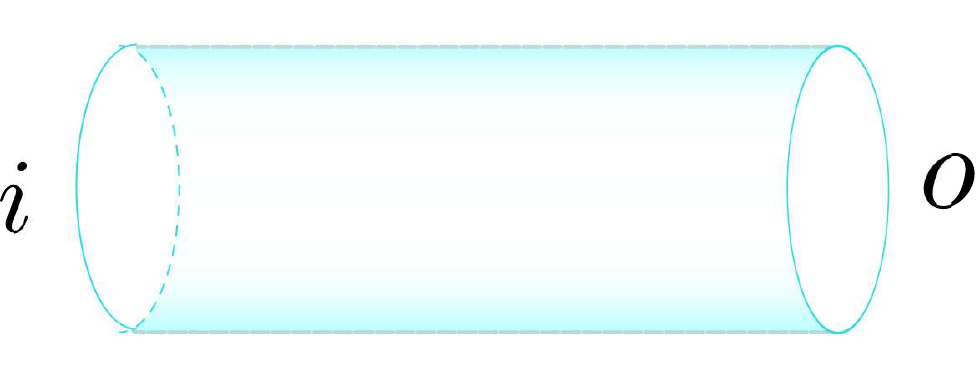}\\
&\Yr_{\Rr16,l}^d=\ \figbox{0.3}{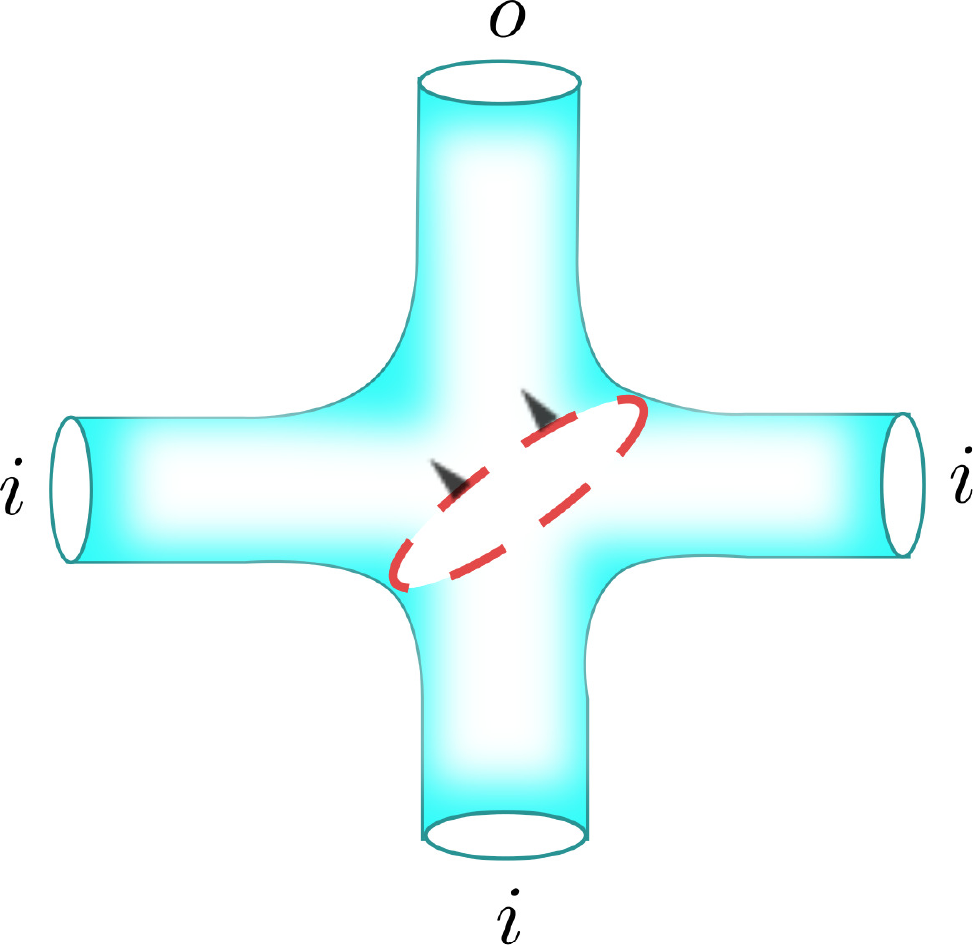}\hspace{27mm}\Yr_{\Rr16,r}^d=\ \figbox{0.3}{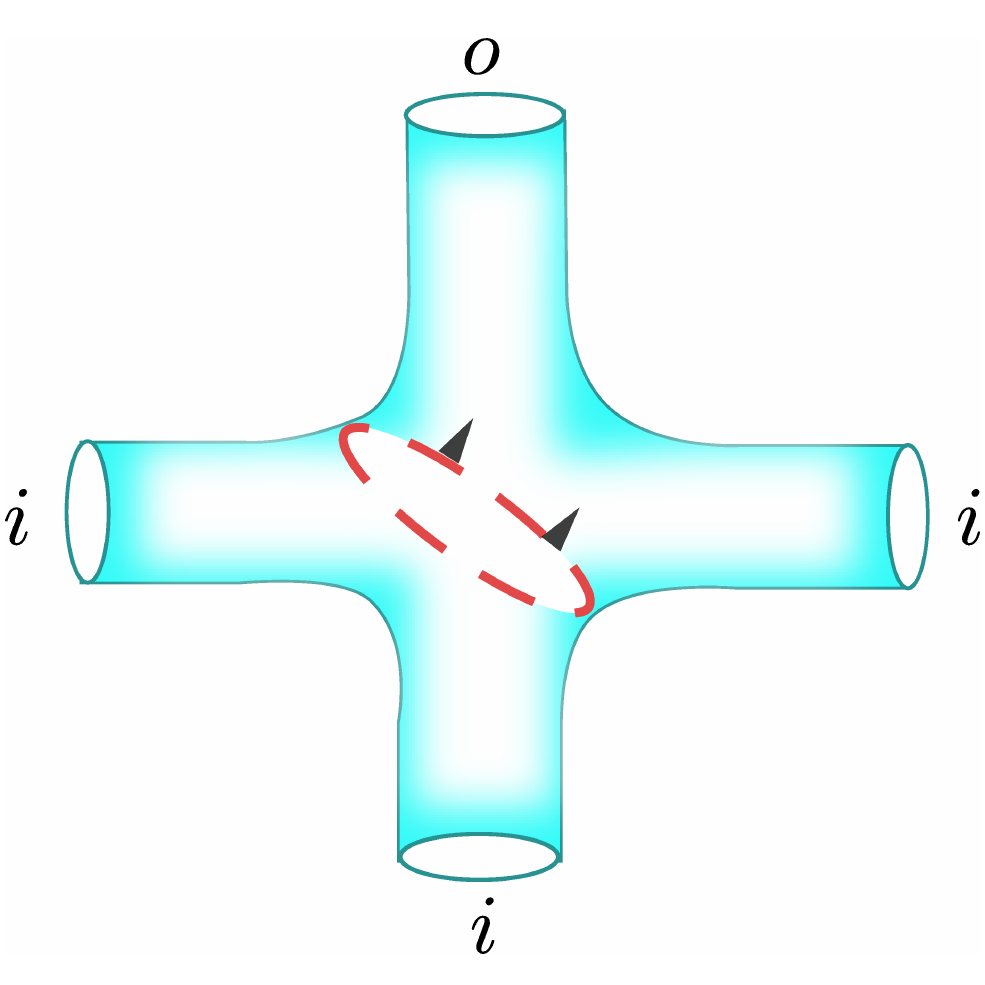}\\
&\Yr_{\Rr17,l}^d=\ \figbox{0.3}{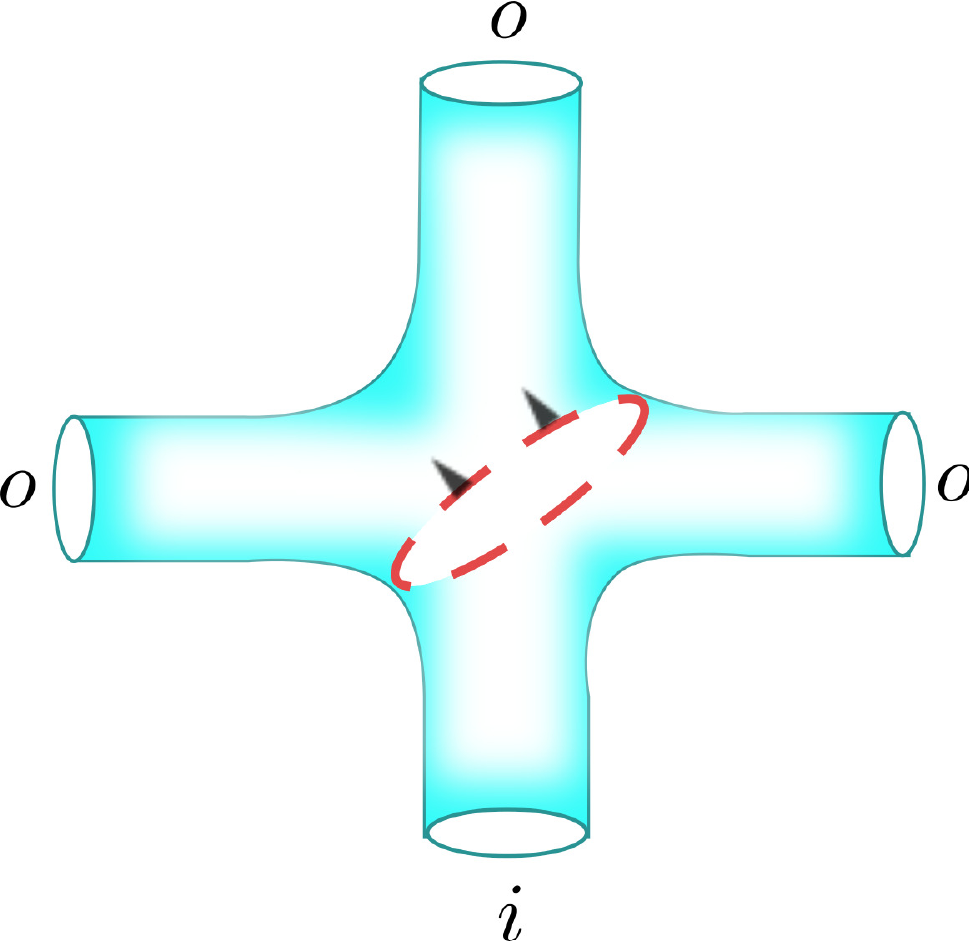}\hspace{27mm}\Yr_{\Rr17,r}^d=\ \figbox{0.3}{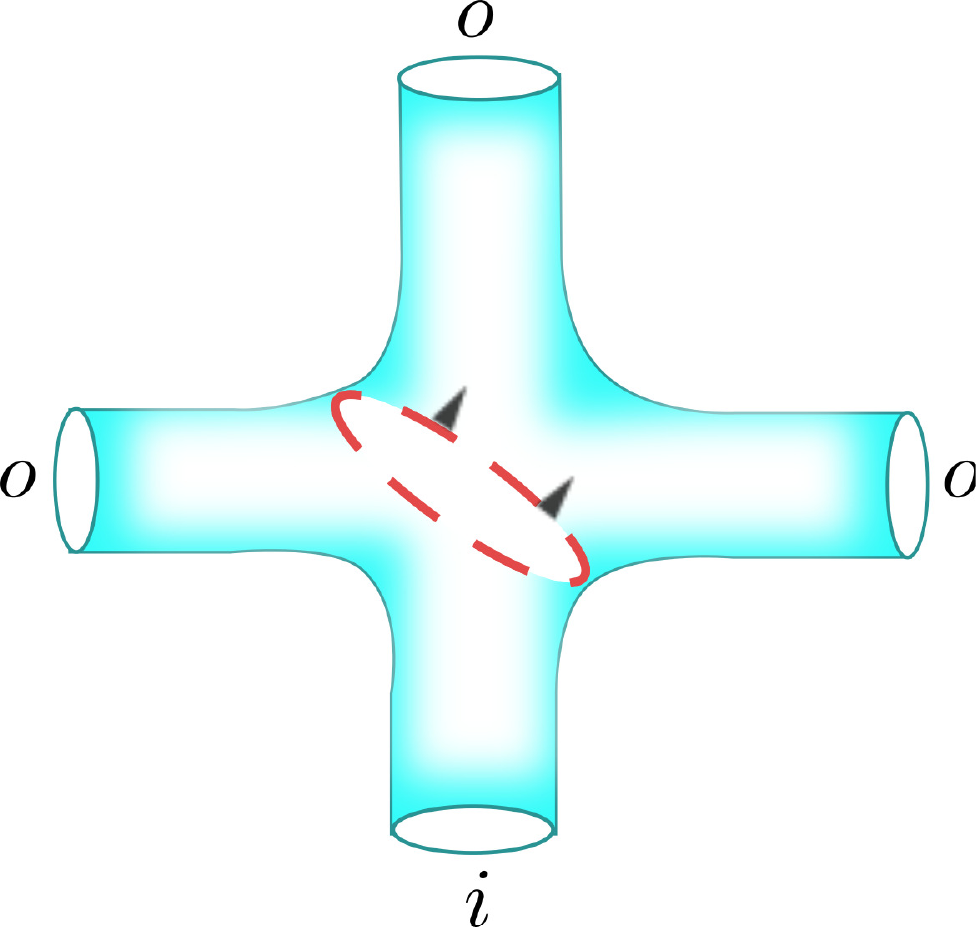}\\
&\Yr_{\Rr18,l}^d=\ \figbox{0.3}{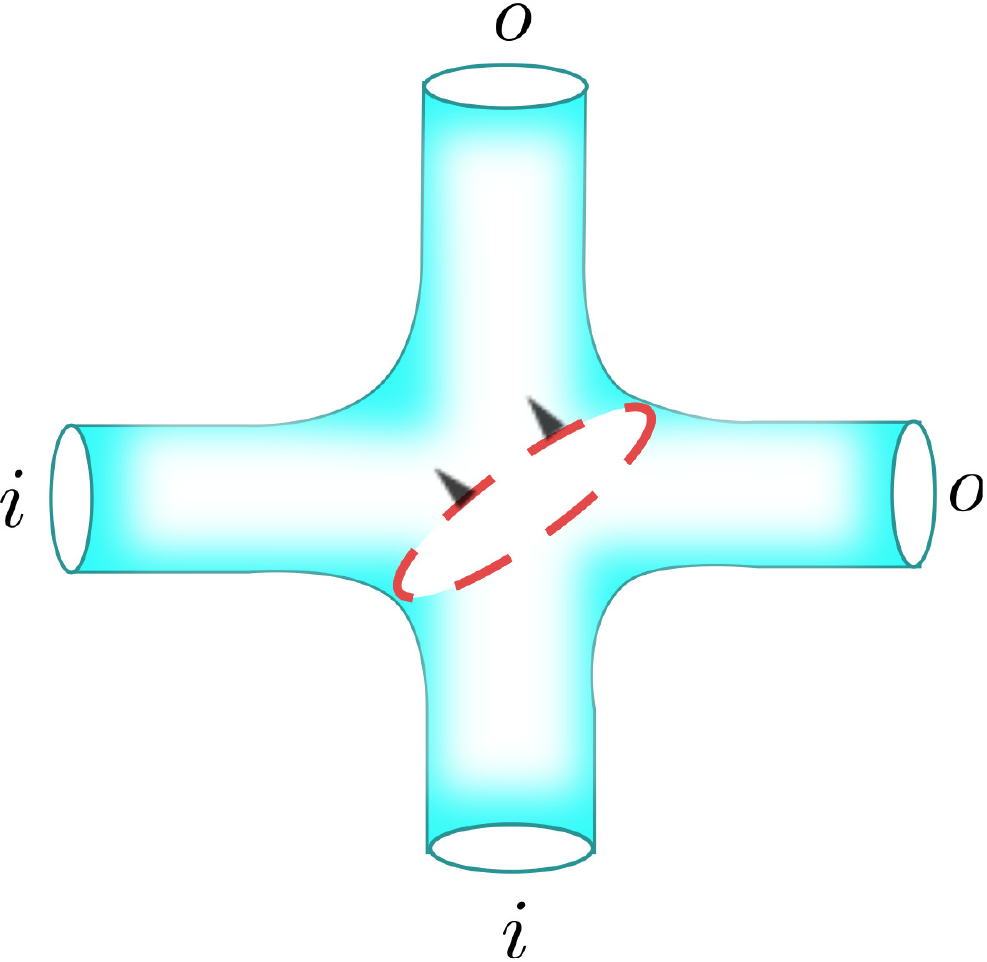}\hspace{27mm}\Yr_{\Rr18,r}^d=\ \figbox{0.3}{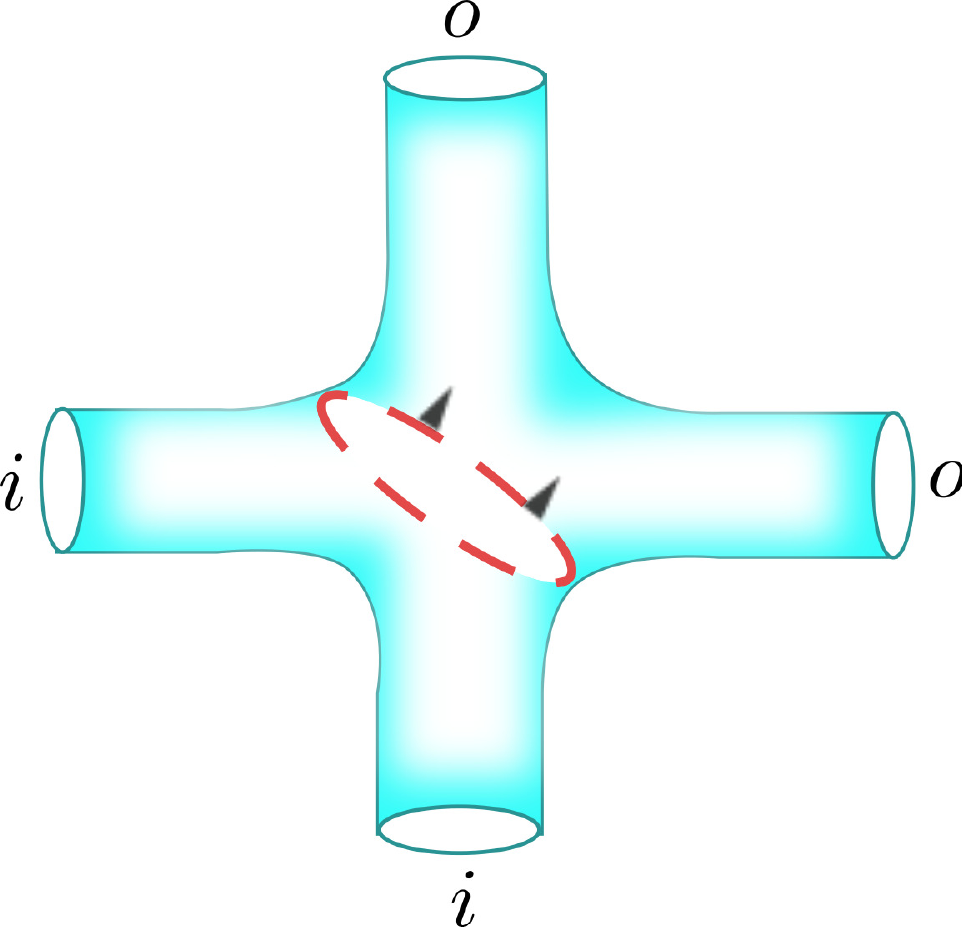}\\
&\Yr_{\Rr19,l}^d=\ \figbox{0.3}{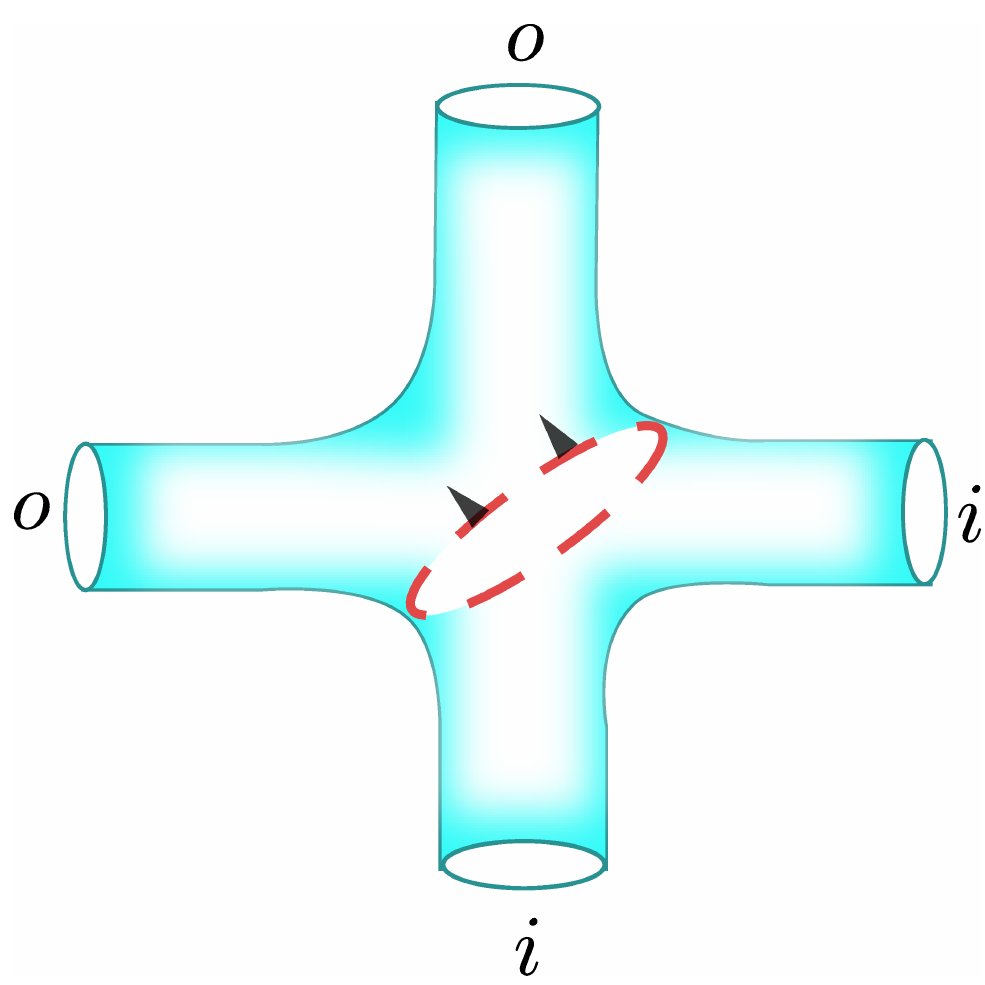}\hspace{27mm}\Yr_{\Rr19,r}^d=\ \figbox{0.3}{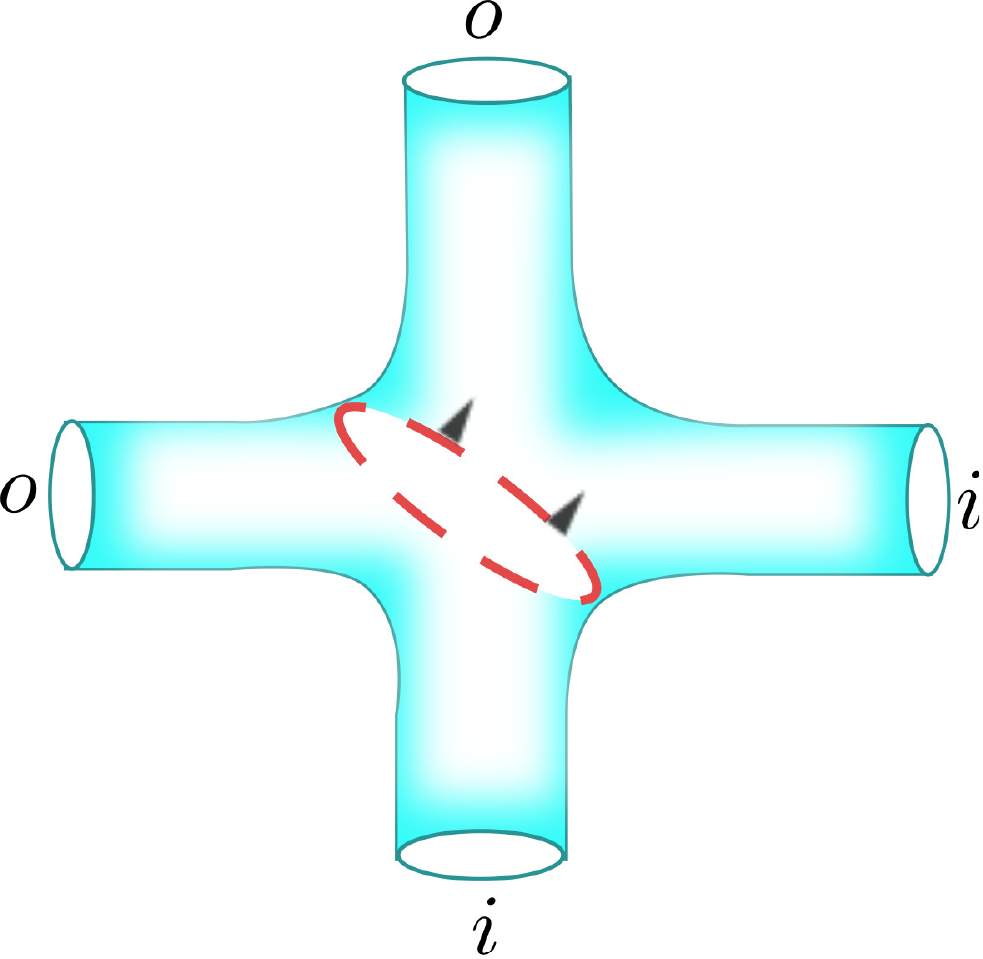}\\
&\Yr_{\Rr20,l}^d=\hspace{4.3mm}\figbox{0.2}{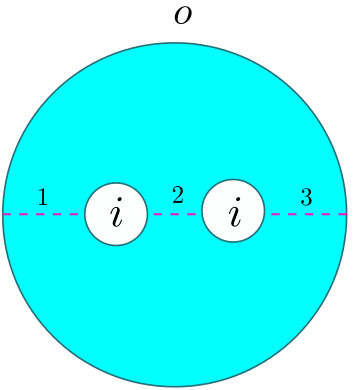}\hspace{29mm}\Yr_{\Rr20,r}^d=\hspace{4.3mm}\figbox{0.2}{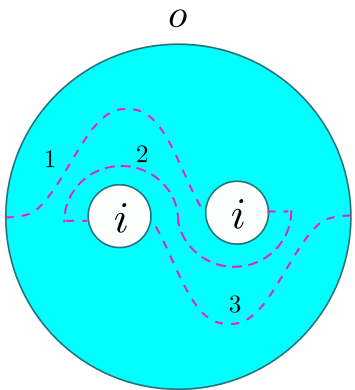}\\
&\Yr_{\Rr21,l}^d=\hspace{5.5mm}\figbox{0.23}{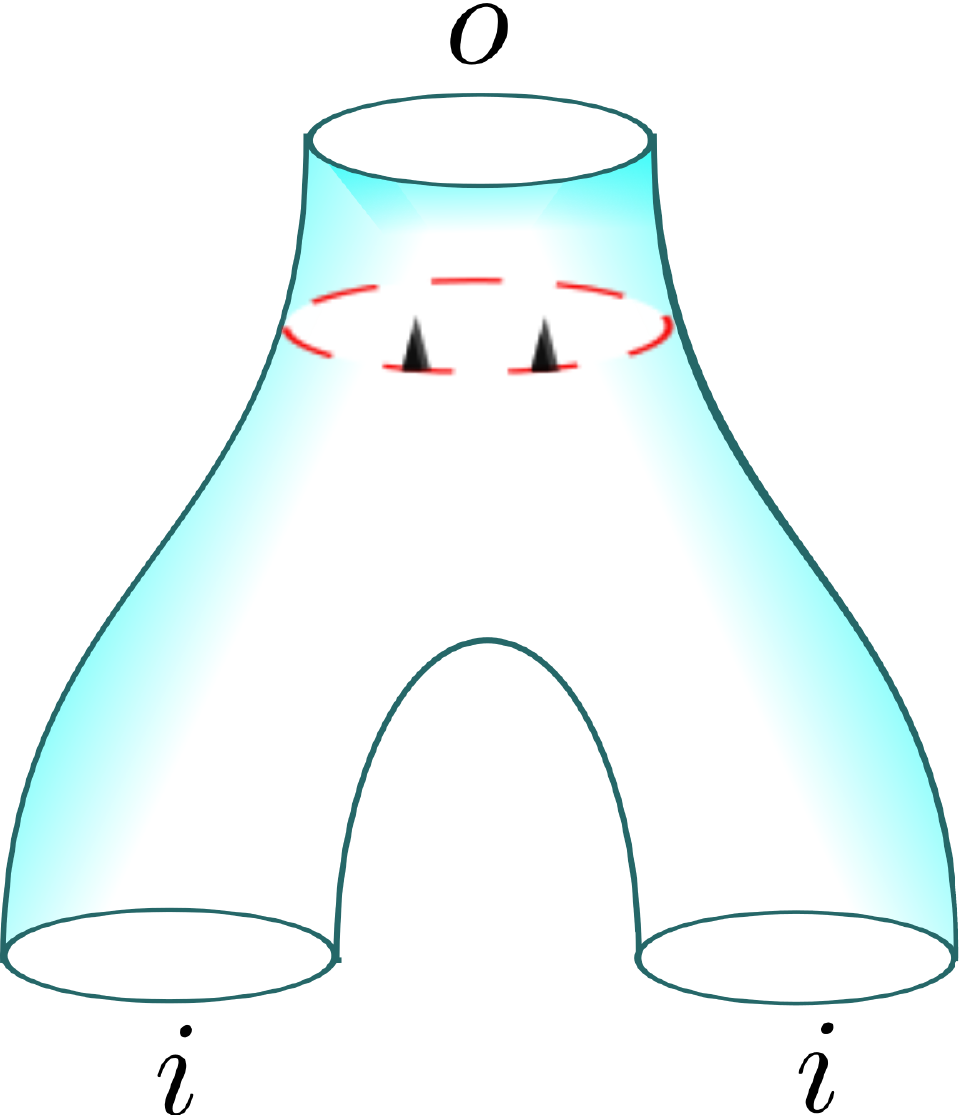}\hspace{30mm}\Yr_{\Rr21,r}^d=\hspace{5.5mm}\figbox{0.23}{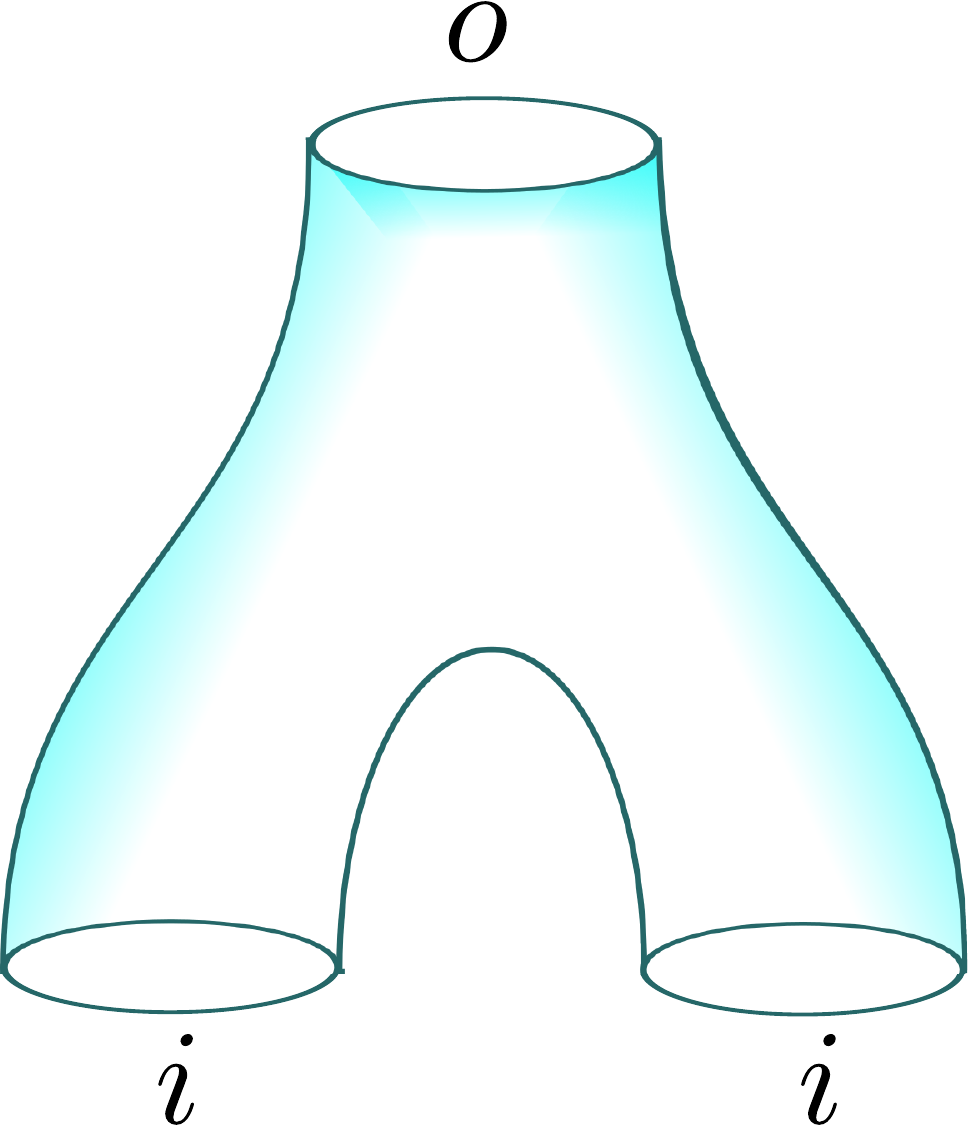}\\
&\Yr_{\Rr22,l}^d=\hspace{4.3mm}\figbox{0.20}{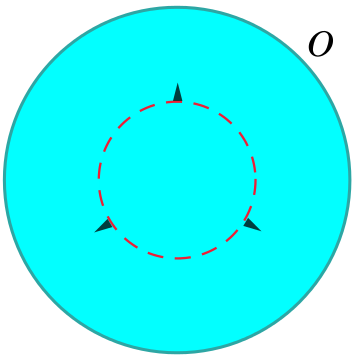}\hspace{29mm}\Yr_{\Rr22,r}^d=\hspace{4.3mm}\figbox{0.20}{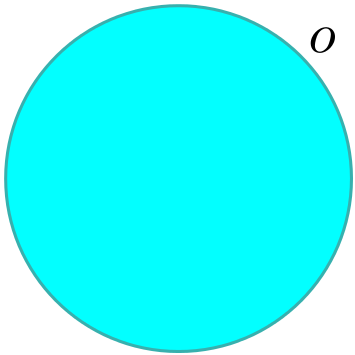}\\
&\Yr_{\Rr23,l}^d=\hspace{5.5mm}\figbox{0.23}{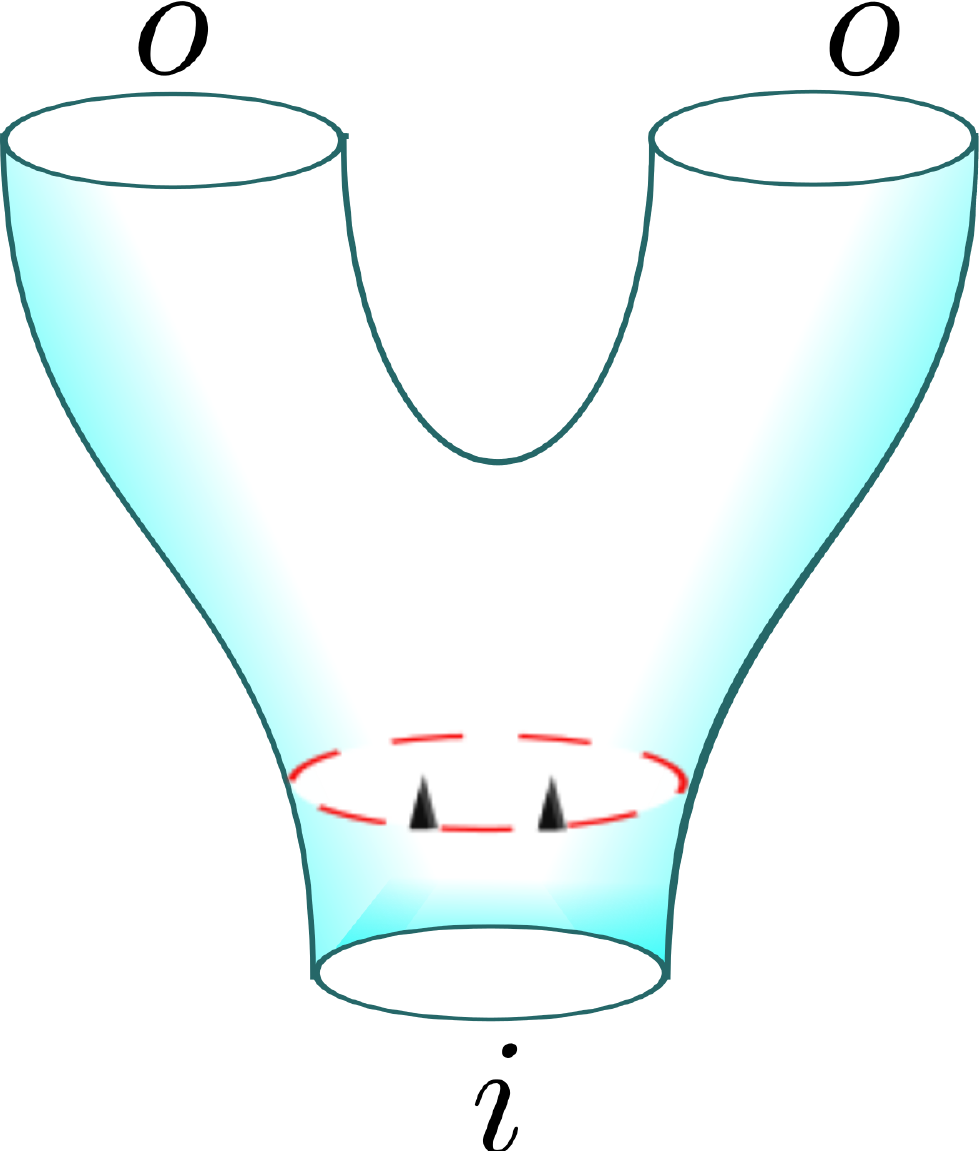}\hspace{30mm}\Yr_{\Rr23,r}^d=\hspace{5.5mm}\figbox{0.23}{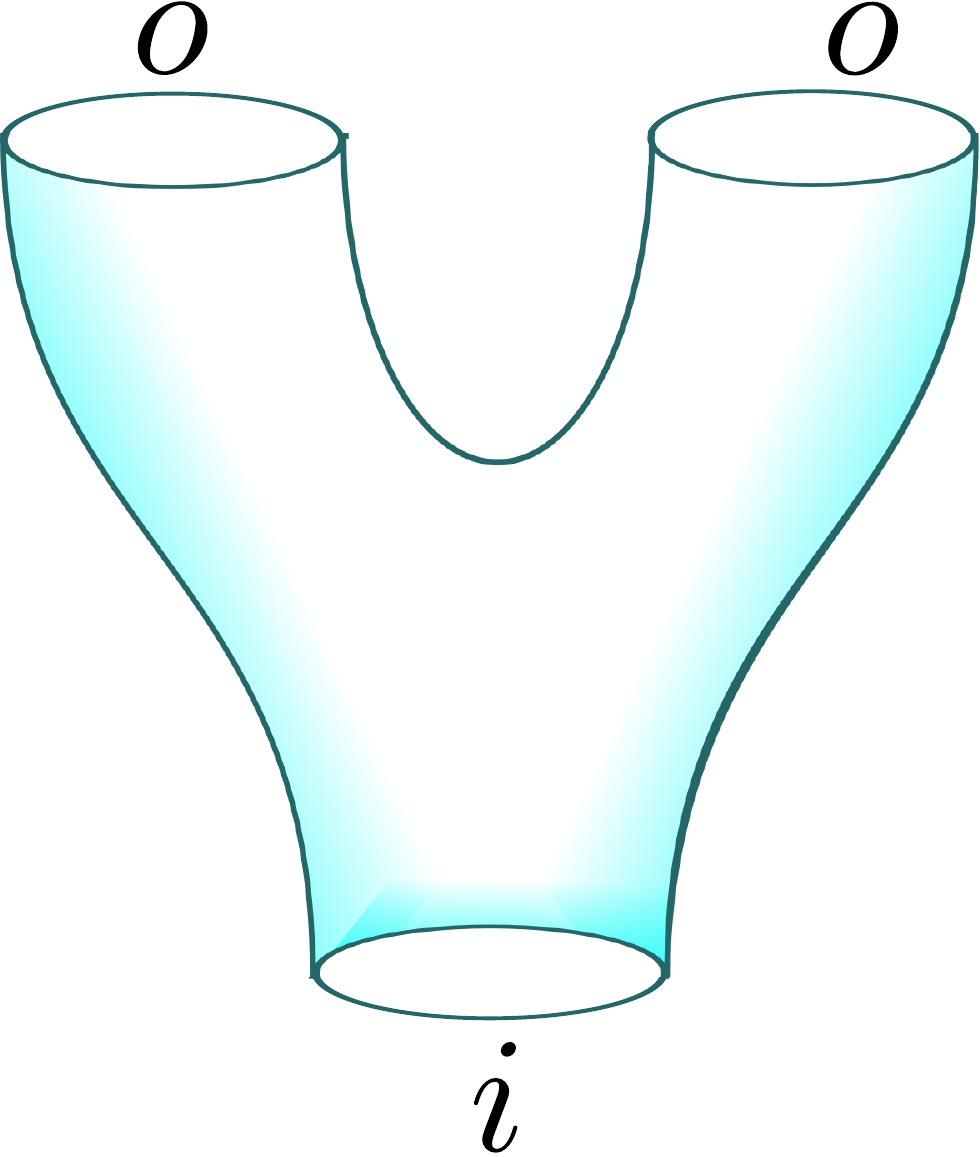}\\
&\Yr_{\Rr24,l}^d=\hspace{4.3mm}\figbox{0.20}{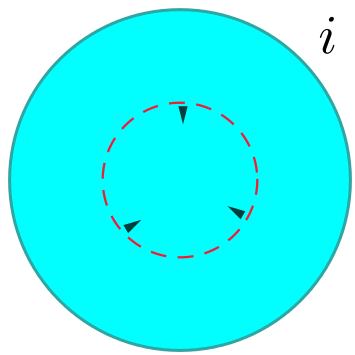}\hspace{30mm}\Yr_{\Rr24,r}^d=\hspace{4.3mm}\figbox{0.20}{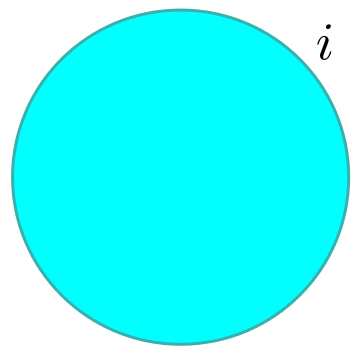}\\
&\Yr_{\Rr25,l}^d=\hspace{9mm}\figbox{0.32}{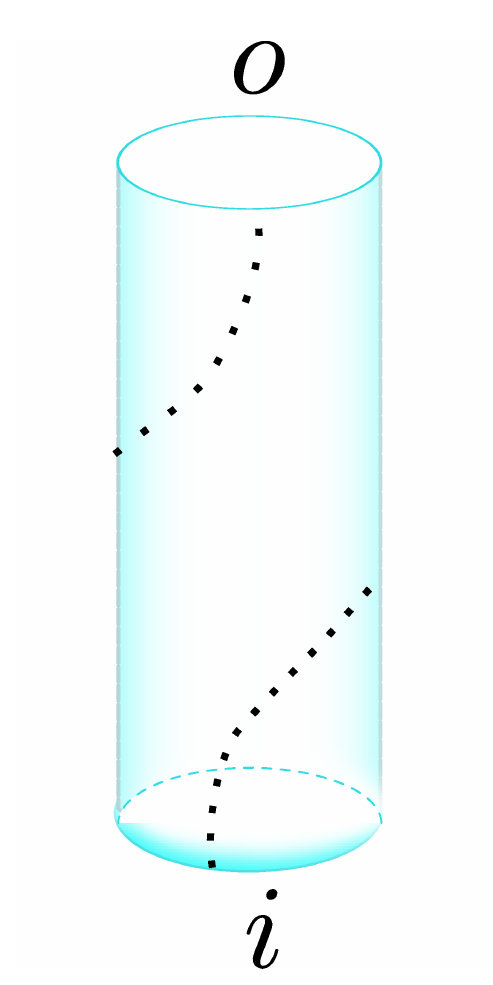}\hspace{34mm}\Yr_{\Rr25,r}^d=\hspace{9mm}\figbox{0.32}{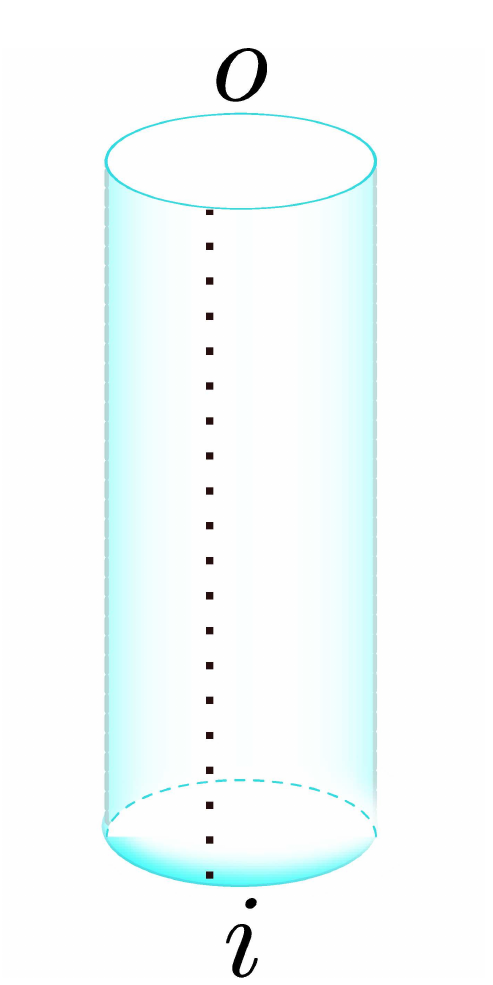}\\
&\Yr_{\Rr26,l}^d=\hspace{4mm} \figbox{0.26}{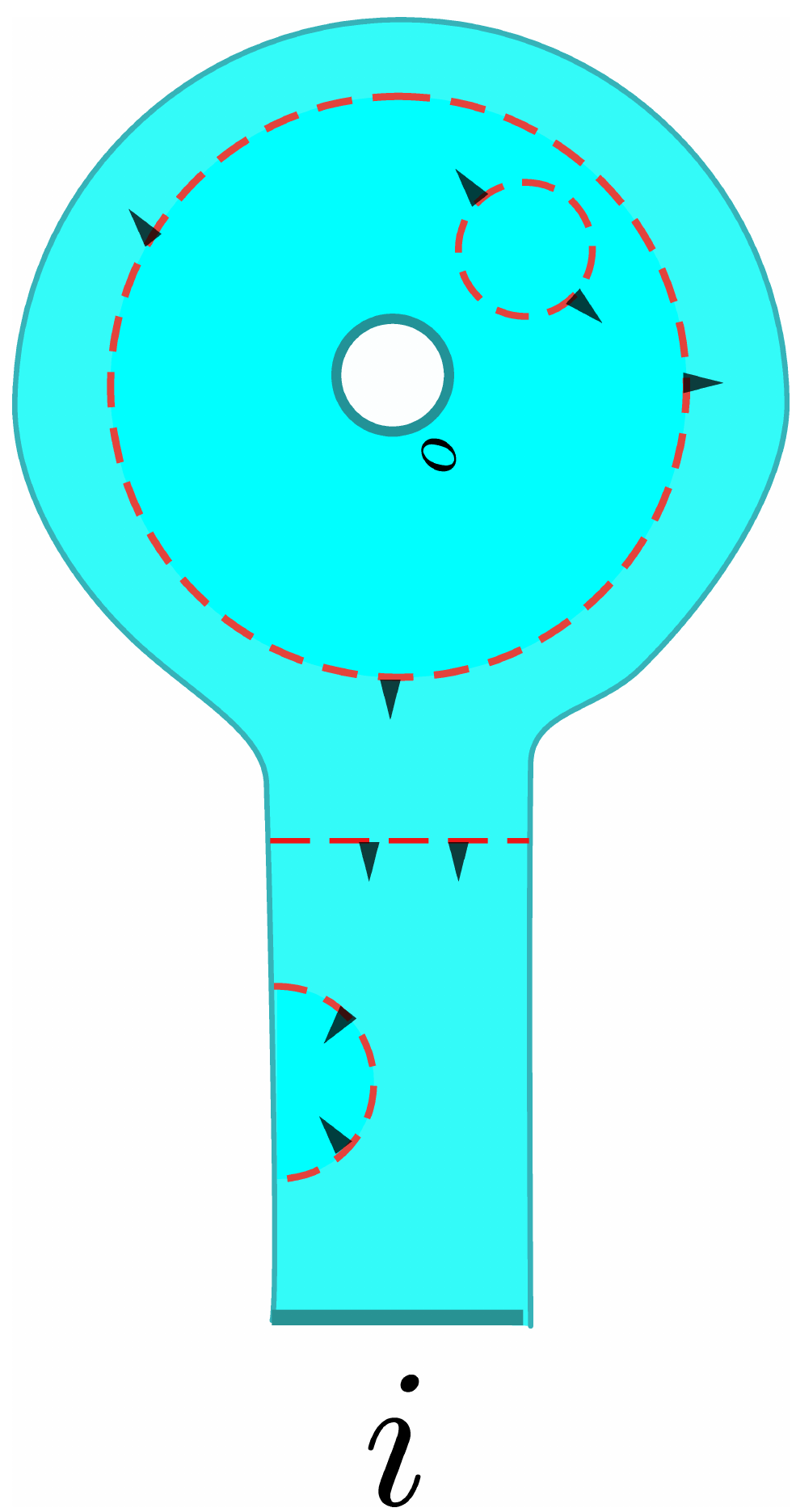}\hspace{29mm}\Yr_{\Rr26,r}^d=\hspace{4mm}\figbox{0.26}{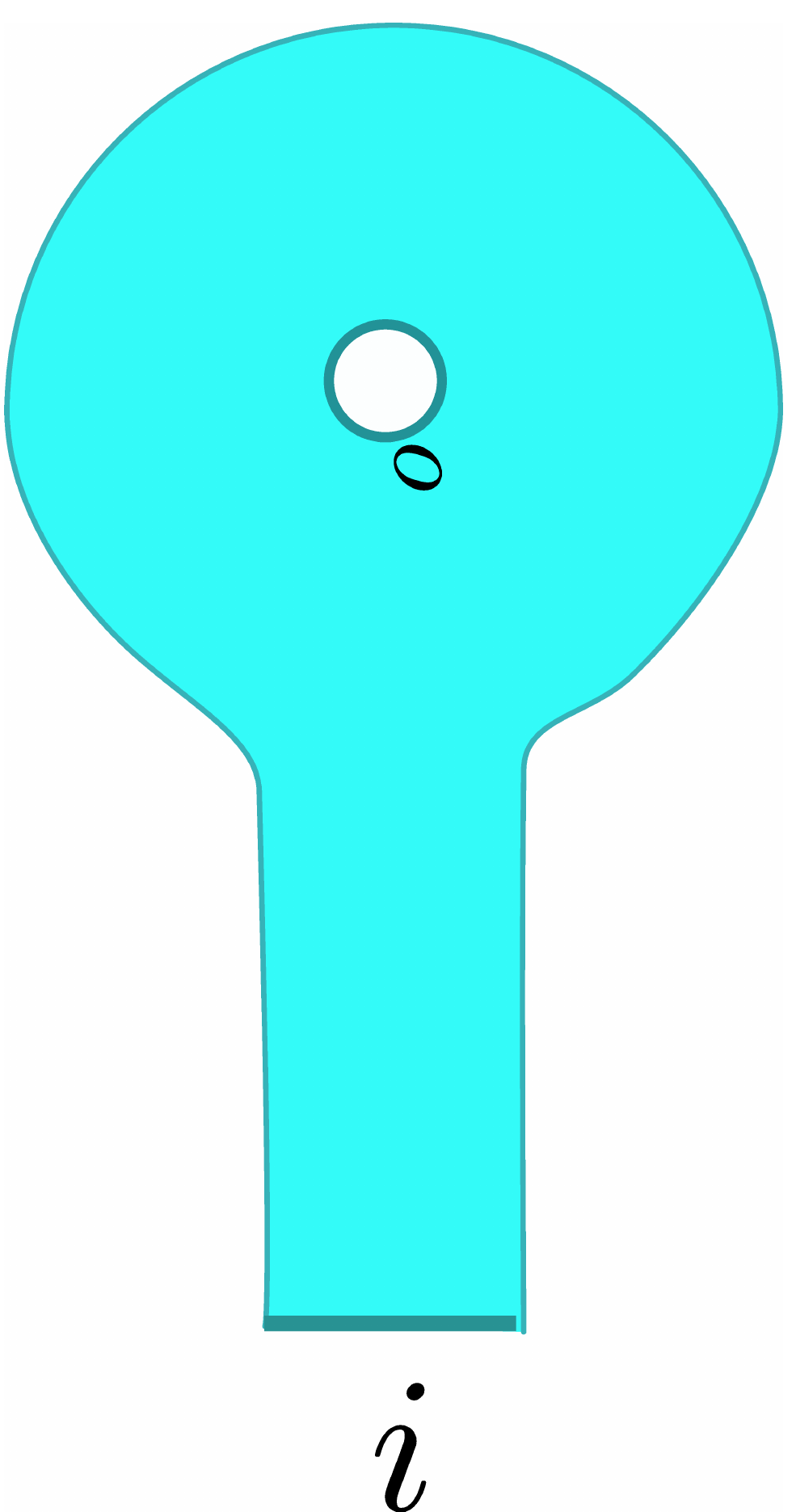}\\
&\Yr_{\Rr27,l}^d=\hspace{4mm} \figbox{0.26}{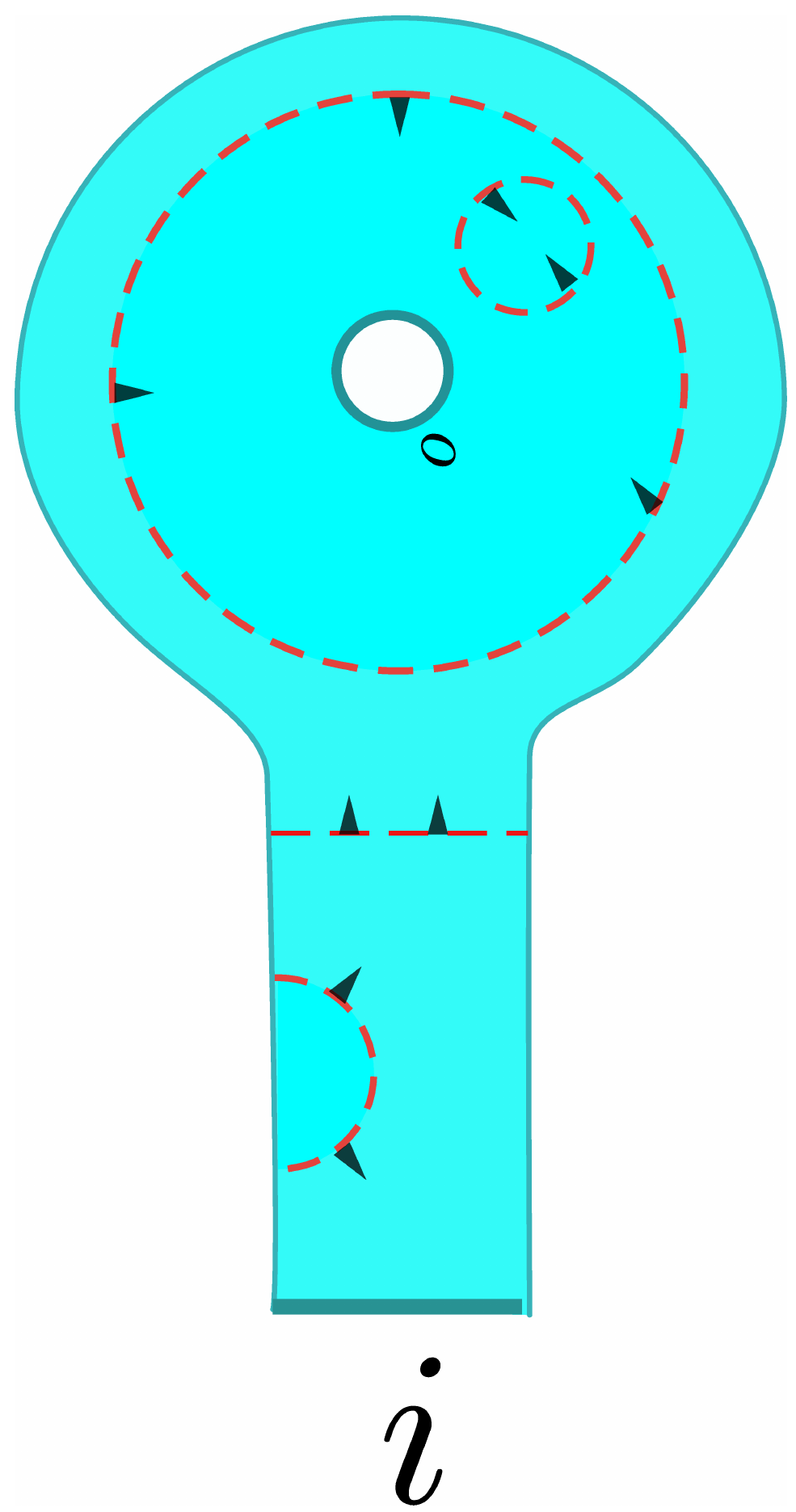}\hspace{29mm}\Yr_{\Rr27,r}^d=\hspace{4mm}\figbox{0.26}{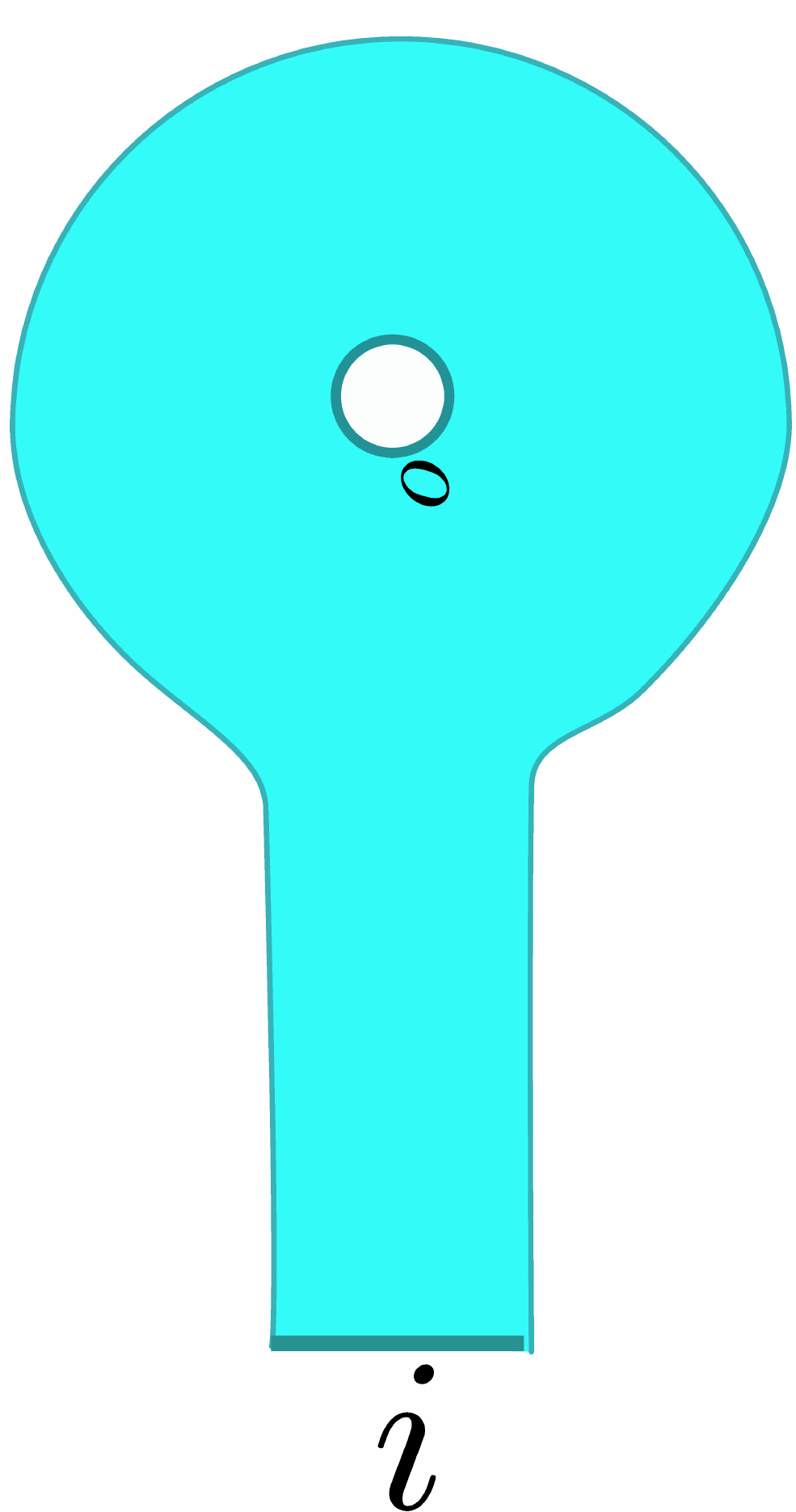}\\
&\Yr_{\Rr28,l}^d=\hspace{6mm} \figbox{0.23}{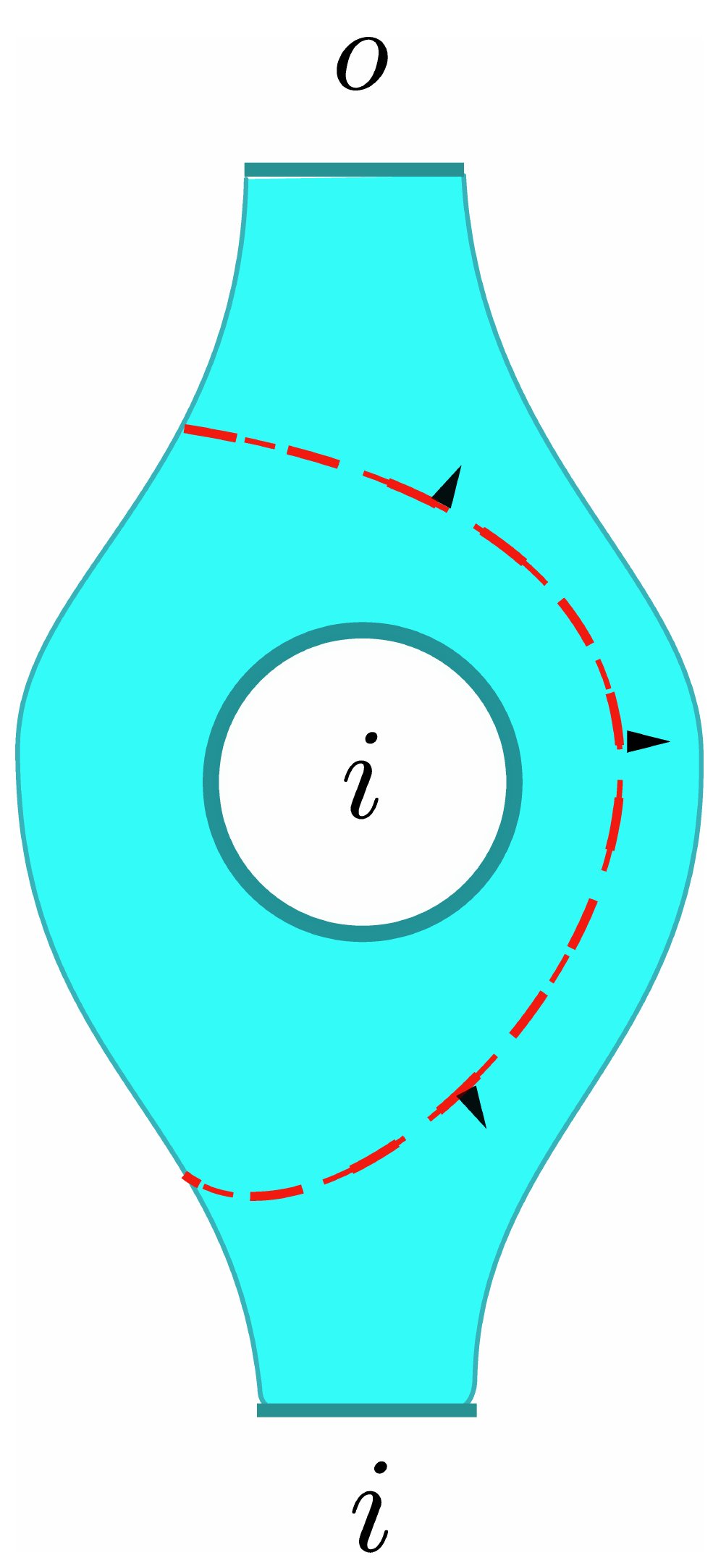}\hspace{31mm}\Yr_{\Rr28,r}^d=\hspace{5mm}\figbox{0.23}{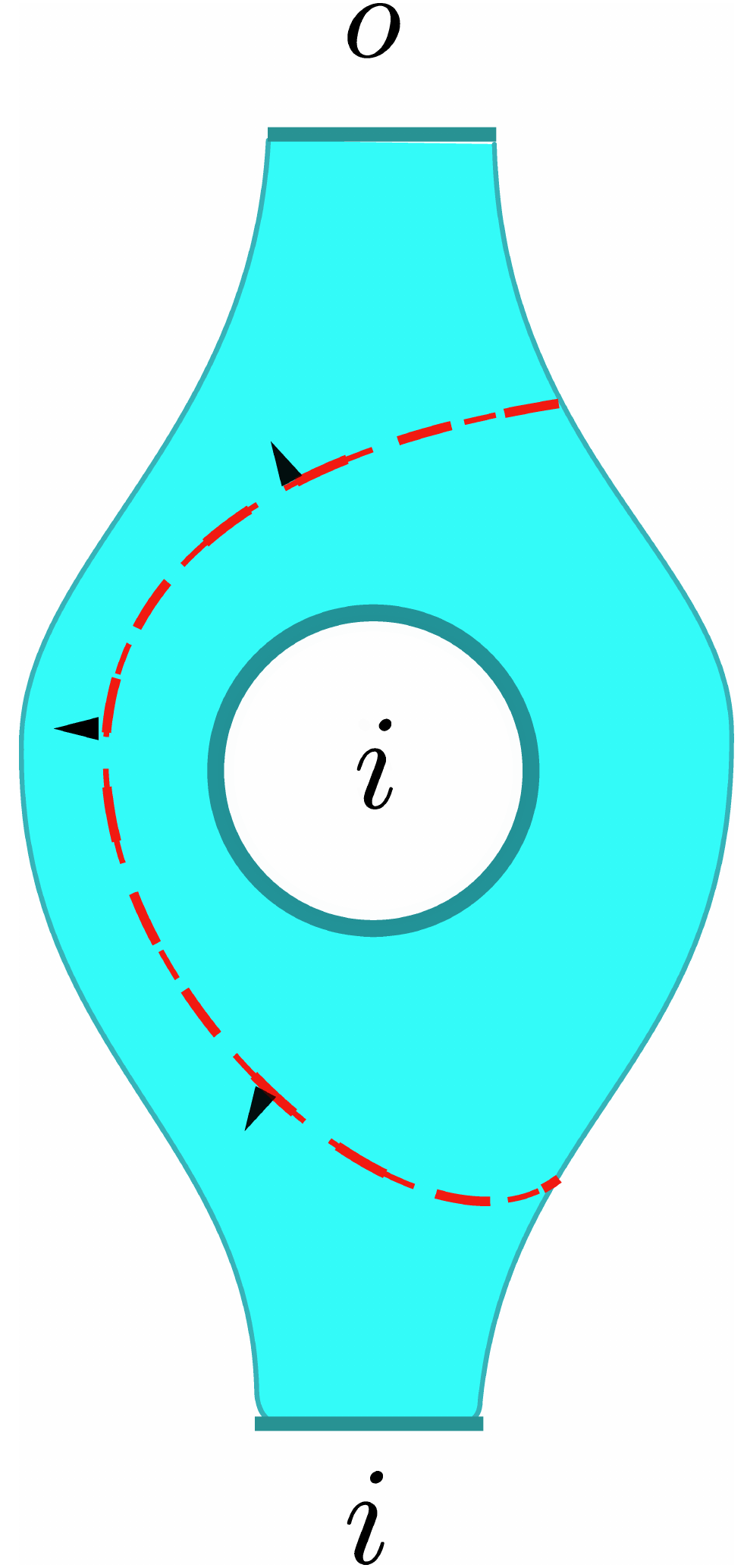}\\
&\Yr_{\Rr29,l}^d=\hspace{5.5mm} \figbox{0.2}{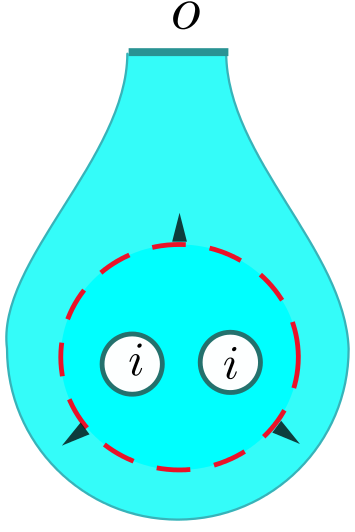}\hspace{29.5mm}\Yr_{\Rr29,r}^d=\hspace{6mm}\figbox{0.2}{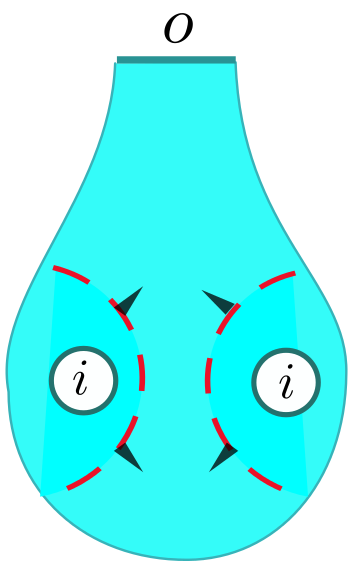}\\
&\Yr_{\Rr30,l}^d=\hspace{4mm} \figbox{0.2}{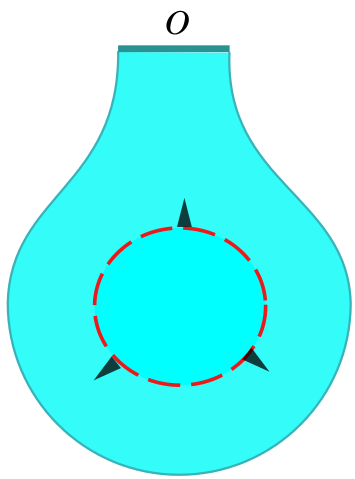}\hspace{30mm}\Yr_{\Rr30,r}^d=\hspace{3mm}\figbox{0.26}{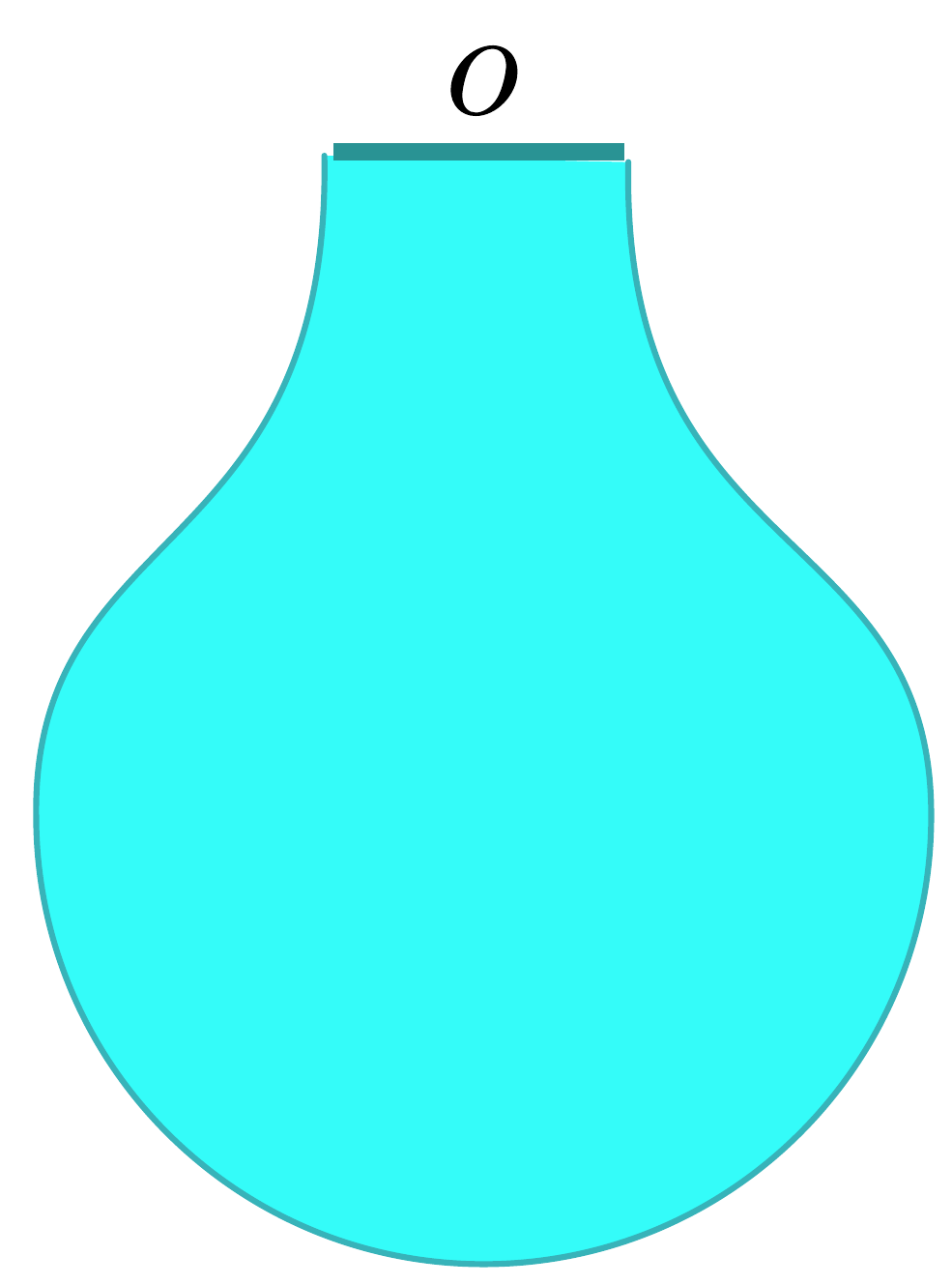}\\
&\Yr_{\Rr31,l}^d=\hspace{2mm} \figbox{0.3}{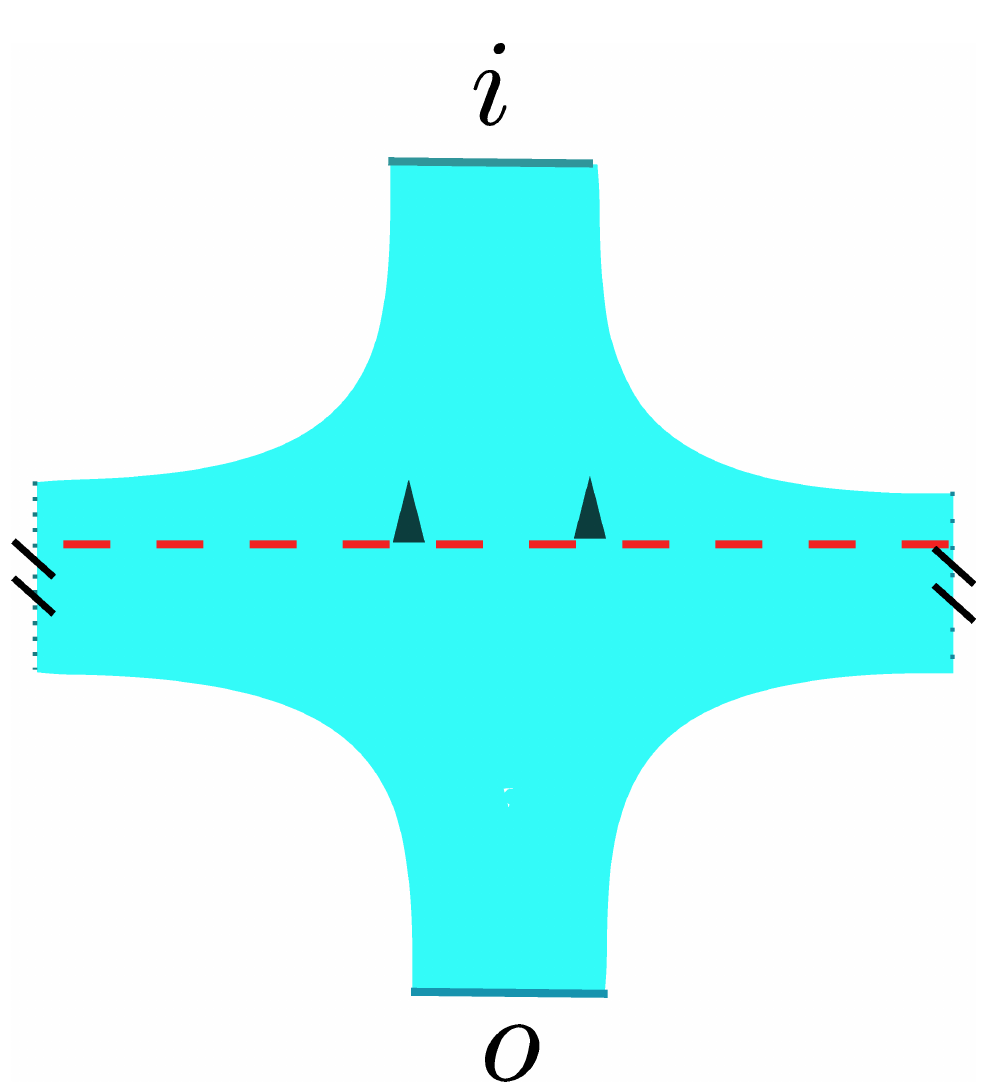}\hspace{27.5mm}\Yr_{\Rr31,r}^d=\hspace{1.5mm}\figbox{0.3}{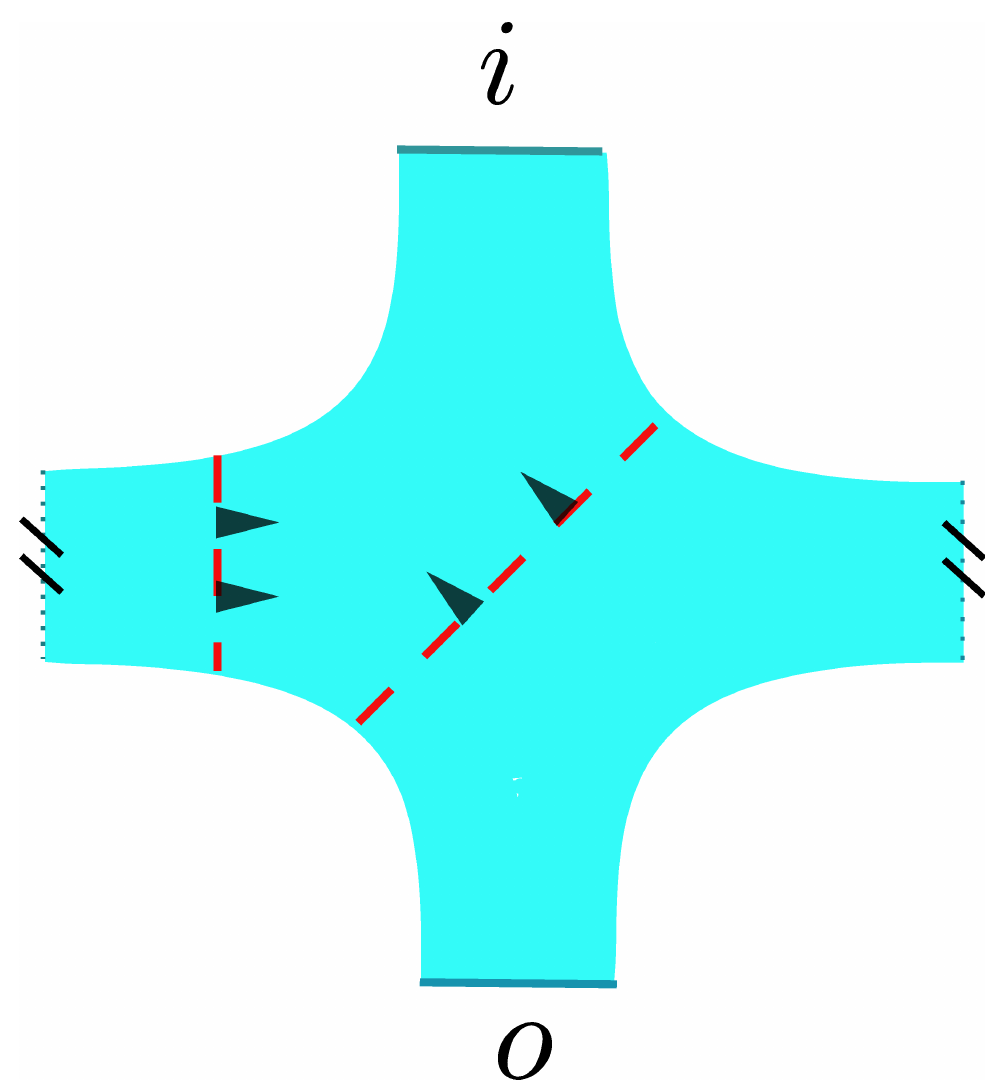}\\
&\Yr_{\Rr32,l}^d=\hspace{6mm} \figbox{0.36}{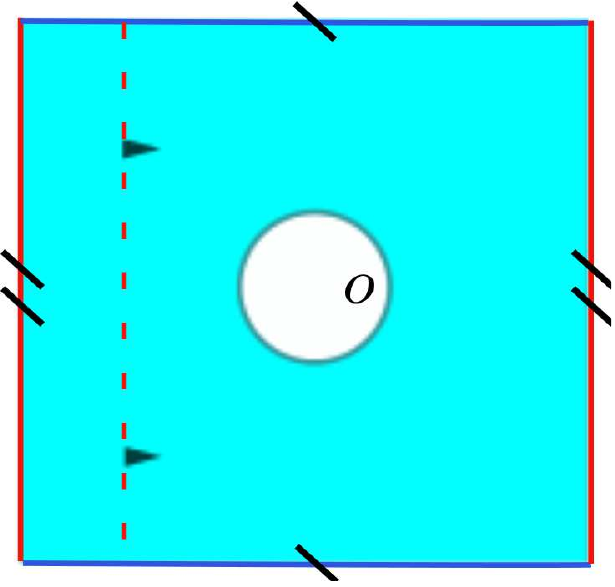}\hspace{30mm}\Yr_{\Rr32,r}^d=\hspace{6mm}\figbox{0.36}{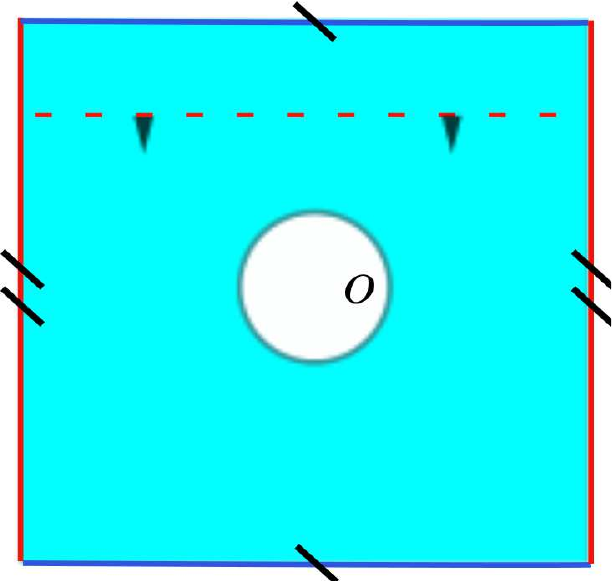}
\end{align*}

These pictures should be interpreted as follows.
For each relation $\Rr n$ there is a fixed world sheet $\Yr_{\Rr n}:= \text{For}^d(\Yr_{\Rr n, l/r})$ in $\WSh$.
For $\Yr_{\Rr n}$ we give two distinct decompositions of $\Yr_{\Rr n}$, namely
\be  \label{eq:Y_Rn-l-r}
  \Yr^d_{\Rr n,l} = ( \Yr_{\Rr n} ,
  ( \alpha_{\Rr n}^1 , \dots , \alpha_{\Rr n}^{l(n)} ) , \varpi_{\Rr n,l} ),
  \quad  \quad
  \Yr^d_{\Rr n,r} = ( \Yr_{\Rr n} ,
  ( \beta_{\Rr n,l}^1 , \dots, \beta_{\Rr n}^{r(n)} ) , \varpi_{\Rr n,r} )   ~.
\ee
The two images $ \Yr^d_{\Rr n,l}$ and $\Yr^d_{\Rr n,r}$ given for each relation $\Rr n$ in the above pictures describe the two different decompositions of $\Yr_{\Rr n}$ into generating world sheets $\Xr_\alpha$, $\alpha \in \Sc$.
The boundary components of $\Yr_{\Rr n}$ labelled $i$/$o$ mark the in-coming/out-going state boundaries, and the dashed lines with arrowheads
 on them show how $\Yr_{\Rr n}$ is to be decomposed into elements of the generating set of world sheets.
For example
\be\begin{array}{rlrl}
 \Yr^d_{\Rr 1,l} \etb=
 ( \Yr_{\Rr 1} , ( mo, \eta o ) , \varpi_{\Rr 1,l} ) ~,
\etb
  \Yr^d_{\Rr 1,r} \etb=
  ( \Yr_{\Rr 1} , ( po ) , \varpi_{\Rr 1,r} ) ~,
\enl
  \Yr^d_{\Rr 26,l} \etb=
  ( \Yr_{\Rr 26} , ( mo, \eps o, \Delta c, \eta c, \iota )
  , \varpi_{\Rr 26,l} ) ~,
\etb
  \Yr^d_{\Rr 26,r} \etb=
  ( \Yr_{\Rr 26} , ( \iota^* ) , \varpi_{\Rr 26,r} ) ~.
\eear\ee

The choice of ordering of the elements of $\Sc$ does not matter, nor does the specific choice of the various morphisms $\varpi_{\Rr n,l\slash r}$, as long as these choices are made once and for all. The only two relations where $\varpi_{\Rr n,l}$ and $\varpi_{\Rr n,r}$ have to be chosen in a correlated way are $\Rr 20$ and $\Rr 25$. In the corresponding pictures the dotted line indicates that $\varpi_{\Rr n,l}$ and $\varpi_{\Rr n,r}$ are related by the action of an element of the mapping class group of $\Yr_{\Rr n}$.
Finally, we will use the convention that if $\Yr_{\Rr n}$ is itself already a generating world sheet
$\Yr_{\Rr n} \equiv \Xr_\alpha$, and $\Yr^d_{\Rr n,r}$  is the trivial decomposition, then we choose
\be
  \Yr^d_{\Rr n,r} = ( \Xr_\alpha , ( \alpha ) , \id)  ~.
\labl{eq:Yn-id-convention}
This is the case for
$\Rr 1$--$\Rr 4$,
$\Rr 10$--$\Rr 15$,
$\Rr 20$--$\Rr 27$ and $\Rr 30$.

\begin{thm}  \label{thm:gen-rel} {\rm
Let $U,F \In \Ob\Fun_\otimes^H(\WSh,\Vect)$ and suppose that $U(f)$ is an invertible linear map for all morphisms $f$ in $\WSh$. Let
$\{ c_\alpha : U(\Xr_\alpha) \rightarrow F(\Xr_\alpha) | \alpha \In \Sc \}$
be a collection of linear maps. There exists a
$\mu \in \text{Nat}_\otimes(U,F)$
with $\mu^\delta = C(\{ c_\alpha | \alpha \In \Sc \} )$ if and only if
\be
  C_{\Yr^d_{\Rr n,l}} = C_{\Yr^d_{\Rr n,r}}
  \qquad \text{for} \quad n=1,\dots,32~,
\labl{eq:sewing-cond}
both sides of which are abbreviations for $C_{\Yr^d_{\Rr n,l}}(\{ c_\alpha | \alpha \In \Sc \} )$ and
$C_{\Yr^d_{\Rr n,r}}(\{ c_\alpha | \alpha \In \Sc \} )$ that are defined in \eqref{eq:nat-C-def}.
}
\end{thm}

This theorem is the key technical result in the present paper. The proof is somewhat lengthy and tedious and has been moved to Appendix \ref{app:proof-nat}.

\subsection{Example: Monoidal natural transformations of $\One$}
\label{sec:Nat-ex}

By Lemma \ref{lem:all-nat-C}, each element of
$\text{Nat}_\otimes(\One^\delta,\One^\delta)$ is of the form
$C(\{ c_\alpha | \alpha \In \Sc \} )$, where in this case
$c_\alpha : \Cb \rightarrow \Cb$ is just a complex number.

We would like to understand the image of $(\,\cdot\,)^\delta :
\text{Nat}_\otimes(\One,\One)
\rightarrow
\text{Nat}_\otimes(\One^\delta,\One^\delta)$.
By Theorem \ref{thm:gen-rel}, $C(\{ c_\alpha | \alpha \In \Sc \} )$
lies in the image of $(\,\cdot\,)^\delta$ if and only if
the relations \erf{eq:sewing-cond} hold. The independent
and non-trivial conditions are
\be
\begin{array}{llllll} \displaystyle
\text{R1} \etb: c_{mo} \, c_{\eta o} = c_{po}
\etb
\text{R3} \etb: c_{\Delta o} \, c_{\eps o} = c_{po}
\etb
\text{R10}\etb: c_{po} \, c_{mo} = c_{mo}
\enl
\text{R11}\etb: c_{po} \, c_{\eta o} = c_{\eta o}
\etb
\text{R12}\etb: c_{po} \, c_{\Delta o} = c_{\Delta o}
\etb
\text{R13}\etb: c_{po} \, c_{\eps o} = c_{\eps o}
\enl
\text{R14}\etb: c_{mc} \, c_{\eta c} = c_{pc}
\etb
\text{R15}\etb: c_{\Delta c} \, c_{\eps c} = c_{pc}
\etb
\text{R21}\etb: c_{pc} \, c_{mc} = c_{mc}
\enl
\text{R22}\etb: c_{pc} \, c_{\eta c} = c_{\eta c}
\etb
\text{R23}\etb: c_{pc} \, c_{\Delta c} = c_{\Delta c}
\etb
\text{R24}\etb: c_{pc} \, c_{\eps c} = c_{\eps c}
\enl
\text{R26}\etb: c_{m o} \, c_{\eps o} \, c_{\Delta c} \,
  c_{\eta c} \, c_{\iota} = c_{\iota^*} \quad
\etb
\text{R27}\etb: c_{\Delta o} \, c_{\eta o} \, c_{m c} \,
  c_{\eps c} \, c_{\iota^*} = c_{\iota} \quad
\etb
\text{R29}\etb: c_{mc} \, c_{\iota} = c_{\iota}^2 \, c_{mo}
\enl
\text{R30}\etb: c_{\eta c} \, c_{\iota} = c_{\eta o}
\etb
\text{R31}\etb: c_{\iota} \, c_{\iota^*} = c_{mo} c_{\Delta o}
\end{array}
\ee

\noindent \nxt
Suppose that $c_{po} \neq 0$. By R1 and R3 this implies
$c_{mo},c_{\eta o},c_{\Delta o},c_{\eps o} \neq 0$.
Then by R31 also $c_{\iota},c_{\iota^*} \neq 0$ and
from R14--R29 we see that in fact $c_\alpha \neq 0$ for all
$\alpha \in \Sc$. R10 and R21 imply $c_{po}=c_{pc}=1$, and
R1, R3, R14, R15, R29 and R30 show
\be
  c_{mo} = \frac{1}{c_{\eta o}} ~~,\quad
  c_{\Delta o} = \frac{1}{c_{\eps o}} ~~,\quad
  c_{mc} = \frac{1}{c_{\eta c}} ~~,\quad
  c_{\Delta c} = \frac{1}{c_{\eps c}} ~~,\quad
  c_{\iota} = \frac{c_{\eta o}}{c_{\eta c}} ~~,\quad
  c_{\iota^*} = \frac{c_{\eps o}}{c_{\eps c}} ~~.
\labl{eq:ex-cond-aux2}
Finally, R31 gives the condition
\be
  c_{\eps c} = \frac{1}{c_{\eta c}} \,
  \big( c_{\eta o} \, c_{\eps o} \big)^2
  ~.
\labl{eq:ex-cond-aux1}
So if $c_{po} \neq 0$, the solutions are parametrised by
$c_{\eta o}, c_{\eps o}, c_{\eta c} \In \Cb^\times$.
The remaining constants are given
by \erf{eq:ex-cond-aux2} and \erf{eq:ex-cond-aux1}.
One checks that this then solves R1--R32.

\medskip\noindent
\nxt
Suppose that $c_{po} = 0$ and $c_{pc} \neq 0$. Then by R10--R13 also
$c_{mo} = c_{\eta o} =c_{\Delta o}=c_{\eps o} = 0$, and by R26 and R27
$c_\iota = c_{\iota^*} = 0$. R14--R24 are solved by $c_{pc}=1$,
$c_{mc} = 1/c_{\eta c}$ and $c_{\Delta c} = 1/c_{\eps c}$ for
any choice of $c_{\eta c},c_{\eps c} \In \Cb^\times$.

\medskip\noindent
\nxt
Suppose that $c_{po} = 0$ and $c_{pc} = 0$. Then R1--R31 force
$c_{\alpha}=0$ for all $\alpha \In \Sc$.

\medskip\noindent
Altogether this gives
$\text{Nat}_\otimes(\One,\One) \cong \{0\} \sqcup
( \Cb^\times \times \Cb^\times ) \sqcup
( \Cb^\times \times \Cb^\times \times \Cb^\times )$ where the element
of $\text{Nat}_\otimes(\One,\One)$ is determined by the linear maps assigned to the generating world sheets,
\be
\begin{array}{rcllllllllllll}
& &
c_{mo} & c_{\Delta o} & c_{\eta o} & c_{\eps o} & c_{mc} &
c_{\Delta c} & c_{\eta c} & c_{\eps c} & c_{po} & c_{pc} & c_{\iota} &
c_{\iota^*}
\\
\{0\} & \mapsto &
0 & 0 & 0 & 0 & 0 & 0 & 0 & 0 & 0 & 0 & 0 & 0
\\
(\alpha,\beta) & \mapsto &
0 & 0 & 0 & 0 & \alpha^{-1} & \beta^{-1} & \alpha & \beta & 0 & 1 & 0 & 0
\\
(\alpha,\beta,\gamma) & \mapsto &
\alpha^{-1} & \beta^{-1} & \alpha & \beta &
  \tfrac{1}{\alpha\beta\gamma} &
\tfrac{\gamma}{\alpha \beta} & \alpha\beta\gamma &
\tfrac{\alpha \beta}{\gamma} & 1 & 1 & \tfrac{1}{\beta\gamma} &
\tfrac{\gamma}{\alpha}
\eear\ee

\sect{Solutions to the sewing constraints}

\subsection{The functor $\Bl$ obtained from a modular tensor category}
\label{sec:Bl-def}

Let $\mathcal{C}$ be an abelian semi-simple finite $\Cb$-linear ribbon category
with simple tensor unit $\one$. Let $\Ic$ be the set of isomorphism classes of simple objects in $\mathcal{C}$ and $\{U_i,i\in \Ic \}$ the chosen representatives. We will adopt the convention that $U_0:=\one$ is the tensor unit in $\mathcal{C}$. For $V \in \Cc$, we choose a basis $\{ b_{V}^{(i;\alpha)} \}$ 
of $\Hom_{\Cc}(V, U_i)$
and the dual basis $\{ b_{(i;\beta)}^{V} \}$ of $\Hom_{\Cc}(U_i, V)$ 
for $i\in \Ic$ such that
$b^{(i;\alpha)}_{V} \circ b_{(i;\beta)}^{V} = 
\delta_{\alpha\beta}\, \id_{U_i}$. We will use the following graphical notation
\be   \label{eq:b-dual-b}
b_V^{(i; \alpha)} = 
  \raisebox{-23pt}{
  \begin{picture}(30,52)
   \put(0,8){\scalebox{.75}{\includegraphics{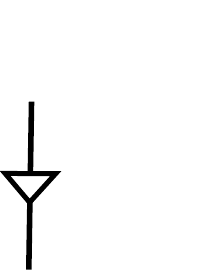}}}
   \put(0,8){
     \setlength{\unitlength}{.75pt}\put(-18,-11){
     \put(39,36)  {\scriptsize $ \alpha $}
     \put(23,65)  {\scriptsize $ U_i $}
     \put(23, 2)  {\scriptsize $ V $}
     }\setlength{\unitlength}{1pt}}
  \end{picture}}
\quad  , \qquad
b_{(i;\alpha)}^V =
  \raisebox{-23pt}{
  \begin{picture}(30,52)
   \put(0,8){\scalebox{.75}{\includegraphics{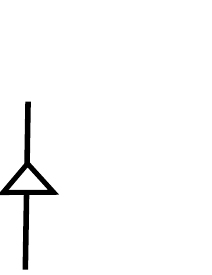}}}
   \put(0,8){
     \setlength{\unitlength}{.75pt}\put(-18,-11){
     \put(37,34)  {\scriptsize $ \alpha $}
     \put(20, 2)  {\scriptsize $ U_i $}
     \put(20,65)  {\scriptsize $ V $}
     }\setlength{\unitlength}{1pt}}
  \end{picture}} 
  ~~.
\ee
for a pair of mutually dual basis in \eqref{eq:mod-inv}, \eqref{eq:cardy-C} and \eqref{Cardy condition}. Let $s_{ij}$ denote following number:
\begin{equation}
s_{ij}=\figbox{0.2}{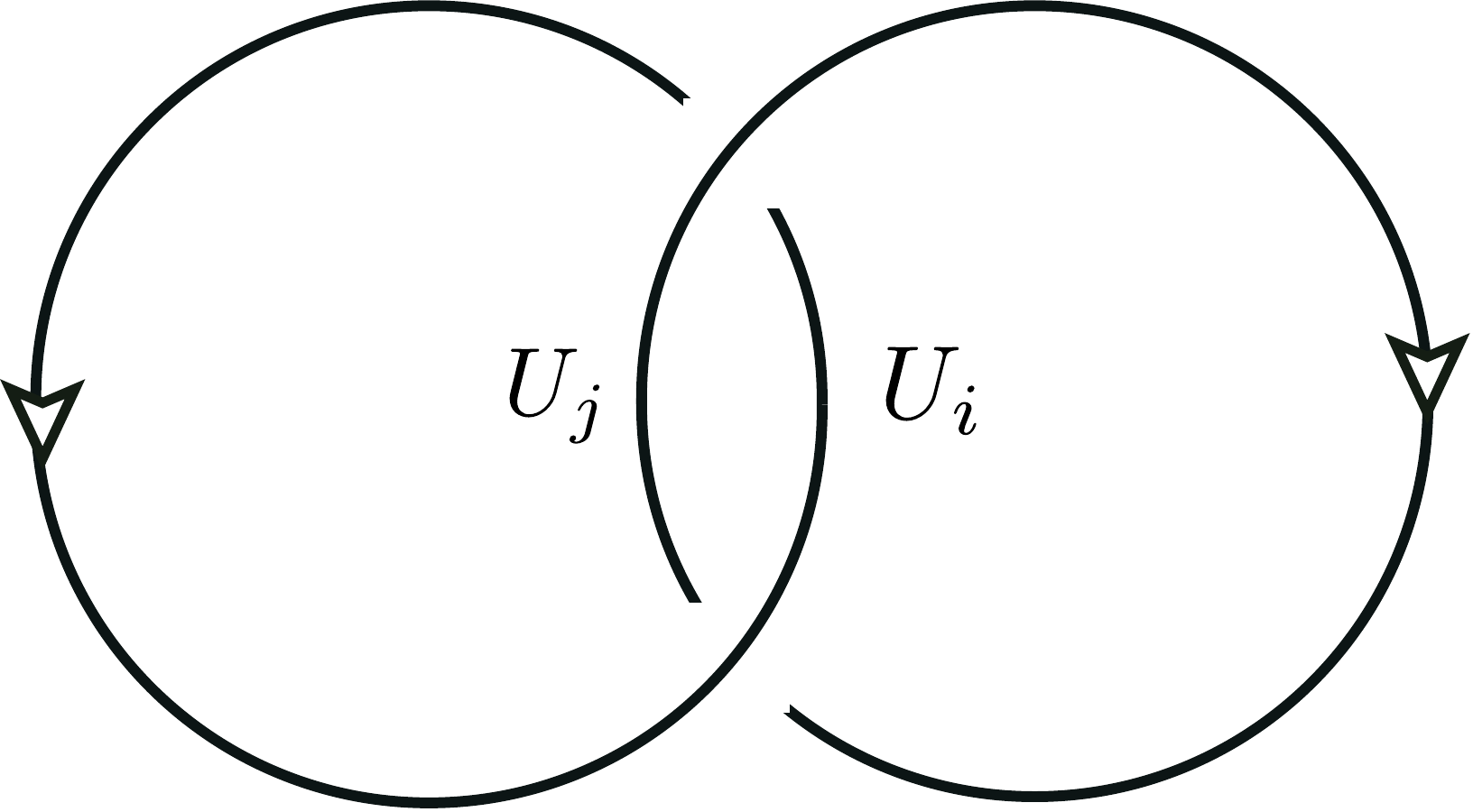}
\end{equation}
Then $\mathcal{C}$ is a modular tensor category if and only if the matrix $(s_{ij})_{i,j \in \Ic}$ is invertible.
The dimension of the modular tensor category $\mathcal{C}$ is defined as $\text{Dim}(\mathcal{C}):=\sum_i (\text{dim} U_i)^2$. We will denote
its square root by $D=\sqrt{\text{Dim}(\mathcal{C})}$. 

\medskip
For the application to rational conformal field theory we need to use a
specific symmetric monoidal functor  $\Bl : \WSh \rightarrow \Vect$.
Below we give a slightly modified (and simplified) version of
the treatment in \cite[Sect.\,3.3]{unique}. The definition of $\Bl \,{\equiv}\, \Bl(\Cc,\Bop,B_l,B_r)$ takes as input
\bnu
\item a modular tensor category $\Cc$,
\item three nonzero objects $\Bop, B_l, B_r \in \Cc$. 
\enu
For the constructions in the next section it is always possible to choose $\Bop = B_l = B_r$. However, these three objects play slightly different roles and so we prefer to distinguish them.

The functor $\Bl$ will be defined with the help of a three-dimensional TFT. Let us briefly
state our conventions for the 3-d TFT; for more details see e.g.\ \cite{turaev-bk,KRT,baki} or \cite[Sect.\,3.2]{unique} which uses the same conventions as here.

Given a modular tensor category $\Cc$, the construction of \cite{turaev-bk} allows
one to construct a 3-d TFT, that is, a symmetric monoidal functor $\tftC$ from a
geometric category $\mathcal{G}_\Cc$ to $\Vect$.
The objects of $\mathcal{G}_\Cc$ are extended surfaces, and we take the morphisms to be homeomorphisms or extended cobordisms. The tensor product is given by the disjoint union.

A {\em extended surface} is an oriented, closed surface
with a finite set of disjoint ordered 
marked arcs labeled by pairs $(U,\epsilon)$, where
$U\In \mathcal{O}bj(\Cc)$ and $\epsilon\in\{+,-\}$, and with a choice of 
Lagrangian subspace $\lambda\,{\subset}\, H^1(\mathrm E,\Rb)$. 
A {\em homeomorphism of extended surfaces} is a homeomorphism of the underlying
surfaces preserving orientation, ordered marked arcs and Lagrangian subspaces.
For an extended surface $\mathrm{E}$ denote by $\overline{\mathrm{E}}$ the extended surface
obtained from $\mathrm{E}$ by inverting the orientation of the surface and of the marked arcs, replacing each label $(U,\epsilon)$ by $(U,-\epsilon)$ 
and keeping the order of the marked arcs. The choice of Lagrangian subspace remains the same.

The boundary $\partial \Mr$ of an oriented three-manifold $\Mr$ will be oriented by the inward pointing normal convention. For example, the boundary of the three-manifold $\{ (x,y,z) | z \ge 0 \}$ with standard orientation is the same
as the one induced by the embedding $(x,y) \mapsto (x,y,0)$ of $\Rb^2$ with its standard orientation.

A {\em extended cobordism} is
a triple $(\Mr, n, h)$ where $\Mr$ is a cobordism of extended surfaces,
$h{:}\ \partial\Mr\To\overline{\mathrm E}\sqcup\Fr$ is a homeomorphism of extended
surfaces, and $n\In\mathbb{Z}$ is a weight which is needed to make $\tftC$
anomaly-free \cite[Sect.\,IV.9]{turaev-bk}. We denote the in-going component (the pre-image of $\overline{\mathrm{E}}$ under $h$) by $\partial_- M$ and the out-going component (the pre-image of $\Fr$) by $\partial_+M$. The cobordism $\Mr$ can
contain ribbons, which are labeled by objects of $\Cc$, and coupons, which are
labeled by morphisms of $\Cc$. Ribbons end on coupons or on the arcs of
$\mathrm E$ and $\mathrm F$.
Two cobordisms $(\Mr,n,h)$ and $(\Mr',n,h')$ from ${\rm E}$ to ${\rm F}$
are {\em equivalent\/} iff there exists a homeomorphism
$\varphi{:}\ \Mr \To \Mr'$ taking ribbons and coupons of $\Mr$ to identically
labeled ribbons and coupons of $\Mr'$ and obeying $h \eq h' \cir \varphi$.
The functor $\tftC$ is constant on the equivalence classes of cobordisms.

\medskip

Next we define the functor $\Bl$. As a start,
from a world sheet $\Xr$ we construct a closed extended surface $\Xhat \;{\equiv}\;
\Xhat(\Bop,B_l,B_r)$,\label{def:Xhat}
which we call the {\em extended double of\/} $\Xr$. It is obtained by gluing
a standard disk with a marked arc to each boundary component of $\Xtil$:
Let $\vec D$ be the unit disk $\{|z|\,{\le}\,1\} \,{\subset}\, \Cb$ with
a small arc embedded on the real axis, centred at $0$ and pointing towards
$+1$. The orientation of $\vec D$ is that induced by $\Cb$. Then we set
  \be
  \Xhat := \Xtil \sqcup \big( \pi_0(\partial\Xtil) \ti \vec D \big) /{\sim}
  \ee
where the equivalence relation is given in terms of the boundary parametrisation,
  \be
  (a,z) \sim \delta_a^{-1}(-\bar{z}) \qquad {\rm for}\quad
  a \In \pi_0(\partial\Xtil),\, z \In \partial \vec D \,.
  \ee
Here the complex conjugation $(-\bar{z})$ is needed for $\Xhat$ to be oriented
and the minus sign is chosen to agree with the sewing procedure in
Definition \ref{def:oc-sewing}(ii).
For $a \In \pi_0(\partial\Xtil)$
the arc on the disk $\{a\} \ti \vec D$ is marked by $(U_a,\eps_a)$,
where $U_a \In \{\Bop,B_l,B_r\}$ and $\eps_a \In \{\pm\}$ are chosen
as follows.
\\[.2em]
\nxt If $a\In b^\text{in}$, then $\eps_a\eq{+}\,$,\, otherwise $\eps_a\eq{-}\,$.
\\[.2em]
\nxt If $\imath_*(a) \eq a$, i.e. $a\in (b^\text{in/out}_\Xr)_\text{op}$ (recall Remark \ref{rem:unor}),
then $U_a \eq \Bop$.
\\[.2em]
\nxt If $\imath_*(a) \,{\neq}\, a$ then $U_a\eq B_l$ if $a\in (b^\text{in}_\Xr)_\text{cl}^+$, and
$U_a\eq B_r$ if $a\in (b^\text{in}_\Xr)_\text{cl}^-$ (recall Remark \ref{rem:unor}). \newline
Note that the involution $\imath{:}\ \Xtil \To \Xtil$ can be extended
to an involution $\hat\imath{:}\ \Xhat \To \Xhat$ by taking it to be
$(a,z) \mapsto (\imath_*(a),\bar{z})$ on $\pi_0(\partial\Xtil) \ti \vec D$. 
To turn $\Xhat$ into an extended surface we still need to
specify a Lagrangian subspace $\lambda \,{\subset}\, H^1(\Xhat,\Rb)$.
To do this we start by taking the {\em connecting manifold\/}
  \be
  \Mr_{\Xr} = \Xhat \ti [-1,1] /{\sim}
  \qquad \text{where for all}~ x\In\Xhat ~,~~ (x,t) \sim (\hat\imath(x),-t) \,,
  \labl{eq:conmf-def}
which has the property that $\partial \textrm{M}_{\Xr} \eq \Xhat$. Then
$\lambda$ is the kernel of the resulting homomorphism
$H^1(\Xhat,\Rb) \To H^1(\textrm{M}_{\Xr},\Rb)$.
We refer to appendix B.1 of \cite{tft5} for more details.
The ordering of the marked arcs on $\Xhat$ is induced from the ordering of the boundary components of the world sheet $\Xr$.

As an example, consider as world sheet $\Xr$ a disk with two open and two closed state boundaries; in this case the extended double is given by a sphere with six marked
arcs:
\be
  \Xr ~=~
  \raisebox{-42pt}{
  \begin{picture}(95,92)
  \put(0,0)     {\scalebox{.28}{\includegraphics{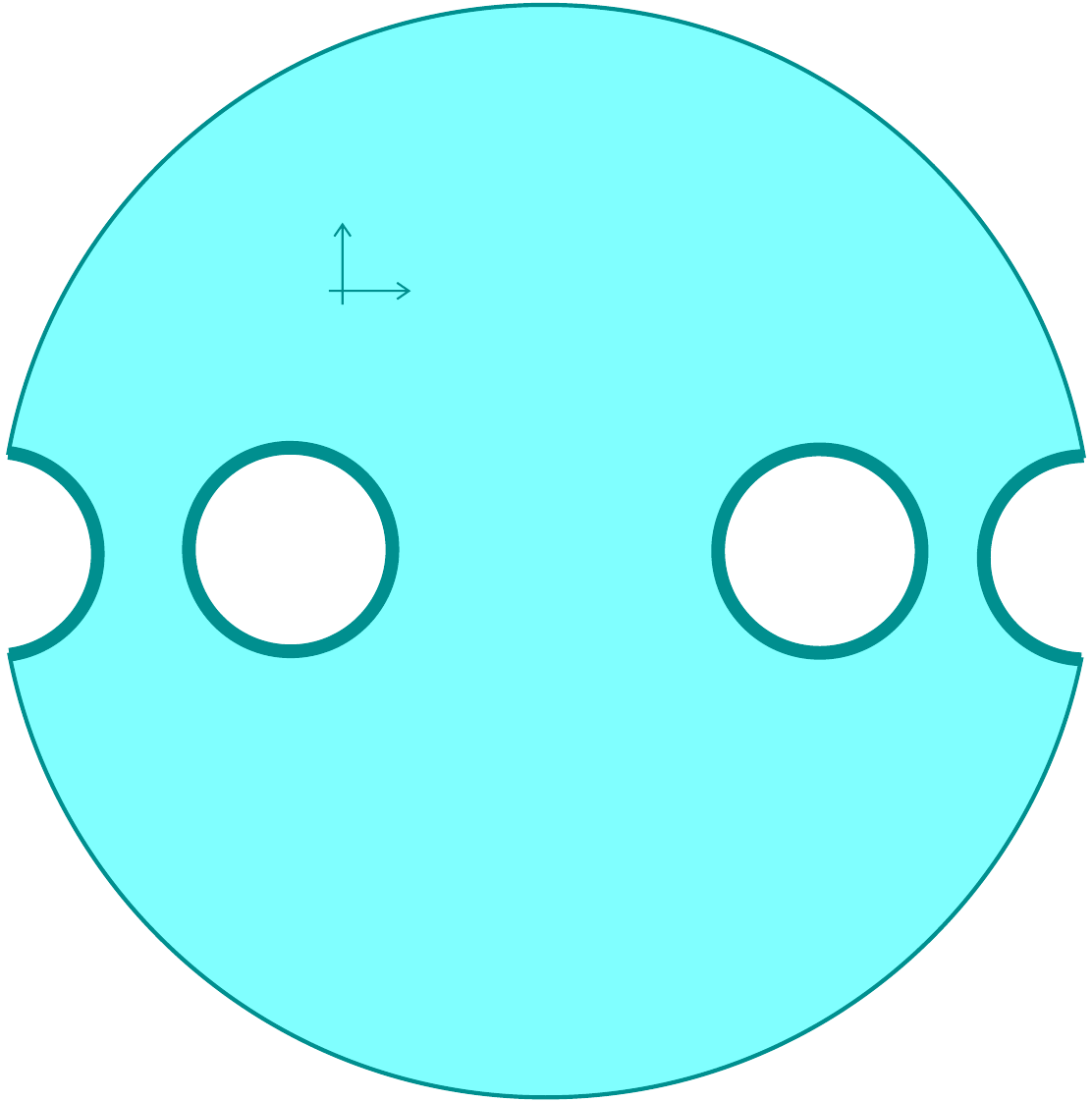}}}
  \put(0,43)    {\small$ i $}
  \put(19,43)   {\small$ i $}
  \put(61,43)   {\small$ o $}
  \put(84,43)   {\small$ o $}
  \put(27,74)   {\tiny{$2$}}
  \put(34,65)   {\tiny{$1$}}
  \end{picture}}
  \qquad , \qquad
  \Xhat ~=~
  \raisebox{-50pt}{
  \begin{picture}(110,105)
  \put(0,0)     {\scalebox{.28}{\includegraphics{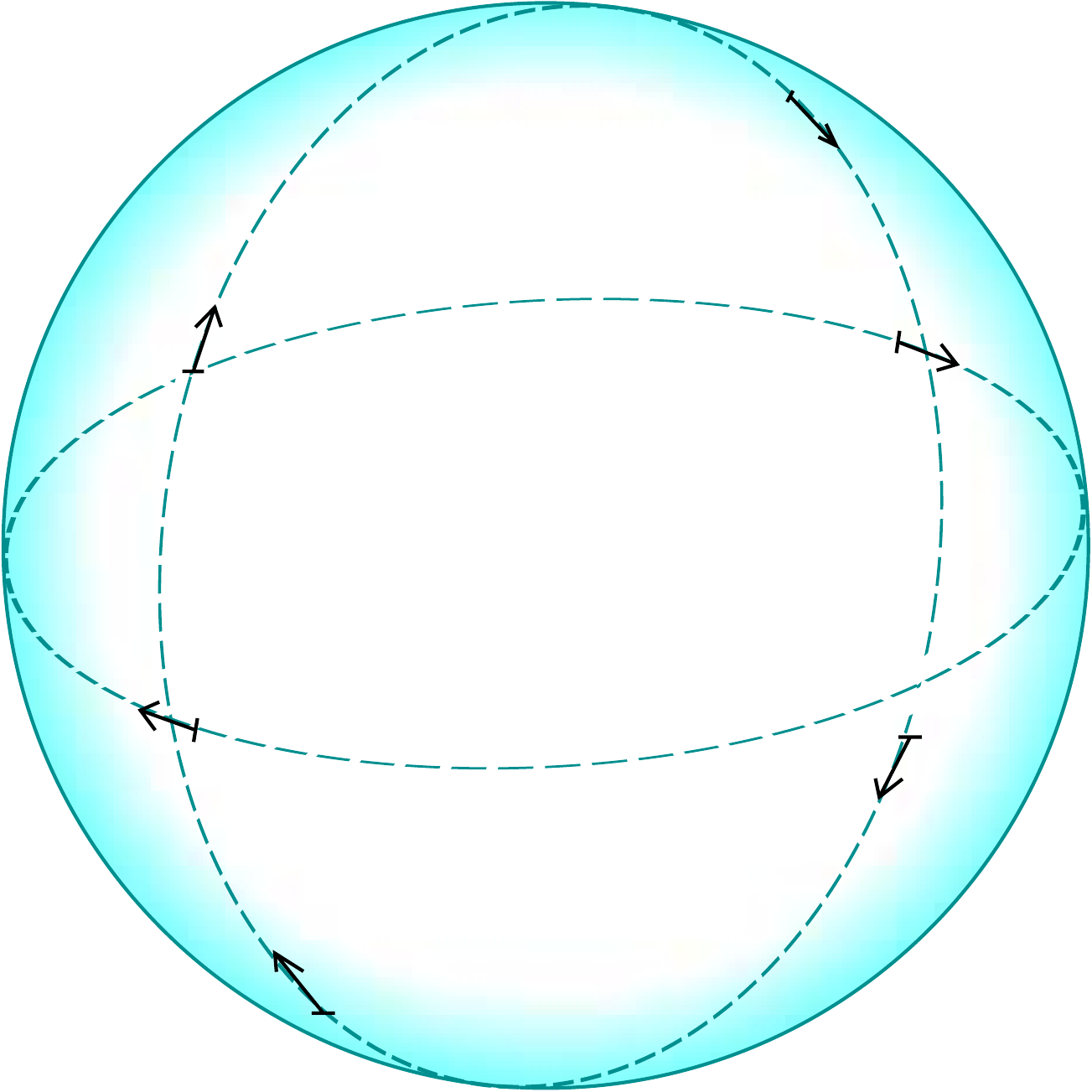}}}
  \put(31,10)   {\tiny{$(B_r{,}+)$} }
  \put(18,36)   {\tiny{$(\Bop{,}+)$} }
  \put(21,69)   {\tiny{$(B_l{,}+)$} }
  \put(60,30)   {\tiny{$(B_r{,}-)$} }
  \put(60,68)   {\tiny{$(\Bop{,}-)$} }
  \put(54,91)   {\tiny{$(B_l{,}-)$} }
  \end{picture}}
\labl{eq:ws-example-1}
The 3-d TFT assigns to the extended surface $\Xhat$ a complex vector space
$\tftC(\Xhat)$. This defines the functor $\Bl$ on objects
(recall that $\Xhat$ is an abbreviation for $\Xhat(\Bop,B_l,B_r)$)
  \be
  \Bl(\Xr) := \tftC( \Xhat ) \,.
  \labl{eq:bl(X)-def}
Now we define the functor $\Bl$ on a morphism $\varpi \eq (\sew,f) \In \Hom(\Xr,\Yr)$. First, we will define an extended cobordism $\Mr_\varpi:=(\Mr_\varpi, 0, h)$ associated to $\varpi$ as follows: 
\bnu
\item Note that we can extend
the isomorphism $f{:}\ \widetilde{\sew(\Xr)} \To \widetilde \Yr$
to an isomorphism $\hat f {:}\ \widehat{\sew(\Xr)} \To \widehat{\Yr}$ between extended surfaces by taking
it to be the identity map on the disks $\vec D$ which are glued to the
boundary components of $\widetilde{\sew(\Xr)}$ and $\widetilde \Yr$. 

\item Next we construct a morphism
$\Mr_\varpi:\, \Xhat \To \widehat{\Yr}$ as an extended cobordism. It is given by
the cylinder over $\Xhat$ modulo an equivalence relation,
  \be
  \Mr_\varpi := \Xhat \times [0,1] / {\sim} \,.
  \labl{eq:shat-cobord}
The equivalence relation identifies points on the boundary $\Xhat\ti
\{1\}$ according to the sewing $\sew$. Namely, for each pair $(a,b) \In \sew$ and
for all $z \In \vec D$ we identify the point $(a,z,1)$ in $\{a\}\ti\vec D
\ti\{1\}$ with the point $(b,-\bar{z},1)$ in $\{b\}\ti\vec D\ti\{1\}$. In other words, $\partial_+\Mr_\varpi = \widehat{\sew(\Xr)}$.
The homeomorphism $h: \partial_- \Mr_\varpi \to \overline{\Xhat}$ is the identity map (which is orientation reversing because the boundary $\partial \Mr$ is oriented by the inward pointing normal) and $h: \partial_+ \Mr \to \widehat{\Yr}$ is given by 
$\hat f$.

\enu
Now we define, for $\varpi \eq (\sew,f)$,
  \be
  \Bl(\varpi) := \tftC(\Mr_\varpi) \,.
  \labl{eq:bl(m)-def}
By \cite[Prop.\,3.8]{unique}, such defined $\Bl: \WSh \to \Vect$ is indeed
a symmetric monoidal functor.

\medskip

We will use the modular tensor category $\Cc_-$, which is obtained from $\Cc$ by replacing the braiding $c$ by the anti-braiding and the twist $\theta$ by its inverse, and the modular tensor category $\CxC:=\Cc_+\boxtimes\Cc_-$, where $\Cc_+ \equiv \Cc$ and $\boxtimes$ denotes the Deligne product of abelian categories, see \cite[Def.\,1.1.15]{baki} and \cite[Sect.\,2.3]{part1}. The objects of  $\CxC$ are pairs $U\ti V$ for $U,V\in \Cc$, as well as direct sums of these, and the morphism spaces are
\be
 \Hom_{\CxC}(V \ti W,V' \ti W') =
 \Hom_{\Cc_+}(V,V')\otimes_\Cb \Hom_{\Cc_-}(W,W')
\ee
and direct sums of these.

\begin{rema}  {\rm
(i) The definition of $\Bl(\Cc,\Bop,B_l,B_r) : \WSh \rightarrow \Vect$ is similar to that of a 
$\Cc$-extended two-dimensional topological modular functor, see \cite[Defs.\,5.7.5 \& 5.7.9]{baki}. 
However, there are also two noteworthy differences. First of all, the surfaces in the source category $\WSh$ are not 
extended surfaces: their boundary components are not labelled by objects of $\Cc$. 
Second, the gluing operation is not required to provide an isomorphism in the target category $\Vect$ (and in general will not do so).
\\[.3em]
(ii) In \cite[Sect.\,3.3]{unique} $\Bl$ was defined in terms of three additional pieces of data: an object $\Acl$ of $\CxC$ and morphisms $e \In
\Hom_{\CxC}(\Acl,\BlxBr)$ and $r \In \Hom_{\CxC}(\BlxBr,\Acl)$
such that $(\Acl,e,r)$ is a retract of $\BlxBr$. The definition given above amounts to the special case $\Acl = \BlxBr$ and $e = r = \id$. We chose the present approach because it simplifies the presentation considerably. This does not result in a loss of generality, instead the realisation of $\Acl$ as a retract of $\BlxBr$ will now appear in the statement of Theorem \ref{thm:cardy-sew} below.
}
\end{rema}

\subsection{Sewing constraints and Cardy algebras}
\label{sec:sew-cardy}

As explained in sections 3.4 and 6.1 of \cite{unique}, the consistency conditions that conformal field theory correlators
have to obey under sewing of world sheets can be understood as the properties of a monoidal natural transformation.
This motivates the following definition.

\begin{defn}  \label{def:sol-sew}  {\rm
Let $\Cc$ be a modular tensor category and let
$\Bop, B_l, B_r \in \text{Ob}(\Cc)$. A
{\em solution to the sewing constraints} is a monoidal natural
transformation $\Cor$ from $\One$ to $\Bl(\Cc,\Bop,B_l,B_r)$.
}
\end{defn}

\begin{rema}  \label{rema:relative-qft} {\rm
In \cite{relative-qft}, Freed and Teleman introduced the notion of a quantum field theory relative to an extended $(n{+}1)$-dimensional quantum field theory $\alpha$. They define a relative quantum field theory to be a homomorphism $F: \mathbf{1} \to \tau_{\leq n} \alpha$, where $\mathbf{1}$ is the trivial theory and $\tau_{\leq n} \alpha$ is the truncation of $\alpha$ to $(n{-}1)-$ and $n$-manifolds. The notion of a solution to the sewing constraints is closely related to that of a relative quantum field theory. Namely, in our case $n=2$ and the $(2{+}1)$-dimensional theory $\alpha$ is the Reshetikhin-Turaev TQFT $\tftC$. The truncated theory $\tau_{\leq 2}\alpha$ corresponds to the symmetric monoidal functor\footnote{
Of course, $\WSh$ is not the same as the category of two-dimensional bordisms (which is the source of $\tau_{\leq 2} \alpha$), but instead (a variant of) $\WSh$ is obtained from two-dimensional bordisms by passing to lists of composable arrows as objects and partial compositions as morphisms.} 
$\Bl : \WSh \rightarrow \Vect$, the trivial theory to the functor $\One$, and the homomorphism $F$ to the monoidal natural transformation $\Cor$. We believe that there is a precise way in which a solution to the sewing constraints for a given rational vertex operator algebra $V$ produces a relative quantum field theory in the sense of \cite{relative-qft}, but we have not worked through the details.
}
\end{rema}

\begin{rema}  {\rm
Theorem \ref{thm:gen-rel} allows us to construct solutions to the sewing
constraints in terms of generators and relations. These relations can
be compared to those in \cite{Le}, which is the original work on
generators and relations for open/closed
conformal field theory correlators, and
to \cite{Lz,LP,MS}, which deal with open/closed
topological field theory.
\\[.3em]
(i)~With only six relations the list of sewing constraints in \cite{Le} is somewhat shorter than the 32 conditions given here. The reason is that in \cite{Le} only in-coming state boundaries are considered, R20 and R25 are implicit because single-valuedness of the three-point function on the sphere and $h\,{-}\,\bar h \in \Zb$ are assumed, R1--5, 10--15, 21--24, 26, 27, 30 do not appear because only decompositions that reduce the complexity of the world sheet are mentioned. The remaining relations are as in \cite[Fig.\,9]{Le}: relation (a) there amounts to R16--19, (b) to R32, (c) to R6--9, (d) to R28, (e) to R29, (f) to R31.
\\[.3em]
(ii)~The relations given in \cite[Prop.\,3.10]{LP} or \cite[Sect.\,2]{MS} for a two-dimensional topological field theory are different from R1--32. This illustrates that even though we work with topological world sheets, the problem of finding a solution of the sewing constraints is generally not the same as giving a two-dimensional topological field theory. Indeed, it should not be, because we want to describe {\em conformal} field theories.
For example, as remarked in \cite[Fig.\,12]{Lz}, in a 2-d TFT relation R32 is trivial as the corresponding cobordisms are equivalent.
\\[.3em]
(iii)~Apart from \cite{Le}, the closest analogue to R1--R32 can be found in \cite{KP}. There a list of combinatorial moves on generalised pants decompositions is given, and it is proved that these act transitively on the set of generalised pants decompositions for a given surface \cite[Thm.\,3.1]{KP}.
}
\end{rema}

\begin{figure}
$$
\begin{array}{llll}
  \Mr_{m o}(f) &=
  \raisebox{-65pt}{\begin{picture}(140,130)
  \put(7,10){\scalebox{0.40}{\includegraphics{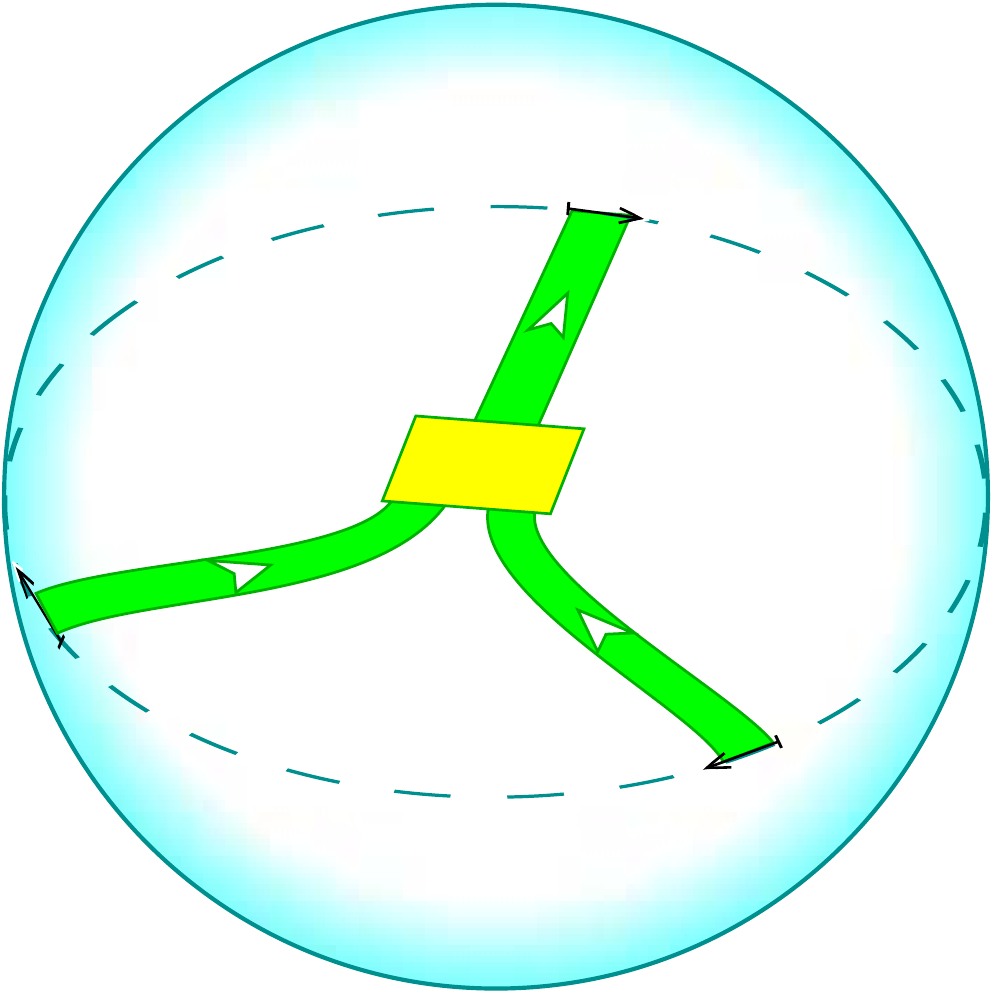}}}
  \put(7,10){
     \setlength{\unitlength}{.40pt}\put(-135,-279){
     \put(265,428)   {\tiny$ f $}
     \put(304,462)   {\tiny$ \Bop $}
     \put(180,414)   {\tiny$ \Bop $}
     \put(314,392)   {\tiny$ \Bop $}
     \put(280,516)   {\tiny$ (\Bop,-) $}
     \put(110,360)   {\tiny$ (\Bop,+) $}
     \put(310,328)   {\tiny$ (\Bop,+) $}
     }\setlength{\unitlength}{1pt}}
  \end{picture}}
  &\Mr_{\eta o}(f) &=
  \raisebox{-65pt}{\begin{picture}(140,130)
  \put(7,10){\scalebox{0.40}{\includegraphics{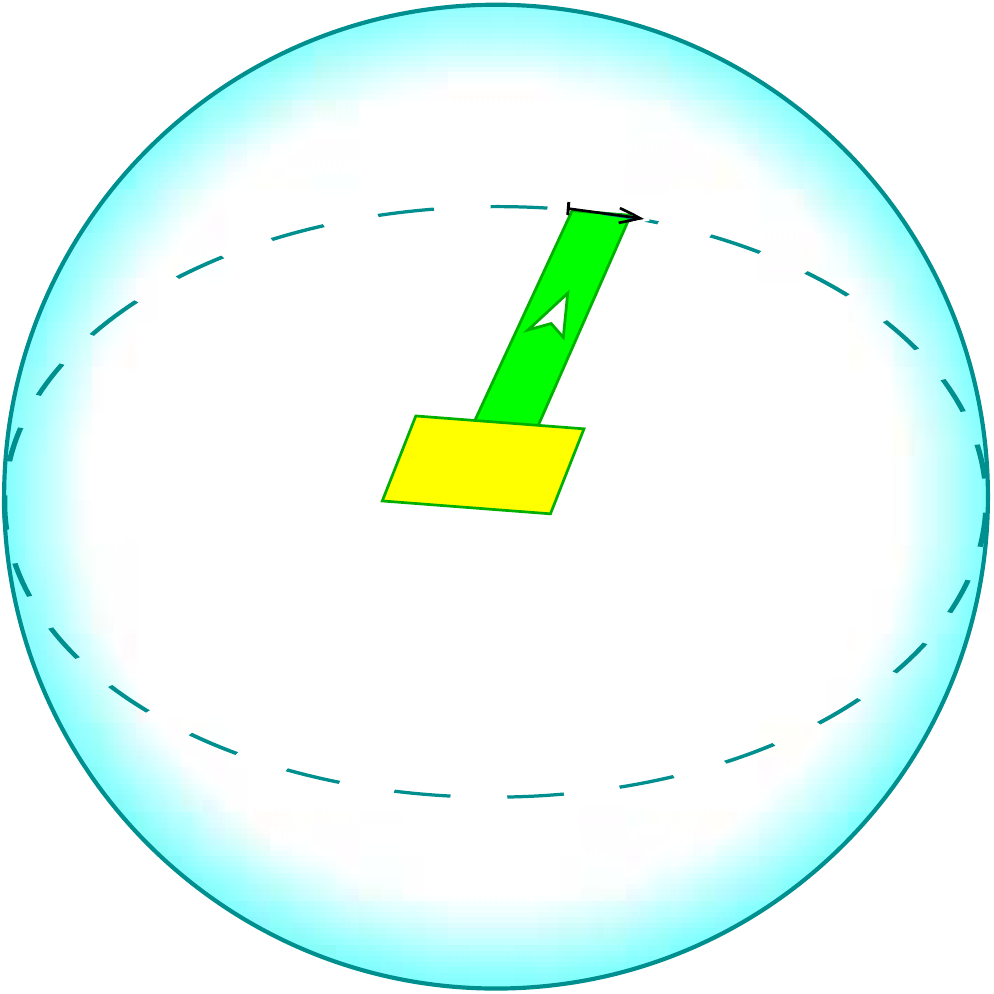}}}
  \put(7,10){
     \setlength{\unitlength}{.40pt}\put(-135,-279){
     \put(265,428)   {\tiny$ f $}
     \put(304,462)   {\tiny$ \Bop $}
     }\setlength{\unitlength}{1pt}}
  \end{picture}}
\\
  \Mr_{\Delta o}(f) &=
  \raisebox{-65pt}{\begin{picture}(140,130)
  \put(7,10){\scalebox{0.40}{\includegraphics{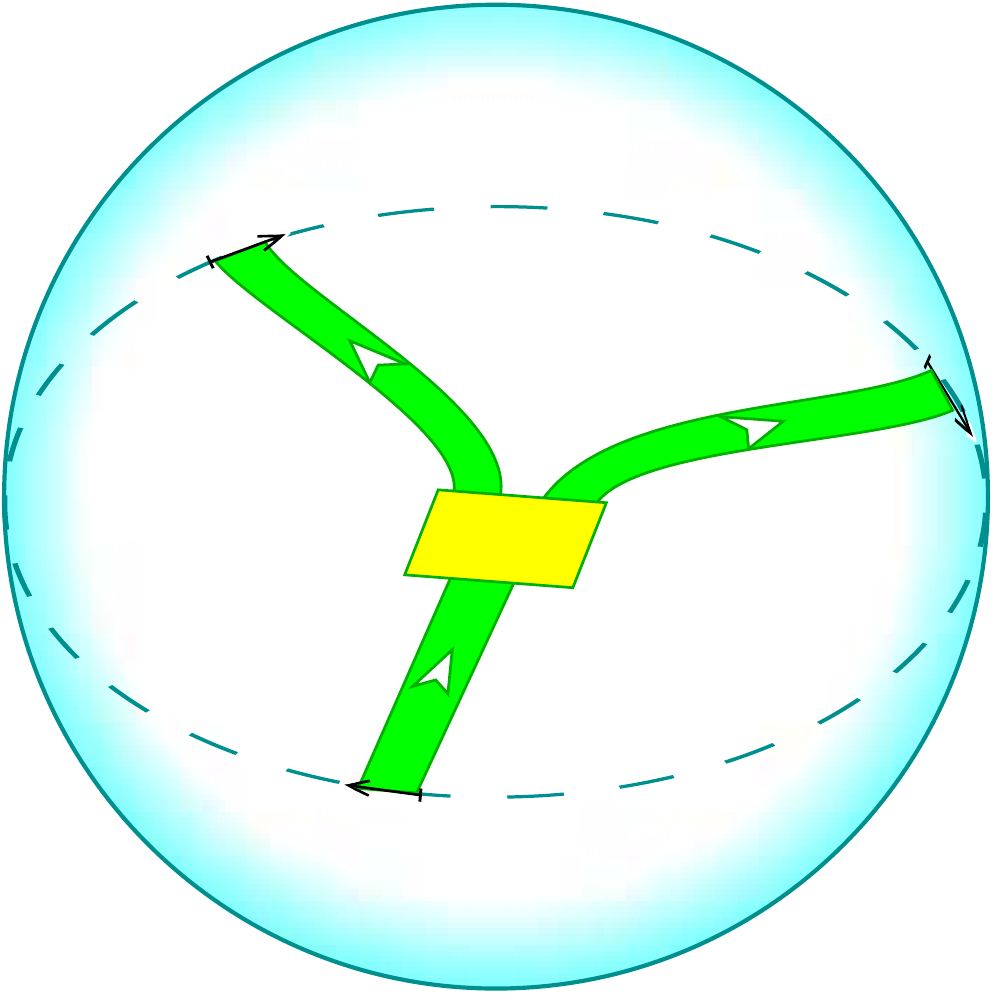}}}
  \put(7,10){
     \setlength{\unitlength}{.40pt}\put(-135,-279){
     \put(273,406)   {\tiny$ f $}
     \put(211,449)   {\tiny$ \Bop $}
     \put(334,457)   {\tiny$ \Bop $}
     \put(272,359)   {\tiny$ \Bop $}
     }\setlength{\unitlength}{1pt}}
  \end{picture}}
  &\Mr_{\eps o}(f) &=
  \raisebox{-65pt}{\begin{picture}(140,130)
  \put(7,10){\scalebox{0.40}{\includegraphics{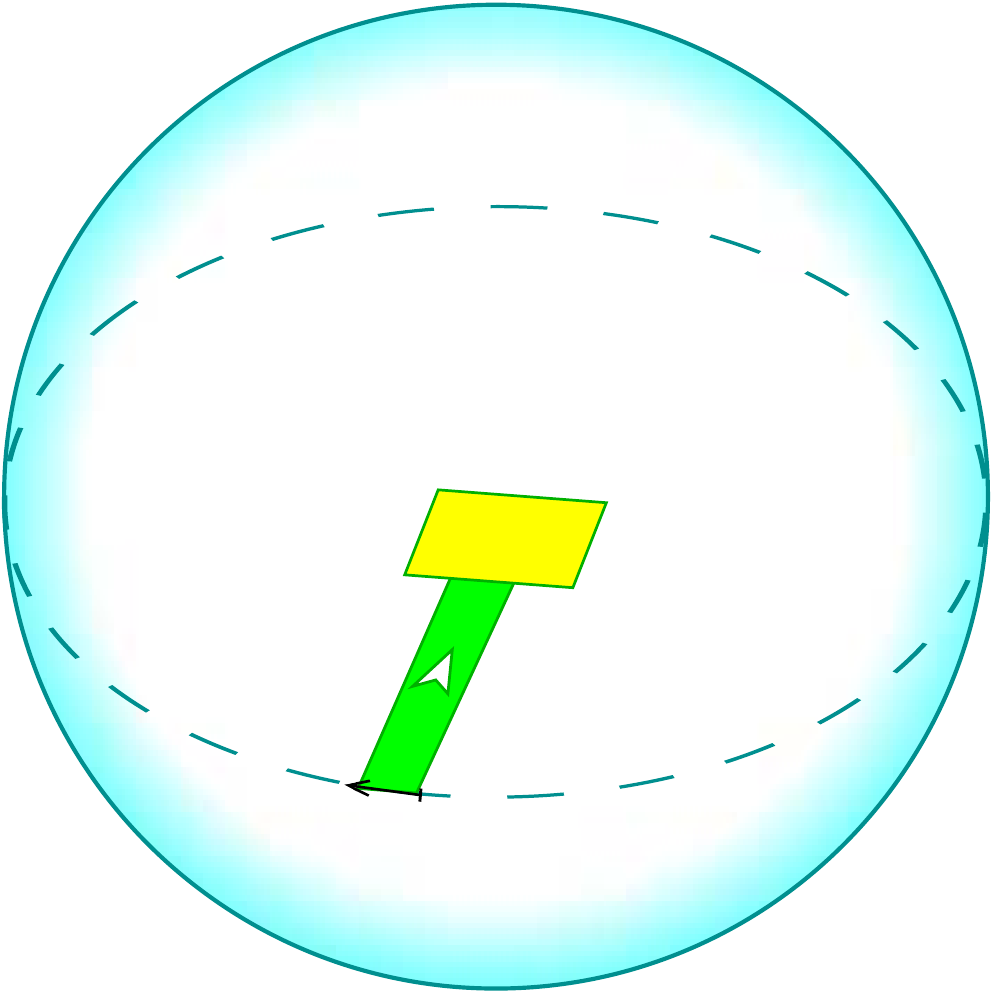}}}
  \put(7,10){
     \setlength{\unitlength}{.40pt}\put(-135,-279){
     \put(273,406)   {\tiny$ f $}
     \put(272,359)   {\tiny$ \Bop $}
     }\setlength{\unitlength}{1pt}}
  \end{picture}}
\\
  \Mr_{\iota}(f) &=
  \raisebox{-65pt}{\begin{picture}(140,130)
  \put(7,10){\scalebox{0.40}{\includegraphics{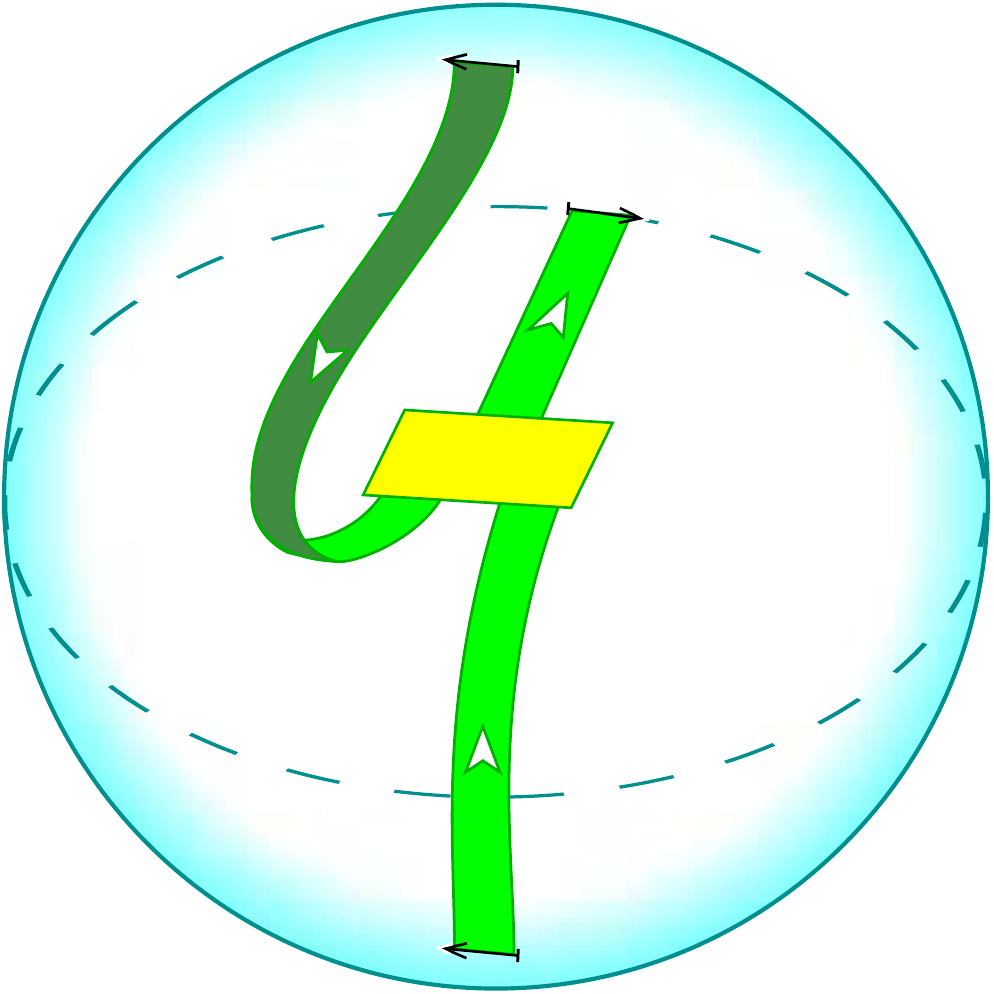}}}
  \put(7,10){
     \setlength{\unitlength}{.40pt}\put(-135,-279){
     \put(273,428)   {\tiny$ f $}
     \put(310,465)   {\tiny$ \Bop $}
     \put(190,454)   {\tiny$ B_l $}
     \put(290,356)   {\tiny$ B_r $}
     }\setlength{\unitlength}{1pt}}
  \end{picture}}
  &\Mr_{\iota^*}(f) &=
  \raisebox{-65pt}{\begin{picture}(140,130)
  \put(7,10){\scalebox{0.40}{\includegraphics{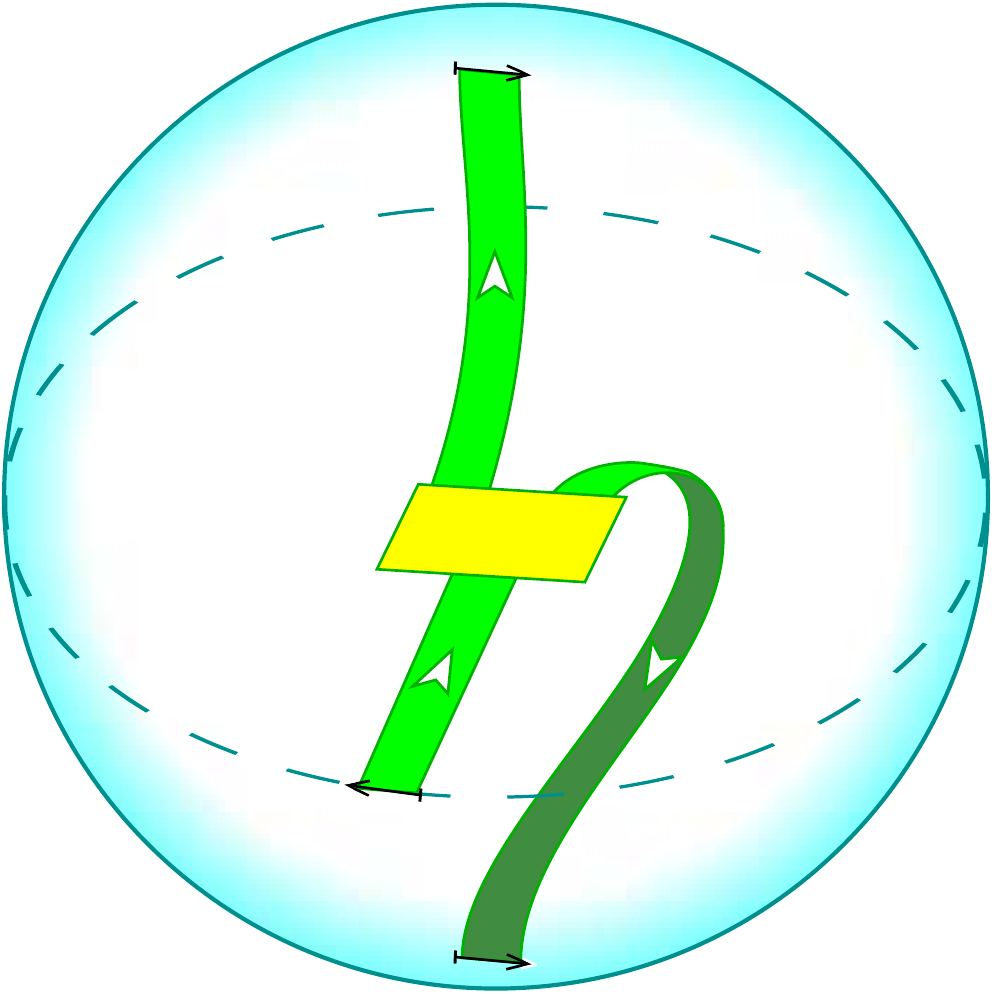}}}
  \put(7,10){
     \setlength{\unitlength}{.40pt}\put(-135,-279){
     \put(274,408)   {\tiny$ f $}
     \put(215,365)   {\tiny$ \Bop $}
     \put(248,478)   {\tiny$ B_l $}
     \put(338,368)   {\tiny$ B_r $}
     }\setlength{\unitlength}{1pt}}
  \end{picture}}
\end{array}
$$
\caption{The cobordisms $\Mr_\alpha : \emptyset  \rightarrow \Xhat_\alpha$
for $\alpha \In \{mo,\eta o,\Delta o,\eps o,\iota,\iota^*\}$.
In $\Mr_{m o}(f)$ we have also included
the labels of the arcs. The solid spheres inherit their orientation from the
embedding in $\Rb^3$, and the ribbons --
with the exception of the darker coloured ribbons in $\Mr_{\iota}$ and
$\Mr_{\iota^*}$ -- are oriented so that their `white' side
faces readers, i.e.\ the ribbon orientation agrees with that of the paper
plane.}
\label{fig:Malpha-cobord}
\end{figure}

Let us for the rest of this paper fix a set of data
$(\Cc,\Bop,B_l,B_r)$. We will abbreviate
\be
  \Bl \equiv \Bl(\Cc,\Bop,B_l,B_r)
  \quad \text{and} \quad
  \Bcl = \BlxBr~.
\ee
The following construction is
basically the one given in \cite[Sect.\,4.1]{unique}. Define the
extended cobordisms
\be
  \Mr_\alpha(\,\cdot\,) :
  \emptyset \rightarrow \Xhat_\alpha
  \qquad \text{for} \quad \alpha \In \Sc
\labl{eq:Malpha-def}
as in Figures \ref{fig:Malpha-cobord}--\ref{fig:Malpha-cobord-3}.
With the help of these we can define linear maps $\Psi_\alpha$ between the morphism
spaces of $\Cc$ and the state spaces of the 3-d TFT.
For generators $\Xr_\alpha$ involving only open state boundaries, we set
\be \label{eq:def-Psi}
\Psi_\alpha(f) = \tftC(\Mr_\alpha(f))1,
\ee
where $1\in \tftC(\emptyset) = \Cb$. The assignment $\Psi_\alpha$
is linear in $f$ because $\tftC$ is multilinear in the labels
of the coupons embedded in the extended cobordisms.
The source and target vector spaces of $\Psi_\alpha$ are explicitly,
\be\begin{array}{llll}\displaystyle
\Psi_{po} \etb: \Hom_\Cc(\Bop,\Bop)
    \longrightarrow \Bl(\Xr_{po}) ~,
\enl
\Psi_{mo} \etb: \Hom_\Cc(\Bop\oti\Bop,\Bop)
    \longrightarrow \Bl(\Xr_{mo}) ~,   \quad
\etb
\Psi_{\eta o} \etb: \Hom_\Cc(\one,\Bop)
    \longrightarrow \Bl(\Xr_{\eta o}) ~,
\enl
\Psi_{\Delta o} \etb: \Hom_\Cc(\Bop,\Bop\oti\Bop)
    \longrightarrow \Bl(\Xr_{\Delta o}) ~,
\etb
\Psi_{\eps o} \etb: \Hom_\Cc(\Bop,\one)
    \longrightarrow \Bl(\Xr_{\eps o}) ~.
\eear\labl{eq:iso-op}
In the same way, the linear maps involving open and closed state boundaries are
\be\begin{array}{llll}\displaystyle
\Psi_{\iota} \etb: \Hom_\Cc(\Bop,B_l\oti B_r)
    \longrightarrow \Bl(\Xr_{\iota}) ~, \quad
\etb
\Psi_{\iota^*} \etb: \Hom_\Cc(B_l\oti B_r,\Bop)
    \longrightarrow \Bl(\Xr_{\iota^*}) ~.
\eear\labl{eq:iso-opcl}
Notice that we choose to put $B_l$ left to $B_r$ according to the second condition in the definition of the order map `$\ord$' in Definition \ref{def:ws}. In order to define the $\Psi_\alpha$ in the case that $\Xr_\alpha$ has only closed state boundaries, we need an intermediate step. Consider $\Psi_{mc}$ as an example. By the definition of $\CxC$, we have
\be
  \Hom_\CxC(\Bcl\oti\Bcl,\Bcl)
  = \Hom_\Cc(B_l \otimes B_l,B_l) \otimes_\Cb
    \Hom_\Cc(B_r \otimes B_r,B_r) ~.
\ee
On elements $f' \otimes_\Cb f''$ of the above space, $\Psi_{mc}$ is defined by
$\Psi_{mc}(f' \otimes f'') = \tftC(\Mr_{mc}(f',f''))1$, 
and on a general element by linear extension. In this way we obtain
\be\begin{array}{llll}\displaystyle
\Psi_{pc} \etb: \Hom_\CxC(\Bcl,\Bcl)
    \longrightarrow \Bl(\Xr_{pc}),
\enl
\Psi_{mc} \etb: \Hom_\CxC(\Bcl\oti\Bcl,\Bcl)
    \longrightarrow \Bl(\Xr_{mo}), \quad
\etb
\Psi_{\eta c} \etb: \Hom_\CxC(\one,\Bcl)
    \longrightarrow \Bl(\Xr_{\eta c}),
\enl
\Psi_{\Delta c} \etb: \Hom_\CxC(\Bcl,\Bcl\oti\Bcl)
    \longrightarrow \Bl(\Xr_{\Delta c}),
\etb
\Psi_{\eps c} \etb: \Hom_\CxC(\Bcl,\one)
    \longrightarrow \Bl(\Xr_{\eps c}).
\eear\labl{eq:iso-cl}
By the construction of the functor $\tftC$ and our choices for $\Mr_\alpha$, the linear maps $\Psi_\alpha$
in \eqref{eq:iso-op}, \eqref{eq:iso-opcl}, and \eqref{eq:iso-cl} are in fact
isomorphisms \cite[Sect.\,IV.2.1]{turaev-bk} (see also \cite[Sect.\,3.2,\,4.1]{unique}).

\begin{figure}
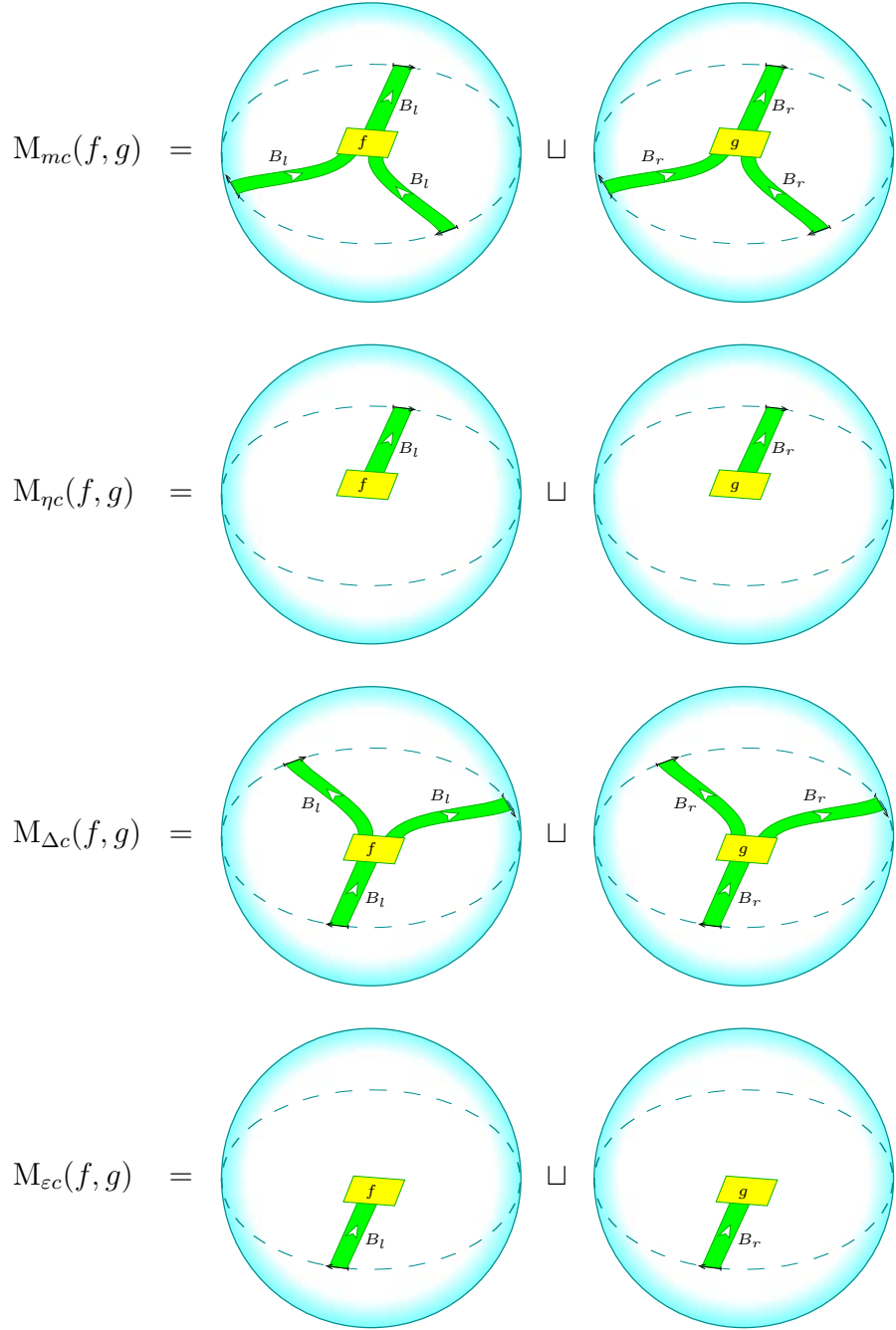

$$
\begin{array}{ll}
  \Mr_{mc}(f,g) &=
  \raisebox{-65pt}{\begin{picture}(140,130)
  \put(7,10){\scalebox{0.40}{\includegraphics{pic-M-m_o}}}
  \put(7,10){
     \setlength{\unitlength}{.40pt}\put(-135,-279){
     \put(265,428)   {\tiny$ f $}
     \put(304,462)   {\tiny$ B_l $}
     \put(180,414)   {\tiny$ B_l $}
     \put(314,392)   {\tiny$ B_l $}
     }\setlength{\unitlength}{1pt}}
  \end{picture}}
  \hspace{-1em}\sqcup
  \raisebox{-65pt}{\begin{picture}(140,130)
  \put(7,10){\scalebox{0.40}{\includegraphics{pic-M-m_o}}}
  \put(7,10){
     \setlength{\unitlength}{.40pt}\put(-135,-279){
     \put(265,428)   {\tiny$ g $}
     \put(304,462)   {\tiny$ B_r $}
     \put(180,414)   {\tiny$ B_r $}
     \put(314,392)   {\tiny$ B_r $}
     }\setlength{\unitlength}{1pt}}
  \end{picture}}
  \\
  \Mr_{\eta c}(f,g) &=
  \raisebox{-65pt}{\begin{picture}(140,130)
  \put(7,10){\scalebox{0.40}{\includegraphics{pic-M-eta_o}}}
  \put(7,10){
     \setlength{\unitlength}{.40pt}\put(-135,-279){
     \put(265,428)   {\tiny$ f $}
     \put(304,462)   {\tiny$ B_l $}
     }\setlength{\unitlength}{1pt}}
  \end{picture}}
  \hspace{-1em}\sqcup
  \raisebox{-65pt}{\begin{picture}(140,130)
  \put(7,10){\scalebox{0.40}{\includegraphics{pic-M-eta_o}}}
  \put(7,10){
     \setlength{\unitlength}{.40pt}\put(-135,-279){
     \put(265,428)   {\tiny$ g $}
     \put(304,462)   {\tiny$ B_r $}
     }\setlength{\unitlength}{1pt}}
  \end{picture}}
  \\
  \Mr_{\Delta c}(f,g) &=
  \raisebox{-65pt}{\begin{picture}(140,130)
  \put(7,10){\scalebox{0.40}{\includegraphics{pic-M-Delta_o}}}
  \put(7,10){
     \setlength{\unitlength}{.40pt}\put(-135,-279){
     \put(273,406)   {\tiny$ f $}
     \put(211,449)   {\tiny$ B_l $}
     \put(334,457)   {\tiny$ B_l $}
     \put(272,359)   {\tiny$ B_l $}
     }\setlength{\unitlength}{1pt}}
  \end{picture}}
  \hspace{-1em}\sqcup
  \raisebox{-65pt}{\begin{picture}(140,130)
  \put(7,10){\scalebox{0.40}{\includegraphics{pic-M-Delta_o}}}
  \put(7,10){
     \setlength{\unitlength}{.40pt}\put(-135,-279){
     \put(273,406)   {\tiny$ g $}
     \put(211,449)   {\tiny$ B_r $}
     \put(334,457)   {\tiny$ B_r $}
     \put(272,359)   {\tiny$ B_r $}
     }\setlength{\unitlength}{1pt}}
  \end{picture}}
  \\
  \Mr_{\eps c}(f,g) &=
  \raisebox{-65pt}{\begin{picture}(140,130)
  \put(7,10){\scalebox{0.40}{\includegraphics{pic-M-eps_o}}}
  \put(7,10){
     \setlength{\unitlength}{.40pt}\put(-135,-279){
     \put(273,406)   {\tiny$ f $}
     \put(272,359)   {\tiny$ B_l $}
     }\setlength{\unitlength}{1pt}}
  \end{picture}}
  \hspace{-1em}\sqcup
  \raisebox{-65pt}{\begin{picture}(140,130)
  \put(7,10){\scalebox{0.40}{\includegraphics{pic-M-eps_o}}}
  \put(7,10){
     \setlength{\unitlength}{.40pt}\put(-135,-279){
     \put(273,406)   {\tiny$ g $}
     \put(272,359)   {\tiny$ B_r $}
     }\setlength{\unitlength}{1pt}}
  \end{picture}}
\end{array}
$$
\caption{The cobordisms $\Mr_\alpha : \emptyset  \rightarrow \Xhat_\alpha$
for $\alpha \In \{mc,\eta c,\Delta c,\eps c,pc\}$.}
\label{fig:Malpha-cobord-2}
\end{figure}

\begin{figure}
$$
\begin{array}{ll}
  \Mr_{p o}(f) &=
  \raisebox{-65pt}{\begin{picture}(140,130)
  \put(7,10){\scalebox{0.40}{\includegraphics{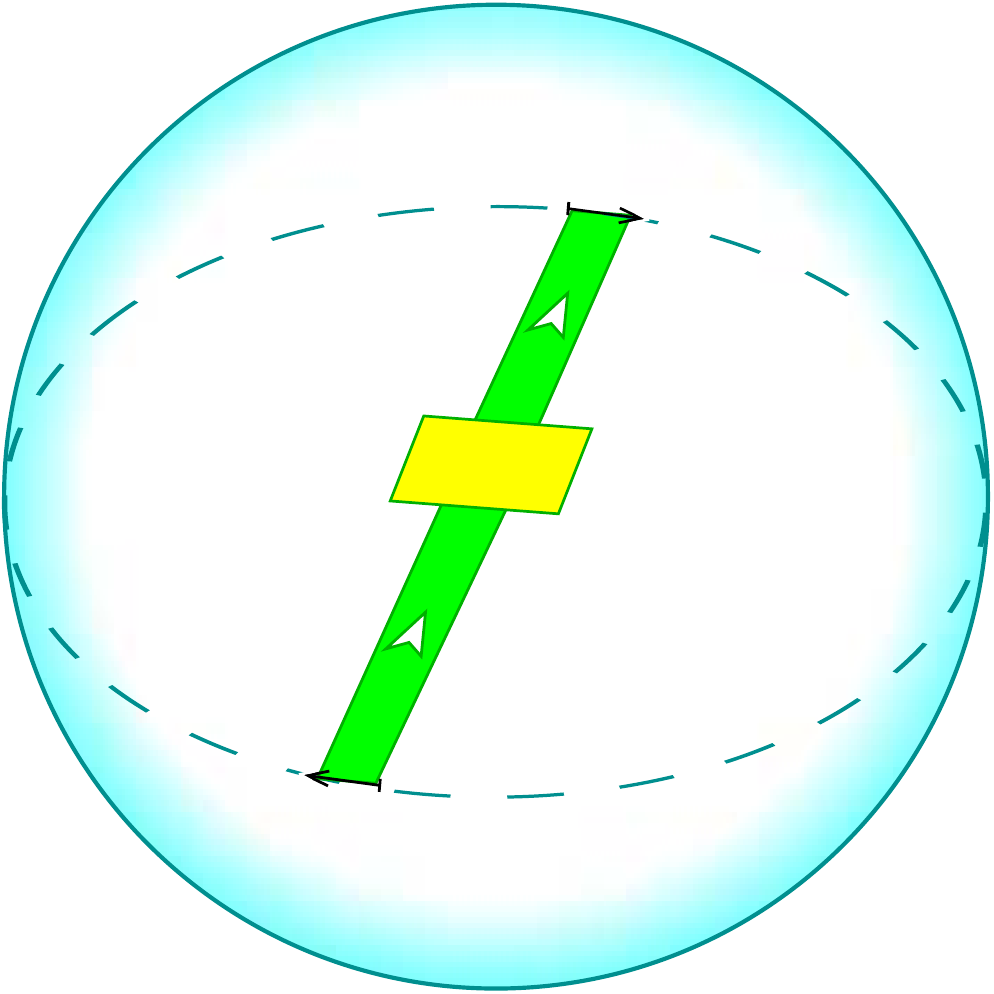}}}
  \put(7,10){
     \setlength{\unitlength}{.40pt}\put(-135,-279){
     \put(268,428)   {\tiny$ f $}
     \put(304,462)   {\tiny$ \Bop $}
     \put(265,370)   {\tiny$ \Bop $}
     }\setlength{\unitlength}{1pt}}
  \end{picture}}
\\
  \Mr_{p c}(f,g) &=
  \raisebox{-65pt}{\begin{picture}(140,130)
  \put(7,10){\scalebox{0.40}{\includegraphics{pic-M-p_o}}}
  \put(7,10){
     \setlength{\unitlength}{.40pt}\put(-135,-279){
     \put(268,428)   {\tiny$ f $}
     \put(304,462)   {\tiny$ B_l $}
     \put(265,370)   {\tiny$ B_l $}
     }\setlength{\unitlength}{1pt}}
  \end{picture}}
  \hspace{-1em}\sqcup
  \raisebox{-65pt}{\begin{picture}(140,130)
  \put(7,10){\scalebox{0.40}{\includegraphics{pic-M-p_o}}}
  \put(7,10){
     \setlength{\unitlength}{.40pt}\put(-135,-279){
     \put(268,428)   {\tiny$ g $}
     \put(304,462)   {\tiny$ B_r $}
     \put(265,370)   {\tiny$ B_r $}
     }\setlength{\unitlength}{1pt}}
  \end{picture}}
\end{array}
$$
\caption{The cobordisms $\Mr_{po} : \emptyset  \rightarrow \Xhat_{po}$
and $\Mr_{pc} : \emptyset  \rightarrow \Xhat_{pc}$.}
\label{fig:Malpha-cobord-3}
\end{figure}

\medskip

After briefly reviewing the definition of Cardy algebras (see \cite[Def.\,5.14]{cardy} and \cite[Sect.\,2.2\,\&\,3.2]{part1}), we can state the relation between a solution to the sewing constraints and a Cardy algebra in Theorems \ref{thm:cardy-determine} and \ref{thm:cardy-sew}. Their proofs will be given in Section \ref{sec:proof-sew} below.

\medskip

We will make use of the graphical notation of morphisms in a modular tensor category, and follow the conventions of \cite[app.\,A.1]{tft5}.

An {\em algebra} in $\Cc$ is a triple $A=(A,m,\eta)$ where $A$ is an object of $\Cc$, $m$ (the multiplication) is a morphism $A \oti A \rightarrow A$ such that $m \cir (m \otimes_\Cb \id_A) = m \cir (\id_A \otimes_\Cb m)$, and $\eta$ (the unit) is a morphism $\one \rightarrow A$ such that $m \cir (\id_A \otimes_\Cb \eta) = \id_A$ and $m \cir (\eta\otimes_\Cb \id_A) = \id_A$. If $m \cir c_{A,A} = m$, then $A$ is called {\em commutative}.
A {\em coalgebra} $A = (A,\Delta,\eps)$ in $\Cc$ is defined analogously to an algebra in $\Cc$, i.e.\ $\Delta : A \rightarrow A \oti A$ and $\eps : A \rightarrow \one$ obey coassociativity and counit conditions.

A {\em Frobenius algebra} $A = (A,m,\eta,\Delta,\eps)$ is an algebra and a coalgebra such that the coproduct is an intertwiner
 of $A$-bimodules,
\be
  (\id_A \otimes_\Cb m) \cir (\Delta \otimes_\Cb \id_A) = \Delta \circ m = (m \otimes_\Cb \id_A) \cir (\id_A \otimes_\Cb \Delta)~.
\ee
We will use the following graphical representation for the morphisms of a Frobenius algebra,
\be
  m = \raisebox{-20pt}{
  \begin{picture}(30,45)
   \put(0,6){\scalebox{.75}{\includegraphics{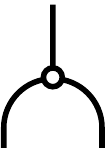}}}
   \put(0,6){
     \setlength{\unitlength}{.75pt}\put(-146,-155){
     \put(143,145)  {\scriptsize $ A $}
     \put(169,145)  {\scriptsize $ A $}
     \put(157,202)  {\scriptsize $ A $}
     }\setlength{\unitlength}{1pt}}
  \end{picture}}
  ~~,\quad
  \eta = \raisebox{-15pt}{
  \begin{picture}(10,30)
   \put(0,6){\scalebox{.75}{\includegraphics{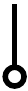}}}
   \put(0,6){
     \setlength{\unitlength}{.75pt}\put(-146,-155){
     \put(146,185)  {\scriptsize $ A $}
     }\setlength{\unitlength}{1pt}}
  \end{picture}}
  ~~,\quad
  \Delta = \raisebox{-20pt}{
  \begin{picture}(30,45)
   \put(0,6){\scalebox{.75}{\includegraphics{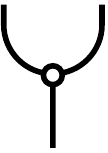}}}
   \put(0,6){
     \setlength{\unitlength}{.75pt}\put(-146,-155){
     \put(143,202)  {\scriptsize $ A $}
     \put(169,202)  {\scriptsize $ A $}
     \put(157,145)  {\scriptsize $ A $}
     }\setlength{\unitlength}{1pt}}
  \end{picture}}
  ~~,\quad
  \eps = \raisebox{-15pt}{
  \begin{picture}(10,30)
   \put(0,10){\scalebox{.75}{\includegraphics{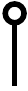}}}
   \put(0,10){
     \setlength{\unitlength}{.75pt}\put(-146,-155){
     \put(146,145)  {\scriptsize $ A $}
     }\setlength{\unitlength}{1pt}}
  \end{picture}}
  ~~.
\ee
Given two Frobenius algebra $(A, m_A, \eta_A, \Delta_A, \eps_A)$ and $(B, m_B, \eta_B, \Delta_B, \eps_B)$
in $\mathcal{C}$. For $f: A\rightarrow B$, we define $f^*: B\rightarrow A$ by
\be  \label{f-star-def}
f^* =  ( (\eps_B \circ m_B) \otimes_\Cb \id_A)
\circ (\id_B \otimes_\Cb f \otimes_\Cb \id_A)
\circ (\id_B \otimes_\Cb (\Delta_A \circ \eta_A)).
\ee
Given a Frobenius algebra in $\Cc$ we can define two isomorphisms $A \rightarrow A^\vee$,
\be
\Phi_A =\,\,
\raisebox{-35pt}{
  \begin{picture}(50,75)
   \put(0,8){\scalebox{.75}{\includegraphics{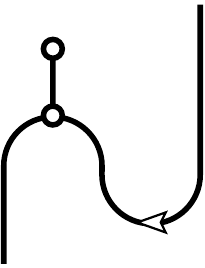}}}
   \put(0,8){
     \setlength{\unitlength}{.75pt}\put(-34,-37){
     \put(31, 28)  {\scriptsize $ A $}
     \put(87,117)  {\scriptsize $ A^\vee $}
     }\setlength{\unitlength}{1pt}}
  \end{picture}},  \quad\quad\quad
\Phi_A' =\,\,
  \raisebox{-35pt}{
  \begin{picture}(50,75)
   \put(0,8){\scalebox{.75}{\includegraphics{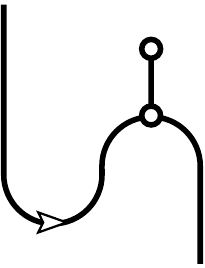}}}
   \put(0,8){
     \setlength{\unitlength}{.75pt}\put(-34,-37){
     \put(87, 28)  {\scriptsize $ A $}
     \put(31,117)  {\scriptsize $ A^\vee $}
     }\setlength{\unitlength}{1pt}}
  \end{picture}}
  ~~.
\labl{eq:Frob-sym-cond}
A Frobenius algebra is called {\em symmetric} if $\Phi_A = \Phi_A'$.

\bigskip

The tensor product of $\Cc$ defines a functor $T : \CxC \rightarrow \Cc$ , i.e. $X \times Y \mapsto X\otimes Y$. Since $\Cc$ is braided, $T$ can be turned into a tensor functor, see \cite[Prop.\,5.2]{js} and \cite[Lem.\,3.5]{op-cl-alg}. More precisely,
For $X^L, X^R, Y^L, Y^R\in \Cc$, we have
$$T(X^L \times X^R) \otimes T(Y^L \times Y^R) = X^L \otimes X^R \otimes Y^L \otimes Y^R$$
$$T((X^L \times X^R) \otimes (Y^L \times Y^R)) = X^L \otimes Y^L \otimes X^R \otimes Y^R$$
We define the natural isomorphism $\phi_2: \otimes \circ (T \times T) \to T\circ \otimes$ by
\be  \label{eq:T-phi-2}
\phi_2 := \id_{X^L} \otimes_\Cb c_{Y^L, X^R}^{-1} \otimes_\Cb \id_{Y^R}
\ee
where $c_{X^R, Y^L}: X^R \otimes Y^L \to Y^L \otimes X^R$ is the braiding isomorphism and all necessary associators in the definition of $\phi_2$ are implicitly assumed.
We also define $\phi_0: \one \to \one\otimes \one$ by $l_\one^{-1}$ where $l_\one: \one\otimes \one \to \one$ is the left unit isomorphism. It is easy to see that $T$ together with $\phi_2$ and $\phi_0$ gives a monoidal functor $\CxC \rightarrow \Cc$. As a consequence, for any Frobenius algebra $A$ in $\CxC$, $T(A)$ is automatically a Frobenius algebra in $\Cc$ \cite[Prop.\,2.13]{part1}.

\begin{defn}  \label{def:cardy-alg}  {\rm
A Cardy $\Cc|\CxC$-algebra is a triple
$(\Aop|\Acl, \ico)$, where
\begin{enumerate}
\item $(\Aop, m_\text{\rm{op}}, \eta_\text{\rm{op}}, \Delta_\text{\rm{op}}, \eps_\text{\rm{op}})$ is
a symmetric Frobenius algebra in $\Cc$,
\item $(\Acl, m_\text{\rm{cl}}, \eta_\text{\rm{cl}}, \Delta_\text{\rm{cl}}, \eps_\text{\rm{cl}})$ is
a commutative symmetric Frobenius algebra in $\CxC$. We also assume that $\Acl=\oplus_{n=1}^N C_n^l \times C_n^r$ for $C_n^l, C_n^R \in \Cc$.
\item $\ico: T(\Acl) \rightarrow \Aop$ an algebra homomorphism. For its restriction on direct summand $C_n^l \otimes C_n^r$ in $T(\Acl)$, we introduce the following graphic notation:
\be   \label{iota-*-def}
\ico^{\, (n)} ~=~
  \raisebox{-30pt}{
  \begin{picture}(54,73)
   \put(0,8){\scalebox{.75}{\includegraphics{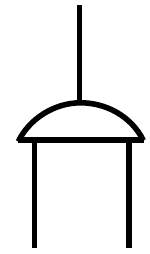}}}
   \put(0,8){
     \setlength{\unitlength}{.75pt}\put(-18,-19){
     \put(19, 8)  {\scriptsize $ C_n^l $}
     \put(51, 8)  {\scriptsize $ C_n^r $}
     \put(35, 98)  {\scriptsize $ \Aop $}
     }\setlength{\unitlength}{1pt}}
  \end{picture}}
~~.
\ee
\end{enumerate}
such that the following conditions are satisfied:
\\[.3em]
(i) Centre condition:
\be   \label{eq:comm-C}
\setlength{\unitlength}{1pt}
\raisebox{-45pt}{
  \begin{picture}(60,90)
   \put(0,8){\scalebox{.75}{\includegraphics{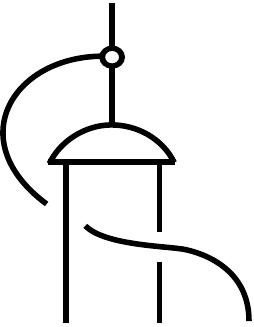}}}
   \put(0,8){
     \setlength{\unitlength}{.75pt}\put(-18,-11){
     \put(33, 2)     {\scriptsize $ C_n^l $}
     \put(63, 2)     {\scriptsize $ C_n^r $}
     \put(45, 110) {\scriptsize $\Aop$}
     \put(85, 2)     {\scriptsize $ \Aop$}
     \put(58, 78)   {\scriptsize $ \Aop$}
     }\setlength{\unitlength}{1pt}}
  \end{picture}}
~=~
\raisebox{-45pt}{
  \begin{picture}(35,90)
  \put(0,8){\scalebox{.75}{\includegraphics{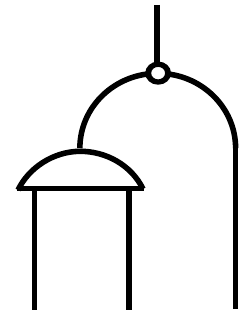}}}
   \put(0,8){
     \setlength{\unitlength}{.75pt}\put(-18,-11){
     \put(23,2)       {\scriptsize $ C_n^l $}
     \put(50,2)       {\scriptsize $ C_n^r $}
     \put(58,110)   {\scriptsize $ \Aop $}
     \put(80,2)       {\scriptsize $ \Aop $}
     \put(23,75)     {\scriptsize $ \Aop $}
  }\setlength{\unitlength}{1pt}}
  \end{picture}}
\qquad \qquad
\text{for}~n=1, \dots, N~~.
\ee
(ii) Modular invariance: for all $i\in \mathcal{I}$, (recall the graphic notations introduced in \eqref{eq:b-dual-b})
\be    \label{eq:mod-inv}
  \frac{\dim U_i \dim U_j}{\mathrm{Dim} \, \Cc} ~~
\raisebox{-55pt}{
  \begin{picture}(110, 140)
   \put(0,8){\scalebox{.75}{\includegraphics{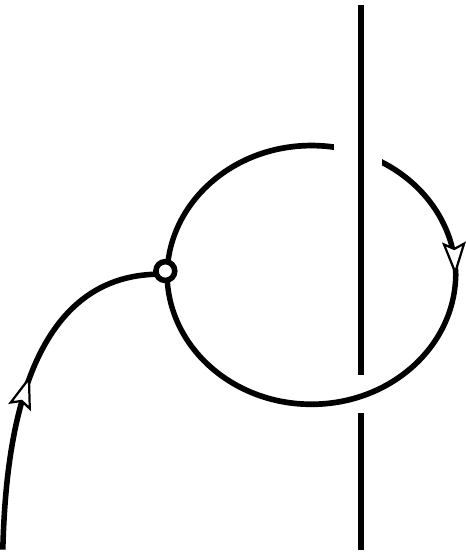}}}
   \put(0,8){
     \setlength{\unitlength}{.75pt}\put(-18,-19){
     \put(14, 8)     {\scriptsize $ \Acl $}
     \put(156,100)  {\scriptsize $ \Acl $}
     \put(118,8)    {\scriptsize $U_i \times U_j$ }
     \put(118,185){\scriptsize $U_i \times U_j $ }
     \put(73, 100)  {\scriptsize $m_{\mathrm{cl}}$}
     }\setlength{\unitlength}{1pt}}
  \end{picture}}
\quad = \quad  \, \sum_{\alpha}
 \raisebox{-55pt}{
  \begin{picture}(90,140)
   \put(8,8){\scalebox{.75}{\includegraphics{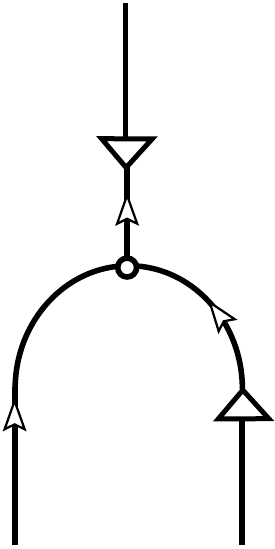}}}
   \put(0,8){
     \setlength{\unitlength}{.75pt}\put(-18,-19){
     \put(28, 8)  {\scriptsize $ \Acl $}
     \put(60, 185)  {\scriptsize $ U_i \times U_j$}
     \put(94, 8)  {\scriptsize $ U_i \times U_j $}
     \put(45,115)  {\scriptsize $ \Acl $}
     \put(98,89)    {\scriptsize $ \Acl $}
     \put(78, 132)   {\scriptsize $\alpha$}
     \put(110, 60)     {\scriptsize $\alpha$}
     \put(60, 89)  {\scriptsize $m_{\mathrm{cl}}$}
}\setlength{\unitlength}{1pt}}
  \end{picture}}
\ee
(iii) Cardy condition:
\be   \label{eq:cardy-C}
\bigoplus_{n=1}^N
      \, \sum_\alpha
\quad
\raisebox{-70pt}{
  \begin{picture}(70,130)
   \put(0,8){\scalebox{.75}{\includegraphics{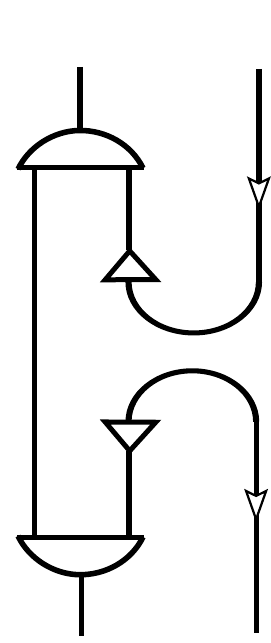}}}
   \put(0,8){
     \setlength{\unitlength}{.75pt}\put(-18,-11){
     \put(9, 95)     {\scriptsize $ C_n^l $}
     \put(60, 54)     {\scriptsize $ C_n^r $}
     \put(60, 134)   {\scriptsize $ C_n^r $}
     \put(40, 118)   {\scriptsize $ \alpha $}
     \put(40, 70)   {\scriptsize $ \alpha $}
     \put(35, 2)       {\scriptsize $\Aop$}
     \put(35, 185)   {\scriptsize $\Aop$}
     \put(85, 185)   {\scriptsize $ U_i^{\vee} $}
     \put(85, 2)       {\scriptsize $ U_i^{\vee} $}
     }\setlength{\unitlength}{1pt}}
  \end{picture}}
=~ \frac{\dim U_i}{\sqrt{\Dim \, \Cc} } \,\,
\raisebox{-55pt}{
  \begin{picture}(85,130)
  \put(0,8){\scalebox{.75}{\includegraphics{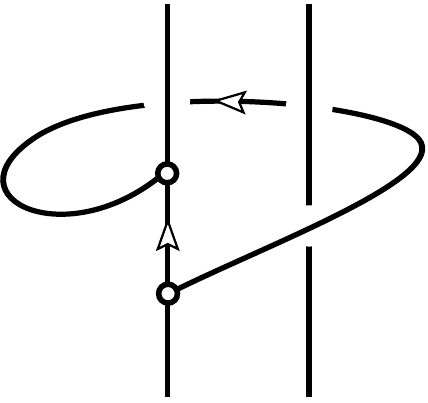}}}
   \put(0,8){
     \setlength{\unitlength}{.75pt}\put(-18,-11){
     \put(140,90)       {\scriptsize $ \Aop $}
     \put(60,2)           {\scriptsize $ \Aop $}
     \put(60,135)       {\scriptsize $ \Aop $}
     \put(100, 2)        {\scriptsize $ U_i^{\vee} $}
     \put(100, 135)    {\scriptsize $ U_i^{\vee} $}
     \put(40, 42)      {\scriptsize $\Delta_\text{\rm{op}}$ }
     \put(72, 77)     {\scriptsize $m_\text{\rm{op}}$ }
  }\setlength{\unitlength}{1pt}}
  \end{picture}}
\qquad \quad \text{for~all}~ i\in \mathcal{I}~,
\ee
where we have used the graphical notation
\be   \label{eq:tilde-iota}
\setlength{\unitlength}{1pt}
\iota_{\mathrm{cl}-\mathrm{op}}^* = \bigoplus_{n=1}^N \,\,
\raisebox{-35pt}{
  \begin{picture}(40,75)
   \put(0,8){\scalebox{.75}{\includegraphics{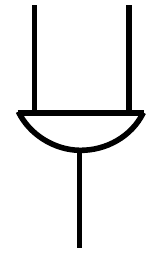}}}
   \put(0,8){
     \setlength{\unitlength}{.75pt}\put(-18,-11){
     \put(20, 90)     {\scriptsize $ C_n^l $}
     \put(50, 90)     {\scriptsize $ C_n^r $}
     \put(35, 2)     {\scriptsize $\Aop$}
     }\setlength{\unitlength}{1pt}}
  \end{picture}}
\ee
and the map
$\iota_{\mathrm{cl}-\mathrm{op}}^\ast: \Aop \to T(\Acl)$ defined via \eqref{f-star-def}.
}
\end{defn}

\begin{rema}  {\rm
(i)
The notion of Cardy $\Cc|\CxC$-algebra was introduced by the first author in \cite[Def.\,5.13]{cardy}. In \cite[Def.\,3.7]{part1}, the first and third author provided the second definition of Cardy algebra by applying the two-sided adjoint functor of $T$. This latter definition is somewhat simpler and convenient to study, and was shown in \cite{part1} to be equivalent to Definition \ref{def:cardy-alg}, but will not be used it in this work. 
\\[.3em]
(ii) The coefficient on the right hand side of the Cardy condition (\ref{eq:cardy-C}) can be modified by an arbitrary non-zero scalar because one can always rescale either $\Delta_{\mathrm{cl}}$ or $\Delta_{\mathrm{op}}$ by a non-zero factor without affecting other axioms but changing the coefficient on the right hand side of (\ref{eq:cardy-C}). In \cite{part1}, a different coefficient was used to in order to make the Cardy condition in the second definition of Cardy algebra free of coefficients (see \cite[Eq. (3.14)]{part1}). In Equation (\ref{eq:cardy-C}), we have used the same coefficient as the one used in \cite{cardy}. It was derived from the natural conventions chosen in the framework of partial conformal field theory introduced in \cite{cardy}. 
}
\end{rema}

\begin{thm}  \label{thm:cardy-determine} {\rm
A solution to the sewing constraints determines a Cardy algebra up to isomorphism.
}
\end{thm}

This theorem will be proved in Section \ref{sec:rel-morph}.

\medskip

It is easy to describe this Cardy algebra explicitly. Suppose that $\Cor$ is a solution to the sewing constraints for $\Bl$. $\Cor$ determines the morphisms
\be
 f_{po} = \Psi_{po}^{-1}(\Cor(\Xr_{po})1) \,:\, \Bop \rightarrow \Bop
 ~~,~~~
 f_{pc} = \Psi_{pc}^{-1}(\Cor(\Xr_{pc})1) \,:\, \Bcl \rightarrow \Bcl ~,
\labl{eq:f-popc-def}
where $1\in \Bl(\emptyset) = \Cb$. We will show in Section \ref{sec:rel-morph} that $f_{po}$ and $f_{pc}$ are idempotents. The Cardy algebra will be defined on the images of $f_{po}$ and $f_{pc}$. By a {\em retract} of an object $V$ with respect to an idempotent $p : V \rightarrow V$ we mean a triple $(U,e,r)$ where $U$ is an object and $e : U \rightarrow V$, $r: V \rightarrow U$ are morphisms such that $e \circ r = p$ and $r \circ e = \id_U$. Choose retracts $(\Aop,\eop,\rop)$ of $\Bop$ with respect to $f_{po}$ and $(\Acl,\ecl,\rcl)$ of $\Bcl$ with respect to $f_{pc}$. Define the remaining morphisms $f_\alpha$ as in \erf{eq:f-popc-def},
\be  \label{eq:f-def}
  f_\alpha = \Psi_{\alpha}^{-1}(\Cor(\Xr_{\alpha})1)
  \quad \text{for} \quad \alpha \in \Sc ~.
\ee
Then it will be proved in Section \ref{sec:rel-morph} that the pair $\Aop,\Acl$ becomes a Cardy algebra when equipped with the structure morphisms
\be\begin{array}{ll}\displaystyle
  m_\text{op} = \rop \circ f_{mo} \circ (\eop \otimes_\Cb \eop)~~,\etb~~~
  \eta_\text{op} = \rop \circ f_{\eta o} ~~,
\enl
  \Delta_\text{op} = (\rop \otimes_\Cb \rop) \circ f_{\Delta o} \circ \eop ~~,\etb~~~
  \eps_\text{op} = f_{\eps o} \circ \eop~~,
\enl
  m_\text{cl} = \rcl \circ f_{mc} \circ (\ecl \otimes_\Cb \ecl) ~~,\etb~~~
  \eta_\text{cl} = \rcl \circ f_{\eta c} ~~,
\enl
  \Delta_\text{cl} = (\rcl \otimes_\Cb \rcl) \circ f_{\Delta c} \circ \ecl ~~,\etb~~~
  \eps_\text{cl} = f_{\eps c}  \circ \ecl ~~,~~~
\enl
  \iota_\text{cl-op} = \rop \cir f_{\iota} \cir T(\ecl) ~~,\etb~~~
  \iota_\text{cl-op}^*=T(\rcl)\circ f_{\iota^*}\circ e_{op}
\eear\labl{eq:falpha-names}

Given another choice $(\Aop',\eop',\rop')$ and $(\Acl',\ecl',\rcl')$ of retracts for $f_{po}$ and $f_{pc}$, we define a morphism between two retracts for $f_{po}$ and $f_{pc}$ by a pair of isomorphisms $\rop' \circ \eop : \Aop \rightarrow \Aop'$ and $\rcl' \circ \ecl : \Acl \rightarrow \Acl'$. It turns out that this pair of isomorphisms gives an isomorphism of Cardy algebras in the sense of \cite[Def.\,3.9]{part1}.

\medskip

However, a Cardy algebra does not determine a unique solution to the sewing constraints. To obtain an one-to-one correspondence, one has to remember the choice of $\Bop, B_l$ and $B_r$, and how the objects $\Aop$ and $\Acl$ underlying the Cardy algebra are embedded in those.

\begin{thm}   \label{thm:cardy-sew}  {\rm
For a given modular tensor category $\Cc$ and three given objects $\Bop,B_l,B_r$ in $\Cc$, 
there is an one-to-one correspondence between the following two notions: 
\bnu
\item A Cardy algebra $(\Aop|\Acl,\iota_\text{cl-op})$ together with a realisation of $\Aop$ as a retract $(\Aop,\eop,\rop)$ of $\Bop$ and a realisation of $\Acl$ as a retract $(\Acl,\ecl,\rcl)$ of $\BlxBr$,
\item A solution $\Cor$ to the sewing constraints for $\Bl(\Cc,\Bop,B_l,B_r)$, together with choices of retracts $(\Aop,\eop,\rop)$ of $\Bop$ and $(\Acl,\ecl,\rcl)$ of $\BlxBr$ such that $\Psi_{po}(\eop \circ \rop) = \Cor(\Xr_{po})$ and $\Psi_{pc}(\ecl \circ \rcl) = \Cor(\Xr_{pc})$.
\enu
Moreover, we can obtain one canonically from the other in an invertible way. In other words, these two notions are equivalent. 
}
\end{thm}

This theorem is proved in Section \ref{sec:rel-morph}.

\begin{rema}  {\rm    
(i) The morphisms $\eop : \Aop \rightarrow \Bop$ and
$\ecl : \Acl \rightarrow \Bcl$ give the images of the morphisms $f_{po}$ and $f_{pc}$ defined in \erf{eq:f-popc-def}. Physically, $f_{po}$ and $f_{pc}$ correspond to the propagator of open and closed states, respectively, and inserting a state lying in the kernel of a propagator into a correlator causes this correlator to vanish. From this point of view, a Cardy algebra precisely captures all the information about non-vanishing correlators.
\\[.3em]
(ii) Suppose we are given another set of objects $\Bop', B_l', B_r'$ and a realisation of $\Bop, B_l, B_r$ as retracts of $\Bop', B_l', B_r'$. Then there is an injective map from solutions of the sewing constraints for $\Bl(\Bop,B_l,B_r)$ to solutions of the sewing constraints for $\Bl(\Bop',B_l',B_r')$. The quickest way to see this is to pass to the Cardy algebra and change the realisations of $\Aop$ and $\Acl$ as retracts. For example, if $(\Bop,e,r)$ is the retract of $\Bop'$, then replace the retract $(\Aop,\eop,\rop)$ of $\Bop$ by the retract $(\Aop,e\circ\eop,\rop\circ r)$ of $\Bop'$. This shows that we can always pass to larger objects $\Bop', B_l', B_r'$ without loosing any solutions.
\\[.3em]
(iii) Theorem \ref{thm:cardy-sew} is an extension of Theorem 4.26 in \cite{unique}. The latter theorem gives a correspondence
(between a special type of Cardy algebra, namely those of the form $(A|Z(A),e)$ as defined in \cite[Thm.\,3.18]{part1}, and a special class of solutions to the sewing constraints, namely those satisfying the conditions stated in \cite[Thm.\,4.26]{unique}. These conditions include in particular that $\Acl$ is haploid. 
As an aside, the reason that in \cite[Thm.\,4.26]{unique} there is no mention of retracts is that in \cite{unique} they are included in the definition of the functor $\Bl$, and it is imposed that $f_{po}$ and $f_{pc}$ (see point (i)) are identities.
}
\end{rema}

\sect{Proofs of Theorem \ref{thm:cardy-determine} and \ref{thm:cardy-sew}}
\label{sec:rel-morph}
The main goal of this section is to provide proofs of Theorem \ref{thm:cardy-determine} and \ref{thm:cardy-sew}.
The main body of the proof contains explicit illustrations of how to reduce the relations R1-R32 given in Section \ref{sec:gen-rel} to the conditions on morphisms in the categories $\Cc$ and $\CxC$ and vice versa.

More precisely, we need to explicitly spell out the condition (\ref{eq:sewing-cond}) in the case of $U=\One$ and $F=\Bl$.
We will set up our problem in the following way. Given a set of maps $c_\alpha: \Cb \to \Bl(\Xr_\alpha)$ for $\alpha\in \Sc$, we define the morphisms $f_\alpha$ (recall (\ref{eq:f-def})) in $\Cc$ or $\CxC$ as follows:
\be   \label{eq:f-alpha-def}
f_\alpha = \Psi_\alpha^{-1}(c_\alpha(1)), \quad \quad \forall \alpha \in \Sc,
\ee
Recall (\ref{eq:Y_Rn-l-r}).
Using the maps:
$$
\Bl(\varpi_{\Rr n,l}): \Bl(\Xr_{\alpha_{\Rr n}^1} \otimes \cdots \otimes \Xr_{\alpha_{\Rr n}^{l(n)}}) \to \Bl(\Yr_{\Rr n}),
$$
$$
\Bl(\varpi_{\Rr n,r}): \Bl(\Xr_{\beta_{\Rr n}^1} \otimes \cdots \otimes \Xr_{\beta_{\Rr n}^{r(n)}}) \to \Bl(\Yr_{\Rr n}),
$$
and the fact that $\Bl$ is a monoidal functor,
we obtain 32 pairs of morphisms $C_{\Yr_{\Rr n, l\slash r}^d}: \Cb \to \Bl(\Yr_{\Rr n})$ defined by
$\Bl(\varpi_{\Rr n,l}) \circ (c_{\alpha_{\Rr n}^1} \otimes_\Cb \cdots \otimes_\Cb c_{\alpha_{\Rr n}^{l(n)}})$ and $\Bl(\varpi_{\Rr n,r}) \circ (c_{\beta_{\Rr n}^1} \otimes_\Cb \cdots \otimes_\Cb c_{\beta_{\Rr n}^{r(n)}})$, respectively. This is nothing but the equation (\ref{eq:nat-C-def}) rewritten in the case $U=\One$ and $F=\Bl$. Therefore, the 32 relations $C_{\Yr_{\Rr n,l}^d} = C_{\Yr_{\Rr n,r}^d}$ give 32 relations on $c_\alpha$ and equivalently on $f_\alpha, \alpha\in \Sc$, i.e. 
\begin{align}  \label{eq:sew-eq-f-alpha}
&\Bl(\varpi_{\Rr n,l}) \circ ( \Psi_{\alpha_{\Rr n}^1}(f_{\alpha_{\Rr n}^1}) \otimes_\Cb \cdots \otimes_\Cb 
\Psi_{\alpha_{\Rr n}^{l(n)}}(f_{\alpha_{\Rr n}^{l(n)}}) )   \nn
&\hspace{3cm} = \Bl(\varpi_{\Rr n,r}) \circ ( \Psi_{\beta_{\Rr n}^1}(f_{\beta_{\Rr n}^1}) \otimes_\Cb \cdots \otimes_\Cb 
\Psi_{\beta_{\Rr n}^{r(n)}}(f_{\beta_{\Rr n}^{r(n)}}) )  \ .
\end{align}
By (\ref{eq:def-Psi}), the two sides of above equation can be expressed by two 3-bordisms equipped with ribbon graphs, their comparison leads to algebraic relations among $f_\alpha$. They are listed in Table 1 to 4. We will derive these relations on $f_\alpha$ in five groups in Section \ref{sec:R1-R13}-\ref{sec:R32} before we prove Theorem \ref{thm:cardy-determine} and \ref{thm:cardy-sew} in Section \ref{sec:proof-sew}.

 \medskip
 Notice that the 32 relations naturally split into the following five groups:
\begin{enumerate}
\item R1-R13 are relations involving only open state boundaries at genus zero;
\item R14-R25 are relations involving only closed state boundaries at genus zero;
\item R26-R30 are relations involving both open and closed state boundaries at genus zero;
\item R31 is the Cardy condition;
\item R32 is the modular invariant condition for genus-one surface.
\end{enumerate}
We will derive some important relations on $f_\alpha$ in each group in full details to illustrate the idea and leave the rest to readers.
To derive the relations on $f_\alpha$ from $\Rr1$-$\Rr32$, we need to set up a convention to move the marked arcs in a consistent
way to a circle on the boundary components of extended cobordisms. Let $\Mr$ be an extended 3-cobordism with at least one marked arc
labelled by $\Bop$, then the component of the boundary containing this marked arc must have an involution, whose fixed point circle contains
this marked arc. We move the arcs on this boundary component of $M$ according to the following rule:
\begin{enumerate}
 \item Ordering: The ordering of the marked arcs is cyclic on the fixed points circle of the involution, and the cyclic orientation induced
 by the ordering of the arcs are opposite to the orientation of the fixed point circle. 
 \item Braiding: A ribbon ending on an arc labelled by $B_l$ can only upper-cross a ribbon ending on an arc labeled by $\Bop$ or $B_r$. A ribbon ending on an arc labelled by $B_{op}$ can only upper-cross a ribbon ending on an arc labeled $B_r$.  
\end{enumerate}

An illustration of the rule of moving the arcs is the following:

\begin{equation}
\figbox{0.45}{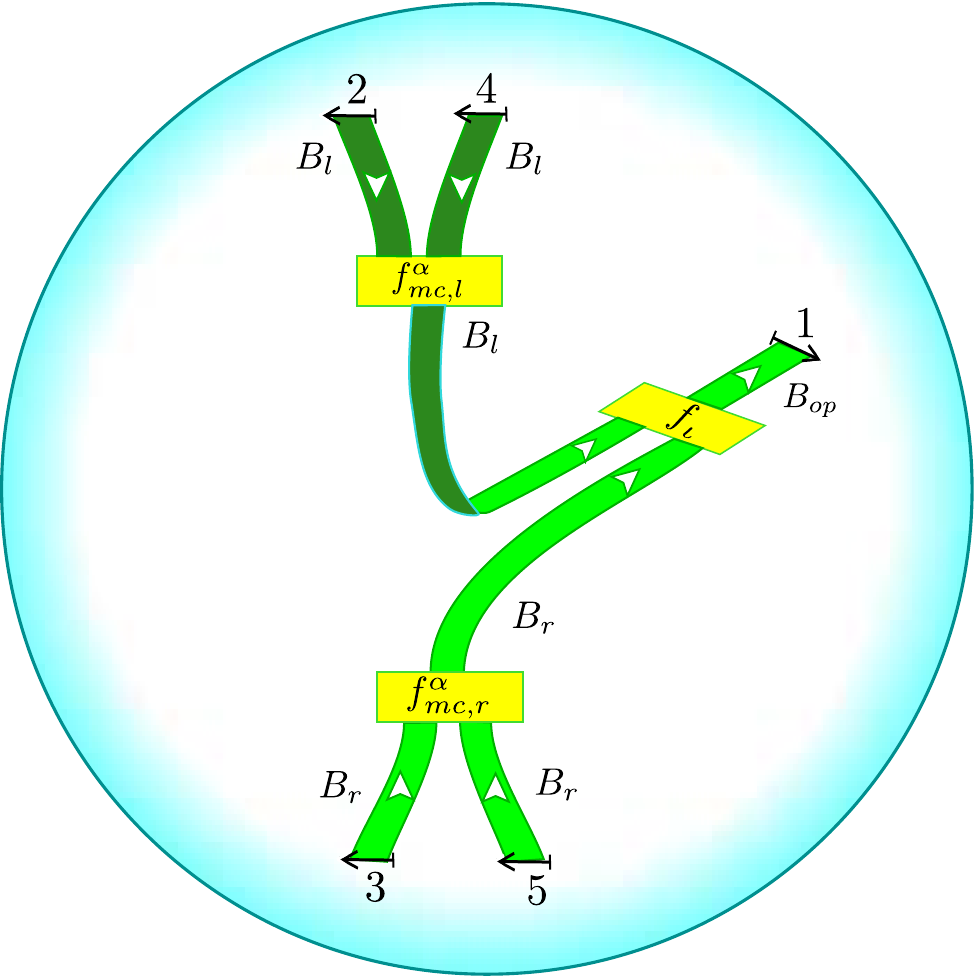}\hspace{6mm}\rightsquigarrow\hspace{6mm}\figbox{0.45}{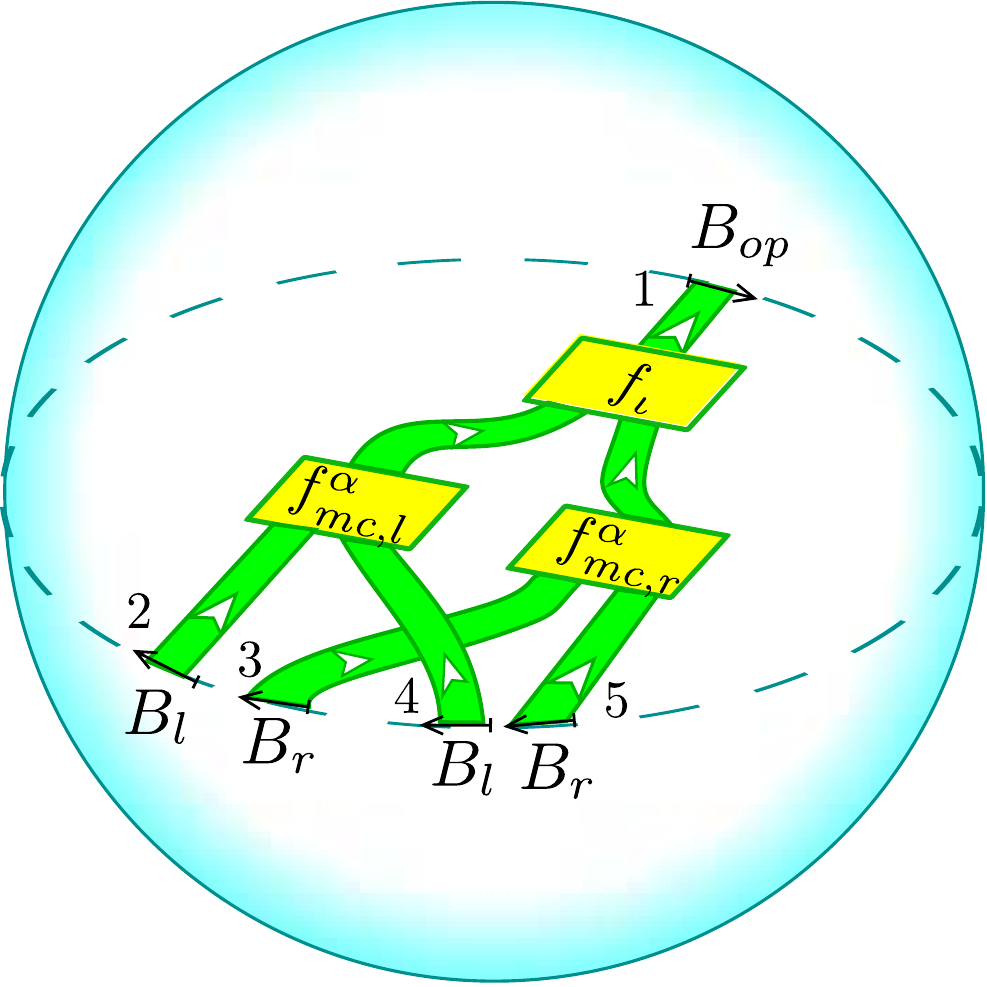} 
\end{equation}

\subsection{Relations R1-R13}\label{sec:R1-R13}

The full list of the relations of $f_\alpha$ obtained from the relations R1-R13 is shown in Table \ref{tab:glue-open}.
They have been studied in \cite{unique}. For the completeness, we derive the relations of $f_\alpha$ obtained from R1 and R9 in details here. 

We start with an explicit choice for the world sheets $\Xr_{\eta o}$, $\Xr_{po}$,
$\Xr_{mo}$ and $\Xr_{\Delta o}$. Then we construct an extended cobordism $\Mr_{\varpi_{\Rr 1,l}}=(\Mr_{\varpi{\Rr 1,l}},0,h)$ associated to the sewing morphism $\varpi_{\Rr 1,l} : \Xr_{mo} \oti \Xr_{\eta o} \rightarrow \Yr_{\Rr 1}$ (recall (\ref{eq:shat-cobord}), (\ref{eq:bl(m)-def})) and an extended cobordism that determines $C_{\Yr^d_{\Rr1,l}}$ as the image of $\tftC$. Repeat this for $C_{\Yr^d_{\Rr1,r}}$. Then we obtain the identity determined by relation R1 listed in Table \ref{tab:glue-open}. Then we go through a similar steps to prove the identity determined by relation R9 in Table \ref{tab:glue-open}.

\subsubsection*{The world sheet $\Xr_{\eta o}$}

Let $D_w = \{ z \,|\, |z{-}w| < 1 \}$ be the open unit disk shifted by $w \in \Cb$ and
let $S_w^\pm =  \{ z \,|\, |z{-}w| = 1 \}$ be the shifted unit circle with anti-clockwise (for label `$+$') or clockwise (for label `$-$') orientation.
Let $\overline \Cb = \Cb \cup \{ \infty \}$ be the Riemann sphere with its standard
orientation, and let $\overline \Hb$ be the upper half plane together with the real line and infinity.

We set $\Xtil_{\eta o} = \overline\Cb - D_0$. The involution is just complex conjugation, $\imath_{\eta o}(z) = \bar{z}$. We identify the quotient surface $\dot\Xr_{\eta o}$ with $\Xtil_{\eta o} \cap \overline\Hb$ 
and choose $\text{or}_{\eta o}(z) = z$. The image of $\text{or}_{\eta o}$ is thus in the upper half plane.

The boundary of $\Xtil_{\eta o}$ is
$\partial \Xtil_{\eta o} = S^-_0$. We take
$\delta_{\eta o}(e^{i \theta}) = - e^{-i\theta} \In S^1$. The minus in the exponent is necessary to make $\delta_{\eta o}$ orientation preserving, and the minus in front of the exponential makes sure that the image of $\text{or}_{\eta o}$ in $S^-_0$ gets mapped to the upper half circle $S^1 \cap \overline\Hb$ as required by Definition \ref{def:ws}.
The partition of the boundary into in-coming and out-going components is $b^\text{in}_{\eta o} = \emptyset$ and $b^\text{out}_{\eta o} = \{ S^-_0 \}$. The order map ord is the unique map. 

We will also need the extended double $\Xhat_{\eta o}$ of $\Xr_{\eta o}$. By definition it is obtained by starting from the disjoint union of $\Xtil_{\eta o}$ and $\vec D$ (recall Section \ref{sec:Bl-def}) and identifying the point $e^{i \theta} \In  S^1 = \partial \vec D$ with
$\delta^{-1}_{\eta o}(-e^{-i \theta}) = e^{i \theta} \In  S^-_0$. This simply means that we can identify $\Xhat_{\eta o}$ with $\overline\Cb$ and a marked arc at zero oriented along the real axis towards $+ \infty$. The arc gets labelled by $(\Bop,-)$.

\subsubsection*{The world sheet $\Xr_{po}$}

We take $\Xtil_{po} = \overline\Cb - (D_0 \cup D_3)$, the involution $\imath_{po}$ to be the complex conjugation and the image of $\text{or}_{po}$ to lie again in the upper half plane. We have $\partial\Xtil_{po} = S^-_0 \cup S^-_3$. Set
$\delta_{po}(e^{i \theta}) = - e^{-i\theta}$ and
$\delta_{po}(3{+}e^{i \theta}) = - e^{-i\theta}$. Finally,
$b^\text{in}_{po} = \{ S^-_0 \}$ and
$b^\text{out}_{po} = \{ S^-_3 \}$. The ordering of the boundary components is given by $ \text{ord}:  S^-_0\mapsto 1,  S^-_3\mapsto 2$.

The extended double $\Xhat_{po}$ is $\overline\Cb$ with an arc at 0 labelled by $(\Bop,+)$ and an arc at 3 labelled by $(\Bop,-)$. Both arcs are oriented towards $+\infty$.

\subsubsection*{The world sheet $\Xr_{mo}$}

We take $\Xtil_{mo} = \overline\Cb - (D_0 \cup D_3 \cup D_6)$, and $\imath_{mo}$ to be the complex conjugation and the image of $\text{or}_{mo}$ to lie in the upper half plane. The boundary is given by $\partial\Xtil_{mo} = S^-_0 \cup S^-_3 \cup S^-_6$ and we parametrise it by 
$\delta_{mo}(e^{i \theta}) =
\delta_{mo}(3{+}e^{i \theta}) =
\delta_{mo}(6{+}e^{i \theta}) = - e^{-i\theta}$.
The in-coming and out-going boundary components are
$b^\text{in}_{mo} = \{ S^-_0 , S^-_3 \}$ and
$b^\text{out}_{mo} = \{ S^-_6 \}$. The ordering of the boundary components is given by$\text{ord}:  S^-_0\mapsto 1,  S^-_3\mapsto 2, S^-_6\mapsto 3$.

The extended double $\Xhat_{mo}$ is $\overline\Cb$ with arcs at 0 and 3 labelled by $(\Bop,+)$ and an arc at 6 labelled by $(\Bop,-)$. All arcs are oriented towards $+\infty$.

\subsubsection*{The world sheet $\Xr_{\Delta o}$}

Let $\Xtil_{\Delta o} = \overline\Cb - (D_0 \cup D_3 \cup D_6)$, and $\imath_{\Delta o}$ be the complex conjugation and the image of $\text{or}_{\Delta o}$ lie in the upper half plane. The boundary is given by $\partial\Xtil_{\Delta o} = S^-_0 \cup S^-_3 \cup S^-_6$,
we parametrise it by
$$\delta_{\Delta o}(e^{i \theta}) =
\delta_{\Delta o}(3{+}e^{i \theta}) =
\delta_{\Delta o}(6{+}e^{i \theta}) = - e^{-i\theta}.$$
The in-coming and out-going boundary components are
$b^\text{in}_{\Delta o} = \{ S^-_0  \}$ and
$b^\text{out}_{\Delta o} = \{ S^-_3 , S^-_6\}$.
 The ordering of the boundary components is given by: $\text{ord}:  S^-_0\mapsto 1,  S^-_3\mapsto 2, S^-_6\mapsto 3$.

The extended double $\Xhat_{\Delta o}$ is $\overline\Cb$ with an arc at 0 labelled by $(\Bop,+)$ and arcs at 3 and 6 labelled by $(\Bop,-)$. All arcs are oriented towards $+\infty$.

\begin{table}[bt]
$$
\begin{array}{llll}
\text{R1} \etb:
f_{mo} \cir (f_{\eta o} \otimes_\Cb \id_{\Bop}) = f_{po}
\qquad \etb
\text{R2} \etb:
f_{mo} \cir (\id_{\Bop} \otimes_\Cb f_{\eta o}) = f_{po}
\enl
\text{R3} \etb:
(f_{\eps o} \otimes_\Cb \id_{\Bop}) \cir f_{\Delta o}  = f_{po}
\qquad \etb
\text{R4} \etb:
(\id_{\Bop} \otimes_\Cb f_{\eps o}) \cir f_{\Delta o}  = f_{po}
\enl
\text{R5} & \multicolumn{3}{l}{\hspace*{-.5em}: \displaystyle
((f_{\eps o} \cir f_{mo}) \otimes_\Cb \id_{\Bop^\vee}) \cir
(\id_{\Bop} \otimes_\Cb b_{\Bop})
=
(\id_{\Bop^\vee} \otimes_\Cb (f_{\eps o} \cir f_{mo})) \cir
(\tilde b_{\Bop} \otimes_\Cb \id_{\Bop}) }
\enl
\text{R6} & \multicolumn{3}{l}{\hspace*{-.5em}: \displaystyle
f_{mo} \cir ( f_{mo} \otimes_\Cb \id_{\Bop} )
= f_{mo} \cir ( \id_{\Bop} \otimes_\Cb f_{mo} ) }
\enl
\text{R7} & \multicolumn{3}{l}{\hspace*{-.5em}: \displaystyle
( f_{\Delta o} \otimes_\Cb \id_{\Bop} ) \cir f_{\Delta o}
=  ( \id_{\Bop} \otimes_\Cb f_{\Delta o} )
\cir f_{\Delta o}}
\enl
\text{R8} & \multicolumn{3}{l}{\hspace*{-.5em}: \displaystyle
(f_{mo} \otimes_\Cb \id_{\Bop}) \cir (\id_{\Bop} \otimes_\Cb f_{\Delta o})
= f_{\Delta o} \cir f_{mo}}
\enl
\text{R9} & \multicolumn{3}{l}{\hspace*{-.5em}: \displaystyle
(\id_{\Bop} \otimes_\Cb f_{mo}) \cir (f_{\Delta o} \otimes_\Cb \id_{\Bop})
= f_{\Delta o} \cir f_{mo}}
\enl
\text{R10} \etb:
f_{po} \cir f_{mo} = f_{mo}
\etb
\text{R11} \etb:
f_{po} \cir f_{\eta o} = f_{\eta o}
\enl
\text{R12} \etb:
f_{\Delta o} \cir f_{po} = f_{\Delta o}
\etb
\text{R13} \etb:
f_{\eps o} \cir f_{po} = f_{\eps o}
\end{array}
$$
\caption{The relations of $f_\alpha$ obtained from the relations R1-R13.}
\label{tab:glue-open}
\end{table}

\subsubsection*{The morphisms $\varpi_{\Rr 1, l\slash r}$}

Since $\Yr_{\Rr 1} = \Xr_{po}$, we take $\varpi_{\Rr 1, r} = \id_{\Xr_{p o}}$ by our convention (recall \eqref{eq:Yn-id-convention}). We need a morphism $\varpi_{\Rr 1, l}=(\sew, f): \Xr_{\eta o} \otimes \Xr_{mo} \rightarrow \Xr_{po}$.
We proceed as follows.

The world sheet $\Xr_{\eta o} \otimes \Xr_{mo}$ is a subset of $\overline\Cb \times \{1,2\}$, namely $(\Xr_{\eta o},1) \cup (\Xr_{mo},2)$. The sewing $\sew$ is the pair $\{ \big( (S^-_0,1), (S^-_0,2) \big) \}$. The sewn world sheet $\sew(\Xr_{\eta o} \otimes \Xr_{mo})$ is given by identifying the points $(-e^{-i\varphi},1) \sim (e^{i\varphi},2)$ of $\Xr_{\eta o} \otimes \Xr_{mo}$. Finally we need a homeomorphism $f : \sew(\Xr_{\eta o} \otimes \Xr_{mo}) \rightarrow \Xr_{po}$. We take $f\big((z,1)\big) = -z^{-1}-3$ and $f\big((z,2)\big)=z-3$. This is compatible with the identification $\sim$ just described, and it maps the boundary circles $S^-_3$ and $S^-_6$ of $\Xr_{mo}$ to the boundary circles $S^-_0$ and $S^-_3$ of $\Xr_{po}$.

\subsubsection*{The cobordism $\Mr_{\varpi_{\Rr 1,l\slash r}}$ associated to $\varpi_{\Rr 1, l\slash r}$}

The cobordism for $\varpi_{\Rr 1, r}$ is simply $\Xhat_{po} \times [0,1]$, where it is understood that the marked arcs on $\Xhat_{po}$ give rise to vertical ribbons inside the cobordism.

The extended cobordism $\Mr_{\varpi_{\Rr 1, l}}$ associated to $\varpi_{\Rr 1, l}$ is defined by
$\Mr_{\varpi_{\Rr 1, l}} = \widehat{\Yr} \times [0,1] / \sim$, where $\Yr = \Xr_{\eta o} \otimes \Xr_{mo}$ and the equivalence relation is given by $(z,1,1) \sim (-\bar{z},2,1)$ for all $|z| \le 1$.
It is understood that the marked arcs on $\widehat\Yr$ give rise to vertical ribbons embedded in $\Mr_{\varpi_{\Rr 1, l}}$. The boundary of $\Mr_{\varpi_{\Rr 1, l}}$ consists of two parts $\partial_- \Mr_{\varpi_{\Rr 1, l}} = \widehat{\Yr} \times \{0\}$ and $\partial_+ \Mr_{\varpi_{\Rr 1, l}} = \widehat\Zr \times \{1\}$, with
$\Zr = \sew(\Xr_{\eta o} \otimes \Xr_{mo})$. The homeomorphism $h\big|_{\partial_- \Mr}$ is given by the identity on $\overline{\widehat\Yr}$ (which is orientation reversing because the boundary $\partial \Mr$ is oriented by the inward pointing normal) and $h\big|_{\partial_+ \Mr_{\varpi_{\Rr 1, l}}}$ is $\hat f$ on $\widehat{\Zr}$. After a little thought, the reader will agree that $\Mr_{\varpi_{\Rr 1,l}}$ homeomorphic to the three-manifold (we nonetheless write `=')
\be
  \Mr_{\varpi_{\Rr 1,l}} = \figbox{0.35}{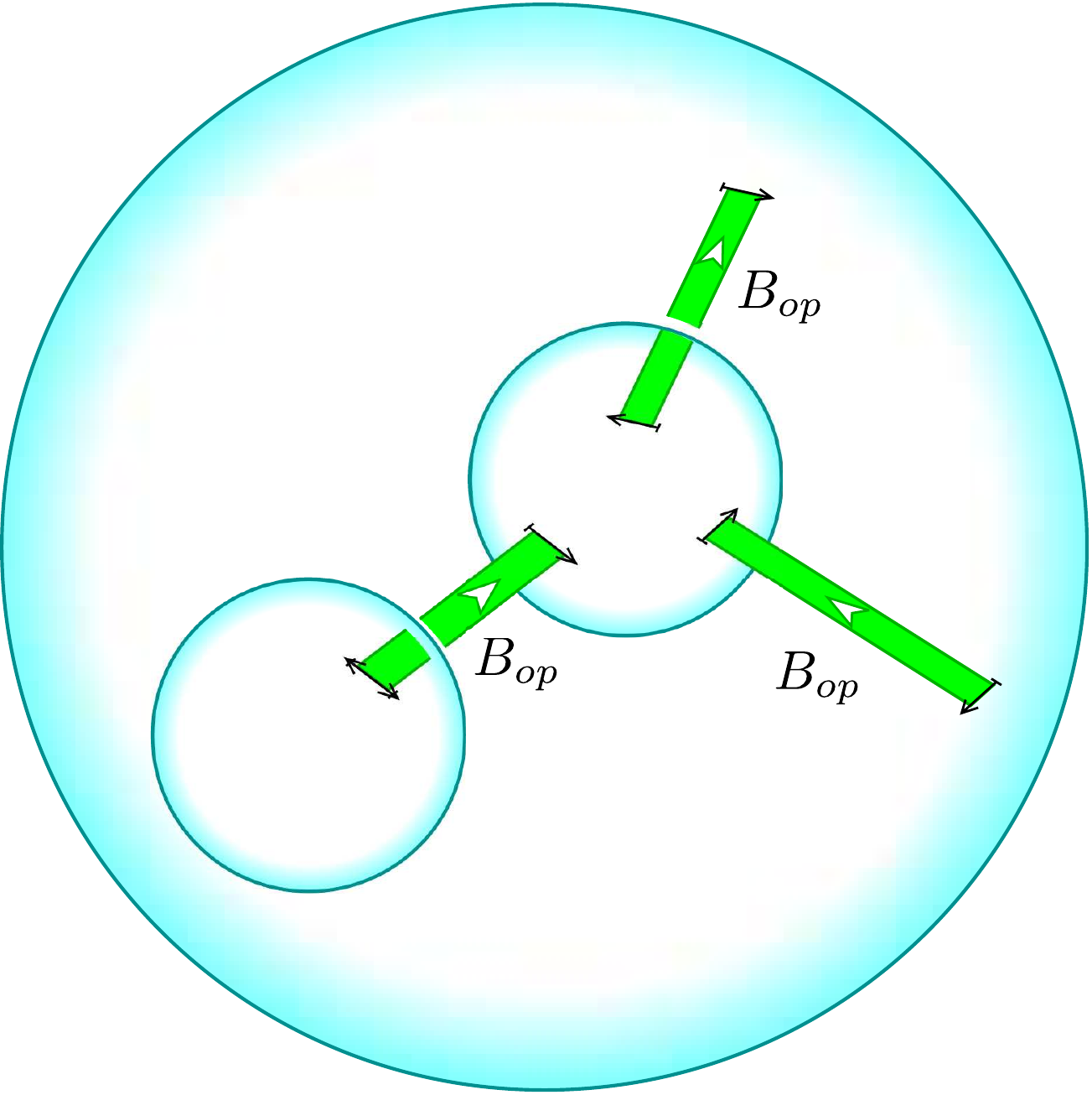}
\ee
i.e.\ a solid three-ball with two solid three-balls cut out and ribbons and marked points as indicated.

\subsubsection*{Relation R1}

Next we need to describe the implication of the condition $C_{\Yr^d_{\Rr1,l}} = C_{\Yr^d_{\Rr1,r}}$.  For the right hand side, we obtain $C_{\Yr^d_{\Rr1,r}} 1= \Bl(\varpi_{\Rr 1,r}) ( \Psi_{po}(f_{po}) )= \Psi_{po}(f_{po})$,
 as $\varpi_{\Rr 1,r}$ was just the identity morphism. As before, the `1' appears, as $C_{\Yr^d_{\Rr1,r}}$ is a linear map
 $\Cb \rightarrow \Bl(\Yr_{\Rr 1})$, and we want to describe a vector in $\Bl(\Yr_{\Rr 1})$. The left hand side of the condition gives
\be
  C_{\Yr^d_{\Rr 1,l}}1 = \Bl(\varpi_{\Rr 1,l})
  \big( \Psi_{\eta o}(f_{\eta o}) \otimes_\Cb
  \Psi_{mo}(f_{mo}) \big)
  = \tftC\Bigg(\ \figbox{0.38}{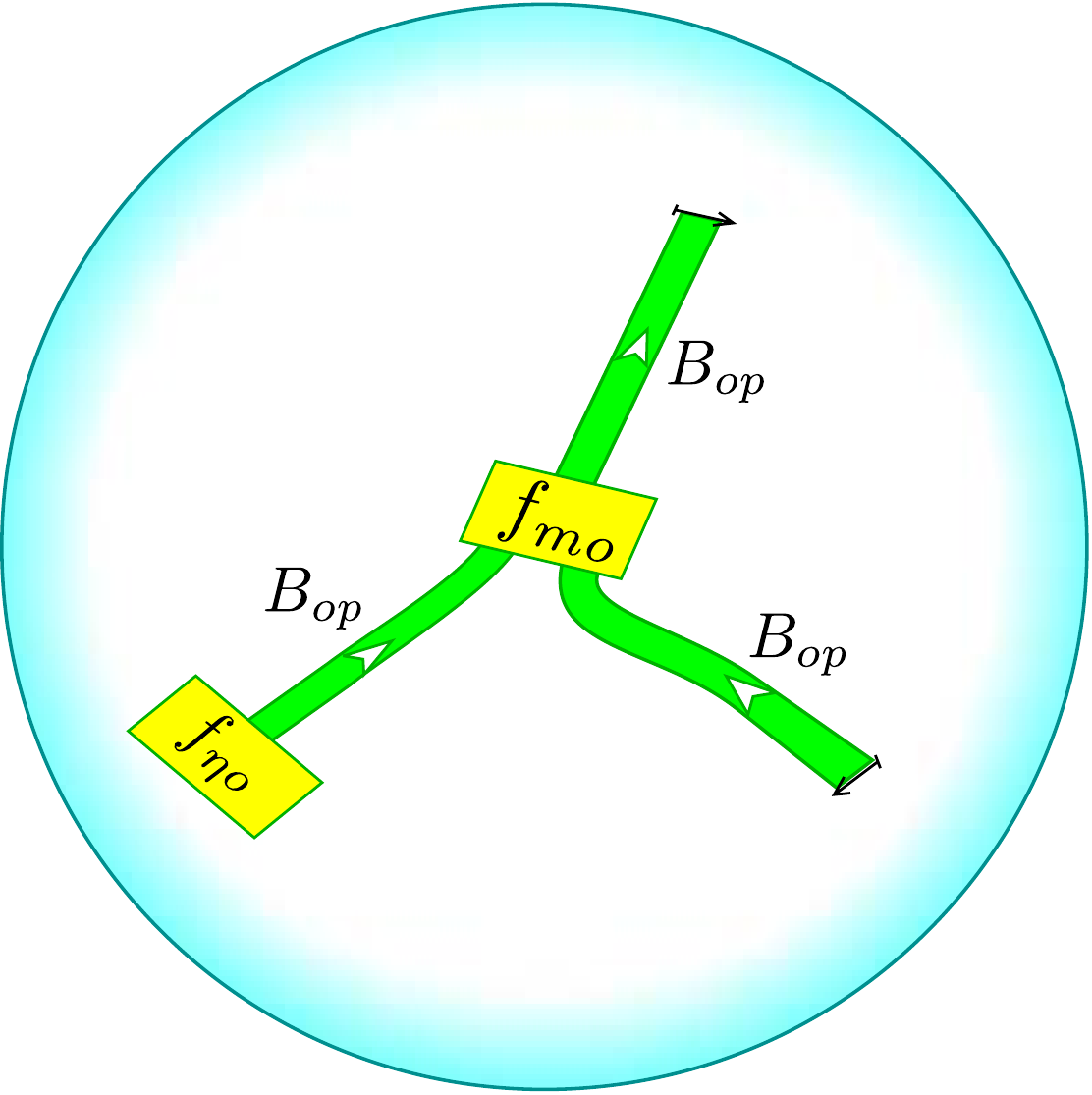}\ \Bigg)1, 
\ee
the cobordism in which is obtained by sewing $\Mr_{mo}(-)$ and $\Mr_{\eta o}(-)$ (defined in Figure \ref{fig:Malpha-cobord}) to $\Mr_{\varpi_{\Rr1,l}}$. Then it is clear that $C_{\Yr^d_{\Rr1,l}} = C_{\Yr^d_{\Rr1,r}}$ if and only if $f_{mo} \cir (f_{\eta o} \otimes_\Cb \id_{\Bop}) = f_{po}$.

\subsubsection*{Relation R9}

The cobordism $\Mr_{\varpi_{\Rr9,l}}$associated to $\varpi_{\Rr9,l}$ is given by $\Mr_{\varpi_{\Rr9,l}}=\hat{\Yr}_l\times[0,1]/\sim$, where $\Yr_l=\Xr_{\Delta o}\otimes \Xr_{mo}$ and the equivalence relation $\sim$ is defined by
$(6+z,1,1)\sim(3-\bar{z},2,1)$, for all $|z|\leqslant 1$, and is given graphically by
$$
\Mr_{\varpi_{\Rr9,l}}=\figbox{0.8}{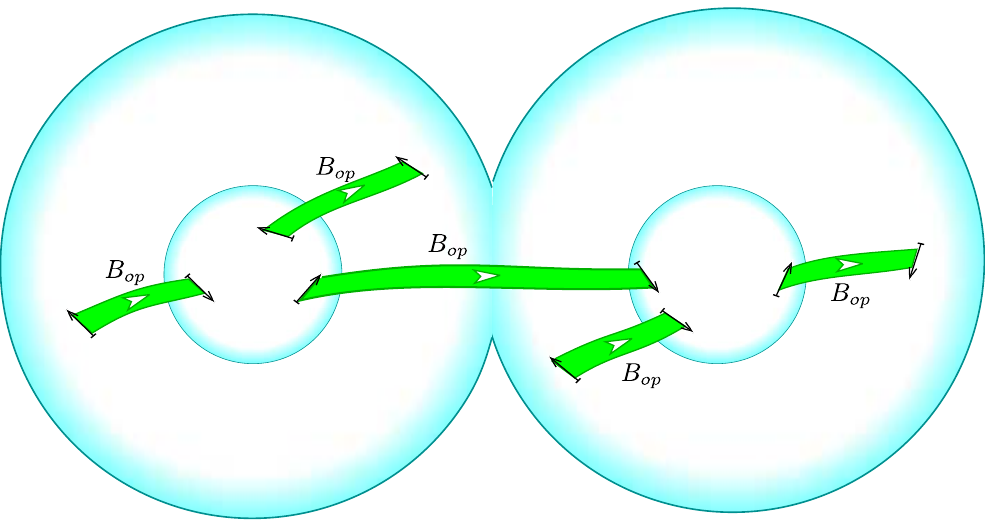}
$$
Consequently, we have 
$$
 C_{\Yr_{\Rr9,l}^d}1=\Bl(\varpi_{\Rr9,l})\circ(\Psi_{\Delta o}(f_{\Delta o})\otimes_\Cb \Psi_{mo}(f_{mo}))=\tftC\Bigg(\figbox{0.45}{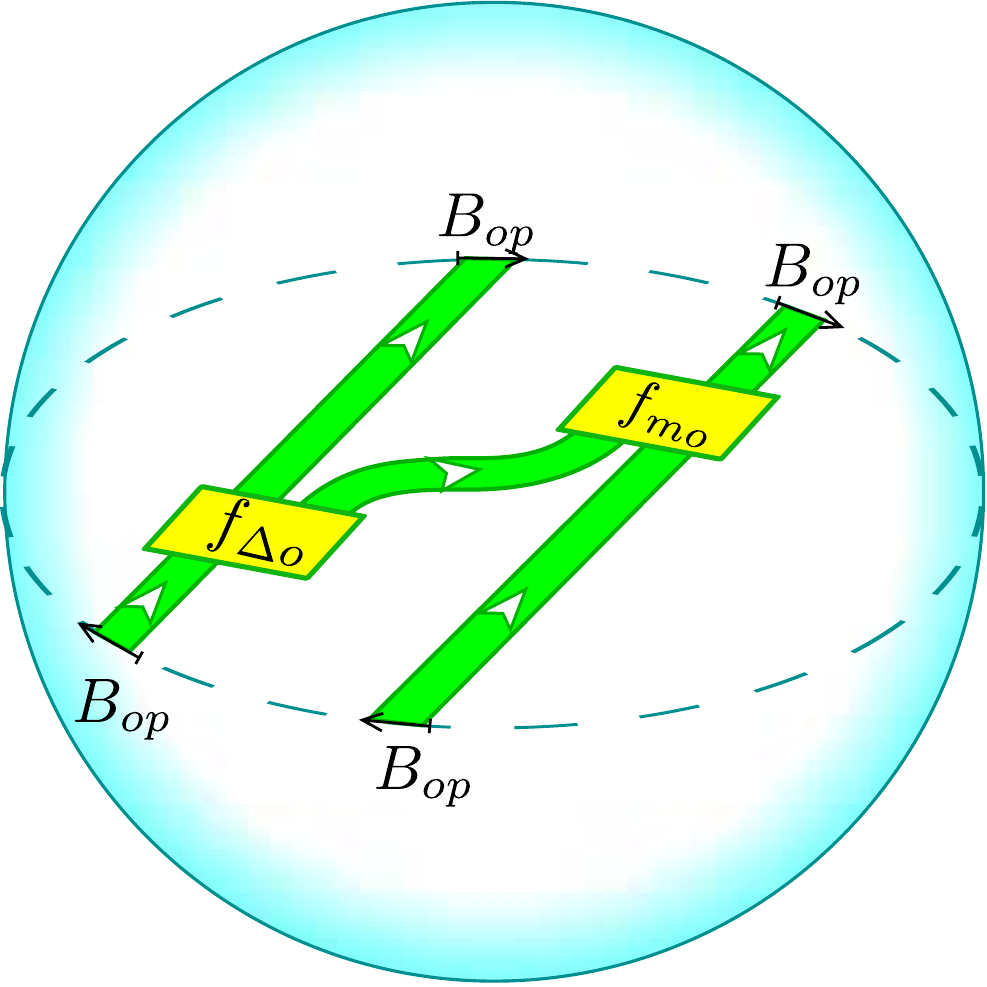}\Bigg)1\, .
$$
The cobordism $\Mr_{\varpi_{\Rr9,r}}$associated to $\varpi_{\Rr9,r}$ is given by $\Mr_{\varpi_{\Rr9,r}}=\hat{\Yr}_r\times[0,1]/\sim$, where $\Yr_r=\Xr_{mo}\otimes \Xr_{\Delta o}$ and the equivalence relation $\sim$ is given by $(6+z,1,1)\sim(-\bar{z},2,1), |z|\leqslant 1$, and is given graphically by
$$
\Mr_{\varpi_{\Rr9,r}}=\figbox{0.8}{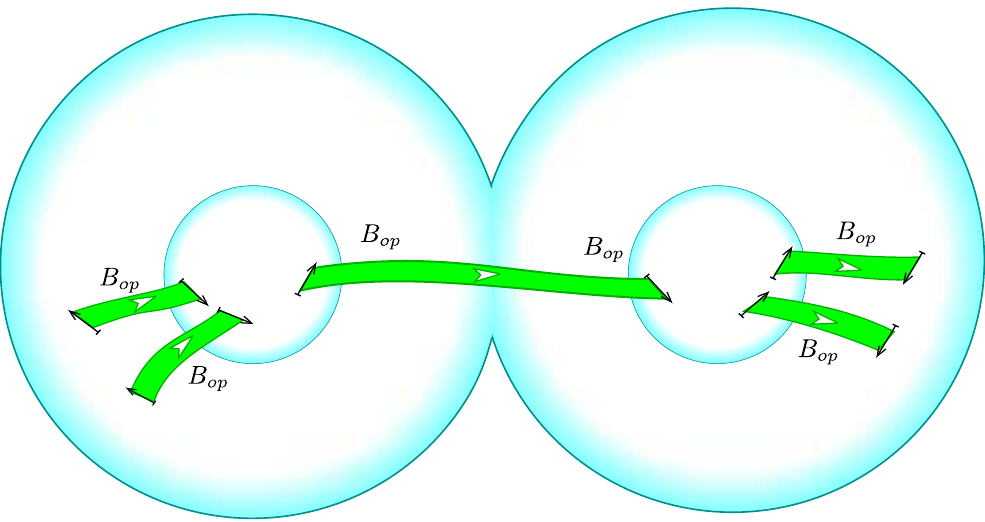}~.
$$
Consequently, we have
$$
 C_{\Yr_{\Rr9,r}^d}1=\Bl(\varpi_{\Rr9,r})\circ(\Psi_{mo}(f_{mo})\otimes_\Cb \Psi_{\Delta o}(f_{\Delta o}))=\tftC\Bigg(\figbox{0.45}{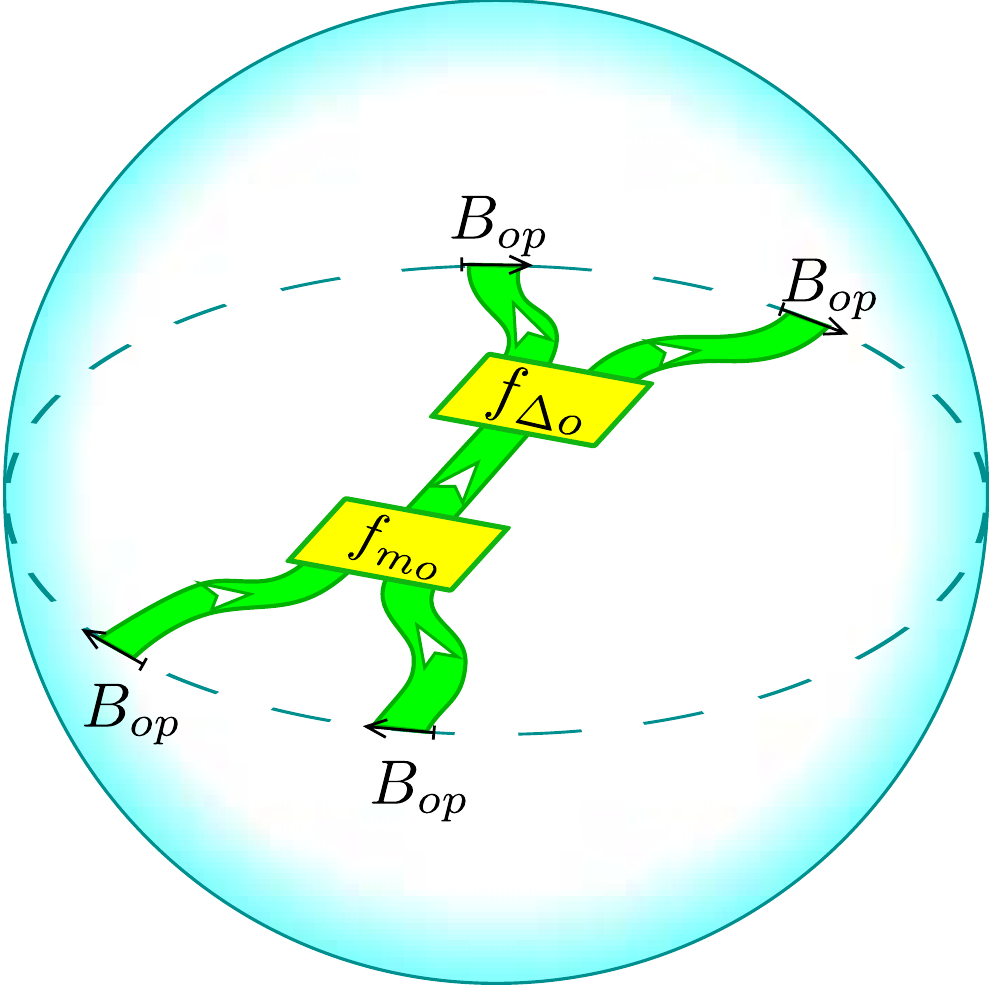}\Bigg)1.
$$
Hence, $C_{\Yr_{\Rr9,l}^d}=C_{\Yr_{\Rr9,r}^d}$ if and only if $(\id_{B_{op}}\otimes_\Cb f_{mo})\circ (f_{\Delta o}\otimes_\Cb \id_{B_{op}} )=f_{\Delta o}\circ f_{mo}$.

\subsection{Relations R14-R25} \label{sec:R14-R25}

Let us briefly describe the implications of relations R14 and R20, the complete list of resulting identities of $f_\alpha$ are given in Table \ref{tab:glue-closed}.

\subsubsection*{The world sheet $\Xr_{\eta c}$}
We take $\tilde{\Xr}_{mc}=(\overline{\mathbb{C}}-D_0)_l \sqcup (\overline{\mathbb{C}}-D_0)_r$. The subscripts $l$ and $r$ here denote the ``left'' and ``right'' component respectively. The involution $\imath_{mc}$ is given by $\imath_{mc}: z_l\mapsto \bar{z}_r, z_r\mapsto \bar{z}_l$. The boundary of $\tilde{\Xr}_{mc}$ is $\partial\tilde{X}_{mc}=(S_0^-)_l\coprod (S_0^-)_r$.
We parametrise the boundary by
$$
\delta_{mc}((e^{i\theta})_l)=-e^{-i\theta}, \quad \delta_{mc}((e^{i\theta})_r)=-e^{-i\theta}.
$$
The ordering of the boundary components is given by $\text{ord}:  ((S^-_0)_l, (S^-_0)_r)\mapsto(1,2)$. The extended double $\hat{X}_{mc}$ is $\overline{\mathbb{C}}_l\cup\overline{\mathbb{C}}_r$ with arcs at $(0)_l$ labelled by $(B_l,-)$, arcs at $(0)_r$ labelled by $(B_r,-)$.

\subsubsection*{The world sheet $\Xr_{pc}$}
We take $\tilde{\Xr}_{mc}=(\overline{\mathbb{C}}-(D_0\cup D_3))_l \sqcup (\overline{\mathbb{C}}-(D_0\cup D_3))_r$. As before, the involution $\imath_{mc}$ is given by 
$\imath_{mc}: z_l\mapsto \bar{z}_r, z_r\mapsto \bar{z}_l$. The boundary of $\tilde{\Xr}_{mc}$ is $\partial\tilde{X}_{mc}=(S_0^-\cup S_3^-)_l\coprod (S_0^-\cup S_3^-)_r$. We parametrise the boundaries by
$$
\delta_{mc}((e^{i\theta})_l)=\delta_{mc}((3+e^{i\theta})_l)=-e^{-i\theta}, \quad \delta_{mc}((e^{i\theta})_r)=\delta_{mc}((3+e^{i\theta})_r)=-e^{-i\theta}.
$$
The ordering of the boundary components is given by
$$ \text{ord}:  ((S^-_0)_l, (S^-_0)_r, (S^-_3)_l,(S^-_3)_r)\mapsto(1,2,3,4).$$
The extended double $\hat{X}_{mc}$ is $\overline{\mathbb{C}}_l\cup\overline{\mathbb{C}}_r$ with arcs at $(0)_l$ labelled by $(B_l,+)$, an arc at $(3)_l$ labelled by $(B_l,-)$, arcs at $(0)_r$ labelled by $(B_r,+)$ and an arc at $(3)_r$ labelled by $(B_r,-)$.

\subsubsection*{The world sheet $\Xr_{mc}$}
We take $\tilde{\Xr}_{mc}=(\overline{\mathbb{C}}-(D_0\cup D_3\cup D_6))_l \sqcup (\overline{\mathbb{C}}-(D_0\cup D_3\cup D_6))_r$. The involution is $\imath_{mc}: z_l\mapsto \bar{z}_r, z_r\mapsto \bar{z}_l$. The boundary of $\tilde{\Xr}_{mc}$ is $\partial\tilde{X}_{mc}=(S_0^-\cup S_3^-\cup S_6^-)_l\coprod (S_0^-\cup S_3^-\cup S_6^-)_r$. We parametrise the boundaries by
\begin{align}
&\delta_{mc}((e^{i\theta})_l)=\delta_{mc}((3+e^{i\theta})_l)=\delta_{mc}((6+e^{i\theta})_l)=-e^{-i\theta}\\
&\delta_{mc}((e^{i\theta})_r)=\delta_{mc}((3+e^{i\theta})_r)=\delta_{mc}((6+e^{i\theta})_r)=-e^{-i\theta}
\end{align}
 The ordering of the boundary components is given by
$$ \text{ord}:  ((S^-_0)_l, (S^-_0)_r, (S^-_3)_l,(S^-_3)_r, (S^-_6)_l,(S^-_6)_r)\mapsto(1,2,3,4,5,6)$$

The extended double $\hat{X}_{mc}$ is $\overline{\mathbb{C}}_l\cup\overline{\mathbb{C}}_r$ with arcs at $(0)_l,(3)_l$ labelled by $(B_l,+)$, an arc at $(6)_l$ labelled by $(B_l,-)$, arcs at $(0)_r,(3)_r$ labelled by $(B_r,+)$ and an arc at $(6)_r$ labelled by $(B_r,-)$.

\begin{table}[bt]
$$
\begin{array}{llll}
\text{R14} \etb:
f_{mc} \cir (\id_{\Bcl} \otimes_\Cb f_{\eta c}) = f_{pc}
\qquad \etb
\text{R15} \etb:
(\id_{\Bcl} \otimes_\Cb f_{\eps c}) \cir f_{\Delta c}  = f_{pc}
\enl
\text{R16} & \multicolumn{3}{l}{\hspace*{-.5em}: \displaystyle
f_{mc} \cir ( f_{mc} \otimes_\Cb \id_{\Bcl} )
= f_{mc} \cir ( \id_{\Bcl} \otimes_\Cb f_{mc} ) }
\enl
\text{R17} & \multicolumn{3}{l}{\hspace*{-.5em}: \displaystyle
( f_{\Delta c} \otimes_\Cb \id_{\Bcl} ) \cir f_{\Delta c}
=  ( \id_{\Bcl} \otimes_\Cb f_{\Delta c} )
\cir f_{\Delta c}}
\enl
\text{R18} & \multicolumn{3}{l}{\hspace*{-.5em}: \displaystyle
(f_{mc} \otimes_\Cb \id_{\Bcl}) \cir (\id_{\Bcl} \otimes_\Cb f_{\Delta c})
= f_{\Delta c} \cir f_{mc}}
\enl
\text{R19} & \multicolumn{3}{l}{\hspace*{-.5em}: \displaystyle
(\id_{\Bcl} \otimes_\Cb f_{mc}) \cir (f_{\Delta c} \otimes_\Cb  \id_{\Bcl})
= f_{\Delta c} \cir f_{mc}}
\enl
\text{R20} \etb:
f_{mc} \cir c_{\Bcl,\Bcl} = f_{mc}
\etb
\text{R21} \etb:
f_{pc} \cir f_{mc} = f_{mc}
\enl
\text{R22} \etb:
f_{pc} \cir f_{\eta c} = f_{\eta c}
\etb
\text{R23} \etb:
f_{\Delta c} \cir f_{pc} = f_{\Delta c}
\enl
\text{R24} \etb:
f_{\eps c} \cir f_{pc} = f_{\eps c}
\etb
\text{R25} \etb:
f_{pc} \cir \theta_{\Bcl} = f_{pc}
\end{array}
$$
\caption{The relations of $f_\alpha$ obtained from the relations R14-R25.}
\label{tab:glue-closed}
\end{table}

\subsubsection*{Relation R14}

The condition we need to impose is $C_{\Yr^d_{\Rr14,l}}= C_{\Yr^d_{\Rr14,r}}$.
In this case $\Yr_{\Rr14,r}=\Xr_{pc}$, so that for $\varpi_{\Rr14,r} : \Xr_{pc} \rightarrow \Yr_{\Rr14,r}$ we can just take the identity. The extended cobordism  $\Mr_{\varpi_{\Rr14,r}}$ associated to $\varpi_{\Rr14,r}$ is just the cylinder over $\Xhat_{pc}$ and so $C_{\Yr^d_{\Rr14,r}}1 = \Psi_{pc}(f_{pc})$, or pictorially,
\be
  C_{\Yr^d_{\Rr14,r}} = \sum_{\alpha} \tftC\Bigg(
  \raisebox{-65pt}{\begin{picture}(140,130)
  \put(7,10){\scalebox{0.40}{\includegraphics{pic-M-p_o}}}
  \put(7,10){
     \setlength{\unitlength}{.40pt}\put(-135,-279){
     \put(260,430)   {\tiny$ f_{pc,l}^\alpha $}
     \put(304,462)   {\tiny$ B_l $}
     \put(265,370)   {\tiny$ B_l $}
     }\setlength{\unitlength}{1pt}}
  \end{picture}}
  \hspace{-1em}\sqcup
  \raisebox{-65pt}{\begin{picture}(132,130)
  \put(7,10){\scalebox{0.40}{\includegraphics{pic-M-p_o}}}
  \put(7,10){
     \setlength{\unitlength}{.40pt}\put(-135,-279){
     \put(260,430)   {\tiny$ f_{pc,r}^\alpha $}
     \put(304,462)   {\tiny$ B_r $}
     \put(265,370)   {\tiny$ B_r $}
     }\setlength{\unitlength}{1pt}}
  \end{picture}}
  \Bigg)
\ee
where we expand $f_{pc}$ as $f_{pc}= \sum_\alpha f_{pc,l}^\alpha \otimes f_{pc,r}^\alpha$
with $f_{pc,l/r}^\alpha \in \Hom_\Cc(B_{l/r},B_{l/r})$. 

\medskip
On the other hand, $\varpi_{\Rr14,l} : \Xr_{\eta c} \otimes \Xr_{mc} \rightarrow \Yr_{\Rr14}$ is a combination of sewing and homeomorphism. We first expand $f_{\eta c}$ and $f_{mc}$ as follows:
\be
  f_{\eta c} = \sum_\alpha f_{\eta c,l}^\alpha \otimes_\Cb f_{\eta c,r}^\alpha
  \qquad \text{and} \qquad
  f_{mc}
  = \sum_\alpha f_{mc,l}^\alpha \otimes_\Cb f_{mc,r}^\alpha ~,
\labl{eq:mcl-decomp}
where $f_{\eta c,l/r}^\alpha \In \Hom_\Cc(\one,B_{l/r})$ and
$f_{mc,l/r}^\alpha \In \Hom_\Cc(B_{l/r} \otimes_\Cb B_{l/r},B_{l/r})$. 
The extended cobordism $\Mr_{\varpi_{\Rr 14, l}}$ associated to $\varpi_{\Rr 14, l}$ is given by $\Mr_l = \widehat{\Yr} \times [0,1] / \sim$, where$\Yr = \Xr_{\eta c} \otimes \Xr_{mc}$ and the equivalence relation $\sim$ is
\begin{align}
((z)_l,1,1) \sim ((-\bar{z})_l,2,1) \text{\ for\ all\ } |z| \le 1, \nonumber  \\
((z)_r,1,1) \sim ((-\bar{z})_r,2,1) \text{\ for\ all\ } |z| \le 1, \nonumber
\end{align}
and graphically by
\be
  \Mr_{\varpi_{\Rr 14, l}} =
\raisebox{-65pt}{\begin{picture}(135,135)
  \put(0,0){\scalebox{0.35}{\includegraphics{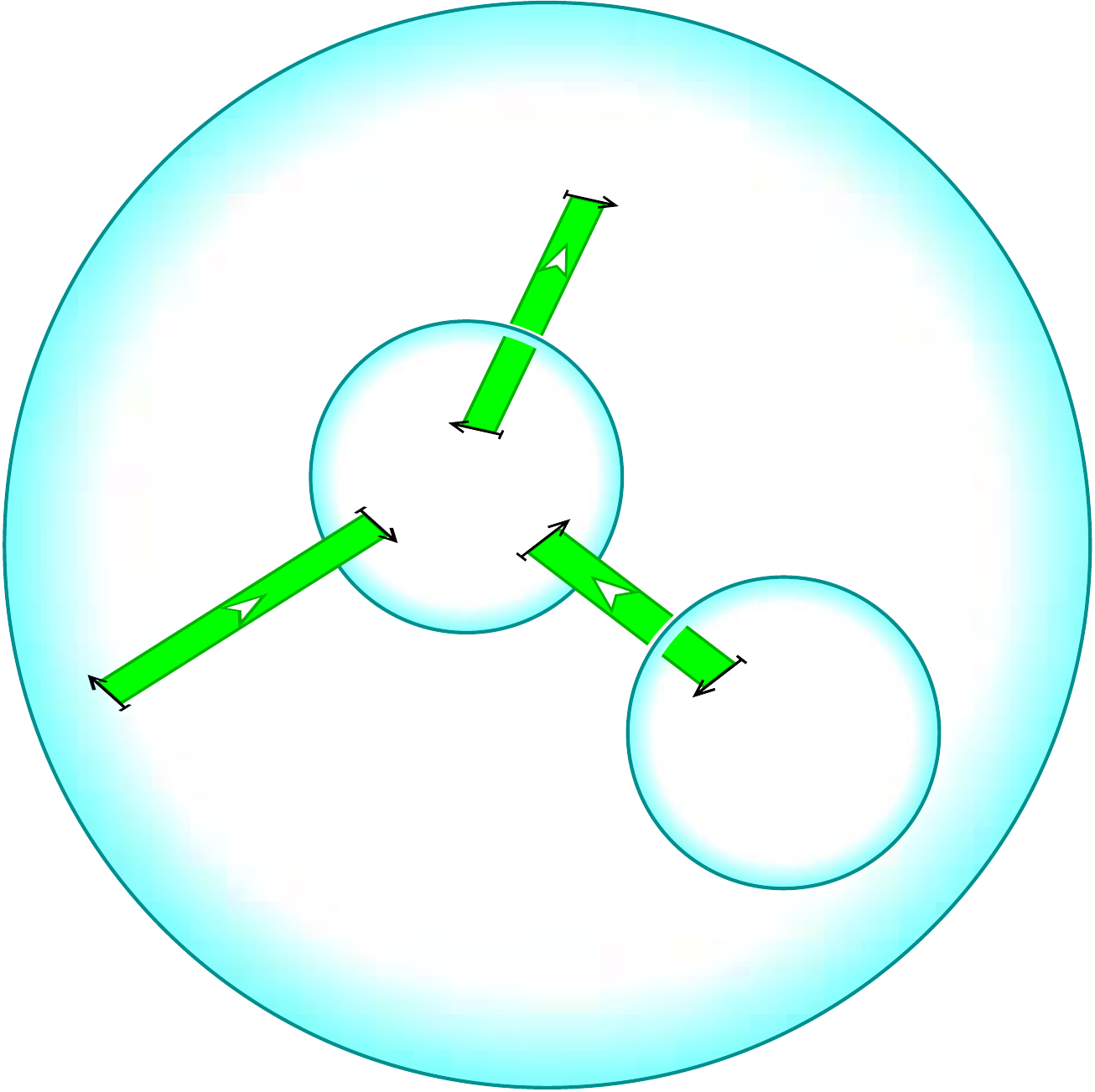}}}
  \put(0,0){
     \setlength{\unitlength}{.35pt}\put(-115,-296){
     \put(155,472)   {\tiny$ B_l $}
     \put(316,567)   {\tiny$ B_l $}
     \put(330,478)   {\tiny$ B_l $}
     }\setlength{\unitlength}{1pt}}
  \end{picture}}
  \sqcup
\raisebox{-65pt}{\begin{picture}(135,135)
  \put(0,0){\scalebox{0.35}{\includegraphics{pic-R1-Ml}}}
  \put(0,0){
     \setlength{\unitlength}{.35pt}\put(-115,-296){
     \put(155,472)   {\tiny$ B_r $}
     \put(316,567)   {\tiny$ B_r $}
     \put(330,478)   {\tiny$ B_r $}
     }\setlength{\unitlength}{1pt}}
  \end{picture}}
\, .
\ee
By the definition $C_{\Yr^d_{\Rr14,l}}1=\Bl(\varpi_{\Rr 14,l})\circ (\Psi_{\eta c}(f_{\eta c}) \otimes_\Cb \Psi_{mc}(f_{mc}))$, we obtain
$$
C_{\Yr^d_{\Rr14,l}} = \sum_{\alpha,\beta}
  \tftC\Bigg(\,
    \raisebox{-65pt}{\begin{picture}(135,135)
  \put(0,0){\scalebox{0.35}{\includegraphics{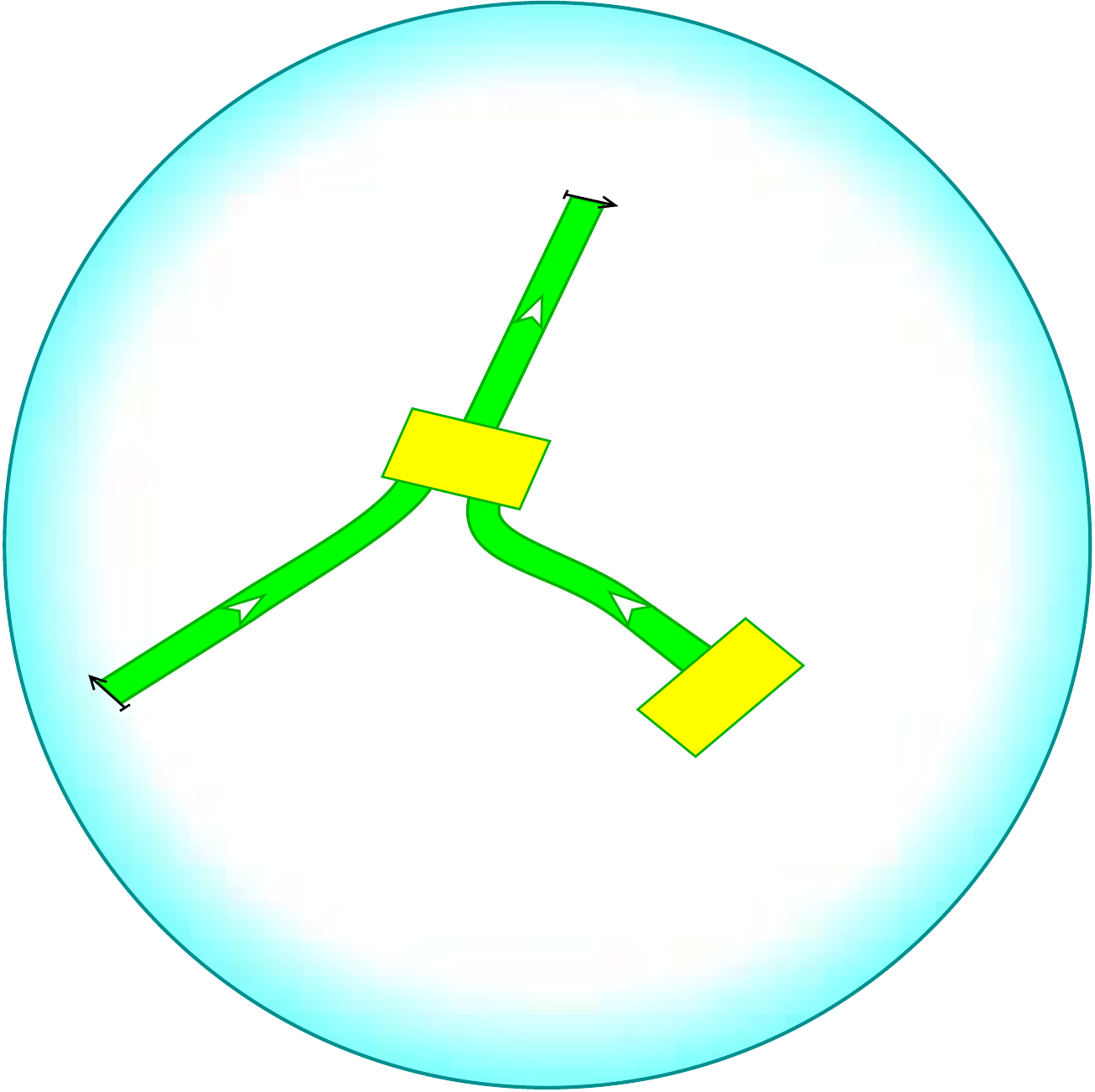}}}
  \put(0,0){
     \setlength{\unitlength}{.35pt}\put(-115,-296){
     \put(155,472)   {\tiny$ B_l $}
     \put(316,567)   {\tiny$ B_l $}
     \put(330,478)   {\tiny$ B_l $}
     \put(346,422)   {\tiny \begin{turn}{+43}$ f_{\eta c,l}^\beta $\end{turn}}
     \put(256,514)   {\tiny \begin{turn}{-12}$ f_{mc,l}^\alpha $\end{turn}}
     }\setlength{\unitlength}{1pt}}
  \end{picture}}
  \sqcup
    \raisebox{-65pt}{\begin{picture}(135,135)
  \put(0,0){\scalebox{0.35}{\includegraphics{pic-R1-Bl}}}
  \put(0,0){
     \setlength{\unitlength}{.35pt}\put(-115,-296){
     \put(155,472)   {\tiny$ B_r $}
     \put(316,567)   {\tiny$ B_r $}
     \put(330,478)   {\tiny$ B_r $}
     \put(346,422)   {\tiny \begin{turn}{+43}$ f_{\eta c,r}^\beta $\end{turn}}
     \put(256,514)   {\tiny \begin{turn}{-12}$ f_{mc,r}^\alpha $\end{turn}}
     }\setlength{\unitlength}{1pt}}
  \end{picture}}
  \Bigg)~.
$$
The condition $C_{\Yr^d_{\Rr14,l}} = C_{\Yr^d_{\Rr14,r}}$ is then equivalent to
the following identity of morphisms in $\Hom_{\CxC}(\BlxBr,\BlxBr)$,
\be
  \sum_{\alpha,\beta}
  \big[ f_{mc,l}^\alpha \circ (\id_{B_l} \otimes_\Cb f_{\eta c,l}^\beta) \big] \otimes_\Cb
  \big[ f_{mc,r}^\alpha \circ (\id_{B_r} \otimes_\Cb f_{\eta c,r}^\beta) \big] =
  \sum_{\gamma} f_{pc,l}^\gamma \otimes_\Cb f_{pc,r}^\gamma ~.
\ee
Carrying out the sums, we obtain the identity of $f_\alpha$ determined by R14 as stated in Table \ref{tab:glue-closed}.

\subsubsection*{Relation R20}

In this case, the extended cobordisms $\Mr_{\varpi_{\Rr20,l/r}}$ associated to $\varpi_{\Rr20,l/r}$ are defined by the same three manifold $\Xr_{mc}\times [0,1]$ but with the homeomorphisms $h\big|_{\partial_\pm \Mr_{\Rr20,r}}$ where $h\big|_{\partial_- \Mr_{\Rr20,l}}$ is given by identity maps on $\Xr_{mc}$
and $h\big|_{\partial_+ \Mr_l}$ is defined pictorially as follows:
$$
\figbox{0.3}{relation20l}\ \  \rightsquigarrow\ \  \figbox{0.3}{relation20r}
$$
Therefore, the extended cobordism $\Mr_{\varpi_{\Rr20,l}}$ can be depicted as follows:
$$
\Mr_{\varpi_{\Rr20,l}}=\figbox{0.45}{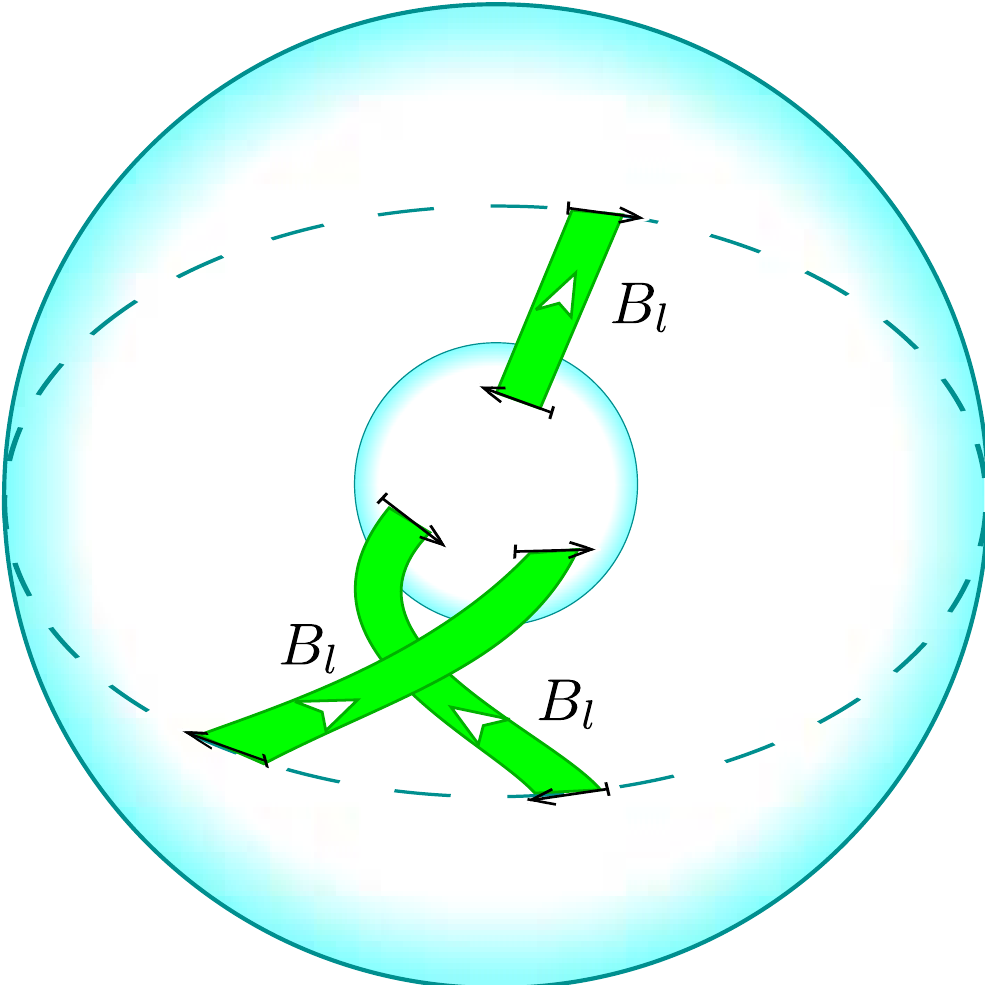}\sqcup \figbox{0.45}{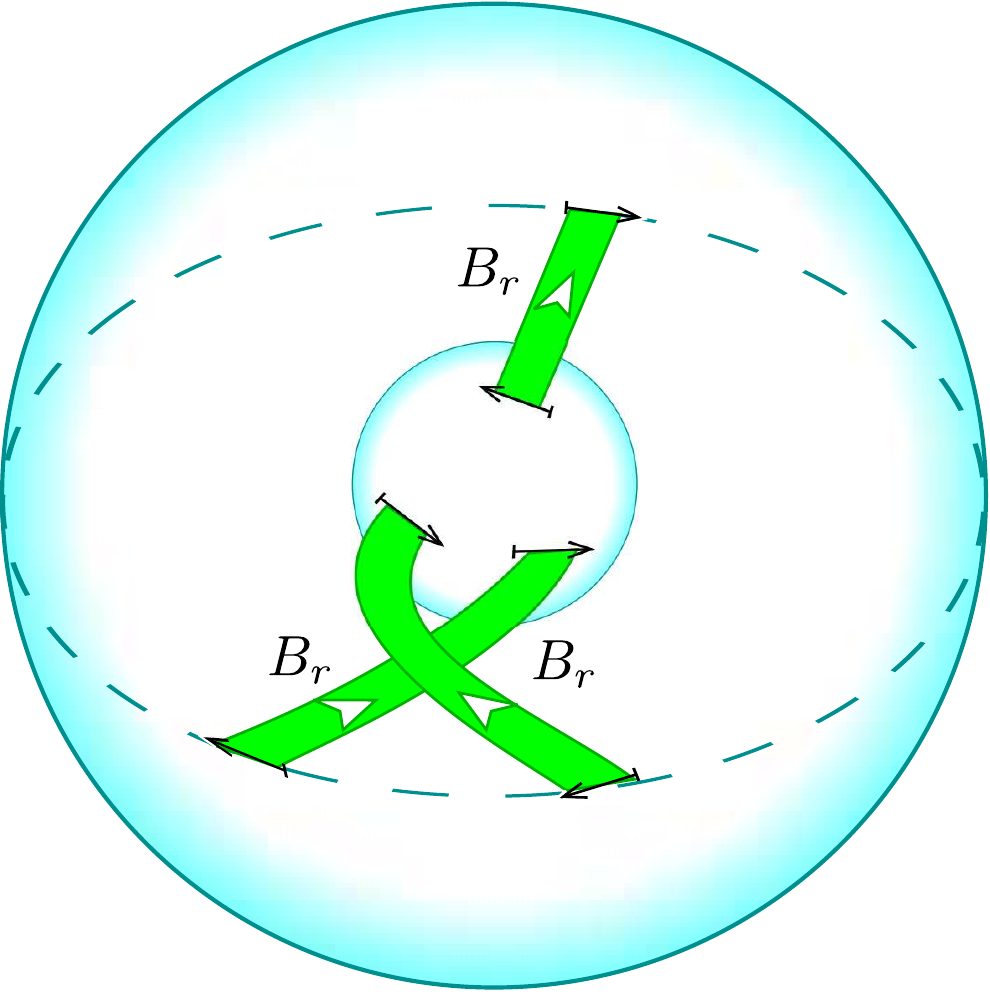}~.
$$
The opposite braiding arises because the involution on $\tilde{\Xr}_{mc}$ is given by complex conjugation, and so the action of the element of the mapping class group on the two connected components of $\tilde{\Xr}_{mc}$ is equally related by complex conjugation.
        
By $C_{\Yr_{\Rr20,l/r}^d}1=\Bl(\varpi_{\Rr20,l/r})\circ(\Psi_{mc}(f_{mc}))$, we obtain
$$
C_{\Yr_{\Rr20,l}^d}=\sum_{\alpha} \tftC\Bigg(\figbox{0.45}{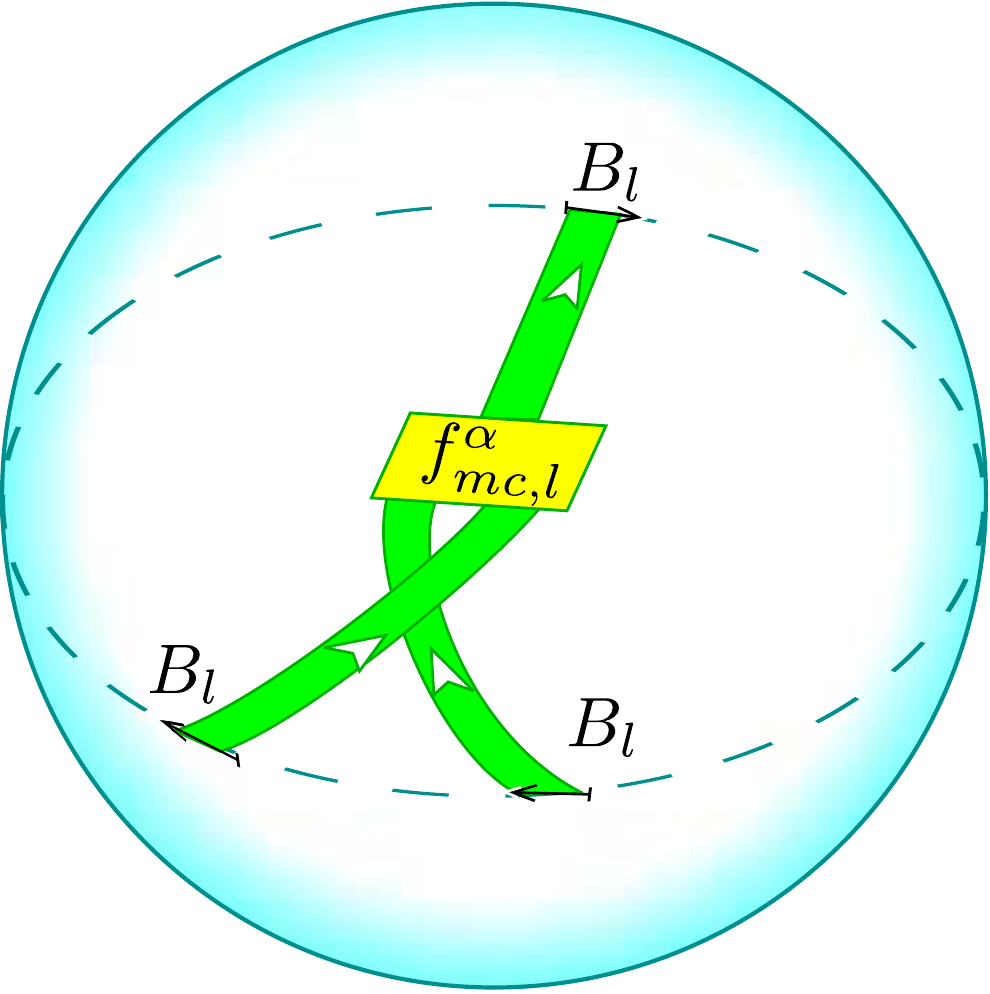}\sqcup\figbox{0.45}{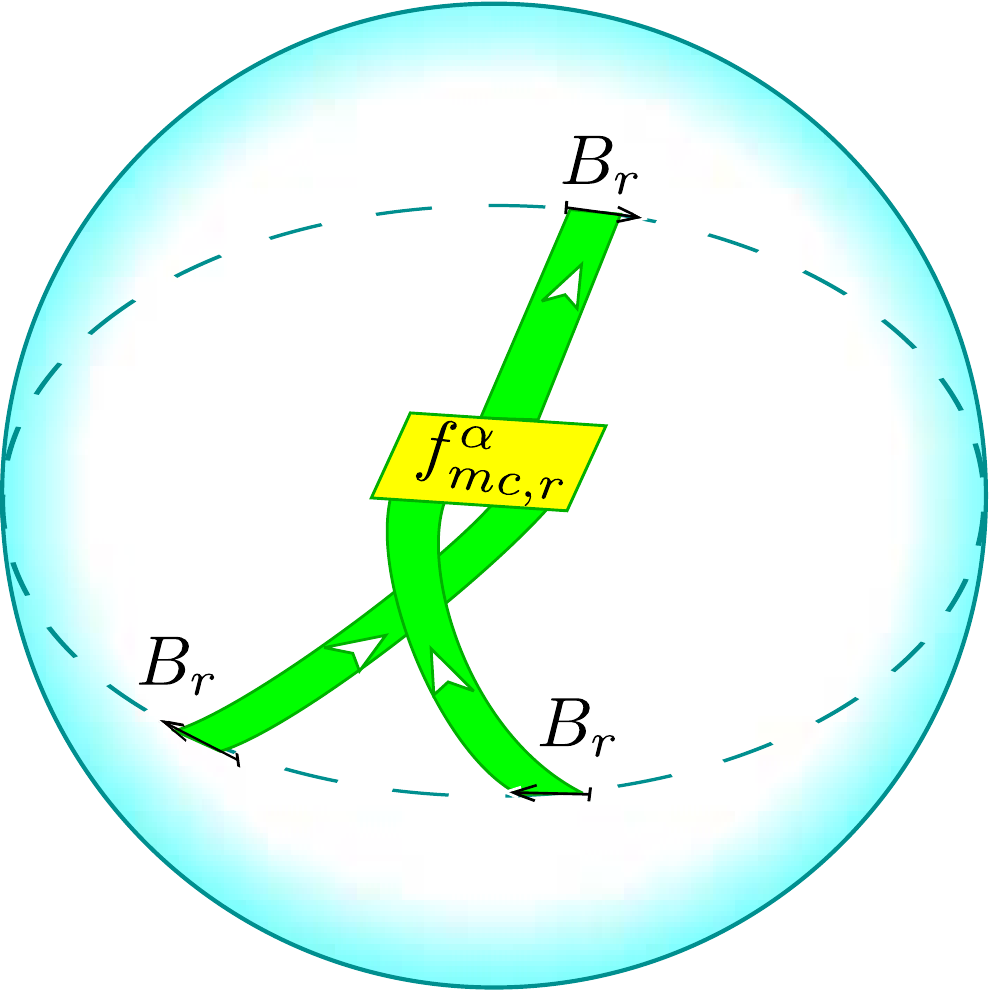}\Bigg)
$$
and
$$
 C_{\Yr_{\Rr20,r}^d}=\sum_{\alpha} \tftC\Bigg(\figbox{0.45}{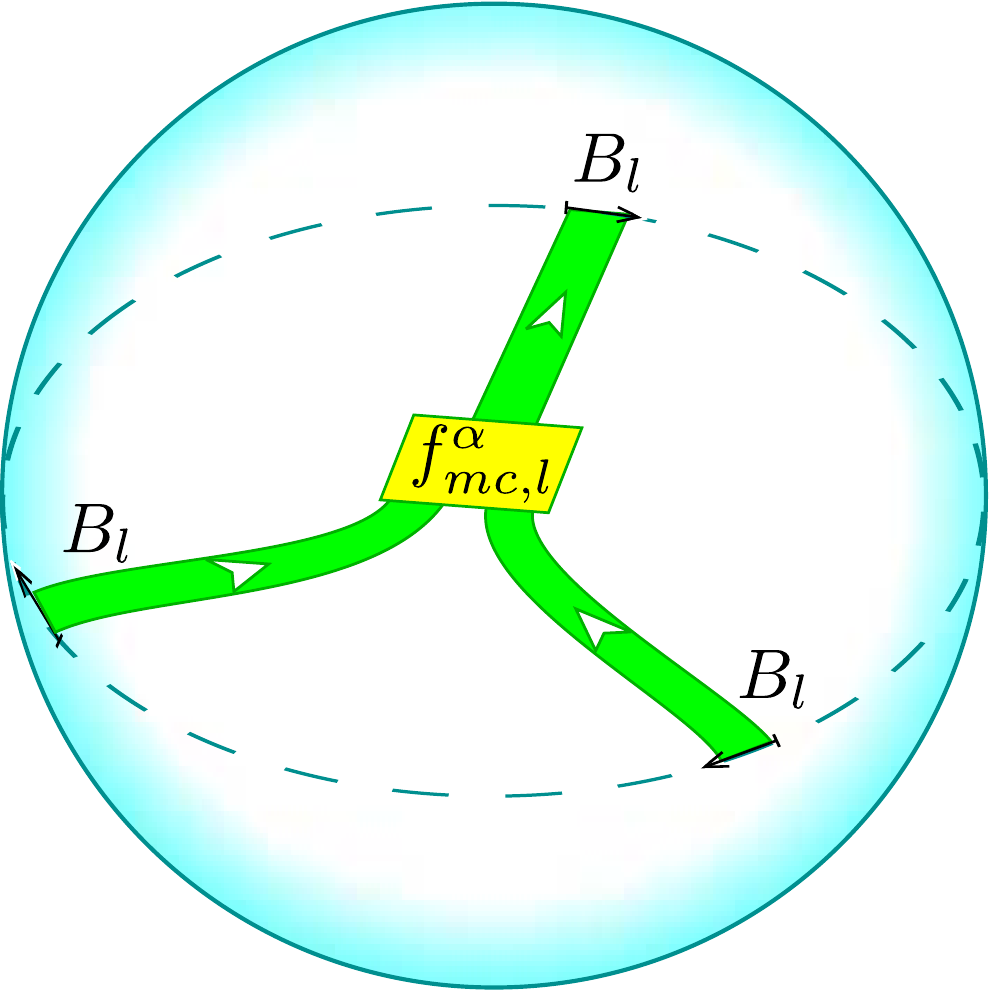}\sqcup\figbox{0.45}{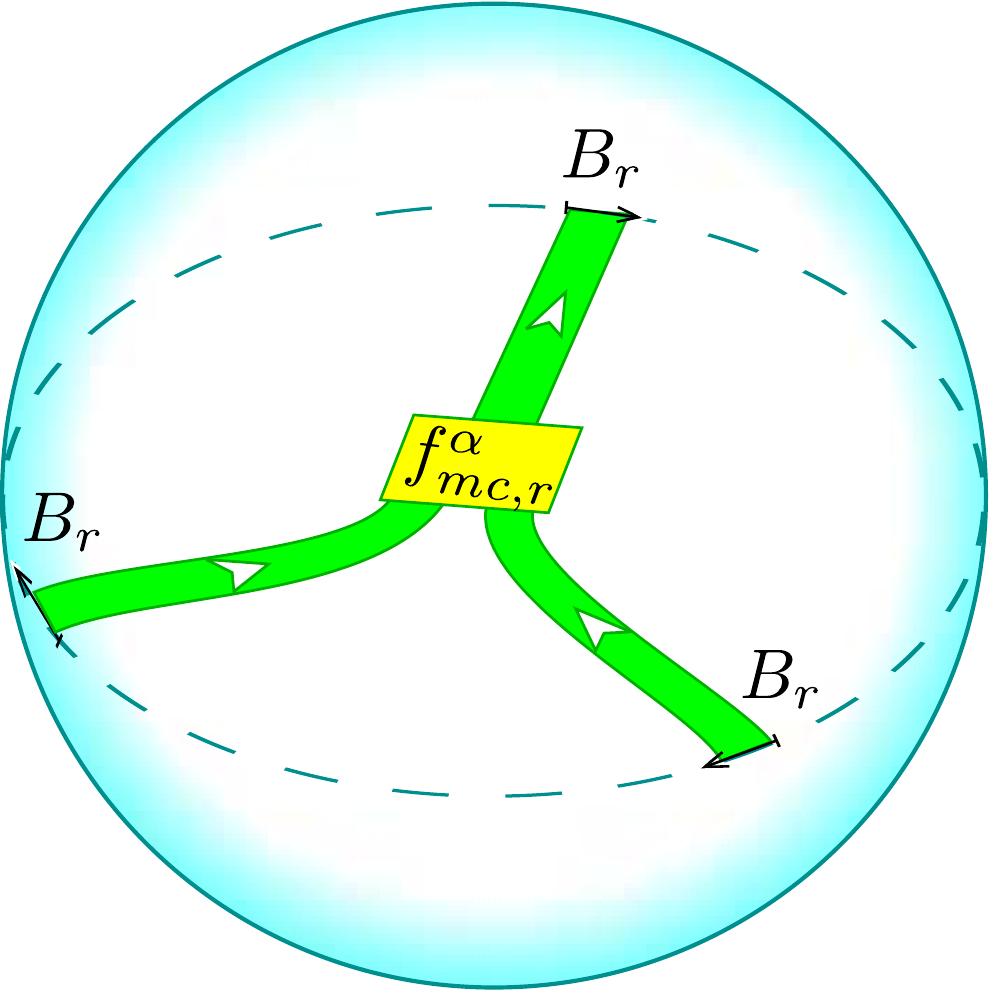}\Bigg).
$$
Hence $C_{\Yr_{\Rr20,l}^d}=C_{\Yr_{\Rr20,r}^d}$ if and only if $f_{mc}\circ c_{B_{cl},B_{cl}}=f_{mc} $

\subsection{Relations R26-R30}  \label{sec:R26-R30}

The complete list of the relations of $f_\alpha$ obtained from the relations R26-R30 is given in Table \ref{tab:glue-mixed}. We will derive those associated to R28 and R29.

\subsubsection*{The world sheet $\Xr_{\iota}$}
We take $\tilde{\Xr}_{\iota}=\overline{\mathbb{C}}-(D_{3i}\cup D_{-3i}\cup D_3)$.
The involution $\imath_{\iota}$ is given by complex conjugation.
The boundary of $\tilde{\Xr}_{\iota}$ is $\partial\tilde{X}_{\iota}=S_{3i}^-\cup S_{-3i}^-\cup S_3^-$.
We parametrise the boundary by
$$\delta_{\iota}(3i+e^{i\theta})=\delta_{\iota}(-3i+e^{i\theta})=\delta_{\iota}(3+e^{i\theta})=-e^{-i\theta}$$
We give the ordering of the boundary components of
$\Xr_{\iota}$ as follows:
$$ \text{ord}: (S_3^-,S_{3i}^-, S_{-3i}^-)\mapsto (1,2,3)$$

The extended double $\hat{\Xr}_{\iota}$ is $\overline{\mathbb{C}}$ with an arc at $3i$ labelled by $(B_l,+)$, an
arc at $-3i$ labelled by $(B_r,+)$ and an arc at $3$ labelled by $(B_{op},-)$.

\subsubsection*{Relation R28}

The extended cobordism $\Mr_{\varpi_{\Rr28,l}}$ associated to $\varpi_{\Rr28,l}$ is given by 
$\Mr_{\varpi_{\Rr28,l}}=\hat{\Yr}\times[0,1]/\sim$, where $\Yr=\Xr_\iota\otimes \Xr_{mo}$ and the equivalence relation $\sim$ is defined by $(3+z,1,1)\sim(-\bar{z},2,1)$ for all $|z|\leqslant 1$, and is given pictorially as follows:
$$
\Mr_{\varpi_{\Rr28,l}}=\figbox{0.8}{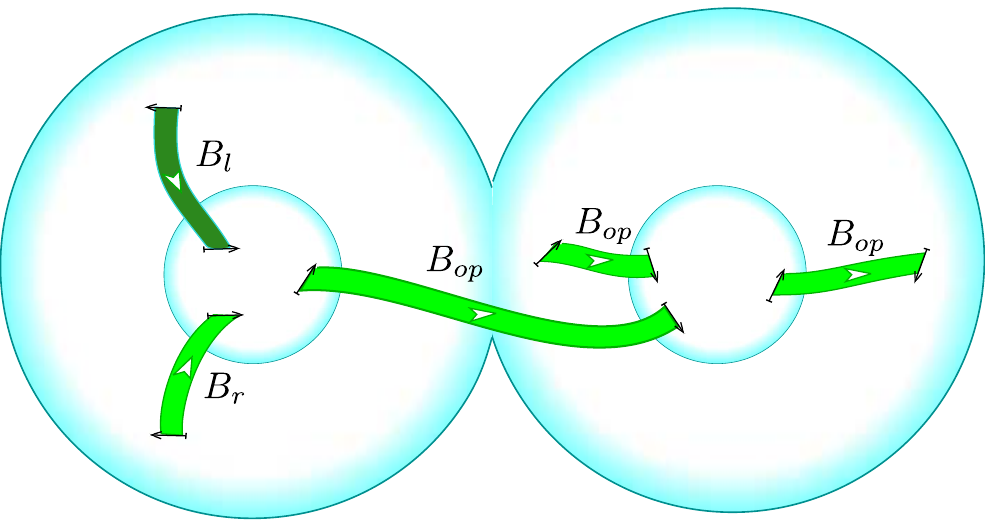}~.
$$
$\Mr_l$ is just the bordism from $\Xr_\iota\otimes \Xr_{mo}$ to $\Yr_{\Rr28}$ that we will apply to define $\Bl(\varpi_{\Rr28,l})$.
By the definition $C_{\Yr_{\Rr28,l}^d}1=\Bl(\varpi_{\Rr28,l})\circ(\Psi_{\iota}(f_{\iota})\otimes_\Cb \Psi_{mo}(f_{mo}))$, 
we obtain
$$
C_{\Yr_{\Rr28,l}^d}= \tftC\Bigg(\figbox{0.45}{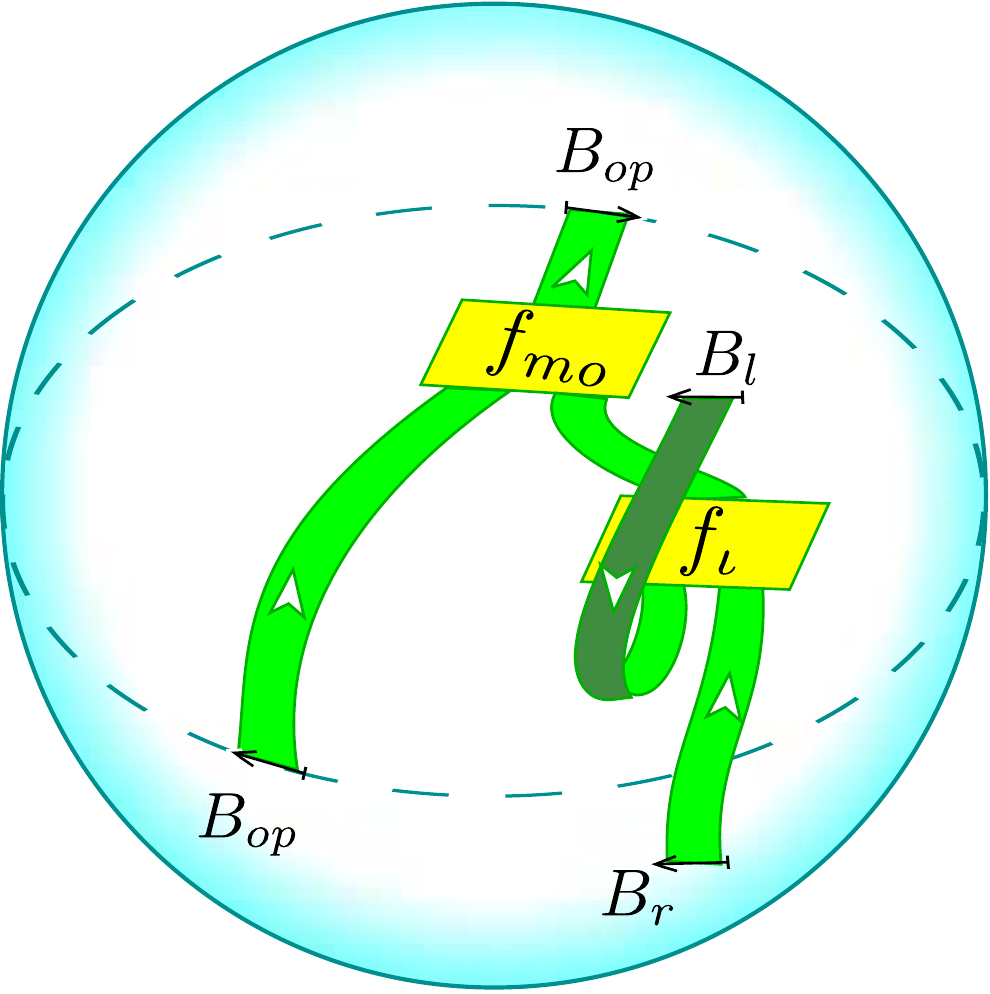}\Bigg)~.
$$

On the other hand, the extended cobordism $\Mr_{\varpi_{\Rr28,r}}$ associated to $\varpi_{\Rr28,r}$ is given by 
$\Mr_{\varpi_{\Rr28,r}}=\hat{\Yr}\times[0,1]/\sim$, where $\Yr=\Xr_\iota\otimes \Xr_{mo}$ and the equivalence relation $\sim$ is defined by $(3+z,1,1)\sim(3-\bar{z},2,1)$ for all $|z|\leqslant 1$, and is given pictorially as follows:
$$
\Mr_{\varpi_{\Rr28,r}}=\figbox{0.8}{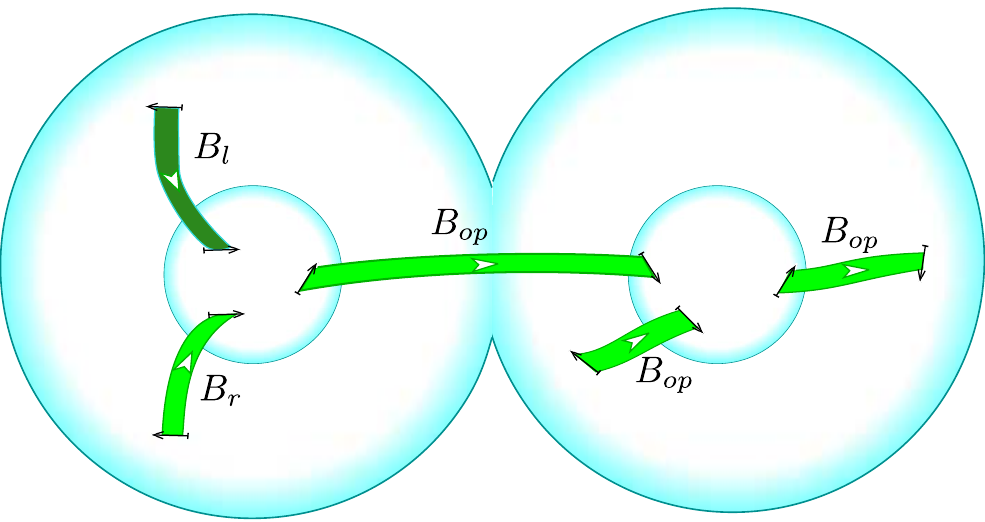}
$$
By the definition $C_{\Yr_{\Rr28,r}^d}1=\Bl(\varpi_{\Rr28,r})\circ(\Psi_{\iota}(f_{\iota})\otimes\Psi_{mo}(f_{mo}))$,
we obtain
$$
C_{\Yr_{\Rr28,r}^d} =\tftC\Bigg(\figbox{0.45}{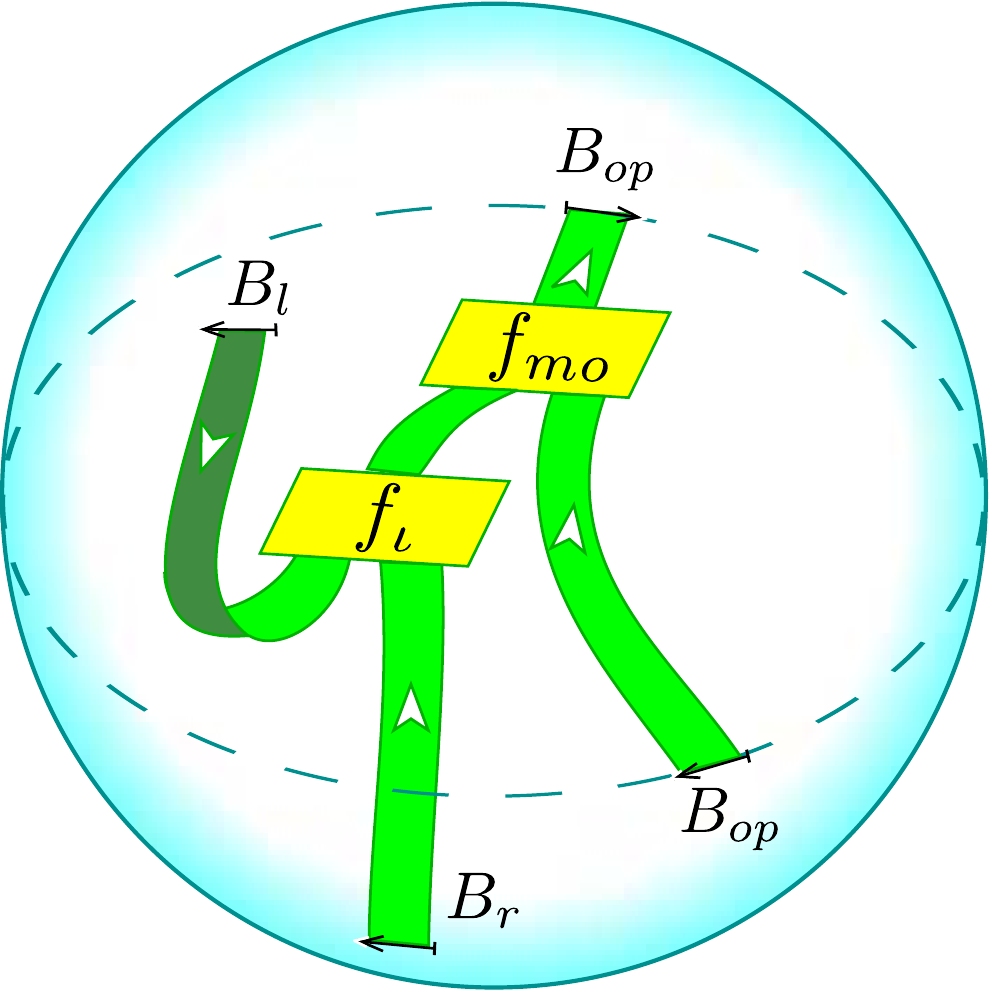}\Bigg)~.
$$

By moving the arcs labeled $B_l$ and $B_r$ to the equator of the sphere according to our convention, we get the following pictures:
$$
 C_{\Yr_{\Rr29,l}^d}=\tftC\Bigg(\figbox{0.45}{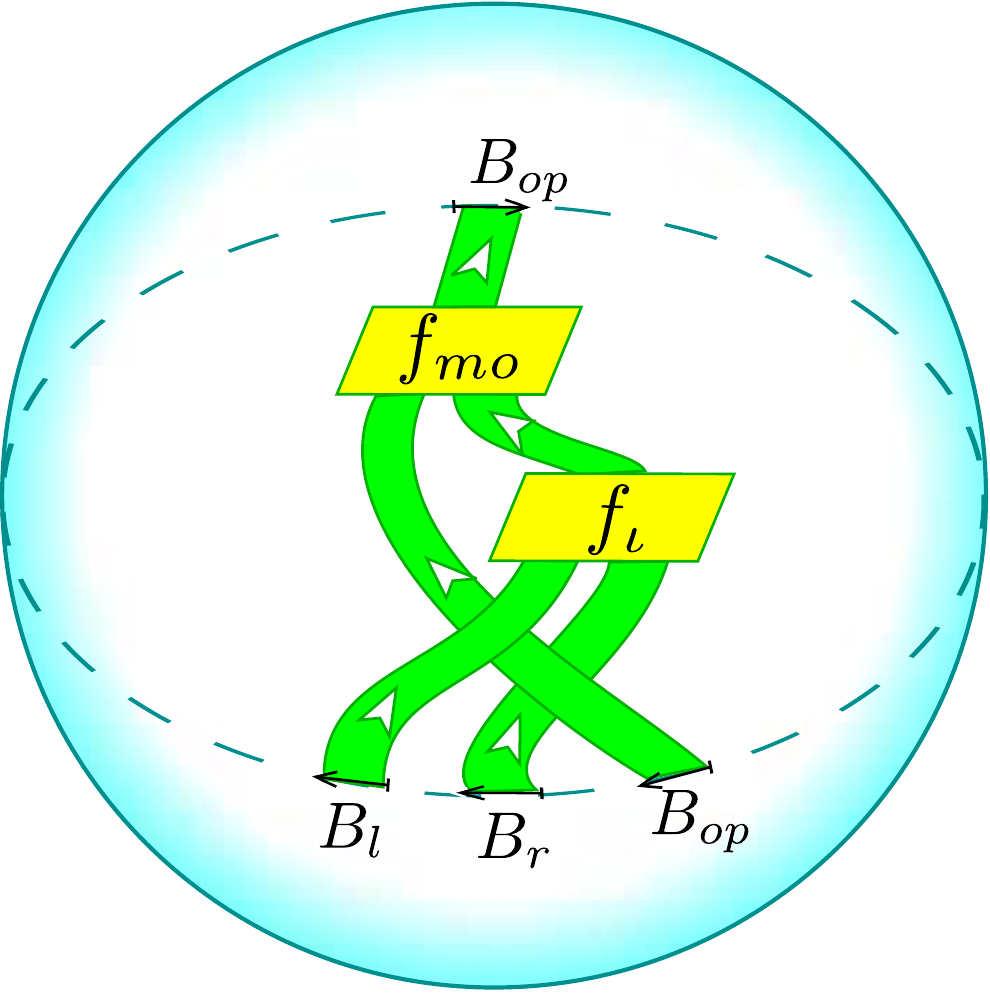}\Bigg),\,\,
 C_{\Yr_{\Rr28,r}^d}=\tftC\Bigg(\figbox{0.45}{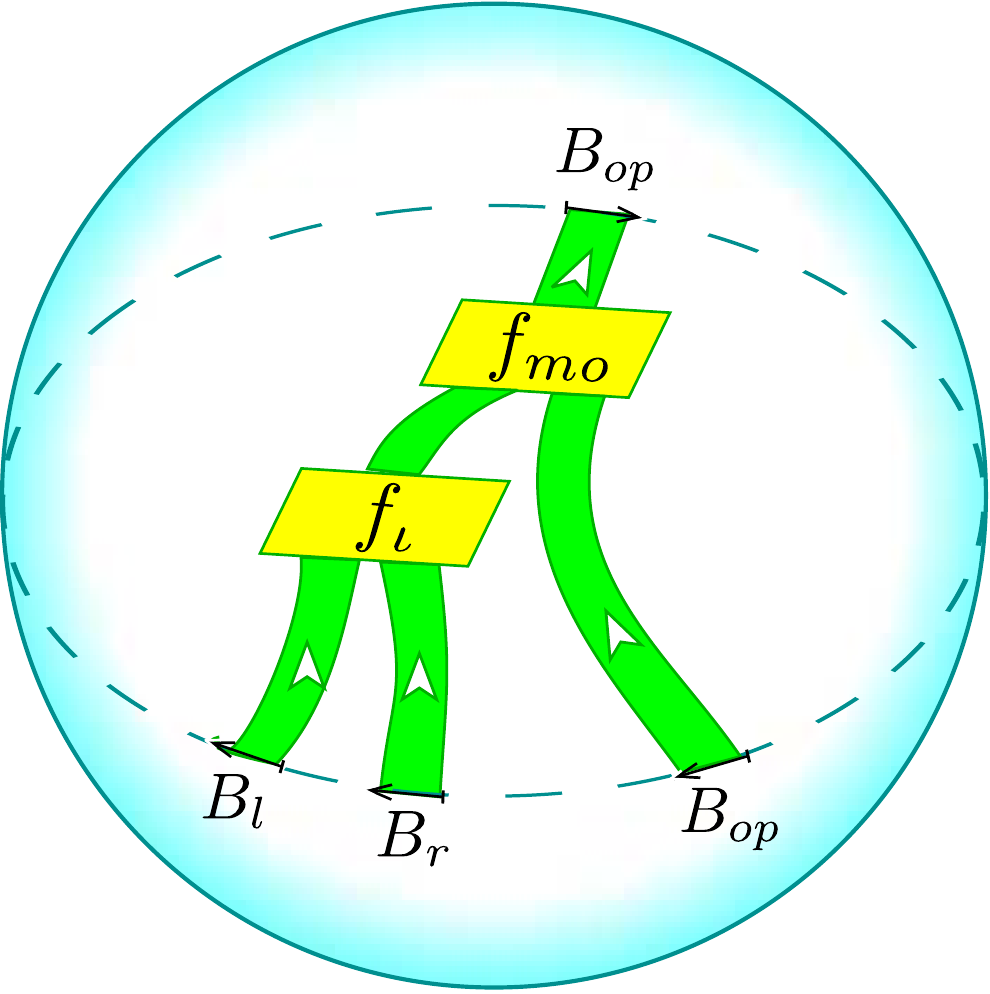}\Bigg).
$$
Hence, $C_{\Yr_{\Rr28,l}^d}=C_{\Yr_{\Rr28,r}^d}$ if and only if the following identity is satisfied:
\begin{equation}\label{center condition}
\figbox{0.3}{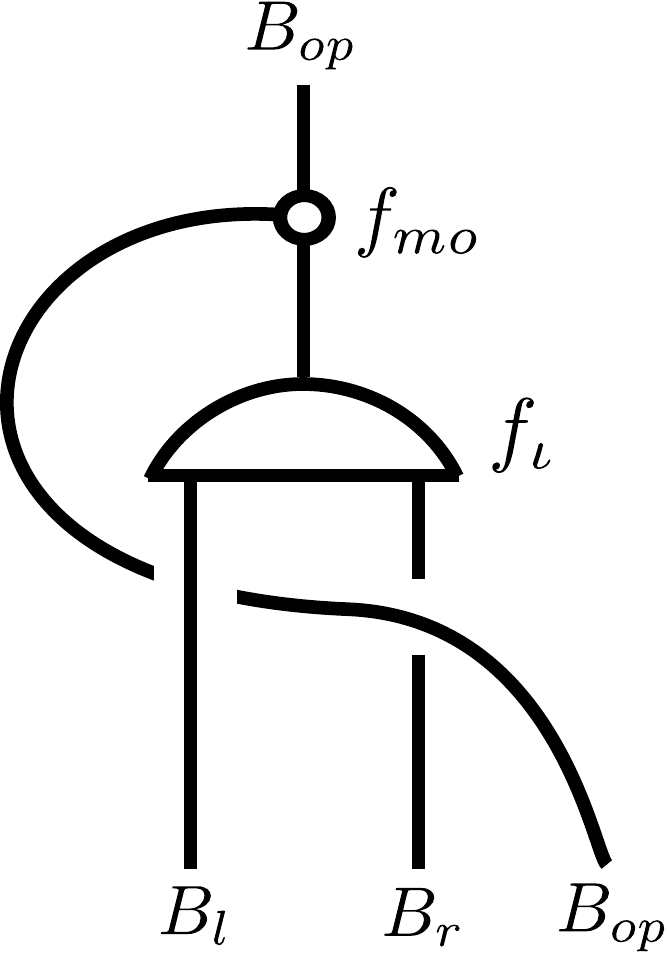}=\figbox{0.3}{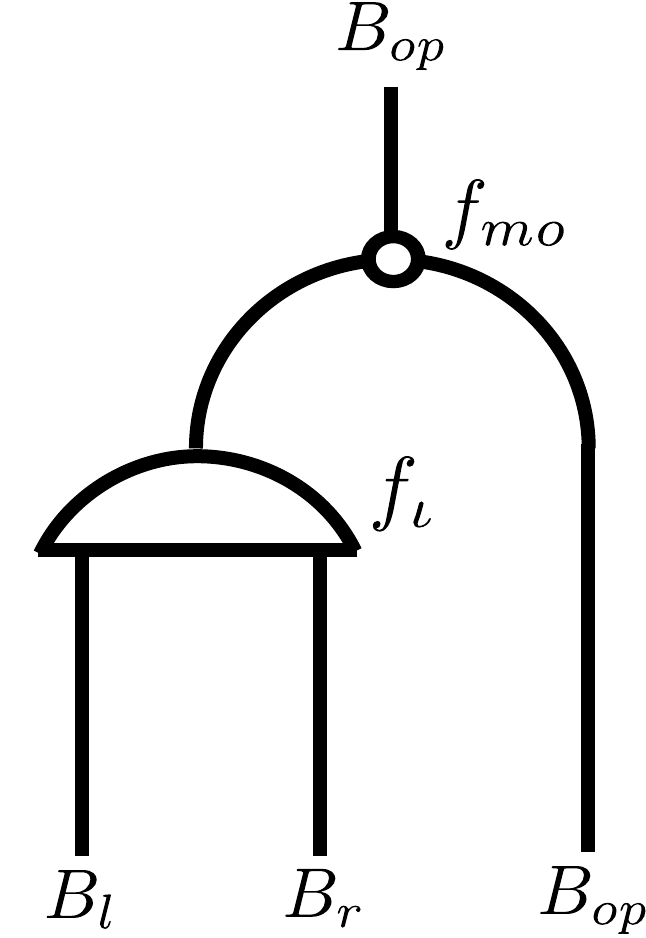} \quad .
\end{equation}

\subsubsection*{Relation R29}

The extended cobordism $\Mr_{\varpi_{\Rr29,l}}$ associated to $\varpi_{\Rr29,l}$ is given by 
$\Mr_{\varpi_{\Rr29,l}}=\hat{\Yr}_l\times[0,1]/\sim$, where
$\Yr_l=\Xr_{mc}\otimes \Xr_{\iota}$ and the equivalence relation $\sim$ is defined by
\begin{align}
((6+z)_l,1,1)\sim(3i-\bar{z},2,1), \text{\ for\ all\ }|z|\leqslant 1,  \nonumber \\
((6+z)_r,1,1)\sim(-3i-\bar{z},2,1), \text{\ for\ all\ }|z|\leqslant 1. \nonumber
\end{align}
The extended cobordism $\Mr_{\varpi_{\Rr29,r}}$ associated to $\varpi_{\Rr29,r}$ is given by 
$\Mr_{\varpi_{\Rr29,r}}=\hat{\Yr}_r\times[0,1]/\sim$, where
$\Yr_r=\Xr_{\iota}\otimes\Xr_{\iota}\otimes\Xr_{mo}$ and the equivalence relation $\sim$ is defined by
\begin{align}
(3+z,1,1)\sim(-\bar{z},3,1),\text{\ for\ all\ }|z|\leqslant 1, \nonumber \\
(3+z,2,1)\sim(3-\bar{z},3,1),\text{\ for\ all\ }|z|\leqslant 1. \nonumber
\end{align}
We are given two morphisms $\varpi_{\Rr29,l} : \Xr_{mc} \otimes \Xr_{\iota} \rightarrow \Yr_{\Rr29}$ and
 $\varpi_{\Rr29,r} : \Xr_\iota \otimes \Xr_\iota\otimes \Xr_{mo}  \rightarrow \Yr_{\Rr29}$. 
Pictorially, we have 
$$
  \Mr_ {\varpi_{\Rr29,l}}=\figbox{0.45}{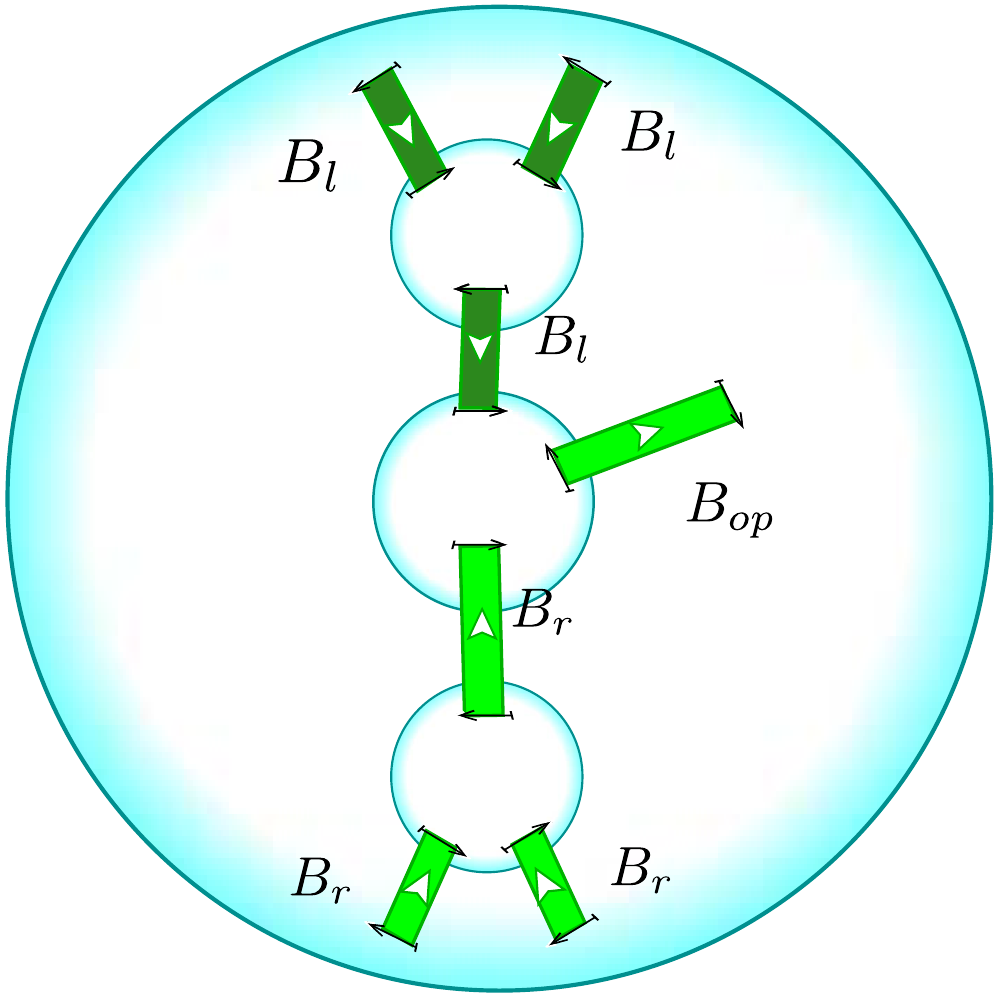}~,
 \quad\quad
  \Mr_{\varpi_{\Rr29,r}} =
\figbox{0.45}{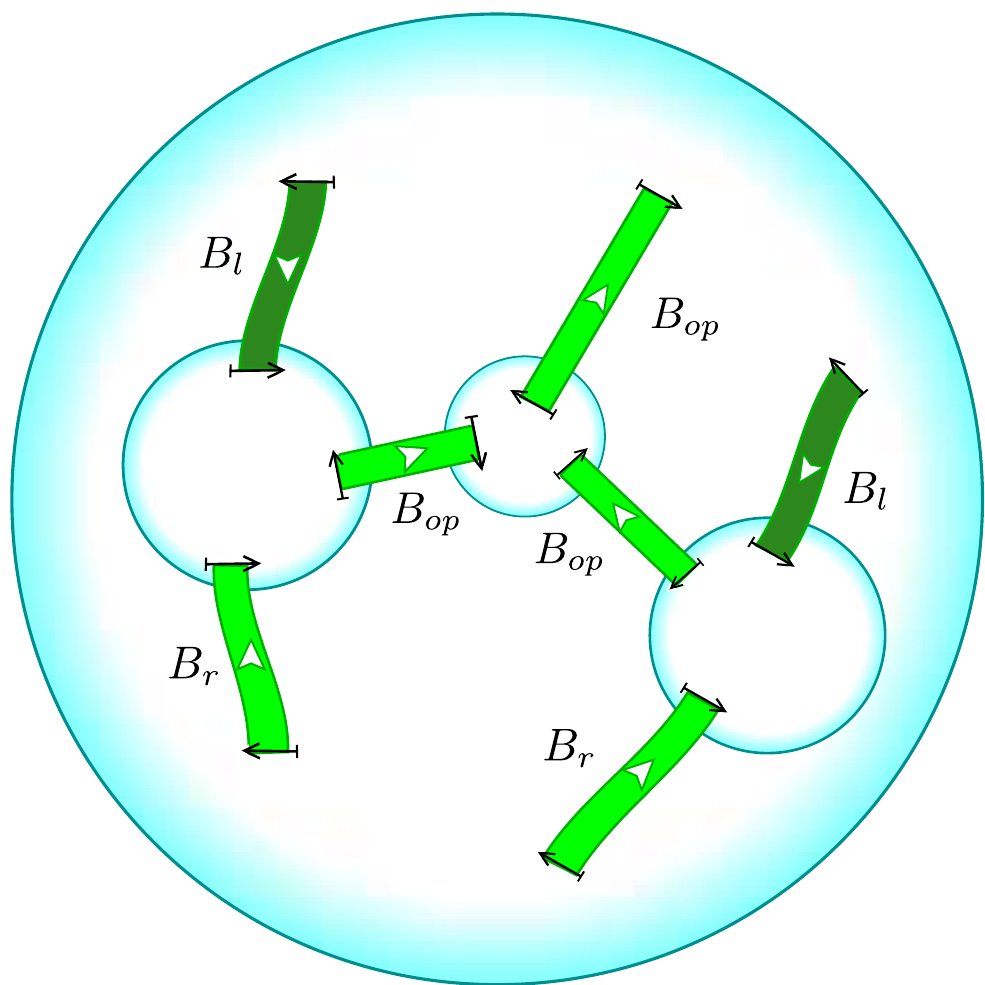}~.
$$
By the definition $C_{\Yr^d_{\Rr29,l}}1 = \Bl(\varpi_{\Rr29,l}) \circ \big( \Psi_{mc}(f_{mc}) \otimes
  \Psi_{\iota}(f_\iota) \big)$, we obtain
$$
C_{\Yr^d_{\Rr29,l}} = \sum_\alpha \tftC\Bigg(\figbox{0.45}{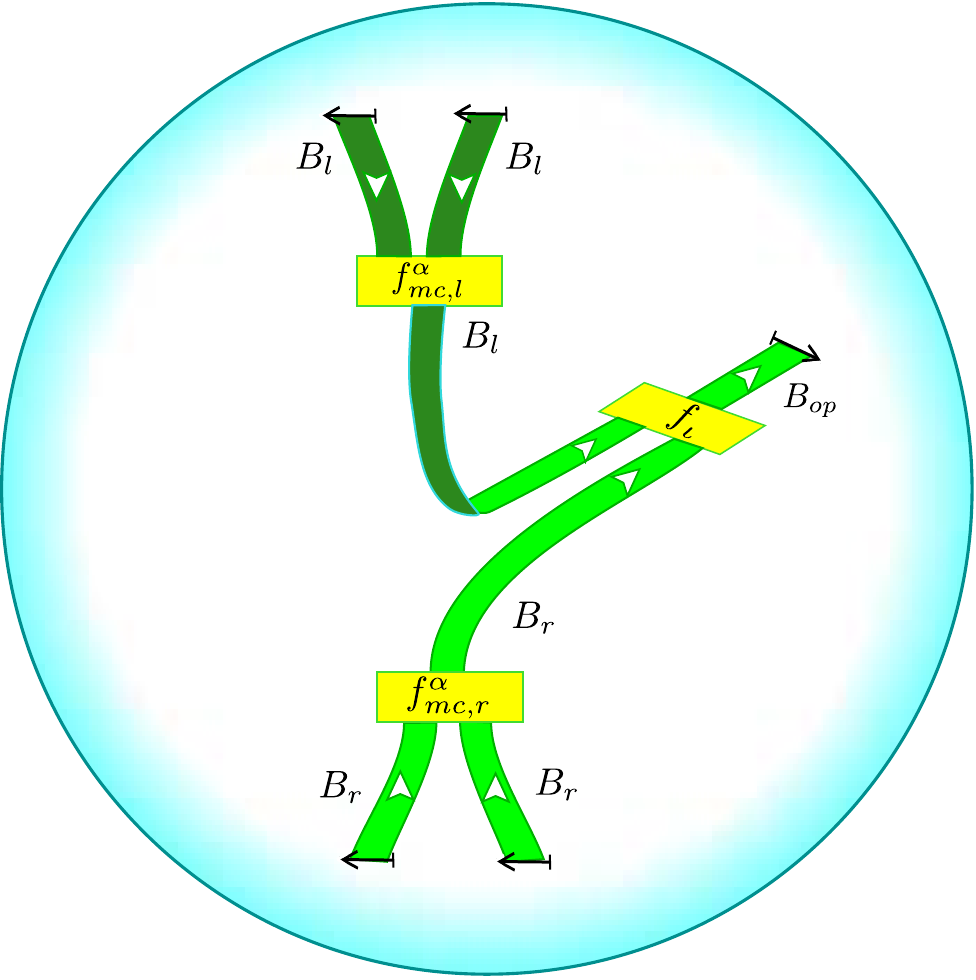}\Bigg) 
= \sum_\alpha \tftC\Bigg(\figbox{0.45}{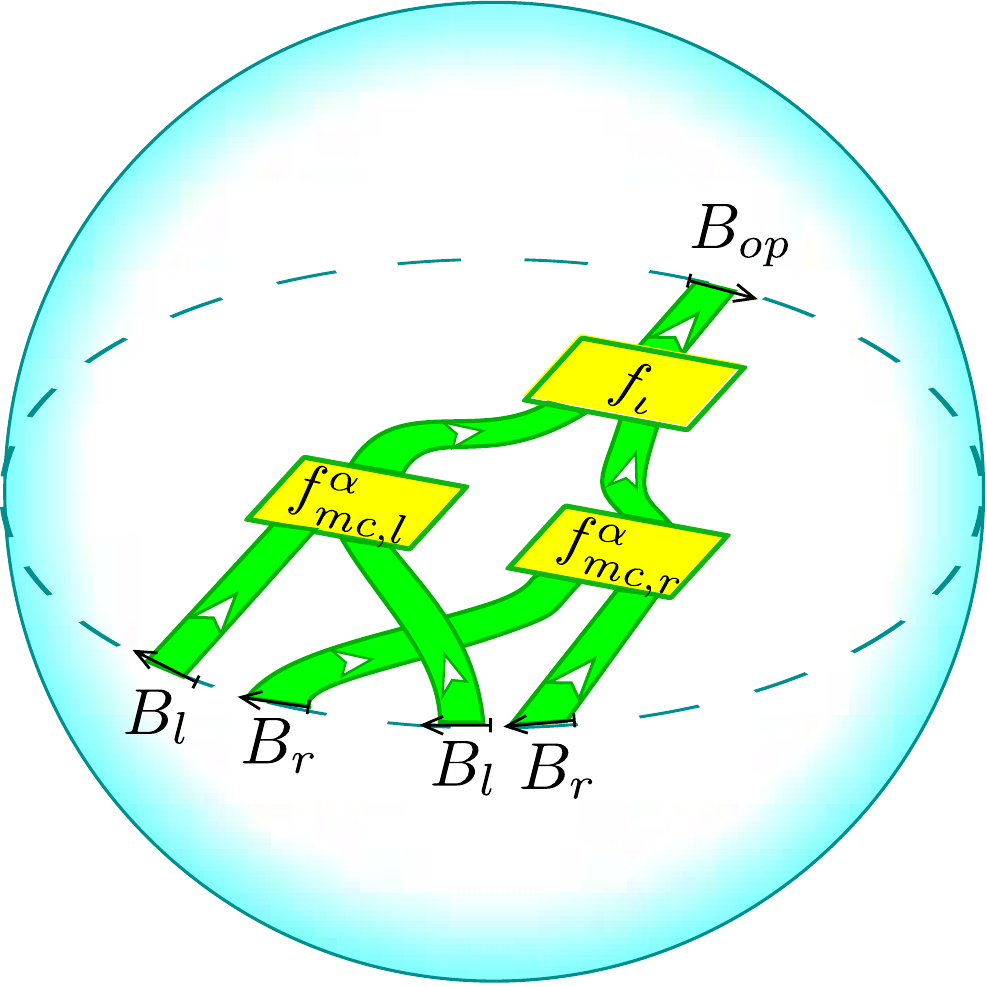}\Bigg)~.
$$
By the definition $C_{\Yr^d_{\Rr29,r}}1 = \Bl(\varpi_{\Rr29,r}) \circ \big( \Psi_{mo}(f_{mo}) \otimes_\Cb \Psi_{\iota}(f_\iota) \otimes_\Cb \Psi_{\iota}(f_\iota) \big)$, we obtain
$$
C_{\Yr^d_{\Rr29,r}} = \tftC\Bigg({\figbox{0.45}{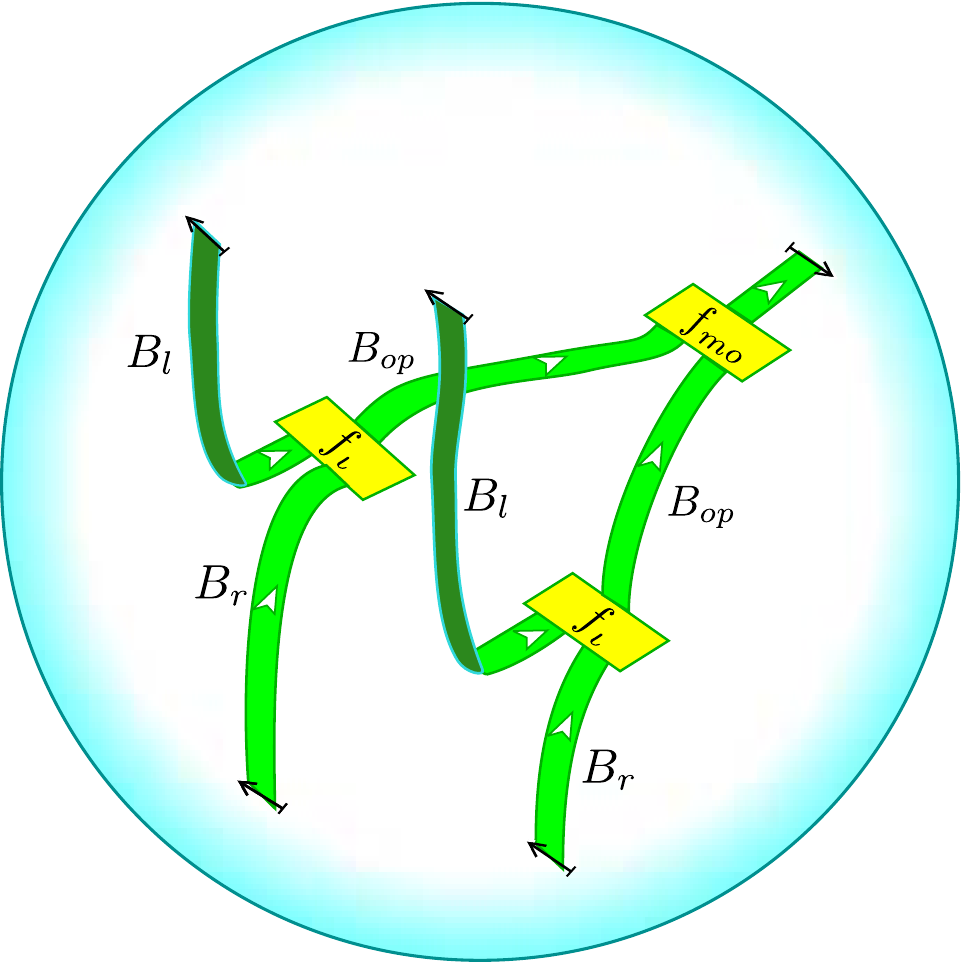}}\Bigg)
= \tftC\Bigg(\figbox{0.45}{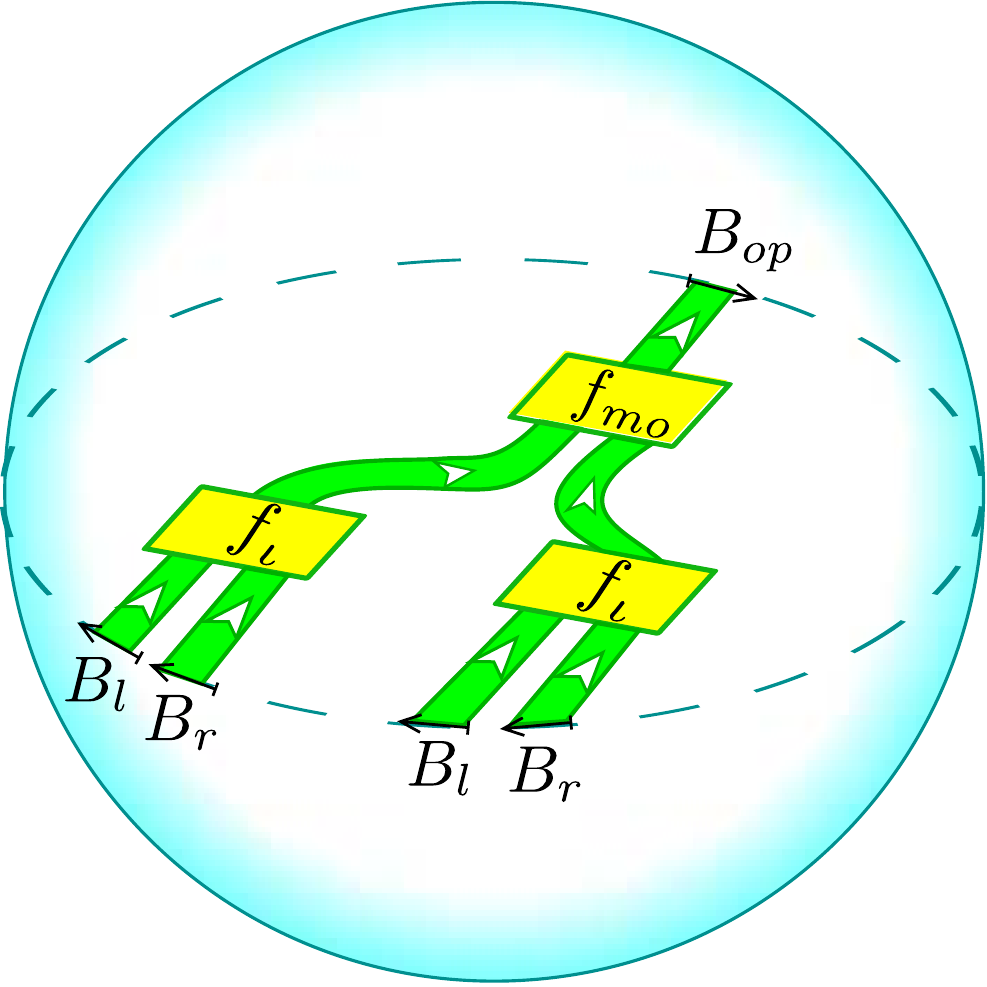}\Bigg)~.
$$
Hence, $C_{\Yr^d_{\Rr29,l}}=C_{\Yr^d_{\Rr29,r}}$ is equivalent to the following identity on morphisms:
$$
  \sum_\alpha f_\iota \circ (f_{mc,l}^\alpha \otimes_\Cb f_{mc,r}^\alpha) \circ
  (\id_{B_l} \otimes_\Cb c_{B_r,B_l}^{-1} \otimes_\Cb \id_{B_r})
  = f_{mo} \circ (f_\iota \otimes_\Cb f_\iota) ~.
$$
By (\ref{eq:T-phi-2}), we obtain the relation of $f_\alpha$ determined by R29 in Table \ref{tab:glue-mixed} immediately.

\begin{table}[bt]
$$
\begin{array}{llll}
\text{R26} & \multicolumn{3}{l}{\hspace*{-.5em}: \displaystyle
\big( (f_{\eps o} \cir f_{mo}) \otimes_\Cb \id_{T(\Bcl)} \big)
\circ \big(\id_{\Bop} \otimes_\Cb f_\iota \otimes_\Cb \id_{T(\Bcl)} \big)
\circ \big(\id_{\Bop} \otimes_\Cb T(f_{\Delta c} \cir f_{\eta c})\big)
= f_{\iota^*} }
\enl
\text{R27} & \multicolumn{3}{l}{\hspace*{-.5em}: \displaystyle
\big( T(f_{\eps c} \cir f_{mc}) \otimes_\Cb \id_{\Bop} \big)
\circ \big(\id_{T(\Bcl)} \otimes_\Cb f_{\iota^*} \otimes_\Cb \id_{\Bop} \big)
\circ \big(\id_{T(\Bcl)} \otimes_\Cb (f_{\Delta o} \cir f_{\eta o})\big)
= f_{\iota} }
\enl
\text{R28} & \text{\ equation\ } (\ref{center condition})
\enl
\text{R29} \etb:
f_\iota \circ T(f_{mc}) \circ \phi^T_{2,\Bcl,\Bcl}
= f_{mo} \circ (f_\iota \otimes_\Cb f_\iota)
\qquad
\etb
\text{R30} \etb:
f_\iota \cir T(f_{\eta c}) = f_{\eta o}
\end{array}
$$
\caption{The relations of $f_\alpha$ obtained from the relations R26-R39.}
\label{tab:glue-mixed}
\end{table}

\subsection{Relation R31} \label{sec:R31}

Let us explain the Cardy condition $\Rr31$ by drawing the world sheet $\Yr_{\Rr31}$ in a slightly deformed way:
\begin{equation}\label{eqn:Cardy}
\figbox{0.6}{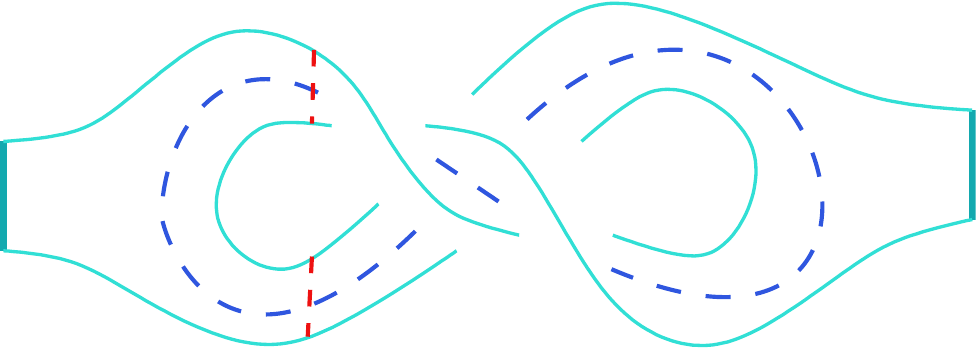}
\end{equation}
The Cardy condition $\Rr31$ arises from two ways of decomposing $\Yr_{\Rr31}$:
one along the blue circle, and the other along the two red intervals. Notice that the pre-image of the red intervals in
$\Yr_{\Rr31}$ consists of two circles which are homologous to each other. The blue circle and one of the red circle in
$\Yr_{\Rr31}$ give rise to the generators of the first homology of $\Yr_{\Rr31}$, and, by tracking their images on the boundaries of 3-dimensional cobordisms,  we will be able to describe boundary parametrisation of 3 cobordisms. 

\medskip

Now let us describe the extended cobordism associated to $\varpi_{\Rr31,l/r}$.
That associated to $\varpi_{\Rr31,l}$ is given by $\Mr_{\varpi_{\Rr31,l}}=\hat{\Yr}_l\times[0,1]/\sim$,
where $\Yr_l=\Xr_{\iota}\otimes\Xr_{\iota^*}$ and the equivalence relation $\sim$ is sketched in \eqref{eq:R31-l-3bord-1}.  By the definition $C_{\Yr^d_{\Rr31,l}}1 = \Bl(\varpi_{\Rr31,l}) \circ \big( \Psi_{\iota}(f_\iota)
\otimes_\Cb
\Psi_{\iota}(f_{\iota^*}) \big)$, we obtain
\begin{eqnarray}
\hspace{-0.3cm}C_{\Yr^d_{\Rr31,l}}&=&\tftC\Bigg(\ \figbox{0.3}{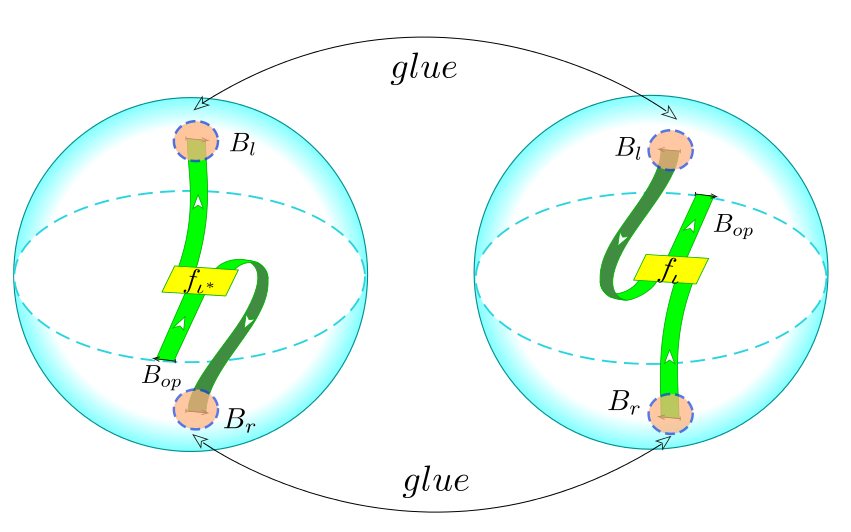}\ \Bigg) \label{eq:R31-l-3bord-1} \\
\hspace{-0.5cm}  &=& \tftC\Bigg(\, \figbox{0.4}{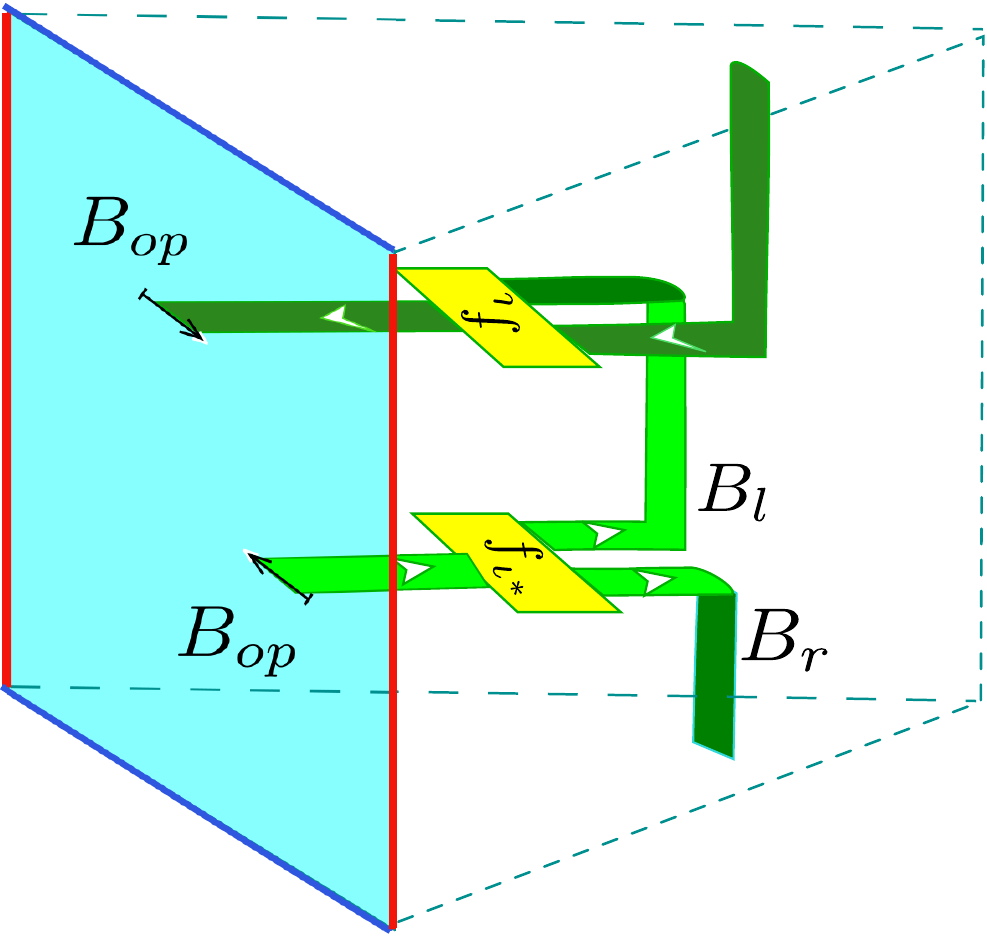} \, \Bigg) 
   = \sum_{k\in I}\sum_{\alpha}\tftC\Bigg(\ \figbox{0.35}{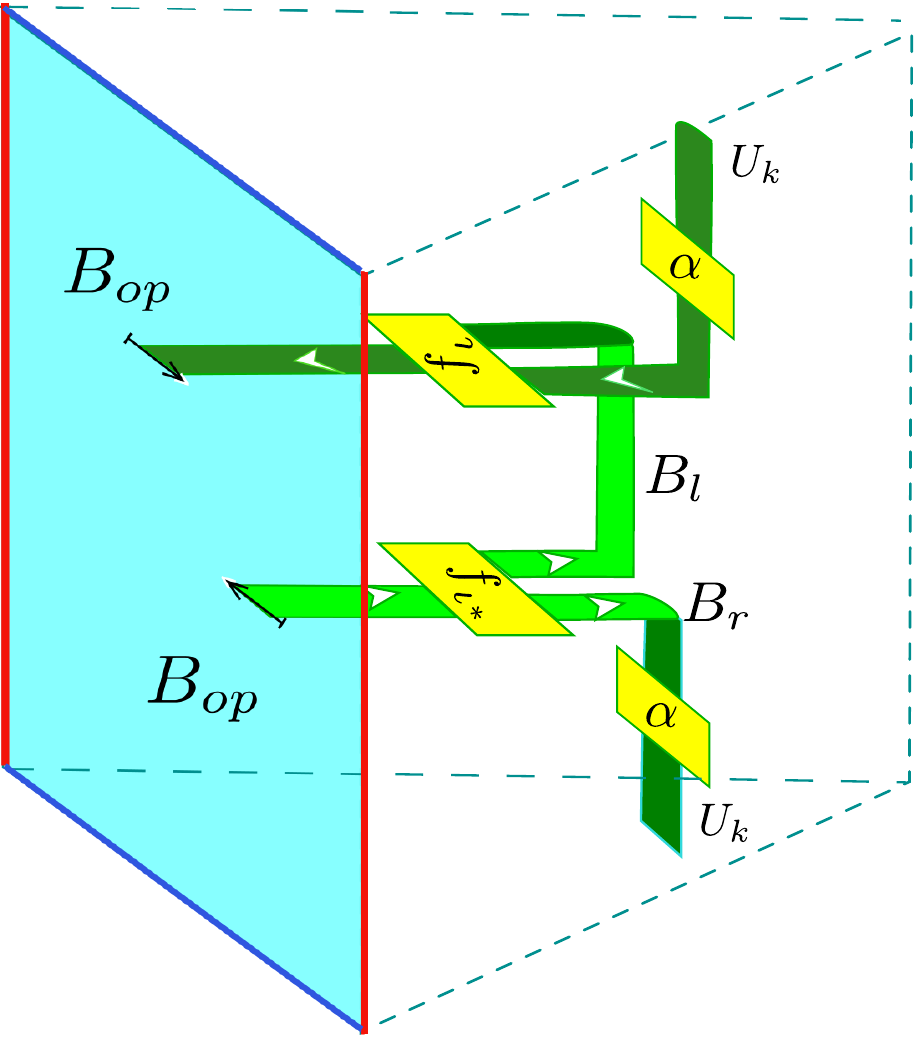}\ \Bigg)
   \label{eq:R31-l-3bord-2}
\end{eqnarray}
where, in the last equality, we have used the decomposition of the identity morphism: $\id_{B_r} = \sum_{k\in \Ic, \alpha}
b_{(k;\alpha)}^{B_r} \circ b_{B_r}^{(k;\alpha)}$ (recall (\ref{eq:b-dual-b})),  and we have used the `wedge-notation' to
represent solid tori. 
Namely, in each wedge, the top and bottom triangular surface are identified, and the dashed squares forming the front and back sides of each wedge are identified, resulting in a solid torus (see \cite[Sect.\,5.1]{tft5} for more explanations of this notation).
After these identifications, the blue and red
lines become circles on the boundary torus. In particular, the blue circles are contractible in solid torus but not
the red ones. We draw the wedge representation in such a way that the red and blue lines on the boundary of the wedge in
(\ref{eq:R31-l-3bord-2}) are the image of the two generators of the first
homology of the standard torus. 

The relation to the blue dashed circle and two red dashed lines in \eqref{eqn:Cardy} is as follows. The dashed blue circle in \eqref{eqn:Cardy} is drawn on the quotient surface and its preimage in the extended double (which is the boundary of the three-manifold shown in \eqref{eq:R31-l-3bord-2}) consists of two disjoint circles, both homologous to the blue circle in  \eqref{eq:R31-l-3bord-2}. Similarly, each of the red dashed lines in \eqref{eqn:Cardy} lifts to a circle in the extended double, both of which are homologous to the red circle in \eqref{eq:R31-l-3bord-2}.
Notice that the boundaries of the disks in (\ref{eq:R31-l-3bord-1}) along which we glue the two
balls are homologous to the blue circles.

Similarly, the extended cobordism associated to $\varpi_{\Rr31,r}$ is sketched pictorially in  
(\ref{eq:R31-r-3bord-1}). By the definition $C_{\Yr^d_{\Rr31,r}}1 = \Bl(\varpi_{\Rr31,r}) \circ \big( \Psi_{\iota}(f_{mo}) \otimes_\Cb \Psi_{\iota}(f_{\Delta o}) \big)$, we obtain
\begin{align}
 C_{\Yr^d_{\Rr31,r}}  &= \tftC\Bigg(\ \figbox{0.32}{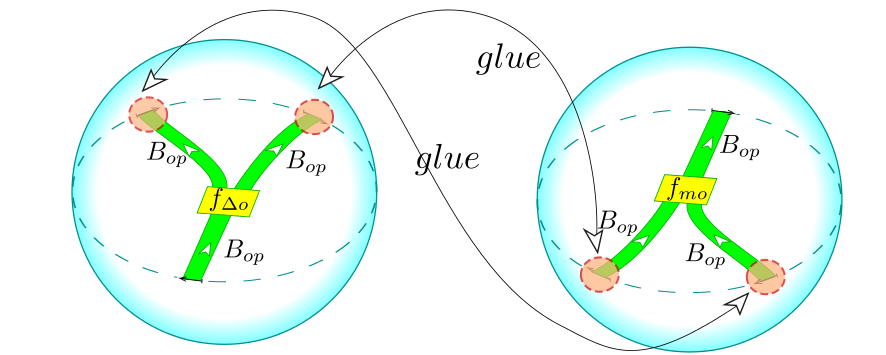}\ \Bigg) \label{eq:R31-r-3bord-1}  \\
  =& \tftC\Bigg(\, \figbox{0.3}{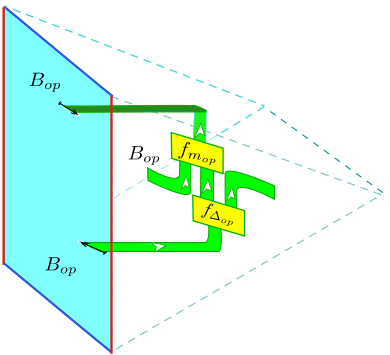}\, \Bigg)  
  \label{eqn:R31-r-3bord-2}
\end{align}

Notice that there is no homeomorphism from the cobordism in (\ref{eqn:R31-r-3bord-2}) to that in (\ref{eq:R31-l-3bord-2}) which is
compatible with the boundary parametrisation since the red circle is contractible in (\ref{eqn:R31-r-3bord-2}) but not in
(\ref{eq:R31-l-3bord-2}). For this, we need to connect the `horizontal wedge' in (\ref{eqn:R31-r-3bord-2}) to a `vertical wedge' so
that we can compare it with the right hand side of (\ref{eq:R31-l-3bord-2}). This is the content of the following lemma:

\begin{lemma}\label{lem:horizontal-to-vertical}{\rm
 \be  \label{eq:S-ident}
\sum_{k\in I}\dfrac{\dim(U_k)}{\sqrt{\text{Dim}\mathcal{C}}}
\ \tftC\Bigg(\ \figbox{0.35}{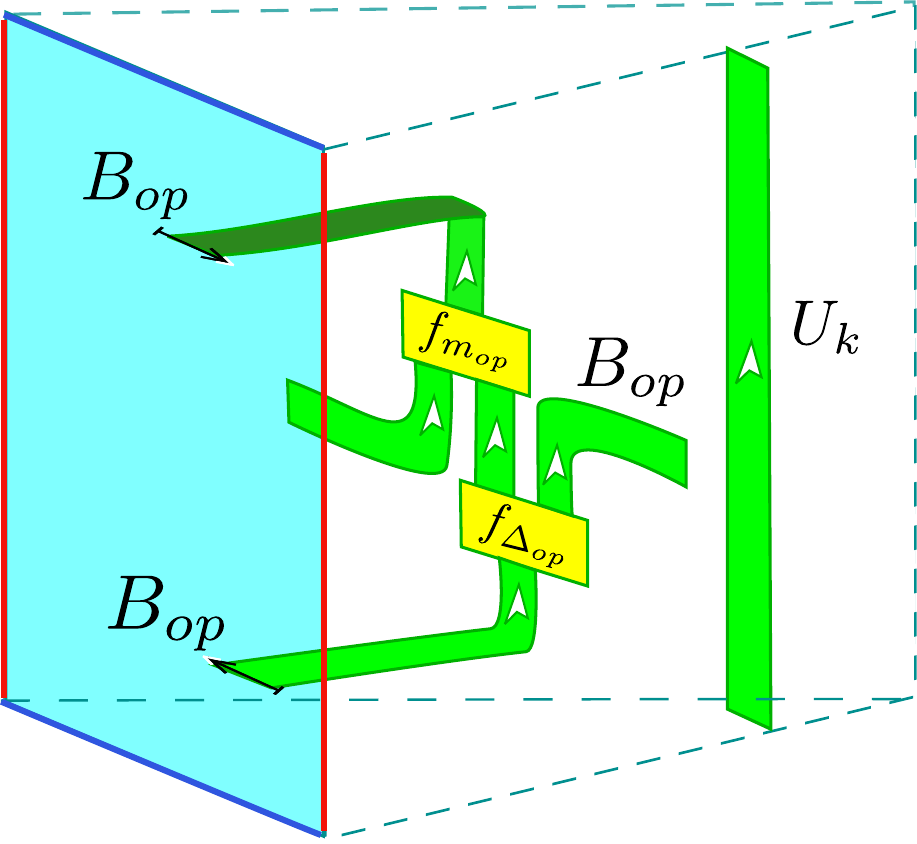}\ \Bigg)
 = \tftC\Bigg(\ \figbox{0.28}{Pic8}\ \Bigg).
\ee}
\end{lemma}

\pf
Let $T=S^1\times S^1$ be a torus with no marked arcs, which is the boundary of a solid torus $D^2\times S^1$. A natural basis of $\tftC(T)$ is given by
\begin{equation}
\Bigg\{e_k=\tftC\Bigg(\figbox{0.4}{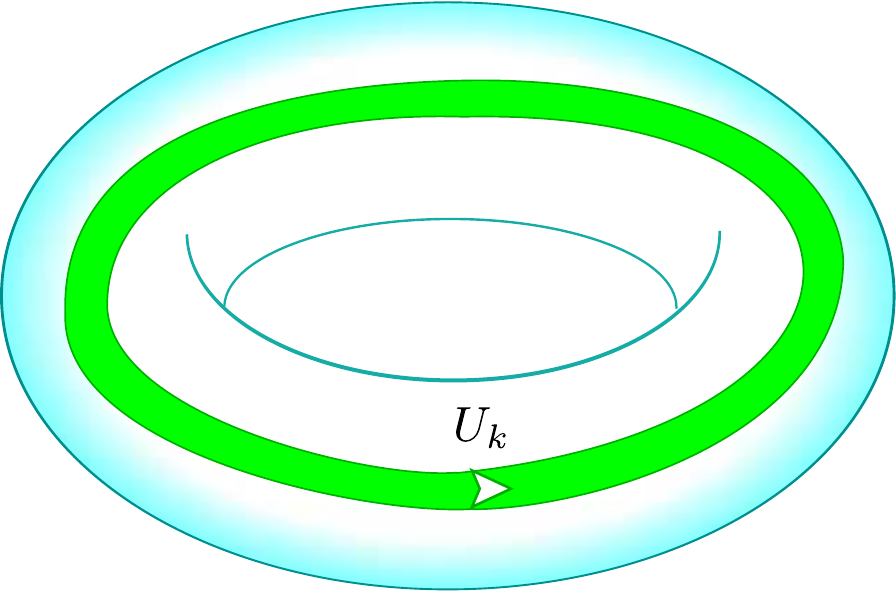}\Bigg)1,k\in I\Bigg\}
\end{equation}
The proof of the lemma is now given by the following calculation:
\begin{align}
& \tftC\Bigg(\, \figbox{0.3}{Pic8}\, \Bigg)  
  \label{eq:R31-r-3bord-2}\\
  =&\tftC\Bigg(\,\figbox{0.35}{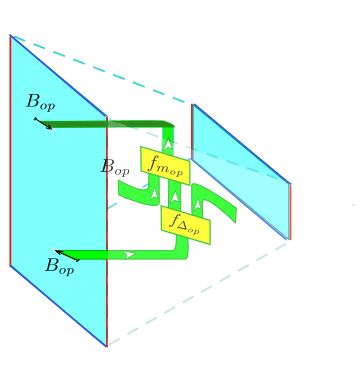}\hspace{-1.5em} \Bigg)\circ \id_{\tftC(T)}\circ\tftC\Bigg(\, \figbox{0.22}{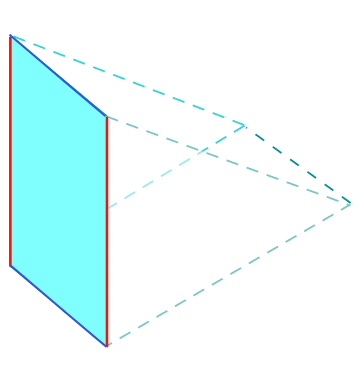}\, \Bigg)  \\
  =&\tftC\Bigg(\,\figbox{0.35}{Pic120}\hspace{-1.5em}  \Bigg)\circ (\sum_{k}(e_k\otimes e_k^\vee))\circ\tftC\Bigg(\,\figbox{0.22}{Pic116}\, \Bigg)\\
  =&\tftC\Bigg(\,\figbox{0.35}{Pic120}\hspace{-1.5em}  \Bigg)\circ \sum_{k\in
I}\left(\tftC(\figbox{0.22}{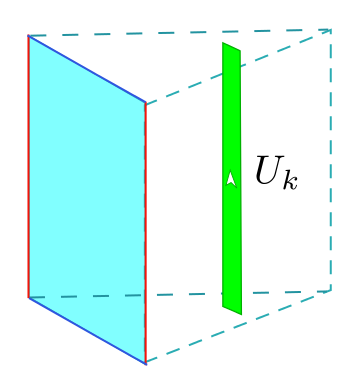}\sqcup\figbox{0.22}{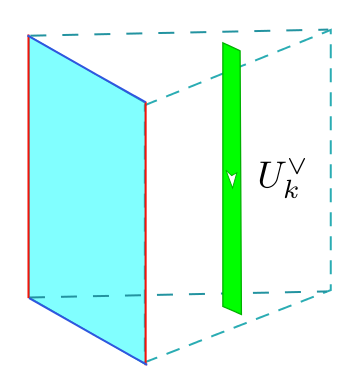})\right) \nonumber \\
    &\circ\tftC\Bigg(\,\figbox{0.22}{Pic116}\, \Bigg)  \\
=&\tftC\Bigg(\,\figbox{0.38}{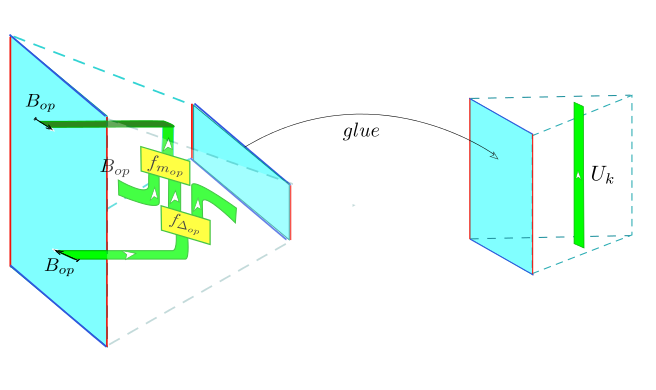}\, \Bigg) \nonumber \\
 &\circ\tftC\Bigg(\,\figbox{0.23}{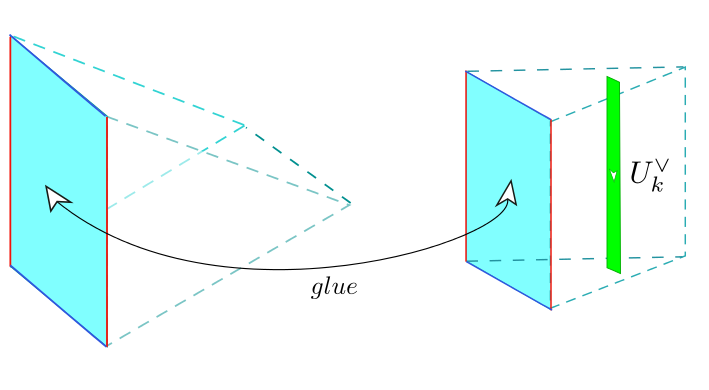}\, \Bigg)\label{eqn: 3-sphere}\\
=&\sum_{k\in I}\tftC\Bigg(\,\figbox{0.36}{117}\, \Bigg)  \cdot\dfrac{\dim U_k}{\sqrt{\text{Dim}\mathcal{C}}}
\end{align}
The only equality worth explanation is the last one: the extended cobordism in (\ref{eqn: 3-sphere}) is an $S^3$ with a
closed ribbon labelled by the simple object $U_k^\vee$ whose value under the functor $\tftC$ is the number $\dfrac{\dim
U_k}{\sqrt{\text{Dim}\mathcal{C}}}$ (see \cite[Sect.\,2.4]{tft1} or \cite[Eqn.\,(4.53)]{cardy}).
\epf

By the construction of 3-d TFT \cite[Sect.\,IV.2.1]{turaev-bk} (see also \cite[Eqn.\,(3.17)]{unique}), we have 
$$
  \tftC( \Mr(f) ) = \tftC( \Mr(g) )
  \quad \Leftrightarrow \quad
  f = g
  \qquad
  \text{for} ~ \Mr(f) =
\figbox{0.28}{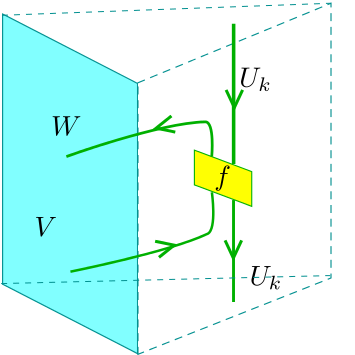}
$$
for objects $V,W$ in $\Cc$, $f,g \In \Hom(V \oti U_k^\vee, W \oti U_k^\vee)$. As a consequence, we see immediately that  the condition $ C_{\Yr^d_{\Rr31,l}}= C_{\Yr^d_{\Rr31,r}}$ is equivalent to the following condition: (recall the graphic notations introduced in (\ref{eq:b-dual-b}))
\begin{equation}\label{Cardy condition}
  \figbox{0.38}{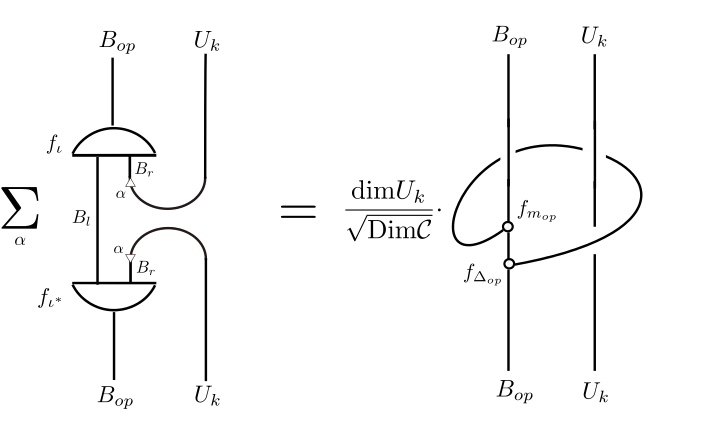} \quad .
\end{equation}
This equation is the condition on the $f_\alpha$ determined by R31, as stated in Table \ref{tab:glue-genus1}. 

\begin{table}[bt]
$$
\begin{array}{llll}
\text{R31} \etb: \text{equation\,(\ref{Cardy condition})}
\qquad
\etb
\text{R32} \etb: \text{equation\,(\ref{eq:R32-via-morph})}
\end{array}
$$
\caption{Relations at genus 1, see R31 and R32 in section \ref{sec:gen-rel}.}
\label{tab:glue-genus1}
\end{table}

\subsection{Relation R32} \label{sec:R32}

In this case, the extended cobordism $\Mr_{\varpi_{\Rr32,l}}$ associated to $\varpi_{\Rr32,l}$ is depicted by the following picture: 
\begin{equation} \label{eqn:R32-l}
\Mr_{\varpi_{\Rr32,l}}=\figbox{0.32}{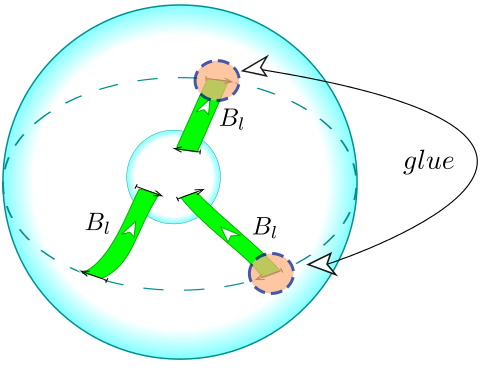}\sqcup\figbox{0.32}{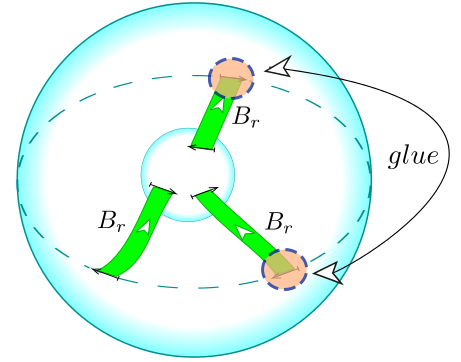}
\end{equation}
where the two red regions on each solid spheres are identified. 

\begin{rema}{\rm
By the definition of the world sheet $\Xr_{mc}$, the involution of the boundaries in \eqref{eqn:R32-l}
 is given by the mapping a point in one component to the point that is symmetric with respect to the equator in the other
component.  
}
\end{rema}

By the definition $C_{\Yr^d_{\Rr32,l}}1 = \Bl(\varpi_{\Rr32,l})
\circ  \big( \Psi_{mc}(f_{mc}) \big)$, we obtain
\begin{align}
C_{\Yr^d_{\Rr32,l}}1=&\sum_{\alpha} \tftC\left(\
\figbox{0.36}{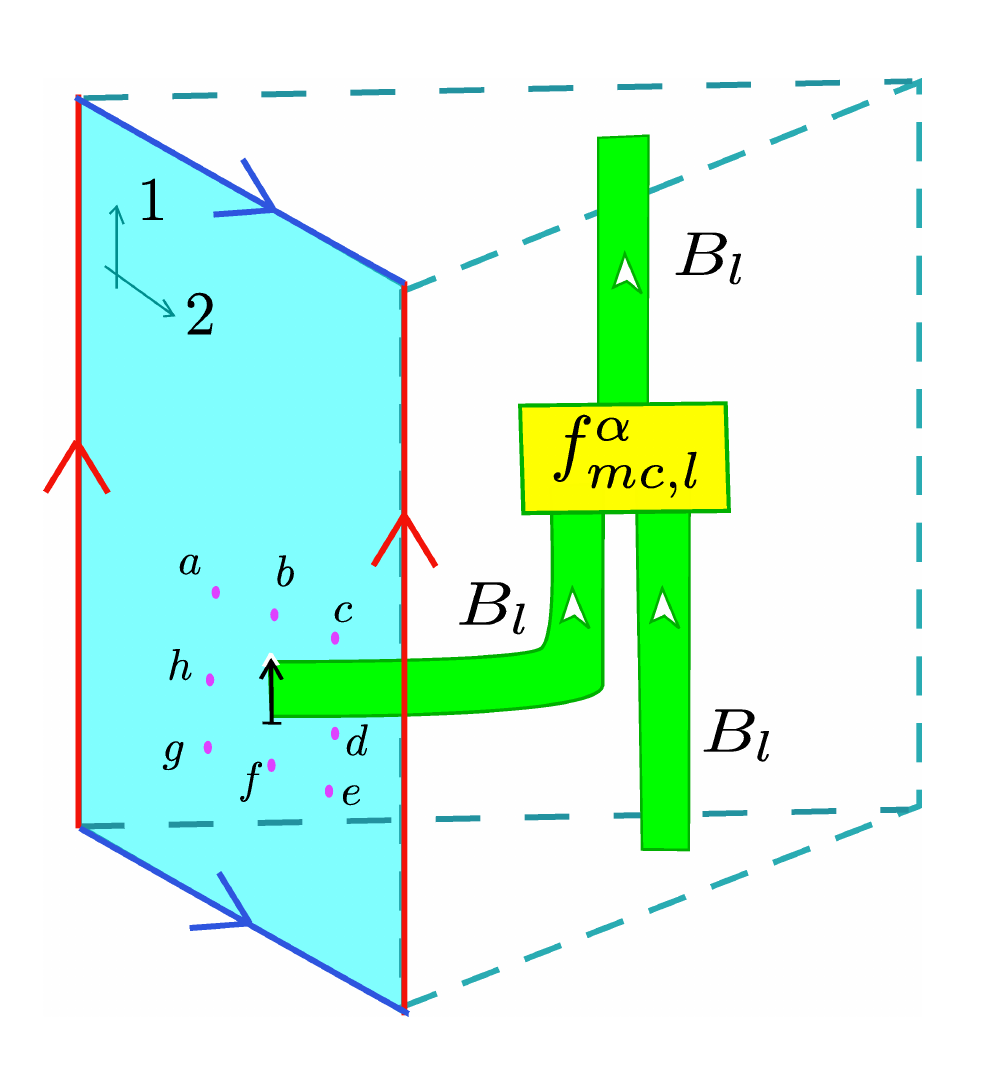}\sqcup\figbox{0.36}{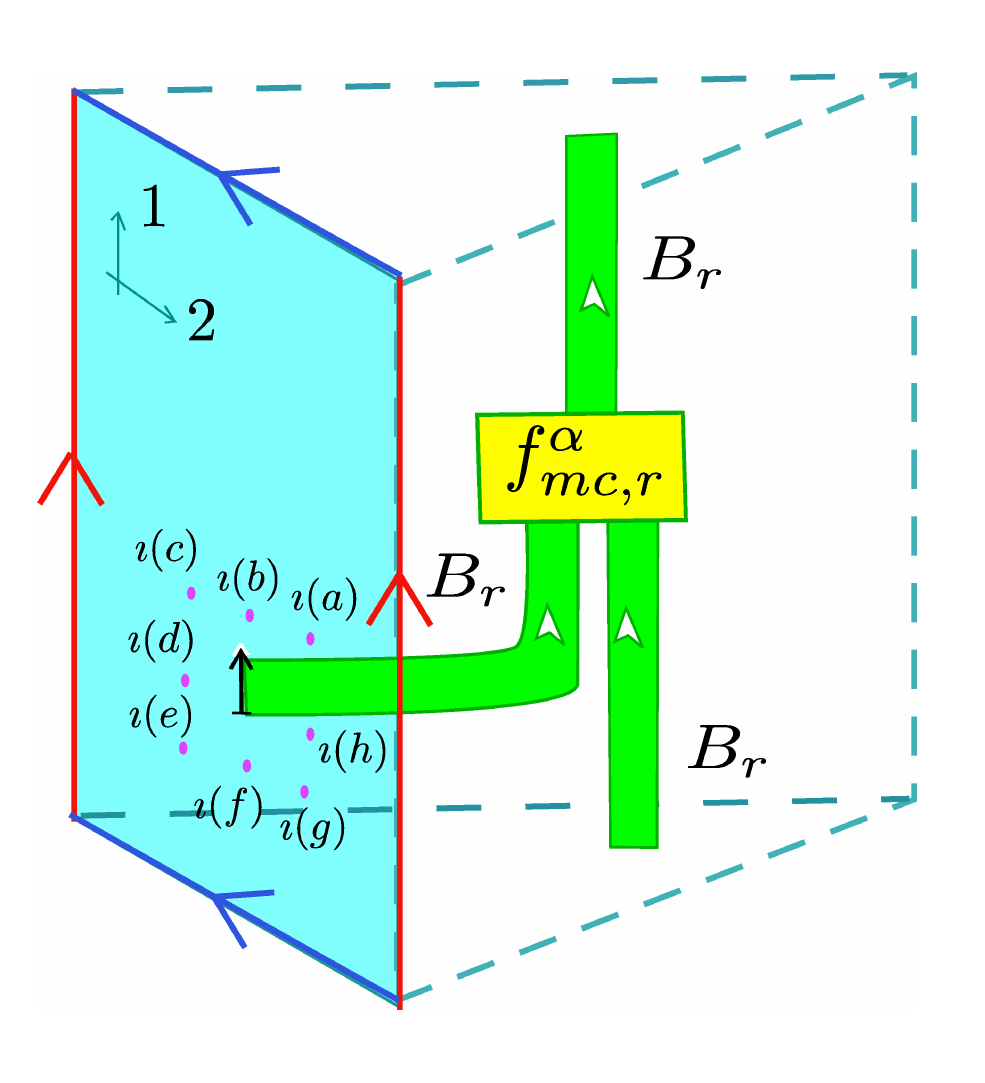}\right)\label{eq:R32-l-c}\\
                    =&\sum_{\alpha,\beta,\gamma}\sum_{i,j\in I}
\tftC\left(\ \figbox{0.36}{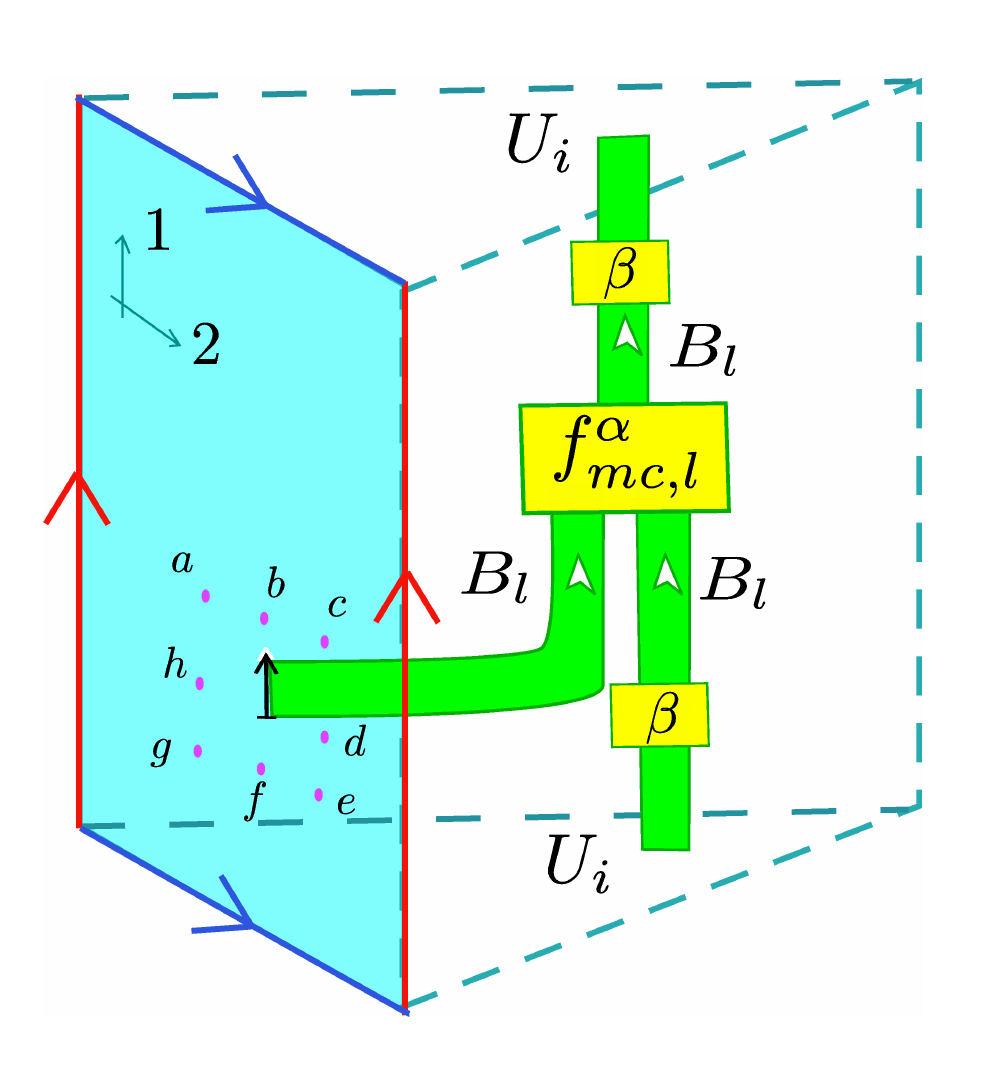}\sqcup\figbox{0.36}{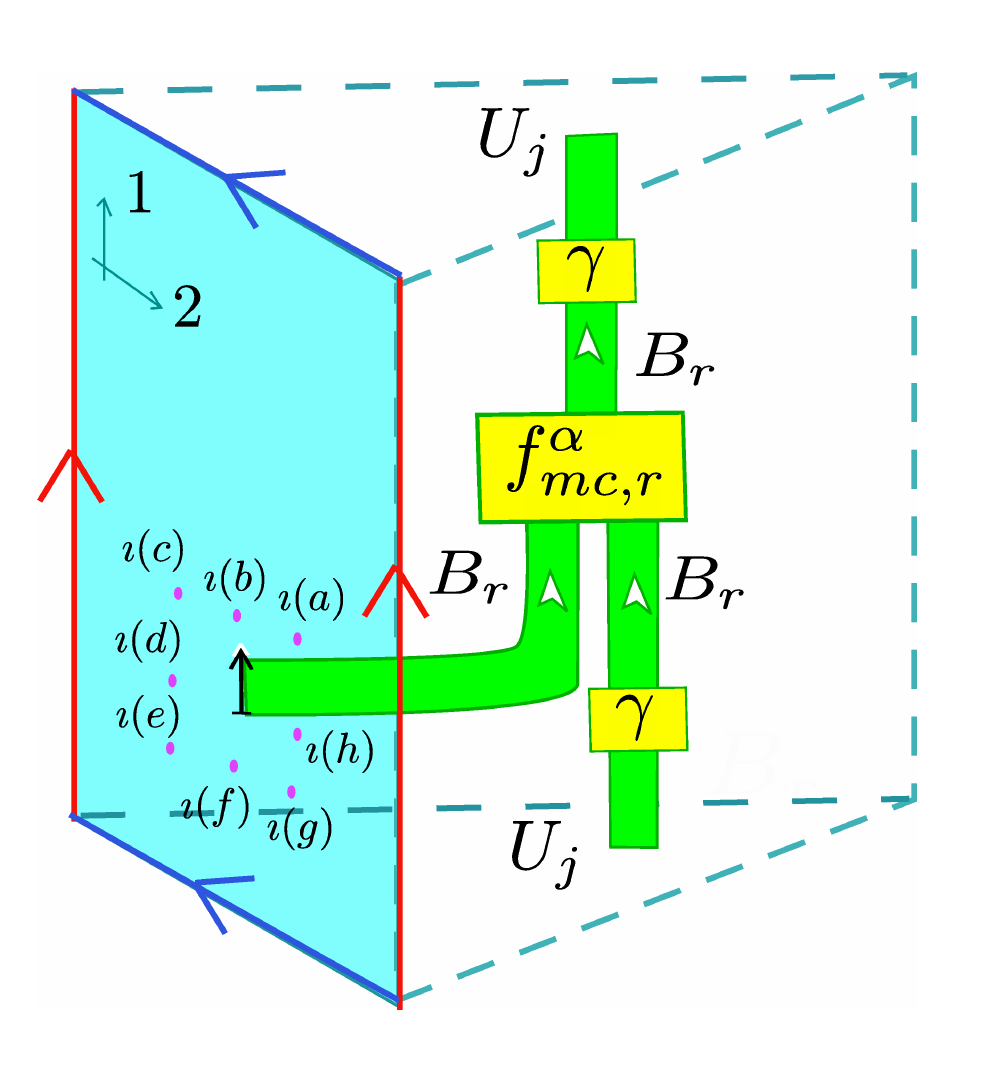}\right)\label{eqn:R32-l-c}
\end{align}
\begin{rema}
{\rm 
(i) In the above wedge representation, the blue and red arrows denote the orientation of the corresponding simple closed curves on
the boundaries. 

(ii)  Notice that it is clockwise from the arrow of the arc to the positive normal direction of the ribbon, an orientation which is
 compatible with that of the boundary surface.

(iii)  We added several purple dots to illustrate the parametrisation of the boundary surfaces and its compatibility with the
involution.  Notice that the involution reverses the orientation of the two boundary surfaces.

}
\end{rema}

Similarly, by $C_{\Yr^d_{\Rr32,r}}1 = \Bl(\varpi_{\Rr32,r}) \circ \big( \Psi_{mc}(f_{mc}) \big)$, we obtain
\begin{align}
C_{\Yr^d_{\Rr32,r}}=&\sum_{\alpha} \tftC\left(\figbox{0.3}{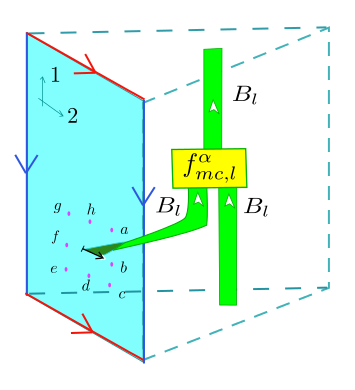}\sqcup
\figbox{0.3}{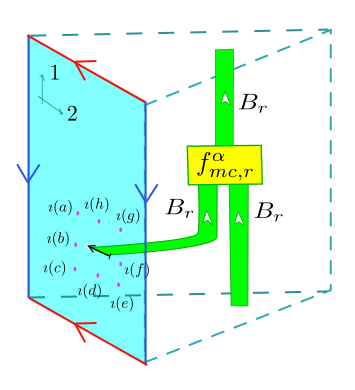}\right)\label{eqn:R32-r-1}\\
&=\sum_{\alpha} \tftC\left(\figbox{0.36}{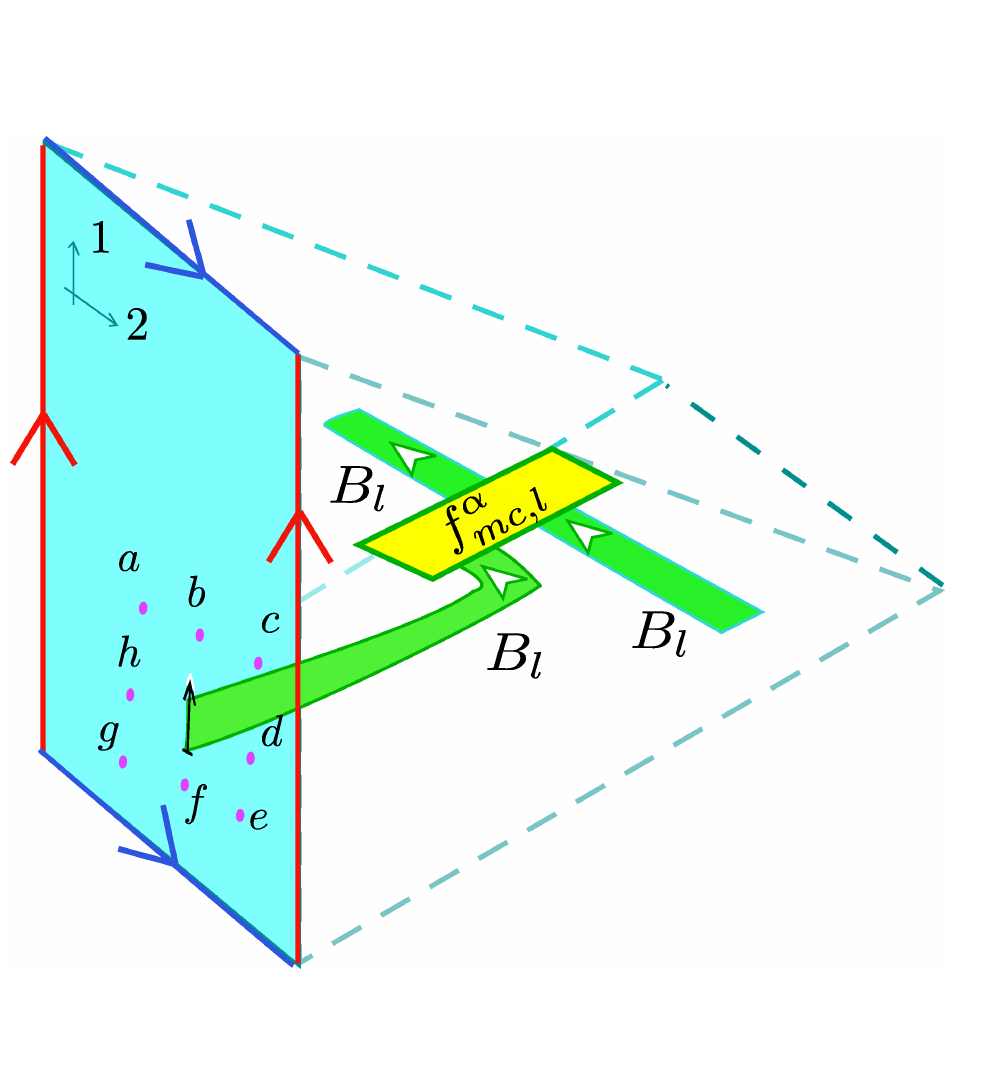}\sqcup \figbox{0.33}{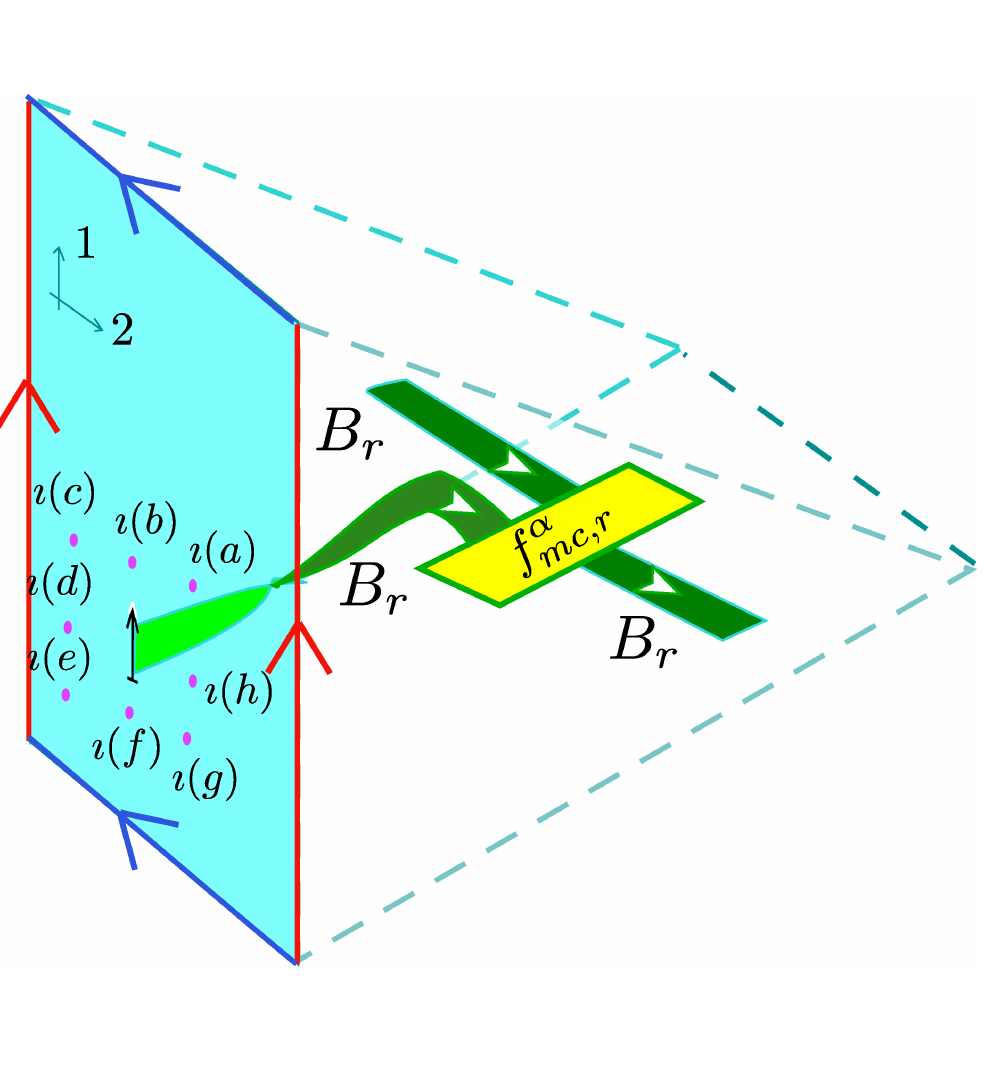}\right)\label{eqn:R32-r-2}
\end{align}
where (\ref{eqn:R32-r-2}) is obtained by rotating the left wedge of (\ref{eqn:R32-r-1}) counter-clockwise by 90 degrees, and
rotating the
right wedge of (\ref{eqn:R32-r-1}) clockwise by 90 degrees.

  In the following pictures, we will omit these purple dots and the orientation of the torus for simplicity. Similar to the proof
of Lemma \ref{lem:horizontal-to-vertical}, we have the following identities:
\begin{align*}
&\sum_{\alpha} \tftC\left(\ \figbox{0.32}{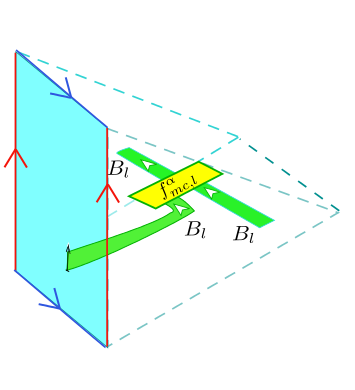}\sqcup \figbox{0.3}{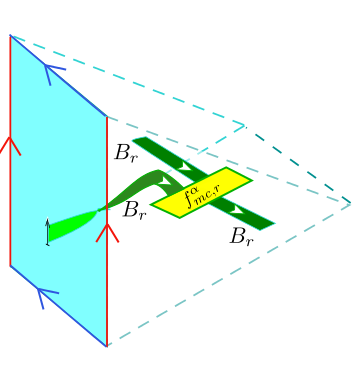}\right)\\
=&\sum_{\alpha} \tftC\left(\ \figbox{0.26}{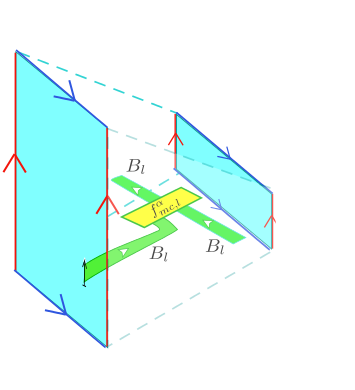}\sqcup \figbox{0.26}{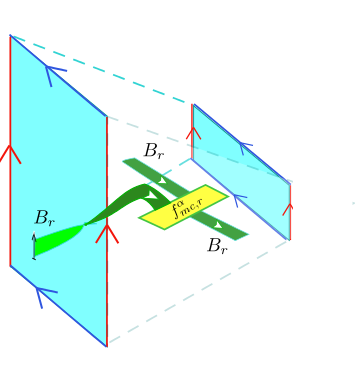}\right)\\
&\hspace{5mm} \circ\id_{\tftC(T^2\sqcup T^2)}\circ \tftC\left(\ \figbox{0.22}{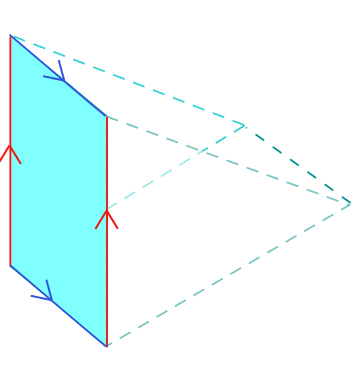}\sqcup
\figbox{0.22}{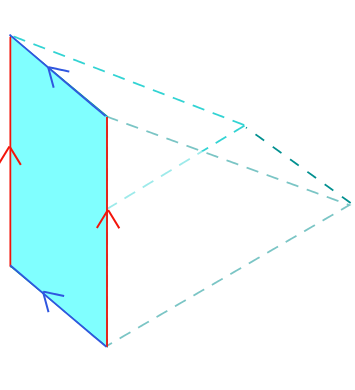}\right)\\
=&\sum_{\alpha} \tftC\left(\ \figbox{0.26}{Pic118r}\sqcup \figbox{0.26}{Pic118l}\right)\\
&\hspace{5mm} \circ \Bigg(\sum_{i,j\in I} 
\tftC\Bigg( \figbox{0.3}{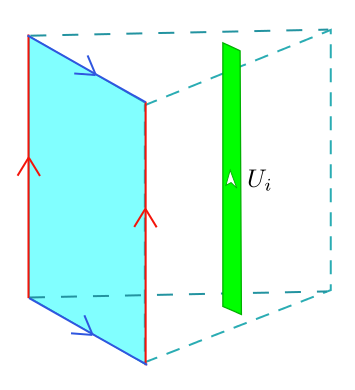}\sqcup\figbox{0.3}{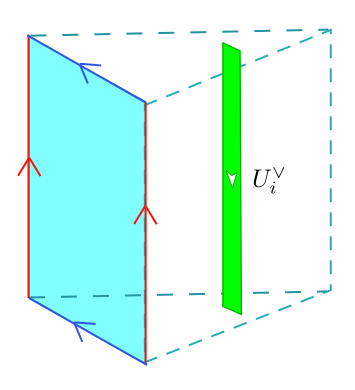}\\
&\hspace{35mm}\sqcup \figbox{0.3}{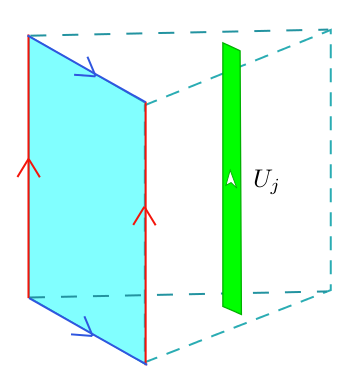}\sqcup\figbox{0.3}{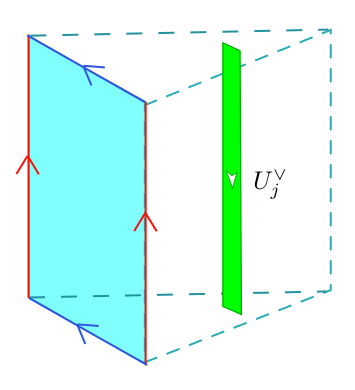}\,\,\Bigg)\Bigg) \\
&\hspace{5mm}\circ \tftC\left(\ \figbox{0.22}{Pic116r}\sqcup
\figbox{0.22}{Pic116l}\right)\\
=& \sum_{\alpha,i,j} \frac{\dim U_i\dim U_j}{\Dim \Cc}\cdot \tftC\Bigg(\ \figbox{0.29}{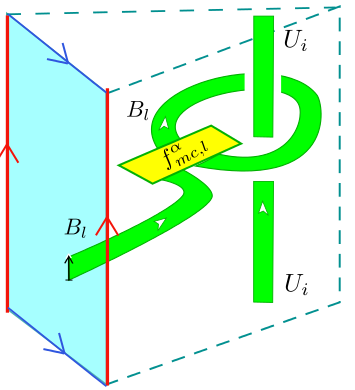}\sqcup
\figbox{0.29}{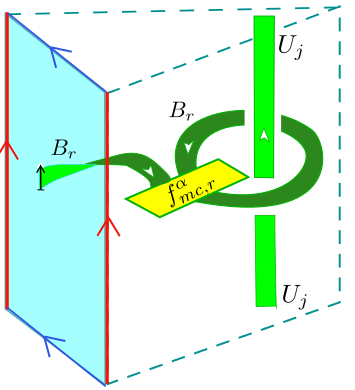}\Bigg)
\end{align*}

By comparing it  with (\ref{eqn:R32-l-c}), we obtain that the condition $C_{\Yr^d_{\Rr32,l}} = C_{\Yr^d_{\Rr32,r}}$
is equivalent to the following identity:
\begin{center}
\includegraphics[width=135mm]{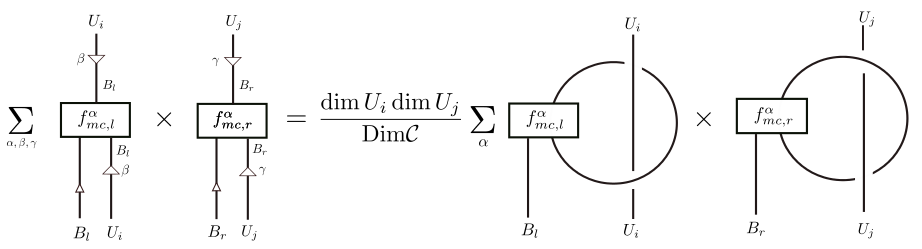}
\end{center}
for all $i,j \in \Ic$. This is nothing but
\be
  \figbox{1.7}{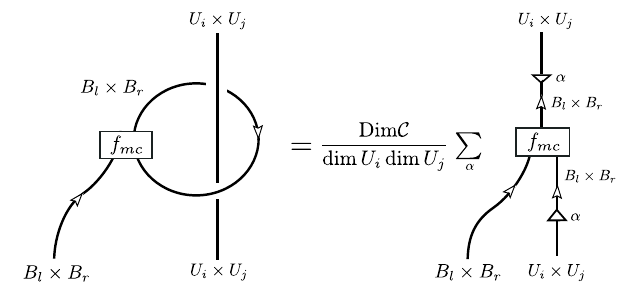} \quad ,
\labl{eq:R32-via-morph}
which is the condition stated in Table \ref{tab:glue-genus1}.

\subsection{Proof of Theorems \ref{thm:cardy-determine} and \ref{thm:cardy-sew}}
\label{sec:proof-sew}

\subsubsection*{Proof of Theorem \ref{thm:cardy-determine}}

Let $\Cor$ be a solution to the sewing constraints for $\Bl(\Cc,\Bop,B_l,B_r)$. Define the morphisms $f_{po}$ and $f_{pc}$ as in \erf{eq:f-popc-def}.
We have
\be
  f_{po} \cir f_{po}
  \overset{R1}{=} f_{po} \cir f_{mo} \cir (f_{\eta o} \otimes \id_{\Bop})
  \overset{R10}{=} f_{mo} \cir (f_{\eta o} \otimes \id_{\Bop})
  \overset{R1}{=} f_{po} ~,
\ee
so $f_{po}$ is an idempotent. In the same way one shows that $f_{pc}$ is an idempotent. Next note that we  have
\begin{eqnarray} \label{fpol}
f_{mo} \cir (\id_{\Bop} \otimes f_{po}) &\overset{R4}{=}& (f_{mo} \otimes f_{\eps o})\cir (\id_{\Bop}
\otimes f_{\Delta o}) \nonumber \\
  &\overset{R8}{=}&
  (\id_{\Bop} \otimes f_{\eps o})\cir  f_{\Delta o} \cir f_{mo}\overset{R4}{=} f_{po} \cir f_{mo}
  \overset{R12}{=} f_{mo}.
\end{eqnarray}
In the same way one shows that
\begin{equation}\label{fpor}
  f_{mo} \cir (f_{po} \otimes \id_{\Bop})
  = f_{mo} \ ,\quad
  (f_{po} \otimes \id_{\Bop}) \cir f_{\Delta o}
  = f_{\Delta o}\ ,\quad
  (\id_{\Bop} \otimes f_{po}) \cir f_{\Delta o}
  = f_{\Delta o}\ ,  
\end{equation}
\be  \label{fpcr}
  f_{mc} \cir (f_{pc} \otimes \id_{B_{cl}})
  = f_{mc}\ ,\quad
  (f_{pc} \otimes \id_{B_{cl}}) \cir f_{\Delta c}
  = f_{\Delta c}\ ,\quad
  (\id_{B_{cl}} \otimes f_{pc}) \cir f_{\Delta c}
  = f_{\Delta c}\ .
\ee
We have
\begin{eqnarray}
f_\iota\circ T(f_{pc}) &\overset{\Rr27}{=}& \big( T(f_{\eps c} \cir f_{mc}) \otimes \id_{\Bop} \big)
\circ \big(\id_{T(\Bcl)} \otimes f_{\iota^*} \otimes \id_{\Bop} \big)  \nn
&& \hspace{3.92cm} \circ \big(\id_{T(\Bcl)} \otimes (f_{\Delta o} \cir f_{\eta o})\big)\circ T(f_{pc}) \nn
&=&\big( T(f_{\eps c} \cir f_{mc}) \otimes \id_{\Bop} \big)
\circ \big(T(f_{pc}) \otimes f_{\iota^*} \otimes \id_{\Bop} \big)
\circ \big(\id_{T(\Bcl)} \otimes (f_{\Delta o} \cir f_{\eta o})\big) \nn
&=&\big( T(f_{\eps c} \cir f_{mc}) \otimes \id_{\Bop} \big)
\circ \big(\id_{T(\Bcl)} \otimes f_{\iota^*} \otimes \id_{\Bop} \big)
\circ \big(\id_{T(\Bcl)} \otimes (f_{\Delta o} \cir f_{\eta o})\big)\nn
&\overset{\Rr27}{=}&f_\iota,  \label{eq:iota-pc}
\end{eqnarray}
where the 3rd equality follows from (\ref{fpcr}). Similarly, we have
\begin{eqnarray}
T(f_{pc})\circ f_{\iota^*}&\overset{\Rr26}{=}&T(f_{pc})\circ\big( (f_{\eps o} \cir f_{mo}) \otimes \id_{T(\Bcl)} \big) \nn
&&\hspace{1.15cm}  \circ \big(\id_{\Bop} \otimes f_\iota \otimes \id_{T(\Bcl)} \big) 
\circ \big(\id_{\Bop} \otimes T(f_{\Delta c} \cir f_{\eta c})\big) \nn
 &=&\big( (f_{\eps o} \cir f_{mo}) \otimes \id_{T(\Bcl)} \big)
\circ \big(\id_{\Bop} \otimes f_\iota \otimes T(f_{pc}) \big)
\circ \big(\id_{\Bop} \otimes T(f_{\Delta c} \cir f_{\eta c})\big)\nn
&=&\big( (f_{\eps o} \cir f_{mo}) \otimes \id_{T(\Bcl)} \big)
\circ \big(\id_{\Bop} \otimes f_\iota \otimes \id_{T(\Bcl)} \big)
\circ \big(\id_{\Bop} \otimes T(f_{\Delta c} \cir f_{\eta c})\big)\nn
&\overset{\Rr26}{=}&f_{\iota^*}  \label{eq:iota*-pc} \ .
\end{eqnarray}

We choose retracts $(\Aop,\eop,\rop)$ and $(\Acl,\ecl,\rcl)$ for $f_{po}$ and $f_{pc}$, respectively.  
We assume that $A_{cl}=\oplus_{n=1}^NC_n^l\times C_n^r$ as in Definition \ref{def:cardy-alg}.
The morphism $e_{cl}:\oplus_{n=1}^NC_n^l\times C_n^r\rightarrow B_l\times B_r$ 
can be expanded as follows: 
$$
e_{cl} = \oplus_{n=1}^N e_{cl,n} = \oplus_{n=1}^N \sum_\alpha e_{cl,n,\alpha}^l\otimes_\Cb e_{cl,n,\alpha}^r \ ,
$$
where $e_{cl,n,\alpha}^l\otimes_\Cb e_{cl,n,\alpha}^r: C_n^l\times C_n^r\rightarrow B_l\times B_r$.  
By the definition of the functor $T$, we have $T(A_{cl})=\oplus_{n=1}^NC_n^l\otimes  C_n^r$. Then we have 
\begin{equation}\label{extension}
T(e_{cl})=\sum_{n=1}^N T(e_{cl,n})=\sum_{n=1}^N\sum_\alpha  e_{cl,n,\alpha}^l\otimes e_{cl,n,\alpha}^r \ .
\end{equation}

We need to show that the morphisms \erf{eq:falpha-names}, assuming relations R1--R32, define a Cardy algebra on $\Aop$ and $\Acl$. That $(\Aop,m_\text{op},\eta_\text{op},\Delta_\text{op},\eps_\text{op})$ is a symmetric Frobenius algebra in $\Cc$ follows from R1--R13 in Table \ref{tab:glue-open}.
In the same way, R14--R25 in Table \ref{tab:glue-closed} guarantee that $(\Acl,m_\text{cl},\eta_\text{cl},\Delta_\text{cl},\eps_\text{cl})$ is a commutative symmetric\footnote{That $\Acl$ is symmetric follows from the fact that  $\Acl$ is commutative (by R20) and has trivial twist (by R25).} Frobenius algebra in $\CxC$. For more details of this part, we refer to \cite[Sect.\,4]{unique}. It remains to check the center condition, the Cardy condition and the modular invariance condition.
\begin{itemize}\setlength{\leftskip}{-1em}
 \item Center condition:
      We have the following sequence of identities. 
      \begin{align}
&\figbox{0.3}{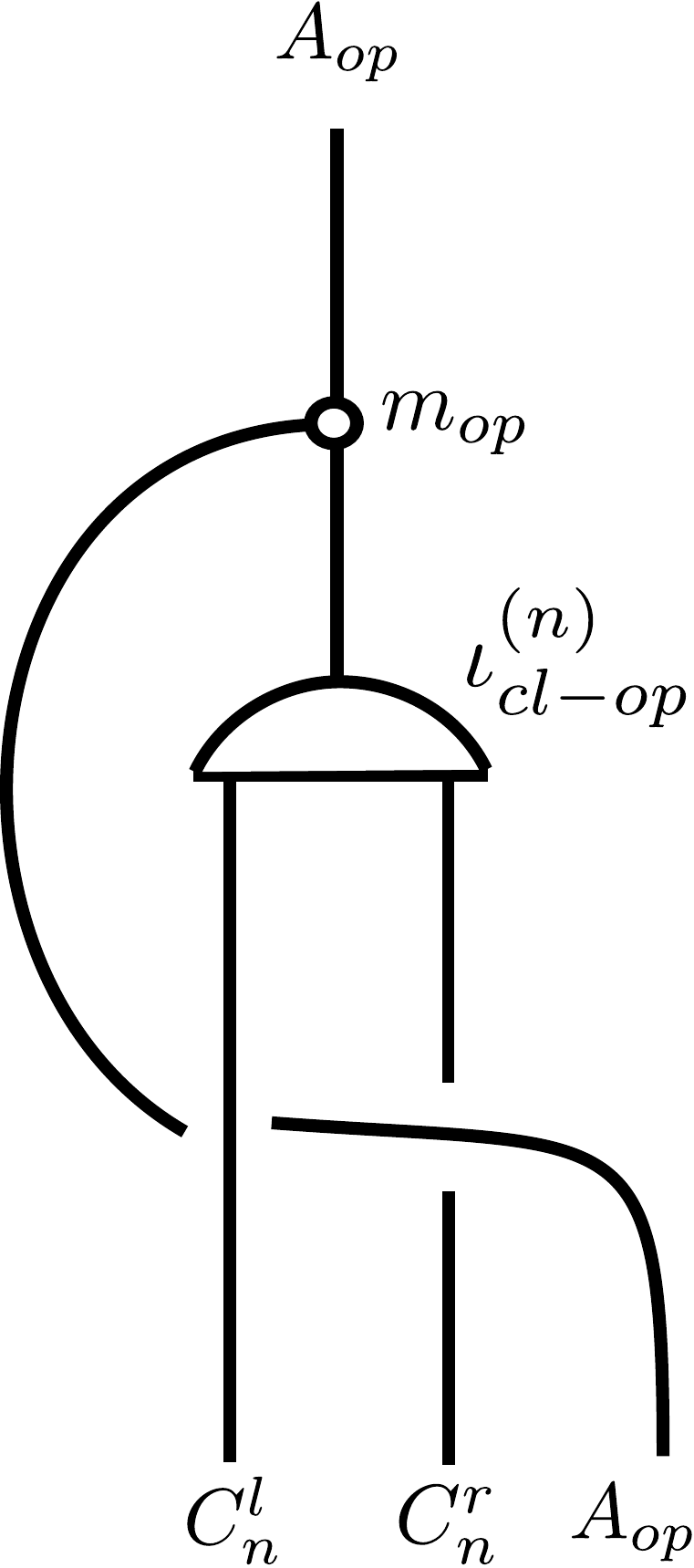}\overset{(1)}{=} \figbox{0.3}{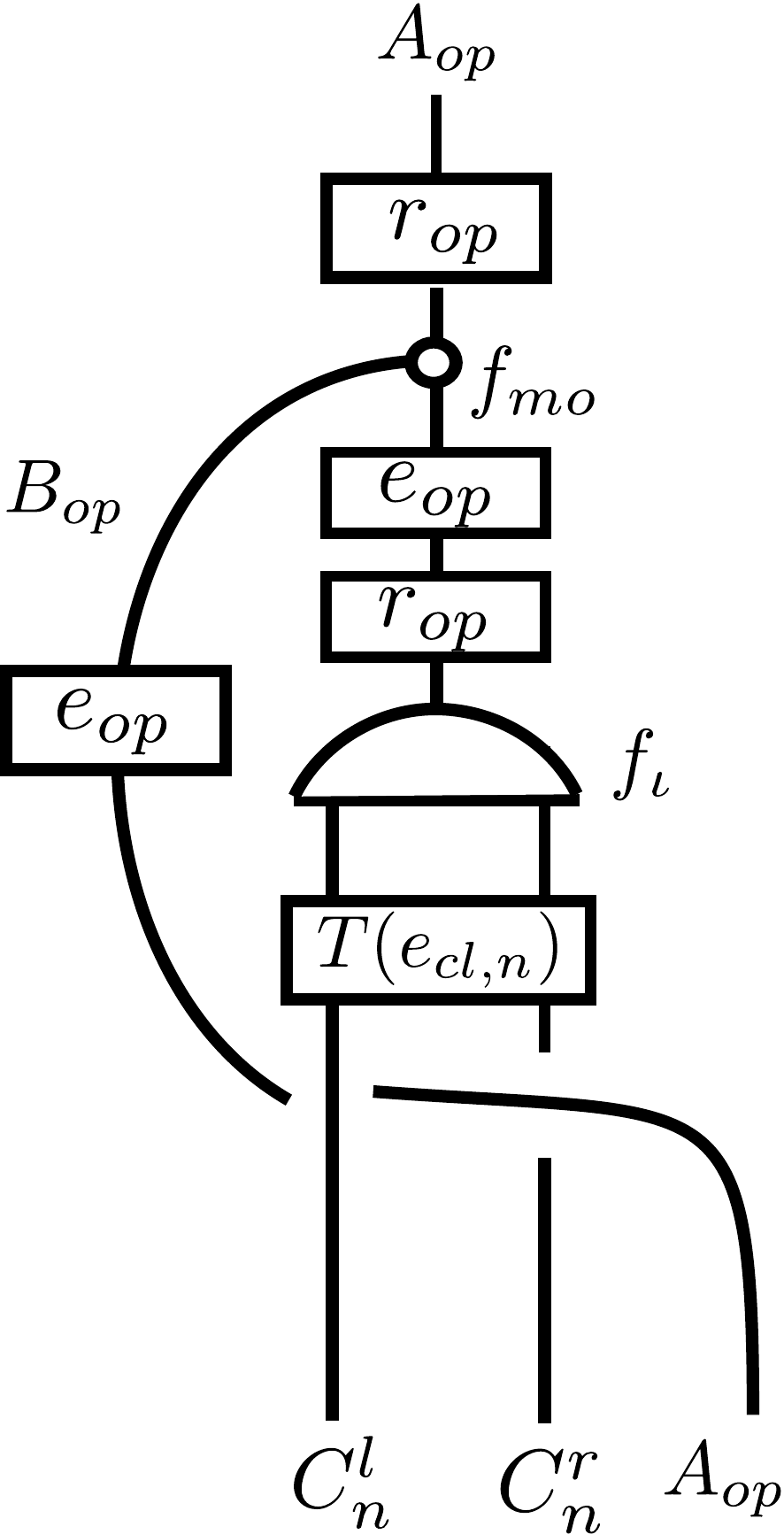}\overset{(2)}{=}
       \sum_\alpha\figbox{0.3}{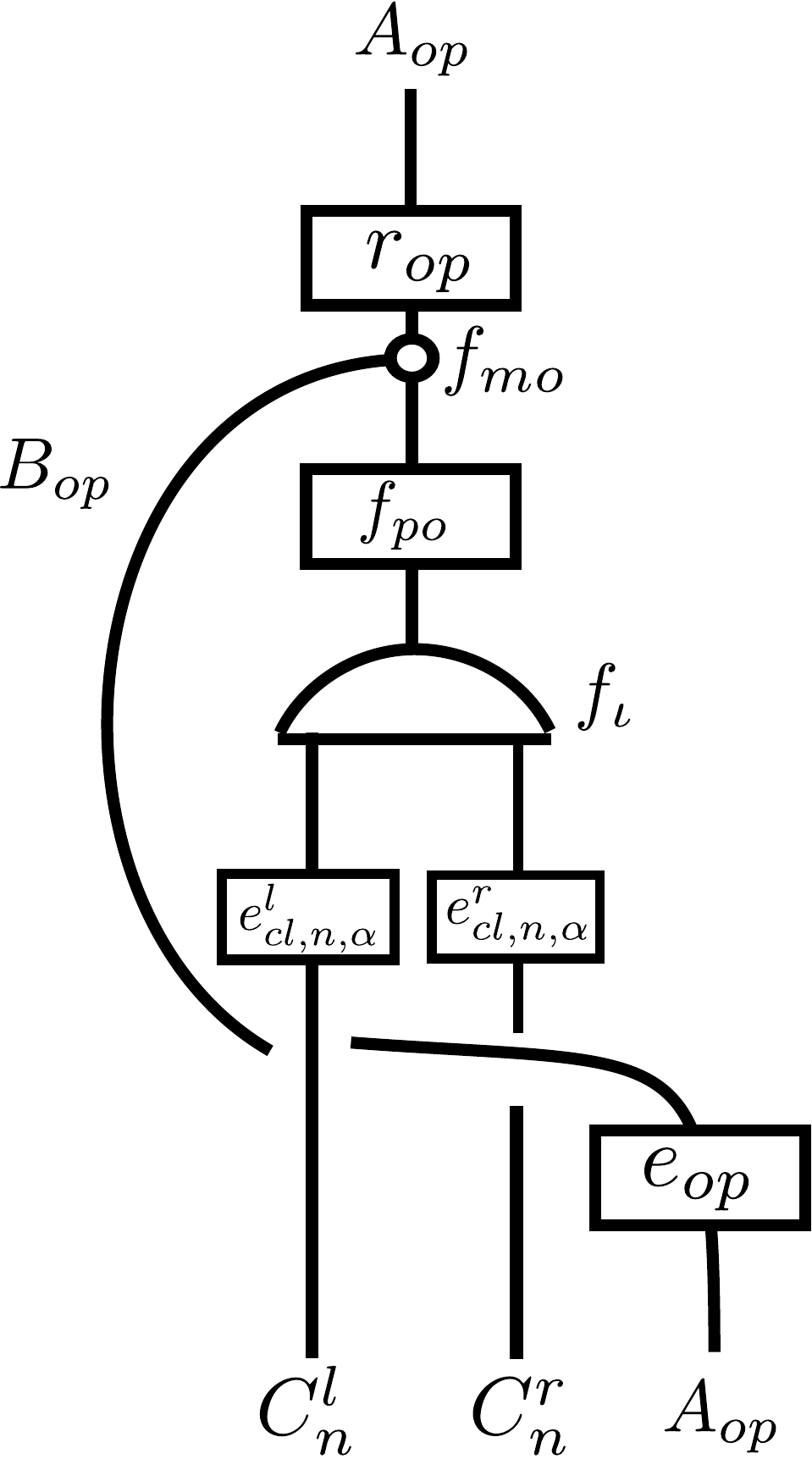}
       \nonumber \\
       & \overset{(3)}{=}\sum_\alpha \figbox{0.3}{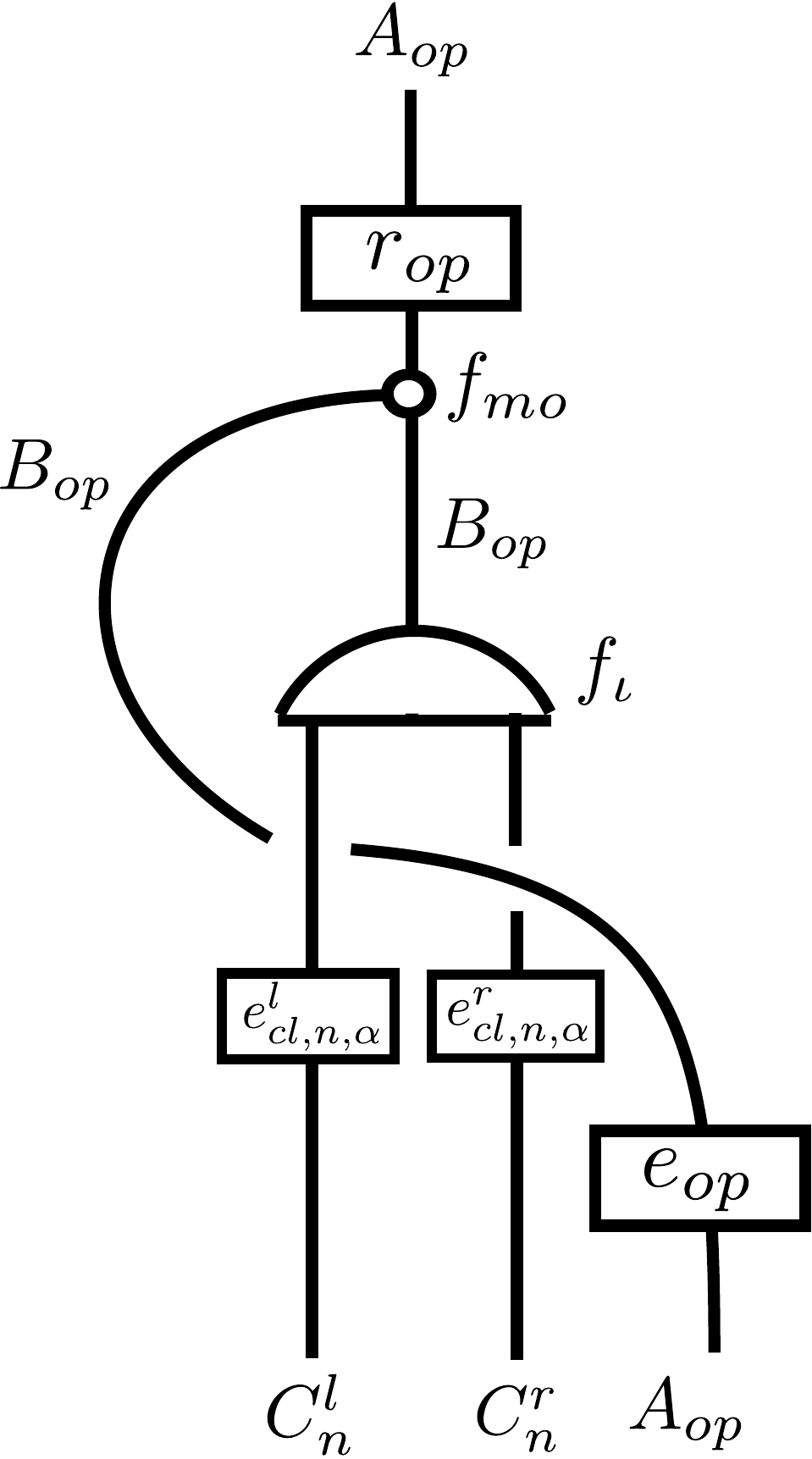}
       \overset{\Rr28}{=}\sum_\alpha\figbox{0.3}{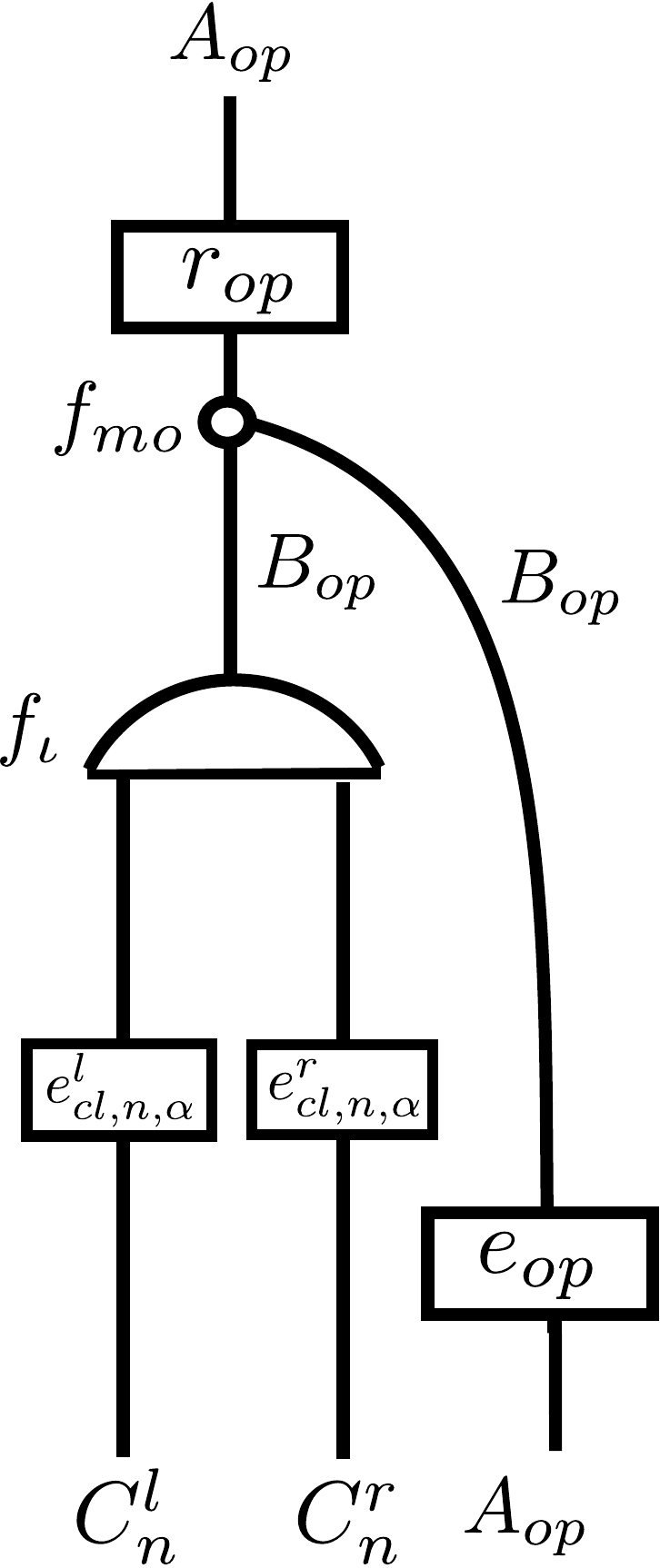}\overset{(5)}{=}\ \figbox{0.3}{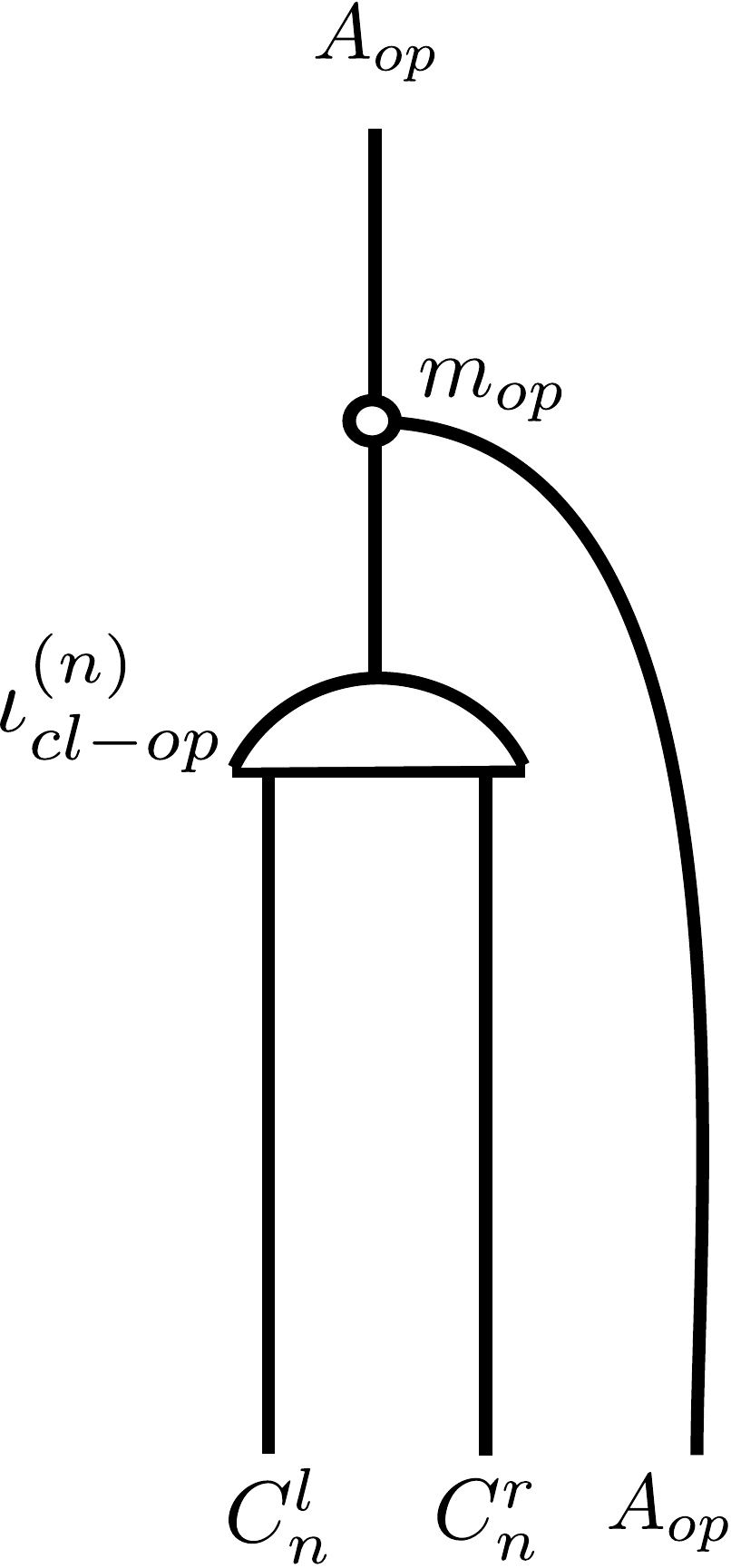}
      \end{align} 
      where step (1) follows from the definition of $m_{op}$ and 
        $\ico^{\, (n)}$; 
      step (2) follows from (\ref{extension}); 
      step (3) follows from the naturality of the braiding and R10; step (5) is the reverse procedure of (1) and (2).
 
 \item Cardy condition:
 
       We start our derivation from the graph on the right hand side of (\ref{eq:cardy-C}). We have the following sequence of identities, 
      \begin{align*}
       &\figbox{0.45}{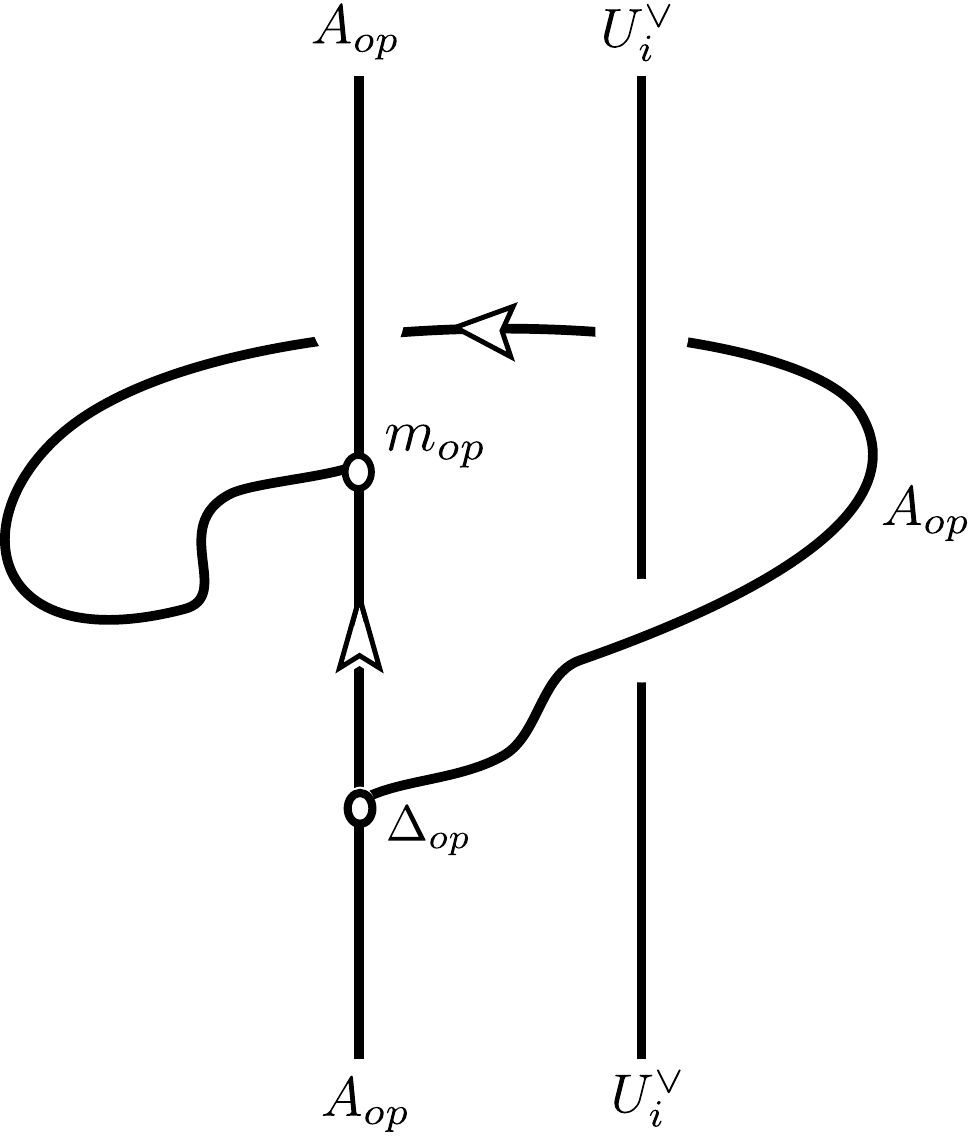} \overset{}{=}\figbox{0.45}{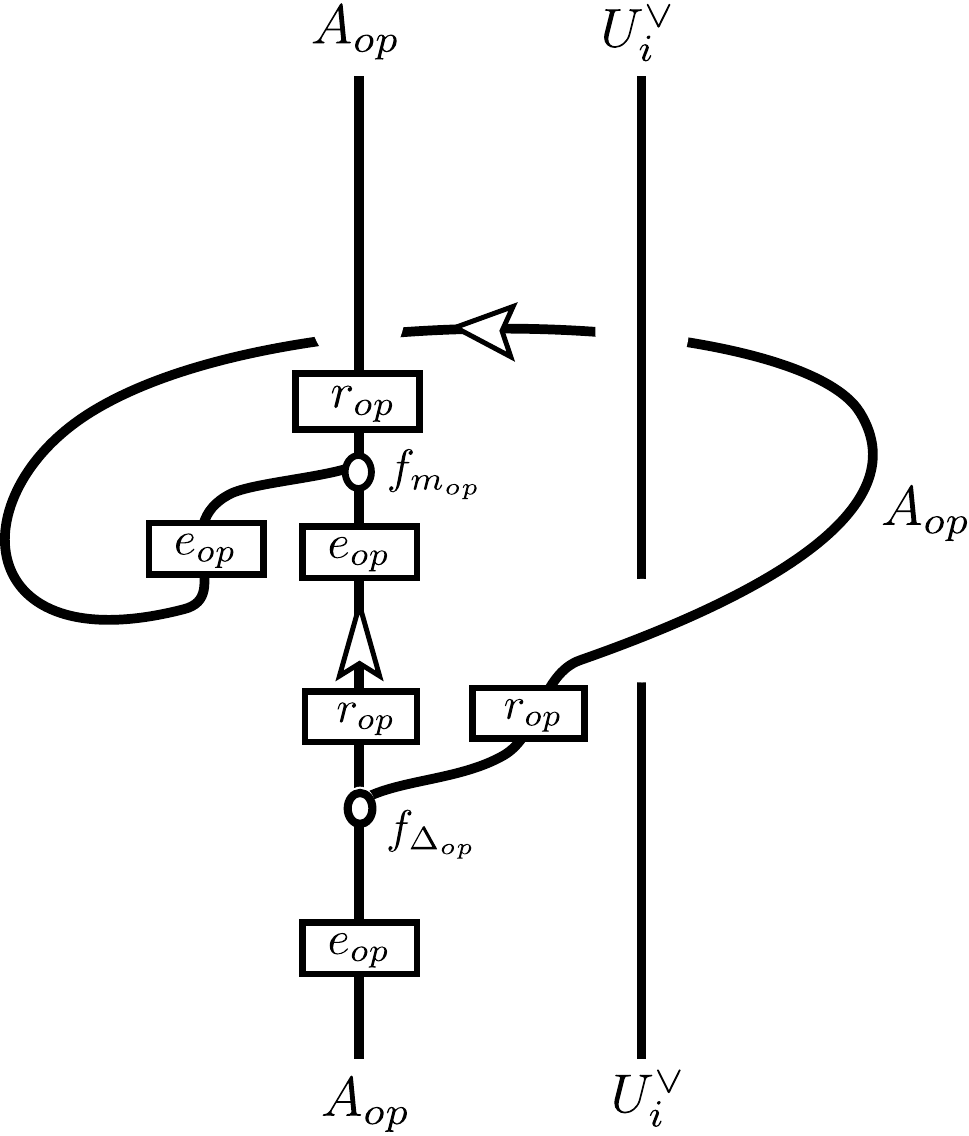}\overset{}{=}\figbox{0.45}{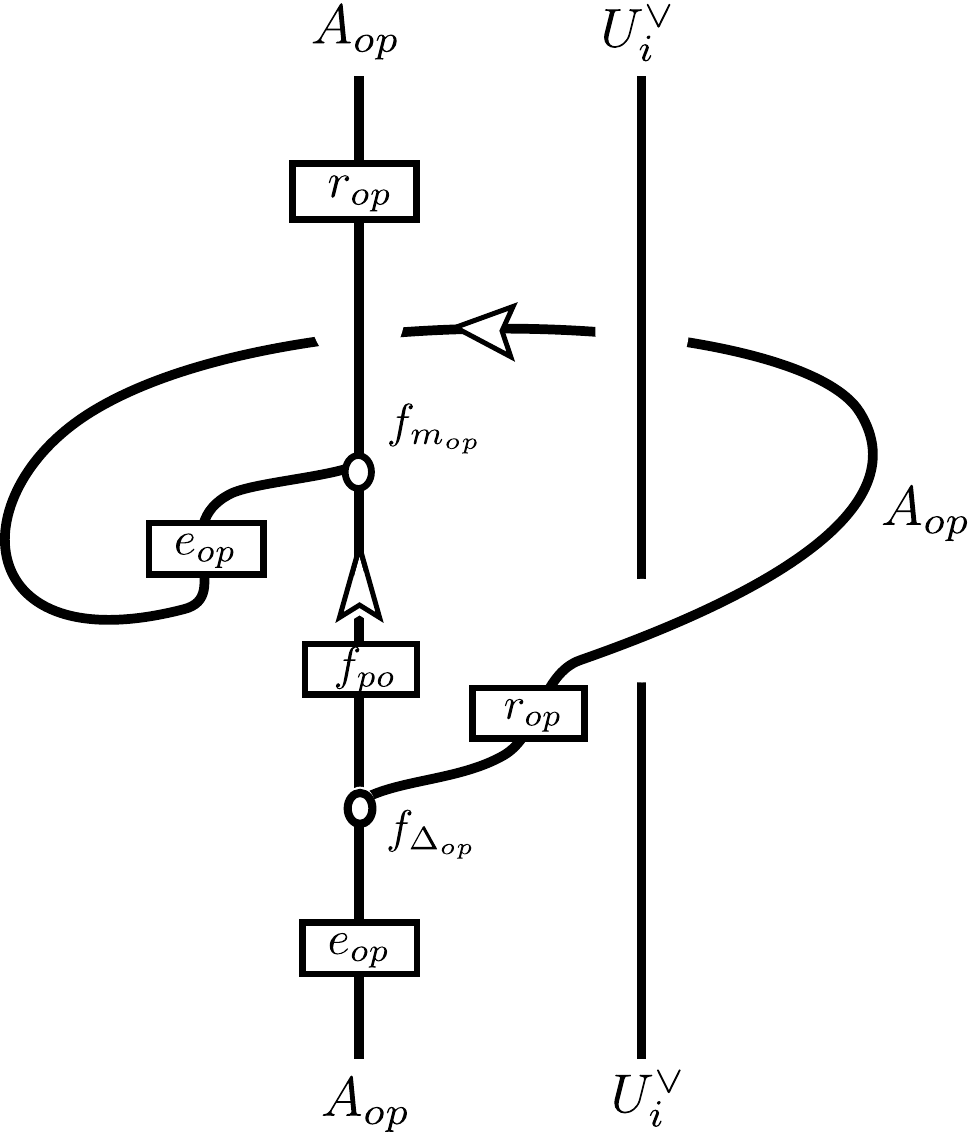}\\
        \overset{(3)}{=}&\figbox{0.45}{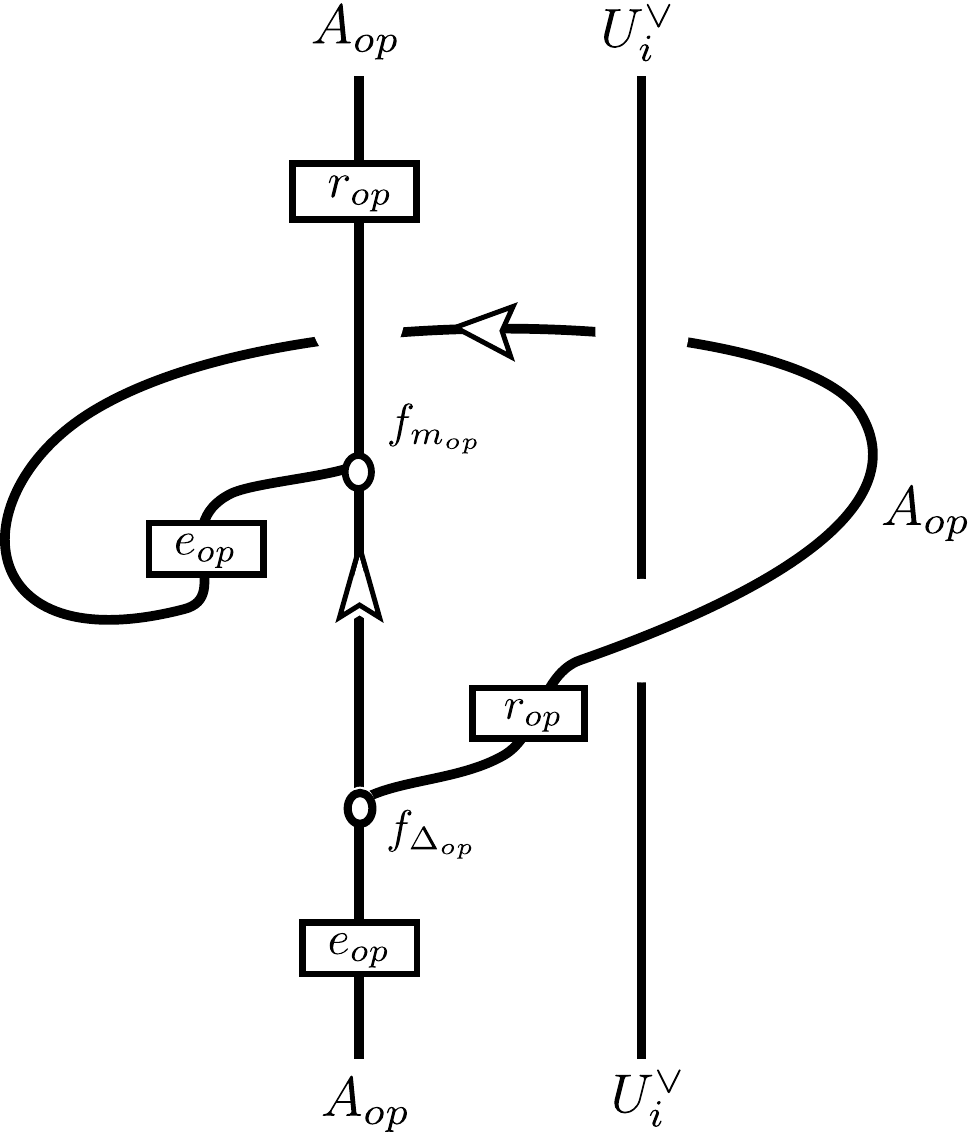}\overset{}{=}\figbox{0.45}{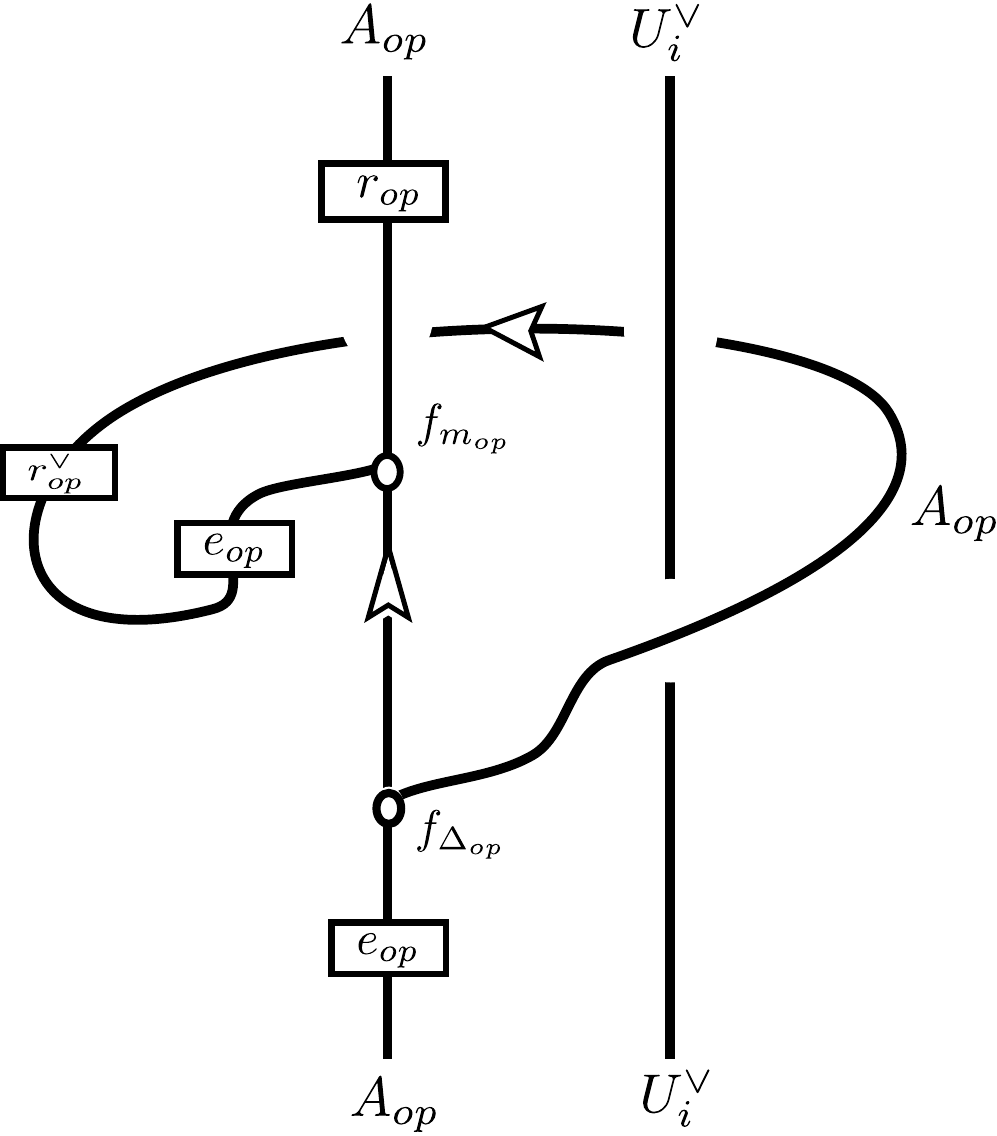}\overset{}{=}\figbox{0.45}{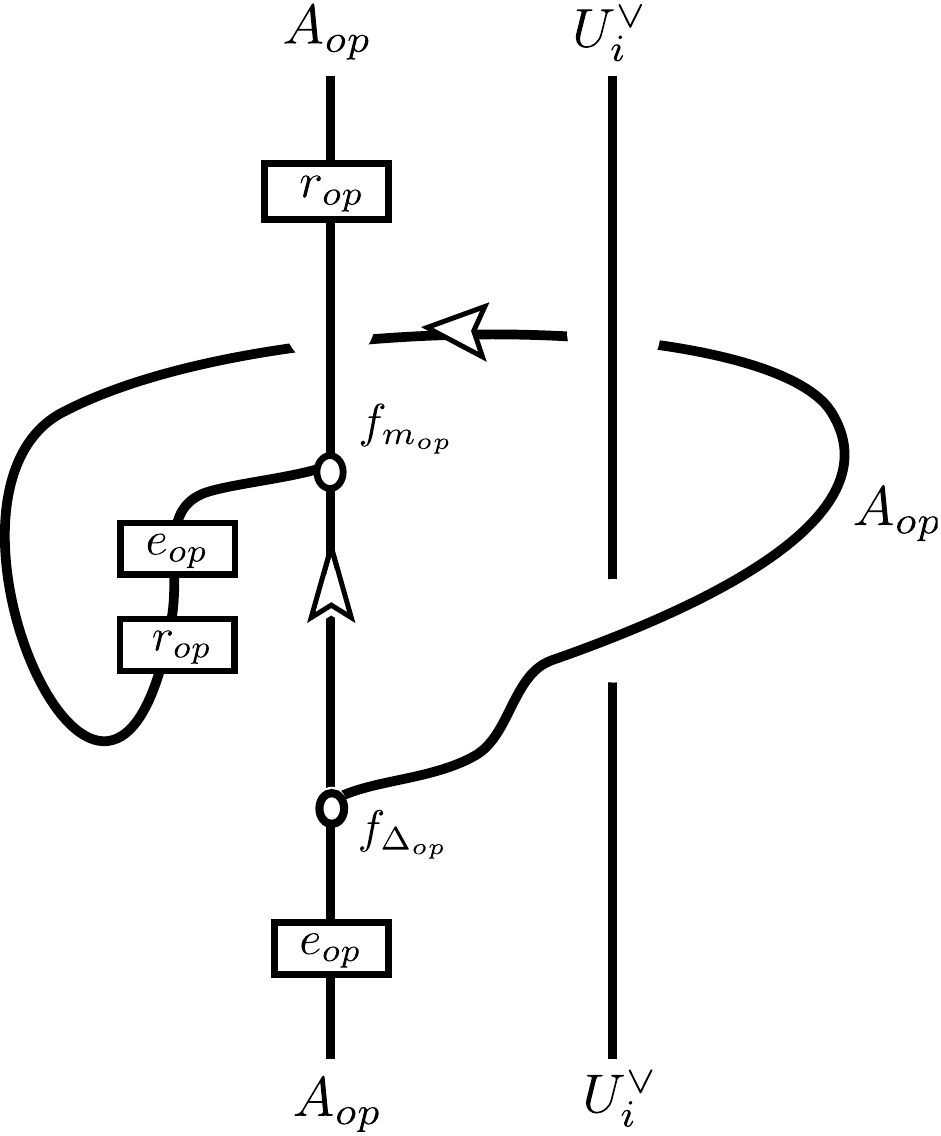}\\
        \overset{}{=}&\figbox{0.45}{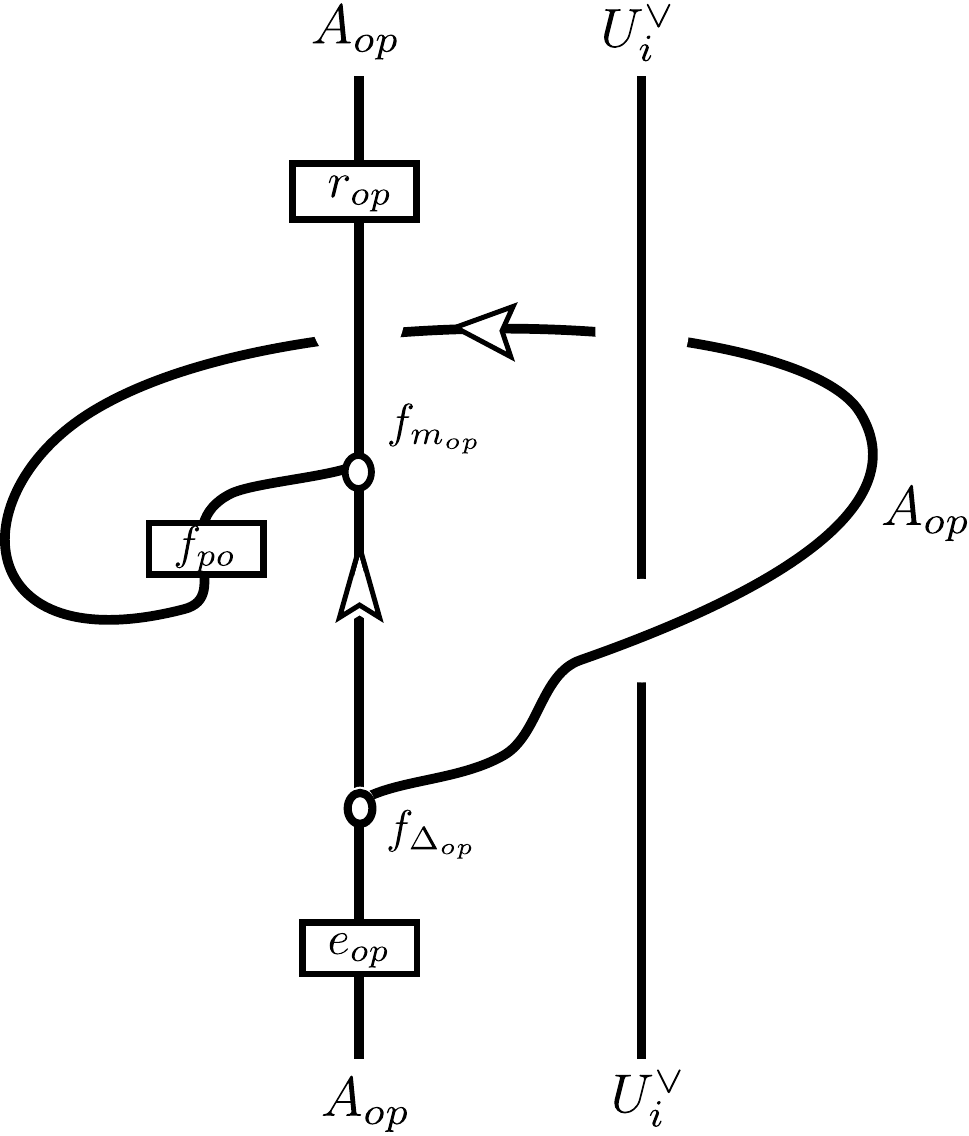}\overset{(7)}{=}\figbox{0.45}{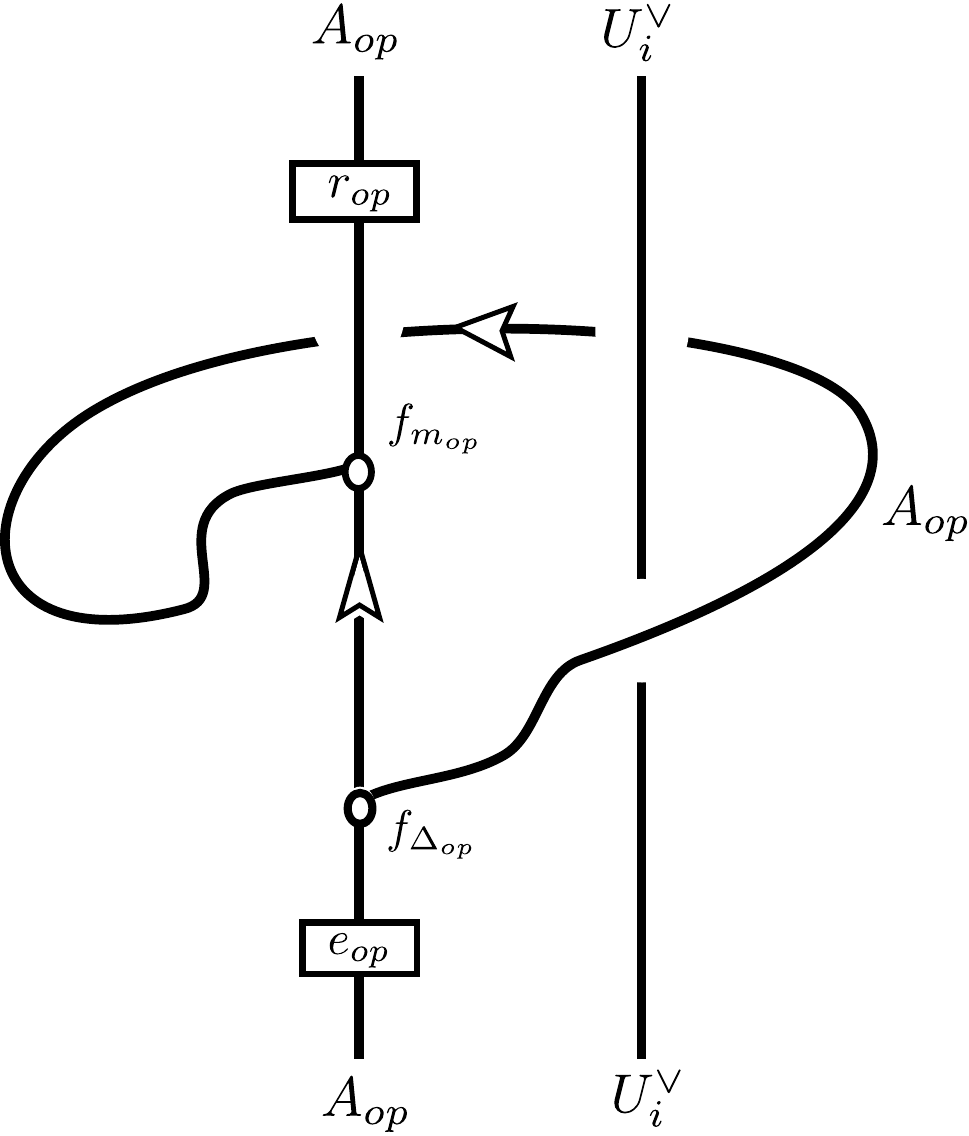} ,\\
       \end{align*}
       where step $(3)$ follows from (\ref{fpol}) and step $(7)$ follows from the first identity in (\ref{fpor}). 
       Continuing from the last graph in above equation, we also have the following sequence of identities: 
        \begin{align*}
         &\frac{\dim(U_i)}{\sqrt{\text{Dim}\mathcal{C}}}\figbox{0.4}{166}\overset{\Rr31}{=}\sum_{\alpha}\figbox{0.6}{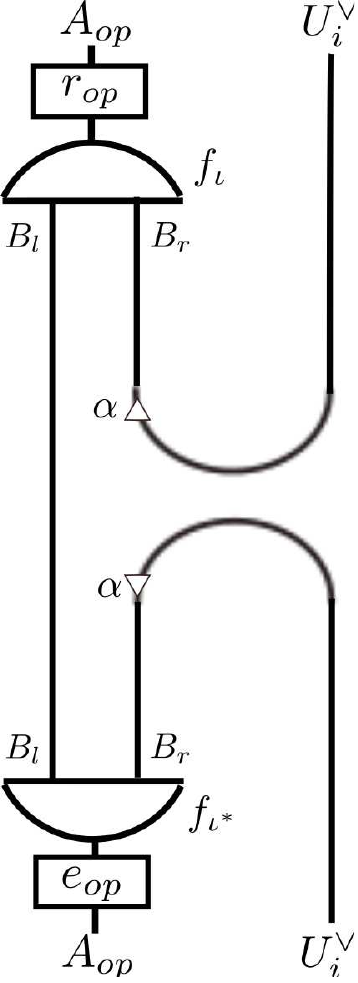} \overset{(2)}{=}\sum_{\alpha}\figbox{0.6}{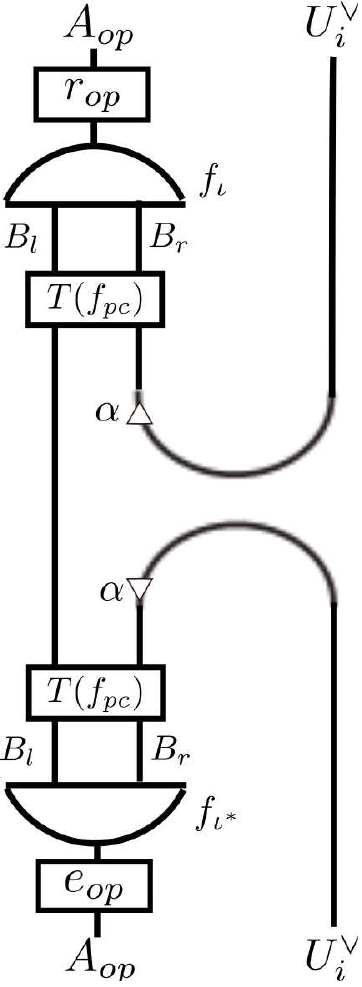}
         \\
        \overset{(3)}{=}&\sum_{\alpha}\figbox{0.6}{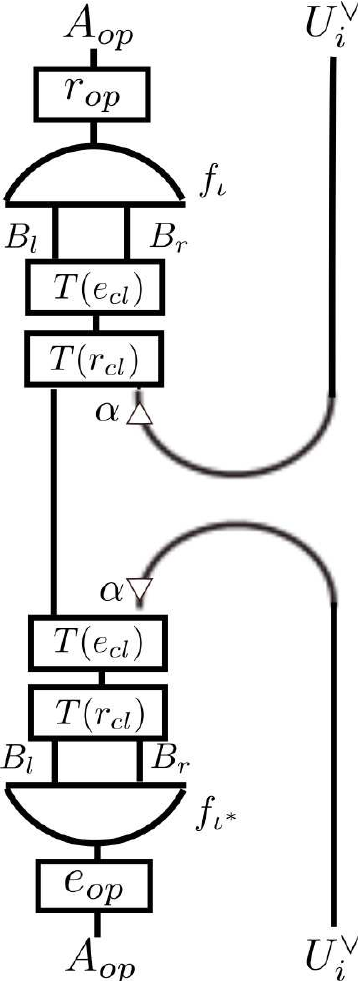}
       \overset{(4)}{=}\sum_{n=1}^{N}\sum_{\beta}\figbox{0.6}{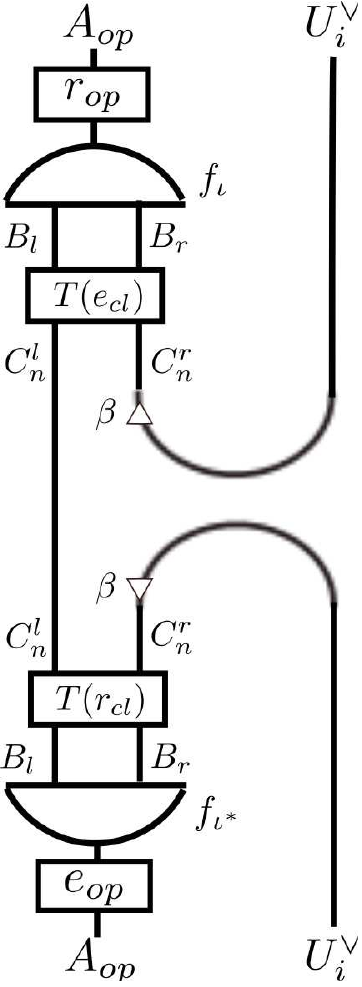} \overset{(5)}{=}\sum_{n=1}^{N}\sum_\beta\figbox{0.6}{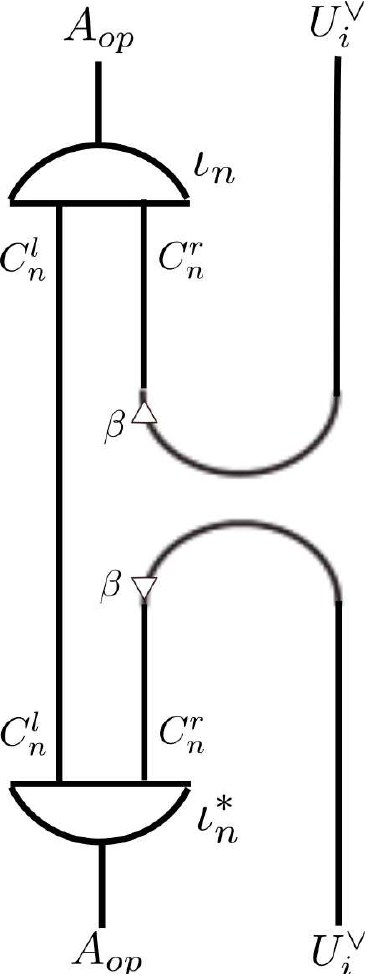}
      \end{align*}
where step (2) follows from (\ref{eq:iota-pc}), (\ref{eq:iota*-pc});
step (3) follows from the definition of $T(f_{pc})$;
step (4) follows from (\ref{eq:BBA}) in the following Lemma \ref{lemma: morphism identity};
step (5) follows from the definition of $\iota_n$ and $\iota^*_n$.
Combining the above two sequence of identities, we obtain the Cardy condition (\ref{eq:cardy-C}).

\end{itemize}

 \begin{lemma}   \label{lemma: morphism identity} {\rm
\begin{equation} \label{eq:BBA}
\begin{aligned}
\sum_{\alpha}\ 
\Bigg(\ \figbox{0.4}{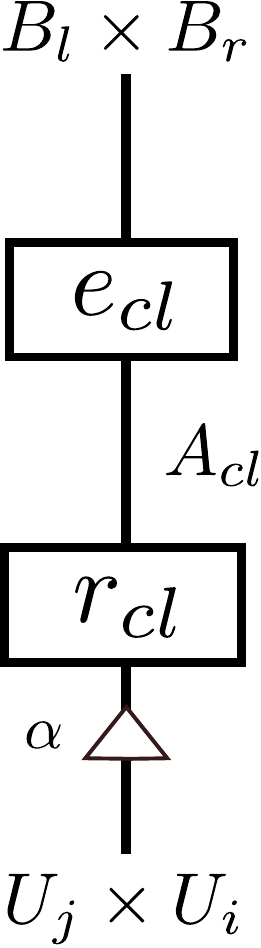}\otimes_\Cb \figbox{0.4}{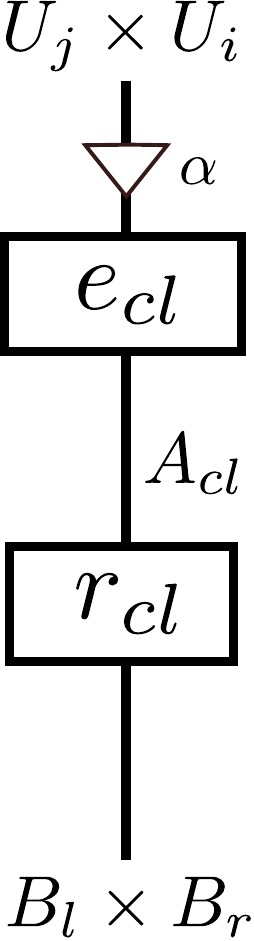}\ \Bigg)
=\sum_\beta\Bigg(\ \figbox{0.4}{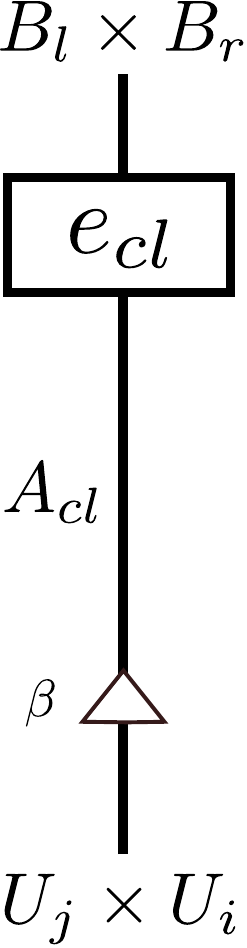}\otimes_\Cb \figbox{0.4}{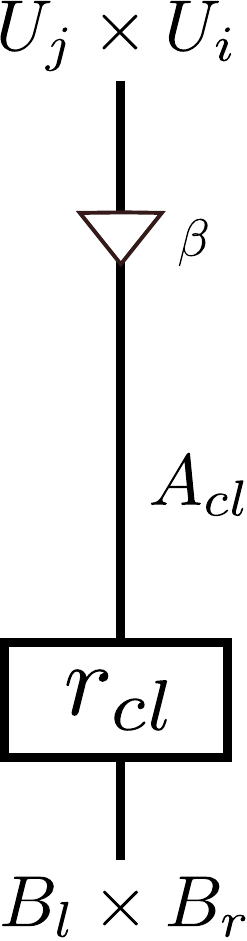}\ \Bigg)
\end{aligned}
\end{equation}
}
\end{lemma}

\pf
Both sides of (\ref{eq:BBA}) belong to the vector space
$$\oplus_{i,j} \Hom_{\mathcal{C}\boxtimes\mathcal{C}}(U_i\times U_j,B_l\times B_r)\otimes_\Cb
 \Hom_{\mathcal{C}\boxtimes\mathcal{C}}(B_l\times B_r,U_i\times U_j),$$
which is isomorphic to the space $\Hom_{\mathcal{C}\boxtimes\mathcal{C}}(B_l\times B_r,B_l\times B_r)$ via the bijective map $\text{ev}: f\otimes g\mapsto f\circ g$. 
Hence, it is enough to show that the image of both sides of (\ref{eq:BBA}) under the map $\text{ev}$ are the same. This follows from the identity $r_{cl}\circ e_{cl}=\id_{A_{cl}}$. 
\epf

\begin{itemize}\setlength{\leftskip}{-1em}
\item Modular invariance: We omit the details as the proof of modular invariance is similar to the proof of Cardy condition.
\end{itemize}

The uniqueness up to isomorphism of Cardy algebras was already demonstrated in the text below Theorem \ref{thm:cardy-determine}.
\epf

\subsubsection*{Proof of Theorem \ref{thm:cardy-sew}}
Suppose that we are given the data of the second notion in Theorem \ref{thm:cardy-sew}. The proof of Theorem \ref{thm:cardy-determine} shows that the morphisms \erf{eq:falpha-names} define a Cardy algebra.

Conversely, assuming that we are given the data of the first notion in Theorem \ref{thm:cardy-sew}, we
define the morphisms $f_\alpha$ as follows:
\be\begin{array}{ll}\displaystyle
   f_{mo}= \eop \circ m_\text{op} \circ (\rop \otimes \rop),\etb~~~
   f_{\eta o} = \eop \circ \eta_\text{op},
\enl
   f_{\Delta o} = (\eop \otimes \eop) \circ\Delta_\text{op} \circ \rop,\etb~~~
   f_{\eps o}= \eps_\text{op} \circ \rop,
\enl
   f_{mc} = \ecl \circ m_\text{cl}  \circ (\rcl \otimes \rcl),\etb~~~
  f_{\eta c} = \ecl \circ \eta_\text{cl},
\enl
  f_{\Delta c} = (\ecl \otimes \ecl) \Delta_\text{cl} \circ \rcl,\etb~~~
  f_{\eps c} = \eps_\text{cl}  \circ \rcl,~~~
\enl
   f_{\iota} = \eop \cir \iota_\text{cl-op}\cir T(\rcl),\etb~~~
  f_{\iota^*}=T(\ecl)\circ \iota_\text{cl-op}^*\circ r_{op}~.
\eear\labl{eq:falpha-names-2}

We define a natural transformation $C$ by $C(\Xr_\alpha)1=\Psi_\alpha(f_\alpha),\alpha\in\mathcal{F}$.
It is straightforward to check that the properties of a Cardy algebra imply the equalities given in Table
\ref{tab:glue-open}--\ref{tab:glue-genus1}. For example, the condition on $f_\alpha$ associated to $\Rr 1$ and $\Rr 9$  are checked below: 
\begin{eqnarray}
\Rr 1:\,\,\,     f_{mo}\circ (f_{\eta o}\otimes \text{id}_{B_{op}})&=&
      e_{op}\circ m_{op}\circ(r_{op}\otimes r_{op})((e_{op}\circ\eta_{op})\otimes \id_{B_{op}})\nn
           &=& e_{op}\circ m_{op}\circ(\eta_{op}\otimes \id)\circ(\id_{\mathbf{1}_\mathcal{C}}\otimes \rop)
           \overset{(1)}{=}\eop\circ \rop = f_{po}, \nn
\Rr 9:\quad       f_{mo}\circ f_{\Delta o}
            &=&(\eop\otimes \eop)\circ\Delta_{op}\circ r_{op}\circ e_{op}\circ m_{op}\circ(r_{op}\otimes r_{op})\nn
            &=&(\eop\otimes \eop)\circ\Delta_{op}\circ  m_{op}\circ(r_{op}\otimes r_{op})\nn
            &\overset{(2)}{=}&(\id_{B_{op}}\otimes e_{op})\circ(\id_{B_{op}}\otimes m_{op})\circ(\id_{B_{op}}\otimes r_{op}\otimes r_{op})  \nn
            &&\hspace{1cm}    \circ   (e_{op}\otimes e_{op}\otimes \id_{B_{op}})\circ(\Delta_{op}\otimes\id_{B_{op}})
                            \circ(\rop\otimes \id_{B_{op}})\nn
            &=&(\id_{B_{op}}\otimes f_{mo})\circ(f_{\Delta}\otimes \id_{B_{op}}),   \nonumber   
\end{eqnarray}
where step (1) follows from the unit property of $\Aop$; step (2) follows from the Frobenius condition of $\Aop$. The proof of the conditions on $f_\alpha$ associated to $\Rr29,\Rr31,\Rr32$ is similar to that of the center condition, the Cardy condition and the modular invariance in the proof of Theorem \ref{thm:cardy-determine}. 

This in turn shows that condition \erf{eq:sewing-cond} in Theorem \ref{thm:gen-rel} is satisfied. Thus $C$ defines a natural transformation from $\One$ to $\Bl$.

It is easy to see that the two constructions are inverse to each other.
\epf

\appendix

\sect{Proof of Theorem \ref{thm:gen-rel}}\label{app:proof-nat}

In this appendix we will abbreviate $C \equiv C(\{ c_\alpha | \alpha \In \Sc \} )$. One direction of the equivalence asserted in Theorem \ref{thm:gen-rel} is trivial, so let us do this one first.

\medskip

\noindent
{\em Proof of Theorem \ref{thm:gen-rel}, part 1:}
\\[.3em]
Suppose that there is a $\mu \in \text{Nat}_\otimes(U,F)$ such that $\mu^\delta = C$. Then \erf{eq:sewing-cond} holds because
 the forgetful functor takes $\Yr^d_{\Rr n,l}$ and $\Yr^d_{\Rr n,r}$ to the same object of $\WSh$, namely to $\Yr_{\Rr n}$, i.e.
$$
  C_{\Yr^d_{\Rr n,l}} = \mu^\delta_{\Yr^d_{\Rr n,l}} = \mu_{\Yr_{\Rr n}}=\mu^\delta_{\Yr^d_{\Rr n,r}}=C_{\Yr^d_{\Rr n,r}}~.
$$
This completes part 1 of the proof of Theorem \ref{thm:gen-rel}.

\medskip

The other direction is more involved and we will prepare the proof with a series of lemmas. We are given $C$ on the generators $\{\Xr_\alpha|\alpha\in\mathcal{F}\}$ such that the condition (\ref{eq:sewing-cond}) is satisfied. We want to define $C$ on a general world sheet $\Xr$ by fixing a
morphism $\varpi : \Xr_{\alpha_1} \otimes \cdots \otimes \Xr_{\alpha_m} \rightarrow \Xr$ and demanding that the diagram
  \be
  \xymatrix{
  U(\Xr_{\alpha_1}) \otimes \cdots \otimes U(\Xr_{\alpha_m})
  \ar[d]^{C_{\Xr_{\alpha_1}}\otimes \cdots \otimes C_{\Xr_{\alpha_m}}}
  \ar[rr]^{\hspace{3em}U(\varpi)} && U(\Xr) \ar[d]^{C_{\Xr}} \\
  F(\Xr_{\alpha_1}) \otimes \cdots \otimes F(\Xr_{\alpha_m}) \ar[rr]^{\hspace{3em}F(\varpi)} && F(\Xr)}
  \ee
commutes. However, for this to be a consistent prescription, the composition
 $F(\varpi) \circ (C_{\Xr_{\alpha_1}}\otimes \cdots \otimes C_{\Xr_{\alpha_m}}) \circ U(\varpi)^{-1}$ has to be
independent of the particular choice of $\alpha_1,\dots,\alpha_m$ and of $\varpi$. This is indeed the case,
provided that \erf{eq:sewing-cond} holds. However, we will need a few lemmas before proving this statement.

\medskip

In a symmetric monoidal category the symmetric braiding $c$ provides, for each permutation $\sigma \in S_n$,
an isomorphism $\Pi_\sigma : V_1 \otimes \cdots \otimes V_n \rightarrow V_{\sigma(1)} \otimes \cdots \otimes V_{\sigma(n)}$.
 A symmetric monoidal functor preserves $\Pi_\sigma$. In particular we have
\be
  U(\Pi_\sigma^{\WSh}) = \Pi_\sigma^{\Vect} \quad \text{and} \quad
  F(\Pi_\sigma^{\WSh}) = \Pi_\sigma^{\Vect} ~,
\labl{eq:UF-permute}
where $\Pi_\sigma^{\WSh}$ and $\Pi_\sigma^{\Vect}$ act by permuting factors in the disjoint union of world sheets and
 in the tensor product of vector spaces, respectively. From hereon we will drop the superscripts $\WSh$ and $\Vect$ on $\Pi_\sigma$.

\begin{lemma}    \label{lem:inv-permute}  {\rm
Let $\sigma \in S_m$ and consider the world sheet $\Xr$ with a decomposition
$$
\Xr_1^d = (\Xr,(\alpha_1,\cdots,\alpha_m),\varpi).
$$ 
If
$\Xr_2^d = (\Xr,(\alpha_{\sigma(1)},\cdots,\alpha_{\sigma(m)}),
\varpi \circ \Pi_\sigma^{-1})$ is another decomposition of $\Xr$, we have $C_{\Xr_1^d}=C_{\Xr_2^d}$.
}
\end{lemma}
\pf
This follows from the definition of $C$ together with \erf{eq:UF-permute}. Write
\begin{equation}
\begin{aligned}
   C_{\Xr_2^d} \overset{(1)}{=}&\
  F(\varpi \circ \Pi_\sigma^{-1}) \circ \big( c_{\alpha_{\sigma(1)}}
    \otimes \cdots \otimes c_{\alpha_{\sigma(m)}}  \big)
    \circ U(\varpi \circ \Pi_\sigma^{-1})^{-1}\\
  \overset{(2)}{=}&\
  F(\varpi) \circ \Pi_\sigma^{-1} \circ \big( c_{\alpha_{\sigma(1)}}
    \otimes \cdots \otimes c_{\alpha_{\sigma(m)}}  \big)
    \circ \Pi_\sigma \circ U(\varpi)^{-1}\\
  \overset{(3)}{=}&\
  F(\varpi) \circ \big( c_{\alpha_{1}}
    \otimes \cdots \otimes c_{\alpha_{m}} \big)
    \circ U(\varpi)^{-1}\\
  \overset{(4)}{=}&\  C_{\Xr_1^d}
\end{aligned}
\end{equation}
Here step (1) is the definition \erf{eq:nat-C-def} of $C_{\Xr^d}$, step (2) is \erf{eq:UF-permute}, in step (3) we permute the factors of the tensor product to move $\Pi_\sigma$ next to $\Pi_\sigma^{-1}$, and step (4) is again \erf{eq:nat-C-def}.
\epf

\medskip

The next lemma is instrumental in using the relations R1--R32 to prove identities among the $C_\Xr$.

\begin{lemma}  \label{lem:change-part}  {\rm
Let $\Xr, \Yr, \Rr \in \WSh$ and let $\xi : \Yr \otimes \Rr \rightarrow \Xr$ be a morphism in $\WSh$.
Let $\Rr^d = (\Rr,(\gamma_1,\dots,\gamma_k), \rho)$ be a decomposition of $\Rr$ and let
\be
  \Yr^d_1 = (\Yr,(\alpha_1,\dots,\alpha_m),\varpi_1)
  \quad , \quad
  \Yr^d_2 = (\Yr,(\beta_1,\dots,\beta_n),\varpi_2)
\ee
be two decompositions of $\Yr$. This gives rise to two decompositions of $\Xr$,
\bea
  \Xr^d_1 = (\Xr, ( \alpha_1,\dots,\alpha_m,\gamma_1,\dots,\gamma_k ), \xi \circ (\varpi_1 \otimes \rho))
  ~,\enl
  \Xr^d_2 = (\Xr, ( \beta_1,\dots,\beta_n,\gamma_1,\dots,\gamma_k ), \xi \circ (\varpi_2 \otimes \rho))~.
\eear\ee
If $C_{\Yr_1^d}=C_{\Yr_2^d}$, then we also have $C_{\Xr_1^d}=C_{\Xr_2^d}$.
}
\end{lemma}
\pf
This is again a simple sequence of identities using relation \erf{eq:sewing-cond} in the middle, 
\begin{eqnarray}
  C_{\Xr_1^d}
  &\overset{(1)}{=}&
  F(\xi \circ (\varpi_1 \otimes \rho)) \circ
  \big( c_{\alpha_1}
    \otimes \cdots \otimes c_{\gamma_k}  \big)
    \circ U(\xi \circ (\varpi_1 \otimes \rho))^{-1}   \nn
  &\overset{(2)}{=}&
  F(\xi \circ (\id_{\Xr} \otimes \rho)) \circ ( C_{\Yr_1^d}
  \otimes
  c_{\gamma_1} \otimes \cdots \otimes c_{\gamma_k} )
  \circ U(\xi \circ (\id_{\Xr} \otimes \rho))^{-1}  \nn
 &\overset{(3)}{=}&
  F(\xi \circ (\id_{\Xr} \otimes \rho)) \circ ( C_{\Yr_2^d}
  \otimes
  c_{\gamma_1} \otimes \cdots \otimes c_{\gamma_k})
  \circ U(\xi \circ (\id_{\Xr} \otimes \rho))^{-1}  \nn
  &\overset{(4)}{=}&
  C_{\Xr_2^d}
\end{eqnarray}
Here step (1) is definition \erf{eq:nat-C-def}, step (2) uses monoidality of $F$ and $U$, step (3) is the identity $C_{\Yr_1^d}=C_{\Yr_2^d}$, and step (4) is just steps (1) and (2) in reverse order.
\epf

\begin{cor}  \label{cor:inv-Rn}  {\rm
Suppose that $C$ obeys \erf{eq:sewing-cond}.
Let $\Xr$ be a world sheet, set $\Zr = \Xr_{\gamma_1} \otimes \cdots \otimes \Xr_{\gamma_k}$, and suppose that there is
a morphism $\xi : \Yr_{\Rr n} \otimes \Zr \rightarrow \Xr$. We can decompose the world sheet $\Xr$ as
$$
  \Xr_1^d = \big(\,\Xr \,,\, (\alpha_1^{(n)},\cdots,\alpha^{(n)}_{l(n)},
  \gamma_1,\dots,\gamma_k),\, \xi \circ (\varpi_l^{(n)} \otimes \id_\Rr)
  \,\big)
$$
and as
$$
  \Xr_2^d = \big(\,\Xr \,,\, (\beta_1^{(n)},\cdots,\beta^{(n)}_{r(n)},
  \gamma_1,\dots,\gamma_k)\,,\, \xi \circ (\varpi_r^{(n)} \otimes \id_\Rr)
  \,\big) ~,
$$
and we have $C_{\Xr_1^d}=C_{\Xr_2^d}$.
}
\end{cor}

Relations R10--13 and R21--24 allow us to omit strips and cylinders attached to state boundaries of generating world sheets
in certain positions. The next two lemmas prove that strips and cylinders can be omitted from generators in arbitrary positions.

\begin{lemma}    \label{lem:inv-omit-cyl}  {\rm
Suppose that $C$ obeys \erf{eq:sewing-cond}.
\\[.3em]
{\it (i) open version:} For $\alpha \in \{mo, \Delta o, \eta o, \eps o, \iota, \iota^*, po\}$ let $\zeta$ be a morphism
 $\Xr_\alpha \otimes \Xr_{po} \rightarrow \Xr_\alpha$. There exists a morphism $\chi : \Xr_\alpha \rightarrow \Xr_\alpha$,
 such that for the two decomposed world sheets
$\Xr_1^d = (\Xr_\alpha,(\alpha,po),\zeta)$
and
$\Xr_2^d = (\Xr_\alpha,(\alpha),\chi)$ we have
$C_{\Xr_1^d}=C_{\Xr_2^d}$.
\\[.3em]
{\it (ii) closed version:}
For $\alpha \in \{mc, \Delta c, \eta c, \eps c, \iota, \iota^*, pc\}$ let $\zeta$ be a morphism
 $\Xr_\alpha \otimes \Xr_{pc} \rightarrow \Xr_\alpha$. There exists a morphism $\chi : \Xr_\alpha \rightarrow \Xr_\alpha$,
such that for the two decomposed world sheets
$\Xr_1^d = (\Xr_\alpha,(\alpha,pc),\zeta)$
and
$\Xr_2^d = (\Xr_\alpha,(\alpha),\chi)$ we have
$C_{\Xr_1^d}=C_{\Xr_2^d}$.
}
\end{lemma}
\pf
We will only treat the case $\alpha = mo$ in some detail. The remaining cases work similarly. There are three possible sewings for $\Xr_{mo} \otimes \Xr_{po}$:
\begin{center}
 \includegraphics[width=90mm]{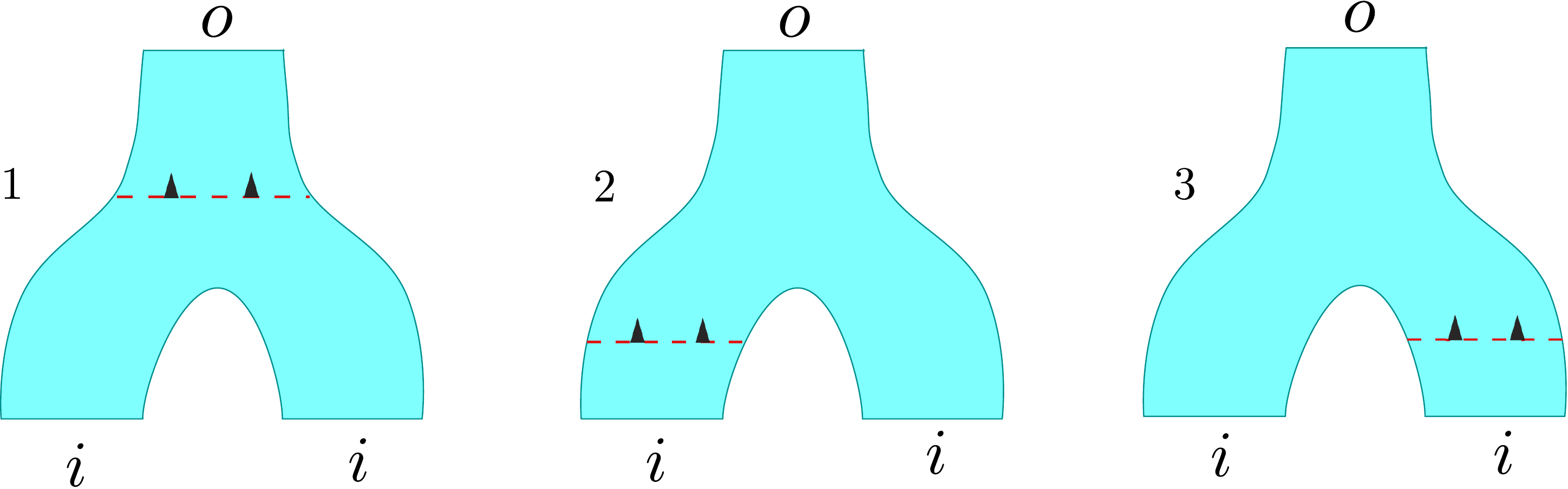}
\end{center}
In sewing 1) we can find a morphism of decomposed world sheets
$\xi : \Yr^d_{\Rr10,l} \rightarrow \Xr_1^d$, in other
words we can write $\zeta = \xi \circ \varpi^{(10)}_l$. The
identity $C_{\Xr_1^d}=C_{\Xr_2^d}$ is then the statement of Corollary
\ref{cor:inv-Rn} if we take $\chi = \xi \circ \varpi^{(10)}_r$
and $\Rr = \emptyset$.

In sewing 2) one considers the following series of decomposed world sheets
\begin{center}
\includegraphics[width=150mm]{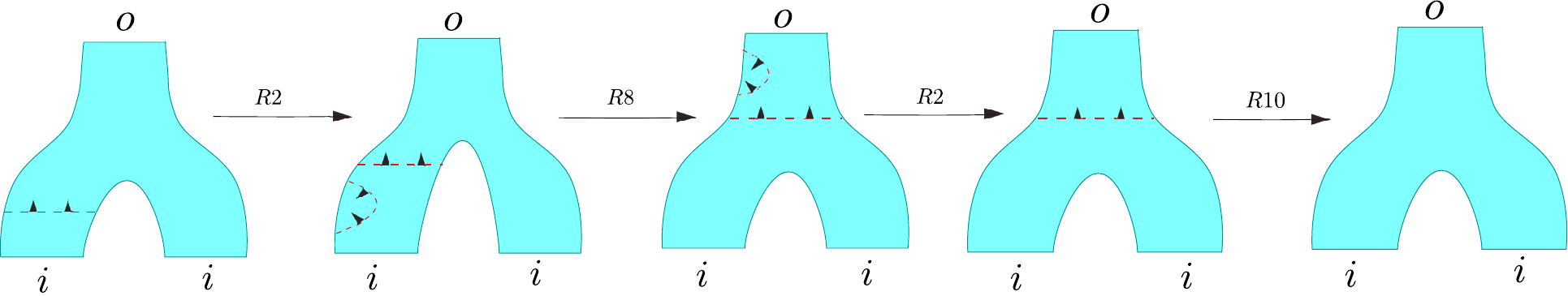}
\end{center}
The first decomposed world sheet is just $\Xr_1^d = (\Xr_{mo},(mo,po),\zeta)$. We would like to apply relation R2 to the generator $\Xr_{po}$ in this decomposition via Corollary \ref{cor:inv-Rn}. Indeed, because by convention $\varpi^{(2)}_r=\id_{\Xr_{po}}$, by Corollary \ref{cor:inv-Rn} we obtain $\Yr_1^d=(\Xr_{mo},(\Delta o, \eps o, mo), \zeta \circ (\varpi^{(2)}_l \otimes \id_{\Xr_{mo}}))$ together with the identity $C_{\Xr^d_1} = C_{\Yr^d_1}$. Then one uses Lemma \ref{lem:inv-permute} to obtain
$\Yr^d_{1b} = (\Xr_{mo},(mo, \Delta o, \eps o), \xi)$ for some $\xi$, together with the identity $C_{\Yr^d_1} = C_{\Yr^d_{1b}}$.
Now apply Corollary \ref{cor:inv-Rn} for relation R8  with auxiliary world sheet $\Rr = \Xr_{\eps o}$. This gives $\Yr^d_2 = (\Xr_{mo},(mo, \Delta o, \eps o), \xi')$ and the identity
$C_{\Yr^d_{1b}} = C_{\Yr^d_{2}}$, and so forth. Altogether this shows that $C_{\Xr_1^d}=C_{\Xr_2^d}$ in sewing 2).

Sewing 3) works in the same way, and the remaining identities in part (i) and (ii) of the statement can be treated analogously. For $\alpha=\iota$ and $\alpha=\iota^*$ one needs to employ R26 and R27 as an intermediate step, for $\alpha=pc$ one uses R1 and one of sewing 1)--3) treated above.
\epf

\medskip

The above lemmas allow us to employ relations R1--R32, and also to remove all cylinders and strips from a decomposed world sheet
(except for one in the case that $\Xr = \Xr_{po}$ or $\Xr = \Xr_{pc}$). The next two ingredients we need are how $C$ behaves
 if two decomposed world sheets are related by a morphism homotopic to the identity, and how the mapping class group of a
generating world sheet acts on $C$.

\begin{lemma}   \label{lem:inv-homotop-id} {\rm
Let $\Xr$ be a world sheet and let $\xi : \Xr \rightarrow \Xr$ be a morphism of world sheets which is homotopic to the identity. Given two decompositions of $\Xr$ of the form $\Xr_1^d = (\Xr,(\alpha_1,\cdots,\alpha_m),\varpi)$ and $\Xr_2^d = (\Xr,(\alpha_1,\cdots,\alpha_m),\xi \circ \varpi)$, we have $C_{\Xr_1^d}=C_{\Xr_2^d}$.
}
\end{lemma}
\pf
Note that $\xi$ is also a morphism from $\Xr_1^d$ to $\Xr_2^d$ in $\WSh^d$. By Lemma \ref{lem:all-nat-C}, $C$ is
a monoidal natural transformation from $U^\delta$ to $F^\delta$. By the definition of the functors $U^\delta$ and
$F^\delta$, commutativity of the square \erf{eq:C-is-natural} implies that
  \be
  \xymatrix{
  U(\Xr) \ar[d]^{C_{\Xr_1^d}} \ar[r]^{U(\xi)} & U(\Xr) \ar[d]^{C_{\Xr_2^d}} \\
  F(\Xr) \ar[r]^{F(\xi)} & F(\Xr)}
  \ee
commutes. Since $U$ and $F$ are constant on homotopy classes of homeomorphisms of world sheets, and since $\xi$
is homotopic to the identity, we have $U(\xi) = \id_{U(\Xr)}$ and $F(\xi) = \id_{F(\Xr)}$. Therefore the above
diagram implies $C_{\Xr_1^d}=C_{\Xr_2^d}$.
\epf

\begin{lemma}   \label{lem:inv-mapgen}  {\rm
Suppose that $C$ obeys \erf{eq:sewing-cond}.
Let $\varpi : \Xr_\alpha \rightarrow \Xr_\alpha$ be a morphism of world sheets, not necessarily homotopic to the
identity. For the two decomposed world sheets
$\Xr_1^d = (\Xr_\alpha,(\alpha),\id)$ and $\Xr_2^d = (\Xr_\alpha,(\alpha),\varpi)$, we have
$C_{\Xr_1^d}=C_{\Xr_2^d}$.
}
\end{lemma}
\pf
Let $E(\alpha) = \Hom(\Xr_\alpha,\Xr_\alpha)$ be the group of automorphisms of $\Xr$ and denote by $E_0(\alpha)$
the normal subgroup of all elements of $E(\alpha)$ which are homotopic to $\id_{\Xr_\alpha}$.

If $\varpi \in E_0(\alpha)$, then $C_{\Xr_1^d}=C_{\Xr_2^d}$ by Lemma \ref{lem:inv-homotop-id}. For $\alpha \in \{\eta o,\eps o,\eta c,\eps c\}$, we have $E_0(\alpha)=E(\alpha)$ because the homeomorphism $h$ in
$\varpi = (\emptyset,h)$ has to leave the unique state boundary of $\dot\Xr$ point-wise fixed. For
$\alpha \in \{mo,\Delta o,po\}$, a nontrivial permutation of the open state boundaries is forbidden because $h$ must preserve the partition into in-coming and out-going state boundaries.
Therefore, we obtain $E_0(\alpha)=E(\alpha)$ in these cases.

For $\alpha \in \{mc, \Delta c, pc, \iota, \iota^*\}$, the mapping class group $E(\alpha)/E_0(\alpha)$ is non-trivial. The classes of the elements $(\emptyset,h)$ of $E(\alpha)$ depicted below (Dehn twists) give a set of generators for each of the $E(\alpha)/E_0(\alpha)$, see e.g.\ \cite[Sect.\,5.1]{baki},
\begin{flalign}
  \Xr_{pc} :
  h(pc) = {\figbox{0.2}{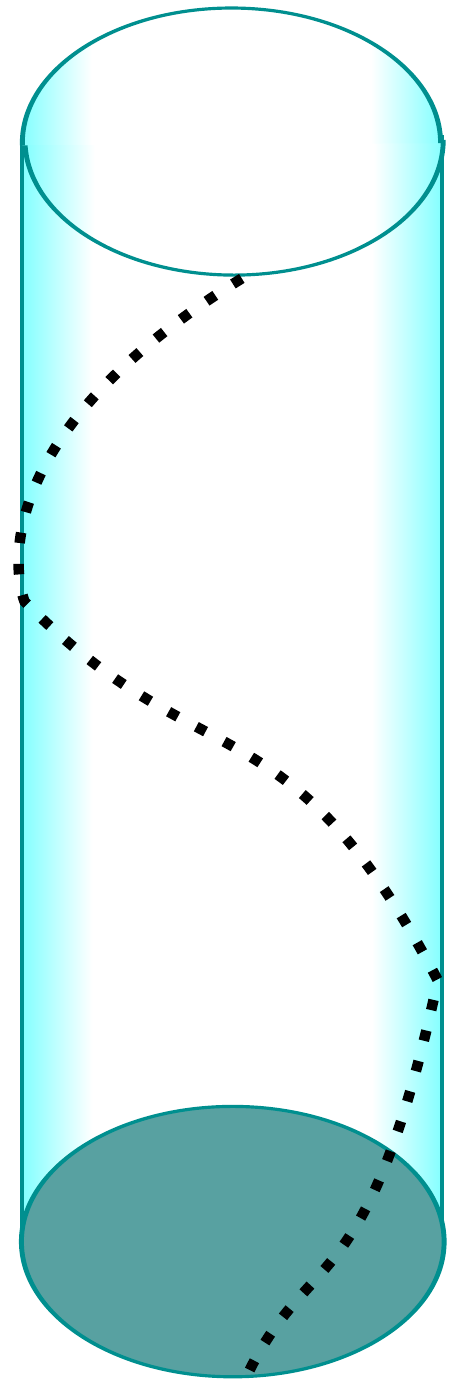}}~,\quad\quad\quad\quad
  h(\iota) ={\figbox{0.2}{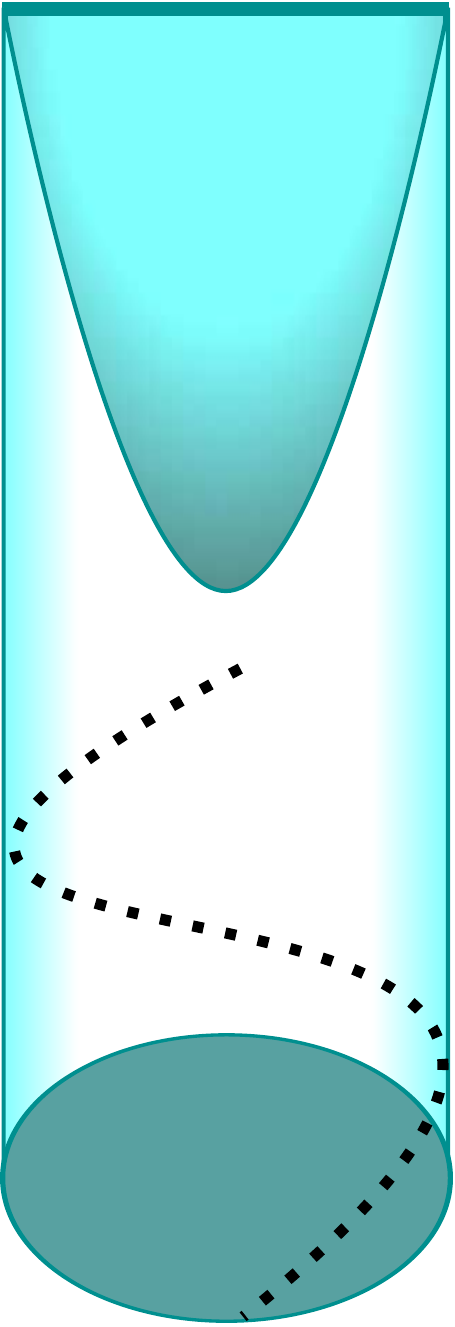}}~,  \nonumber
\end{flalign}
\begin{equation}
\begin{aligned}
  \Xr_{mc} :&h(mc)_1 ={\figbox{0.18}{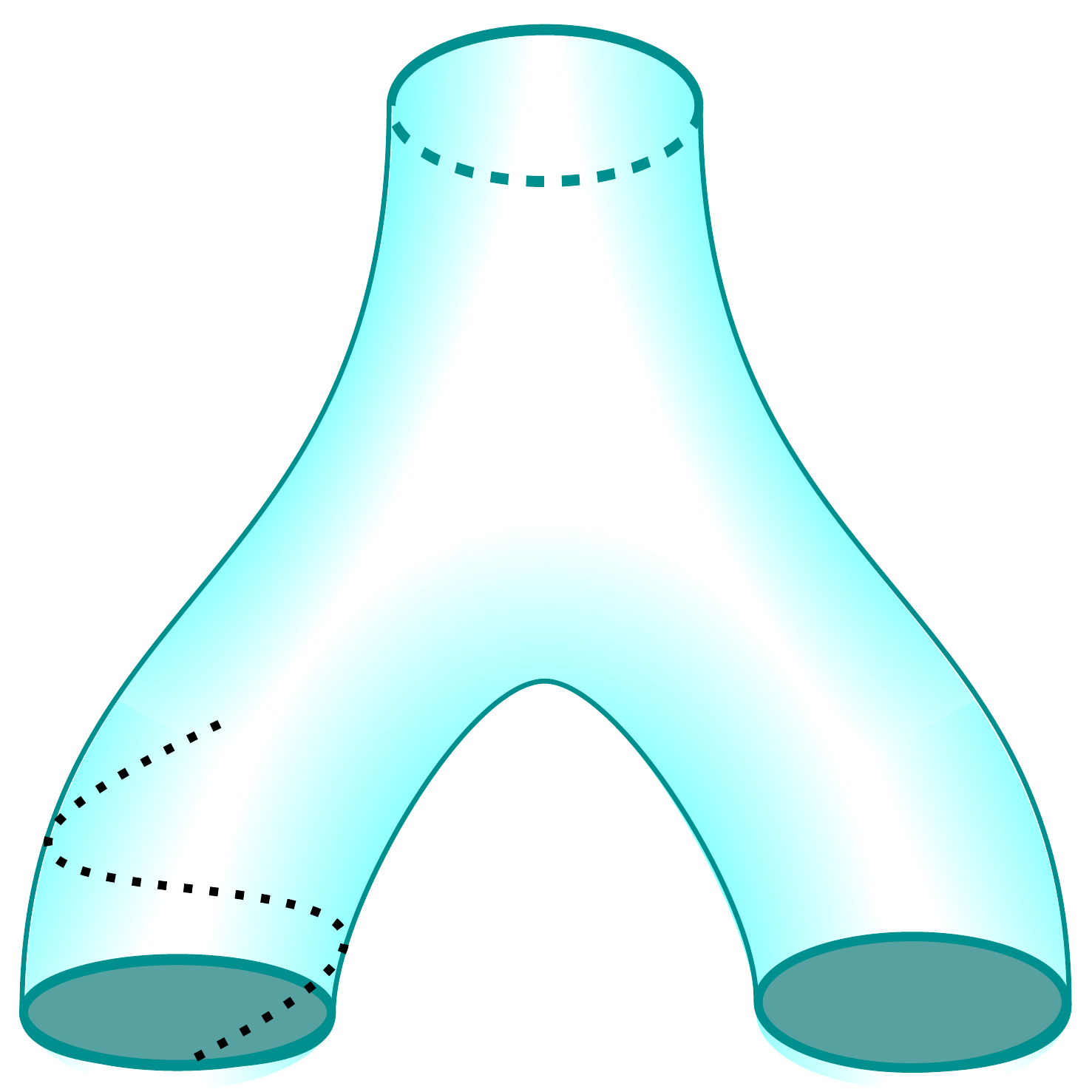}},\
  h(mc)_2 = {\figbox{0.18}{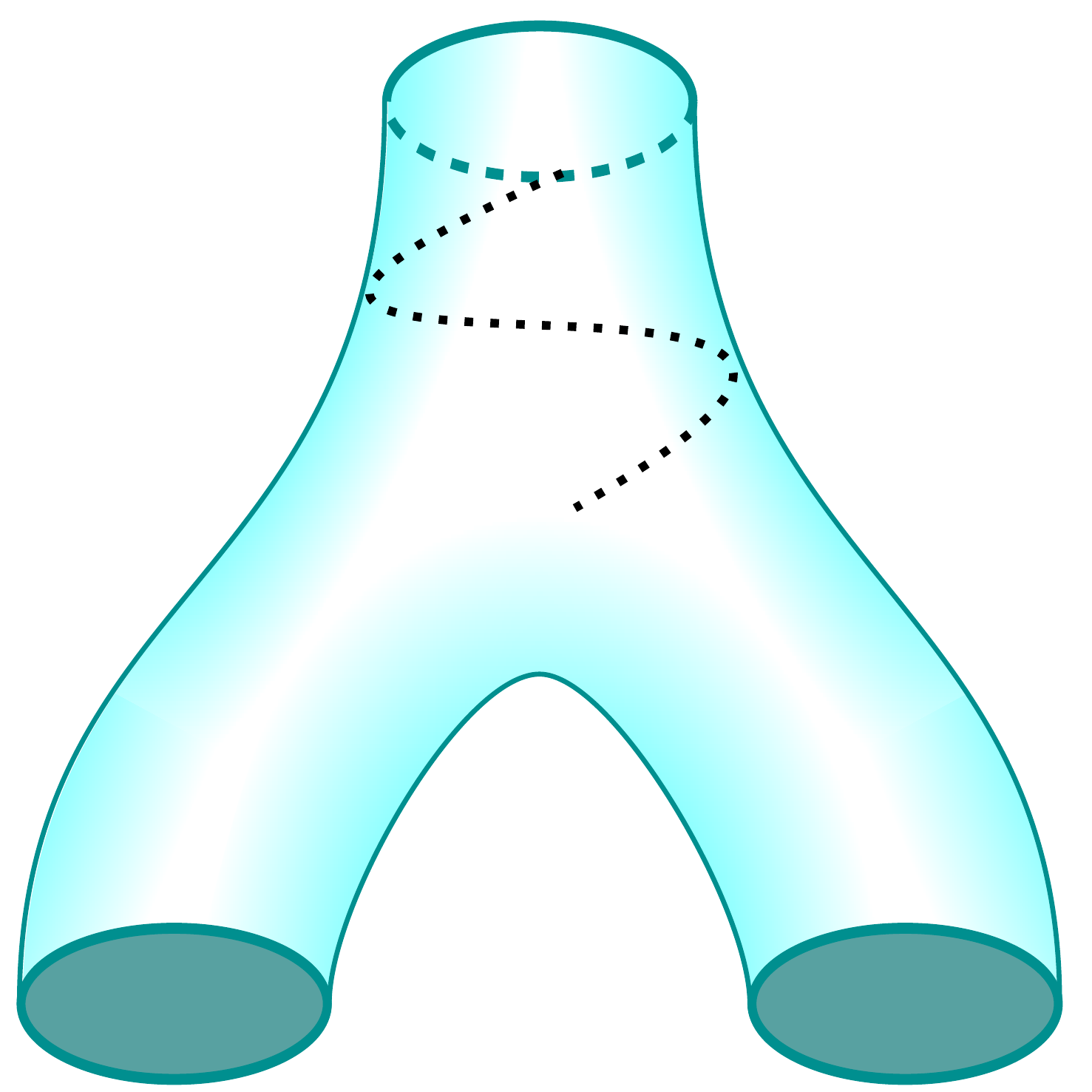}},\
  h(mc)_3 = {\figbox{0.18}{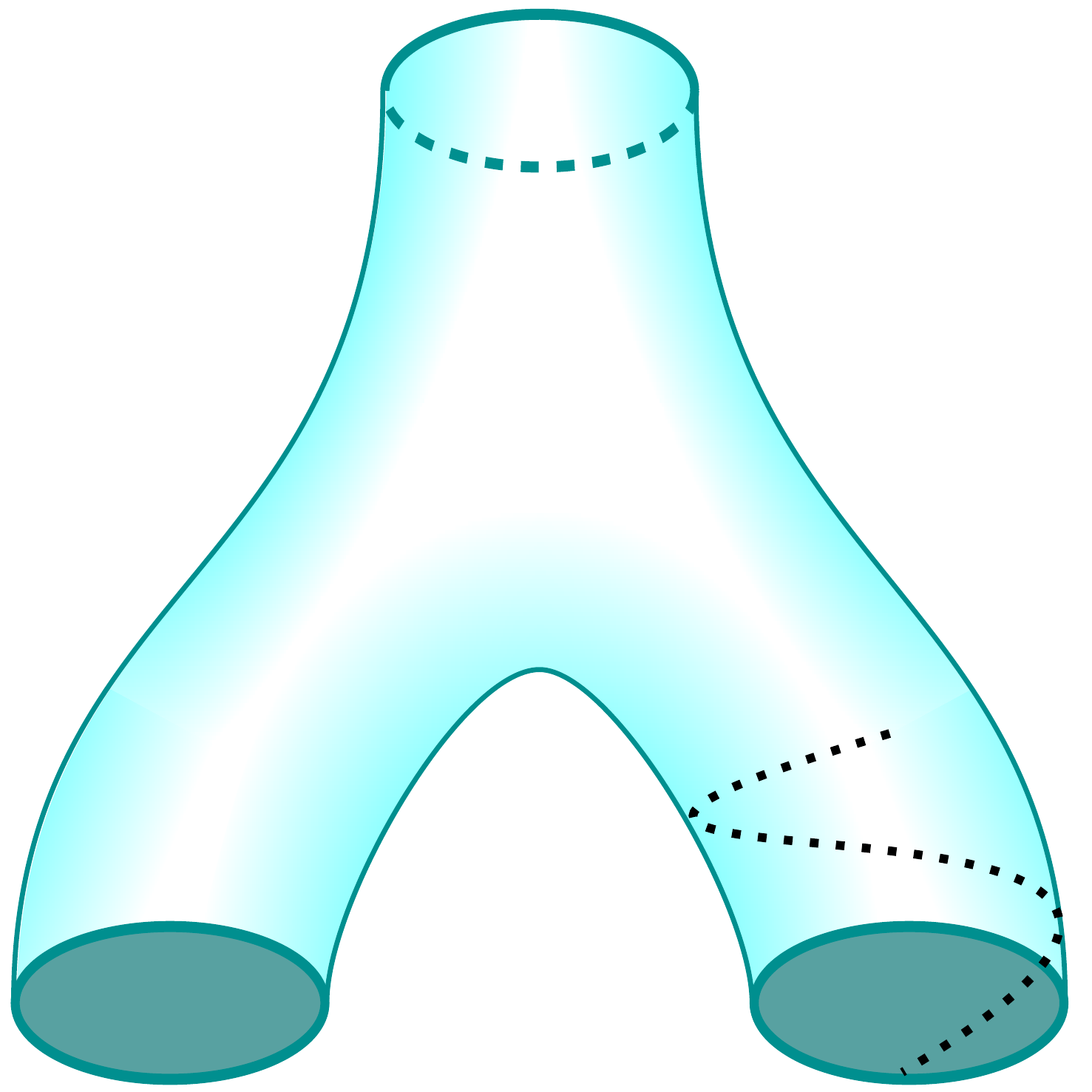}}~.  \nonumber\\
   & h(mc)_4 ={\figbox{0.26}{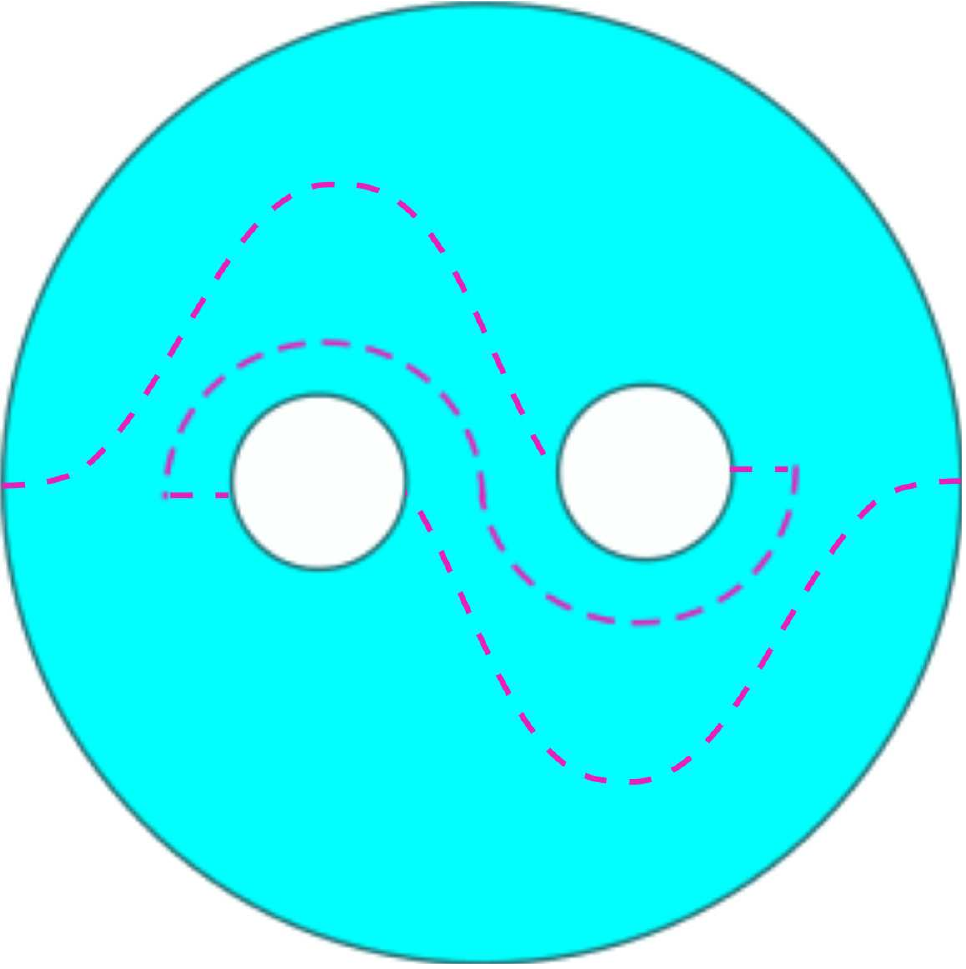}},
\end{aligned}
\end{equation}
For $\Xr_{\Delta c}$ and $\Xr_{\iota^*}$, we can take the same set of generators as for $\Xr_{mc}$ and $\Xr_{\iota}$ but with the
labelling of in-coming and out-going state boundaries interchanged.
The Lemma is proved if we can show that $C_{\Xr_1^d}=C_{\Xr_2^d}$ with $\varpi = (\emptyset,h)$ for every $h$ in the above list.

By assumption $C$ obeys \erf{eq:sewing-cond}. Then for $\alpha=pc$ and $h=h(pc)$, $C_{\Xr_1^d}=C_{\Xr_2^d}$ is implied by R25. 
Similarly, for $\alpha=mc$ and $h=h(mc)_4$ the statement follows from R20. The remaining cases can all be proved by the following
method. We only present the details for $\alpha=\iota$.

Our starting point is the decomposed world sheet
$\Xr_1^d = (\Xr_\iota,(\iota),\id)$, and we want to arrive at $\Xr_2^d = (\Xr_\iota,(\iota),\varpi)$ with $\varpi = (\emptyset,h(\iota))$. 

Pick any decomposed world sheet $\Yr_0^d=(\Xr_\iota,(\iota,pc), \eta)$ of the generator $\Xr_\iota$. By Lemma \ref{lem:inv-omit-cyl}\,(ii)
there exists a decomposed world sheet  $\Xr_0^d=(\Xr_\iota,(\iota),\chi)$, such that $C_{\Xr_0^d}=C_{\Yr_0^d}$. 
Let $\Yr_1^d=(\Xr_\iota,(\iota,pc),\chi^{-1}\circ\eta=\zeta)$. Applying Lemma \ref{lem:change-part} (with $R = \emptyset$ and $\xi = \chi^{-1}$), we obtain that $C_{\Yr_1^d}=C_{\Xr_1^d}$. 
By the same argument we also have $C_{\Xr_2^d}=C_{\Yr_2^d}$ with $\Yr_2^d = (\Xr_\iota,(\iota,pc),\varpi \circ \zeta)$.

Let's write the morphism $\zeta$ in the form $\zeta = (\sew,f)$ and denote by $P \subset \Xr_\iota$ the $f$-image of $\Xtil_{pc}$,
 considered as a subset of $\sew(\Xr_\iota \otimes \Xr_{pc})\widetilde{~}$.  We are free to choose the generator $h(\iota)$ such that it is different from the identity only on $P$. It is therefore possible to find a $\xi : \Xr_{pc} \rightarrow \Xr_{pc}$ such that 
\begin{equation} \label{eq:varpi-zeta}
\varpi \circ \zeta = \zeta \circ (\id_{\Xr_{\iota}} \otimes \xi).
\end{equation}

Now recall that R25 states $C_{\Yr^d_{\Rr25,l}} = C_{\Yr^d_{\Rr25,r}}$ for the decomposed world sheets $\Yr^d_{\Rr25,l} = (\Xr_{pc},(pc),\varpi^{(25)}_l=\id)$ and $\Yr^d_{\Rr25,r} = (\Xr_{pc},(pc),\varpi^{(25)}_r)$, where $\varpi^{(25)}_r : \Xr_{pc} \rightarrow \Xr_{pc}$ is chosen to generate the mapping class group $E(pc)/E_0(pc)$ of $\Xr_{pc}$.
Composing both $\varpi^{(25)}_l$ and $\varpi^{(25)}_r$ with 
$(\varpi^{(25)}_r)^{-1}$, we  obtain that for the decomposed world sheets $\Yr^d_l = (\Xr_{pc},(pc),(\varpi^{(25)}_r)^{-1})$ and $\Yr^d_r = (\Xr_{pc},(pc),\id)$ we have (again as a special case of Lemma \ref{lem:change-part} with $R = \emptyset$)
\be
  C_{\Yr^d_l} = C_{\Yr^d_r} ~.
\labl{eq:gen-mcg-aux1}

Since $\varpi^{(25)}_r$ generates $E(pc)/E_0(pc)$ there exists a morphism $\xi_0$ homotopic to the identity and an integer $m$ such that $\xi = \xi_0 \circ (\varpi^{(25)}_r)^m$ (recall (\ref{eq:varpi-zeta})). Altogether we have rewritten the morphism $\varpi \circ \zeta$ in $\Yr_2^d$ as
\be
\Yr_2^d = (\Xr_\iota,(\iota,pc),\zeta \circ (\id_{\Xr_{\iota}} \otimes \xi_0)
 \circ (\id_{\Xr_{\iota}} \otimes (\varpi^{(25)}_r)^m))
\ee
Applying Corollary \ref{cor:inv-Rn} $|m|$ times in the case of the relation R25 or \erf{eq:gen-mcg-aux1}, we see that $C_{\Yr^d_2} = C_{\Yr^d_3}$ for $\Yr^d_3 = (\Xr_\iota,(\iota,pc),\zeta \circ (\id_{\Xr_{\iota}} \otimes \xi_0))$. But $\zeta \circ (\id_{\Xr_{\iota}} \otimes \xi_0)$ is homotopic to $\zeta$ and so $C_{\Yr^d_2} = C_{\Yr^d_1}$. Altogether we have now proved that $C_{\Yr_1^d}=C_{\Yr_2^d}$ which implies $C_{\Xr_1^d}=C_{\Xr_2^d}$.
\epf

\bigskip
The proof of Theorem \ref{thm:gen-rel} needs Morse functions for open/closed world sheets. 
These functions are described in \cite[app.\,A.1 \& A.2]{MS} and \cite[Sect.\,3.3]{LP}; see also \cite{BNR} and references therein for more details on Morse theory on manifolds with boundary. We include the definition of Morse functions here:

\begin{defn}  \label{def: Morse function} {\rm
A {\em Morse function on a world sheet} $\Xr$ is a smooth function $f:\dot\Xr\rightarrow[0,1]
$ such that
\begin{enumerate}
\item All interior critical points of $f$ are non-degenerate, i.e., the Hessian of $f$ at the critical points are non-degenerate.
\item The restriction of $f$ to the physical boundary of $\dot\Xr$ has only discrete critical points, which are all non-degenerate (as functions on
1-manifolds).
\item $f^{-1}(0)$ is the union of the in-coming state boundaries of $\dot\Xr$ and $f^{-1}(1)$ is the union of the out-going state boundaries.
\end{enumerate}
 A Morse function
is called {\em generic} if for all distinct critical points $c_1\not=c_2$ of $f$, we have $f(c_1)\not=f(c_2)$.
}
\end{defn}

The following lemma describes how two generic Morse functions are related to each other if they lie in the same
connected component of generic Morse functions:

\begin{lemma}\label{lem:Morse-family-isotopy}{\rm
Let $f_t:\dot\Xr\rightarrow [0,1],t\in[0,1]$ be a family of generic Morse functions. Then there exist isotopies $\varphi_t:[0,1]\rightarrow [0,1]$ and 
$h_t:\dot\Xr\rightarrow\dot\Xr$ with $\varphi_0$ and $h_0$ identity maps such that $f_0=\varphi_t\circ f_t\circ h_t$ for
$t\in[0,1]$. }
\end{lemma}

We consider all world sheets that allow a generic Morse function with one critical point. Those generators
$\Xr_\alpha$ for $\alpha\in\mathcal{S}\setminus\{pc,po\}$ are such world sheets but not vice versa.
 Consider the following world sheets:
\begin{equation}  \label{eq:X-o4}
\Xr_{o4}={\figbox{0.4}{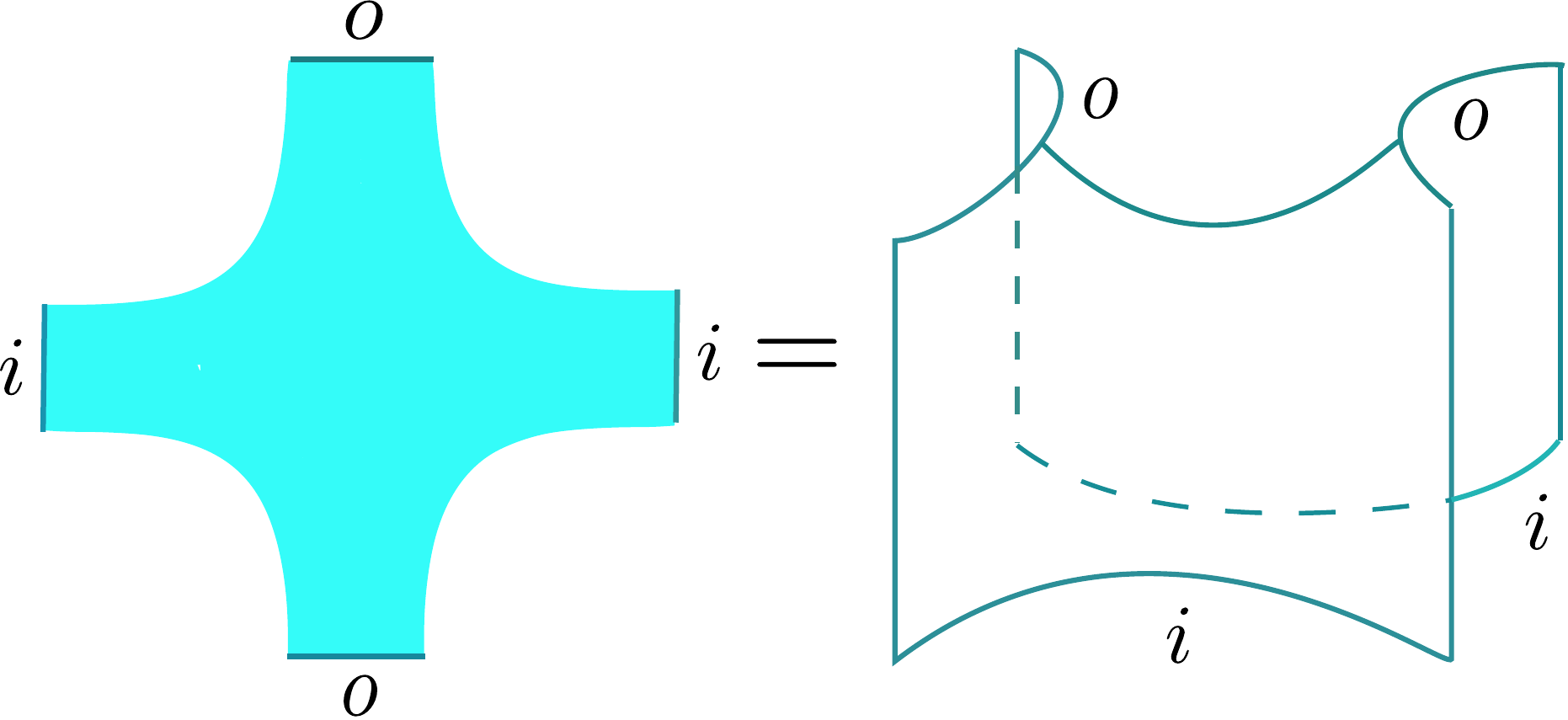}}\ ,\ \ \ \ \ \
\end{equation}
\begin{equation} \label{eq:X-oco}
 \Xr_{oco(\alpha,\beta)}={\figbox{0.33}{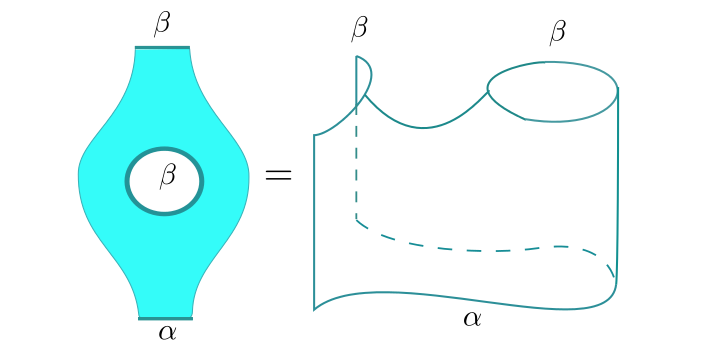}}
\end{equation}
where $\alpha,\beta\in \{i,o\}$ and $\alpha \neq \beta$. The height function in the right picture of (\ref{eq:X-o4}) is
a Morse function on $\Xr_{o4}$ with one critical point. Similarly, if $\alpha=i,\beta=o$, the height function in the right
picture of (\ref{eq:X-oco}) is also a Morse function with one critical point. The set 
\be 
  \mathcal{M} := \{ o4, oco(i,o), oco(o,i) \},
\ee
together with the set $\Sc$, is complete in the following sense. Let $\Xr$ be a world sheet and let $f$ be a generic Morse function on $\Xr$. Let $0<f(c_1)<\cdots<f(c_n)<1$ be the images of the critical points of $f$. If an interval $[s,t]$ contains a single
critical value of $f$, then the preimage $f^{-1}([s,t])$ is a disjoint union of $\Xr_\alpha$ for $\alpha \in \Sc \cup \mathcal{M}$,
where all but one $\alpha$ are in $\{po,pc\}$. This is proved in \cite[app.\,A.2]{MS} and also in \cite[Prop.\,3.8]{LP}. 

\begin{rema} {\rm
In the definition of the set of extended generators $\mathcal{S}\cup\mathcal{M}$, we did not consider the data of ordering of
boundaries of world sheets, since it is not essential in the proof of Theorem \ref{thm:gen-rel}.
}\end{rema}

There are several non-isotopic ways of decomposing  $\Xr_\alpha$'s for $\alpha\in\mathcal{M}$ into two generators from the set $\mathcal{S}$, namely
\begin{equation}\label{pic:decom-extended-gen}
 \figbox{0.3}{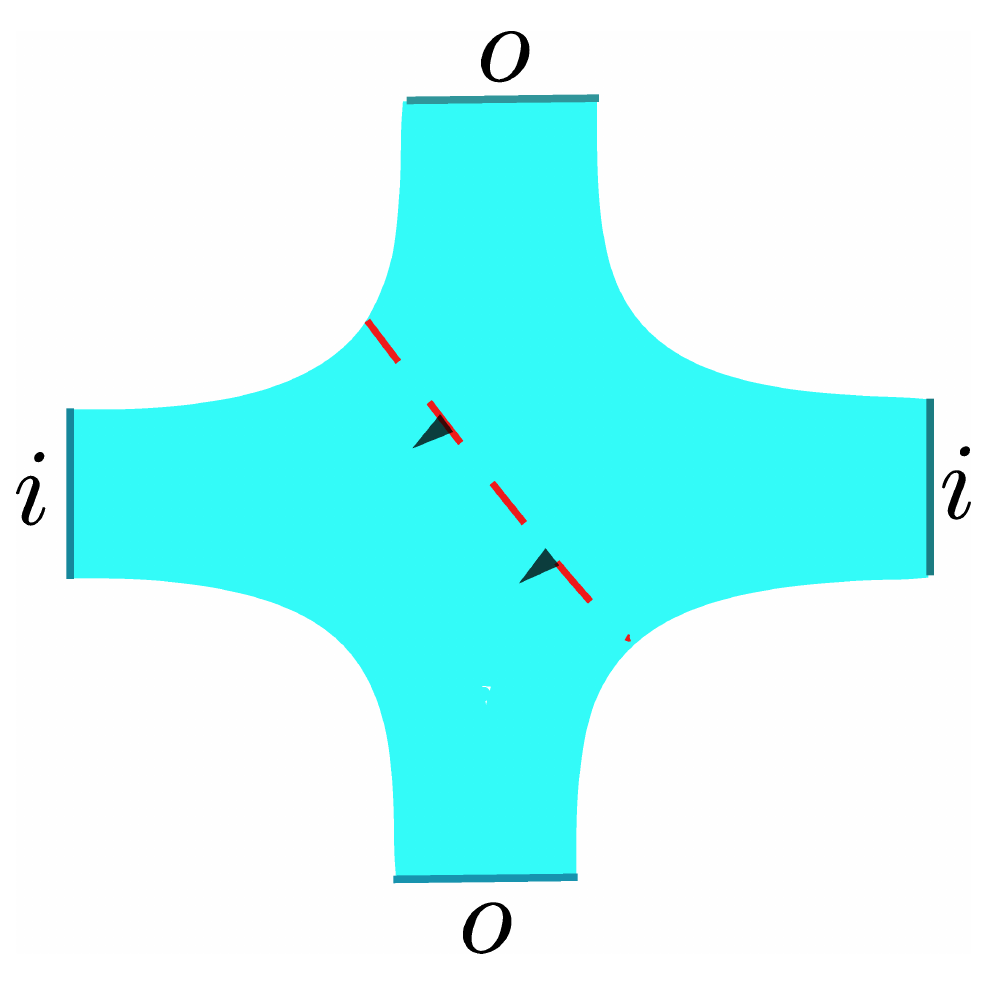}\hspace{5mm}\figbox{0.3}{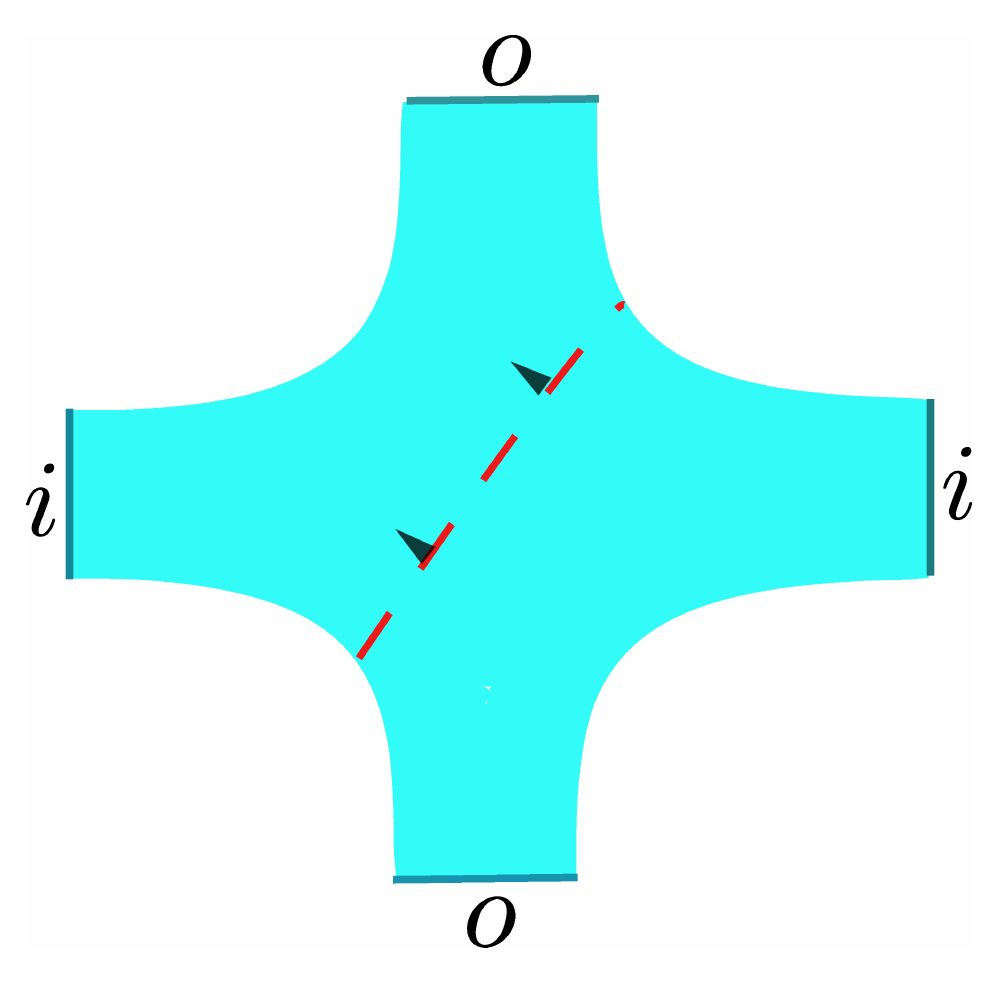}\hspace{5mm} \figbox{0.18}{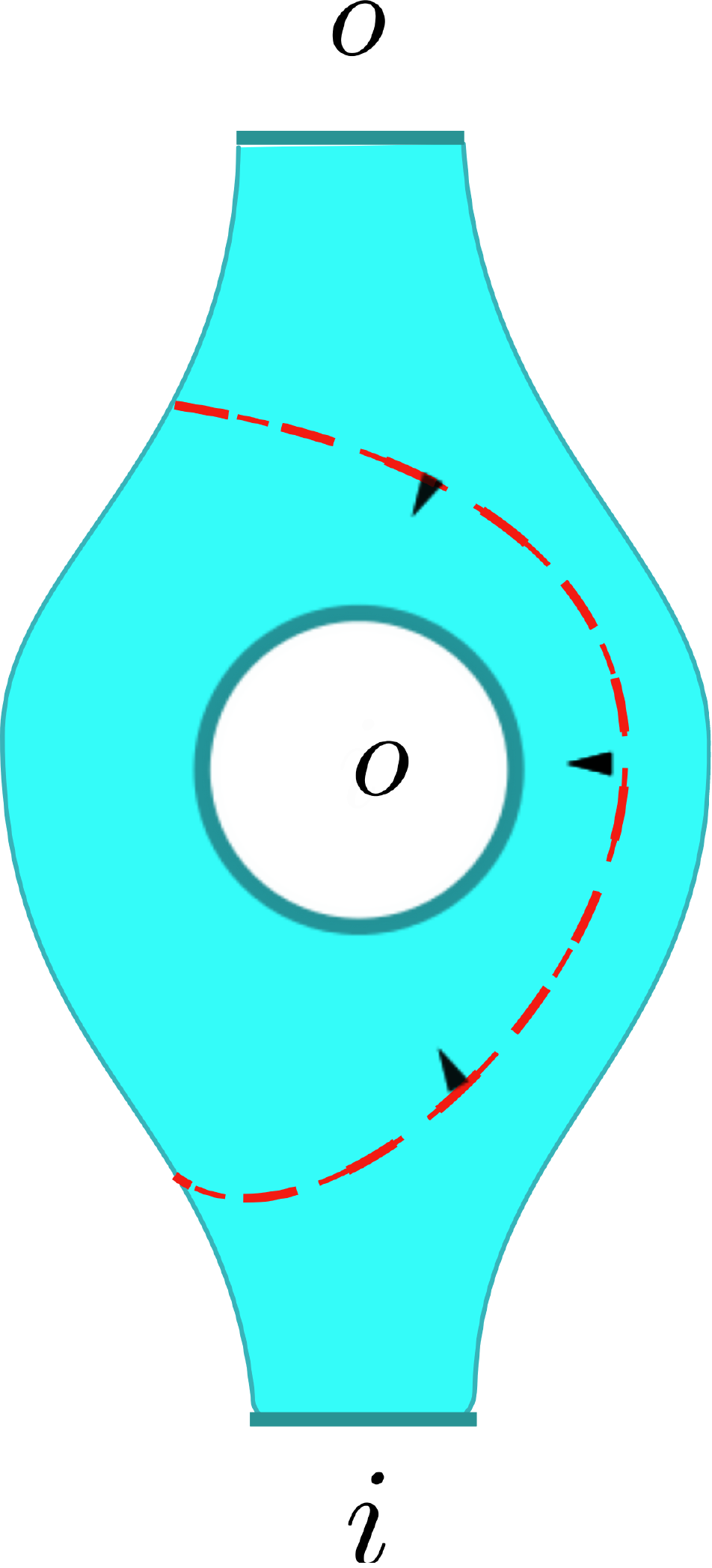} \hspace{5mm}
\figbox{0.18}{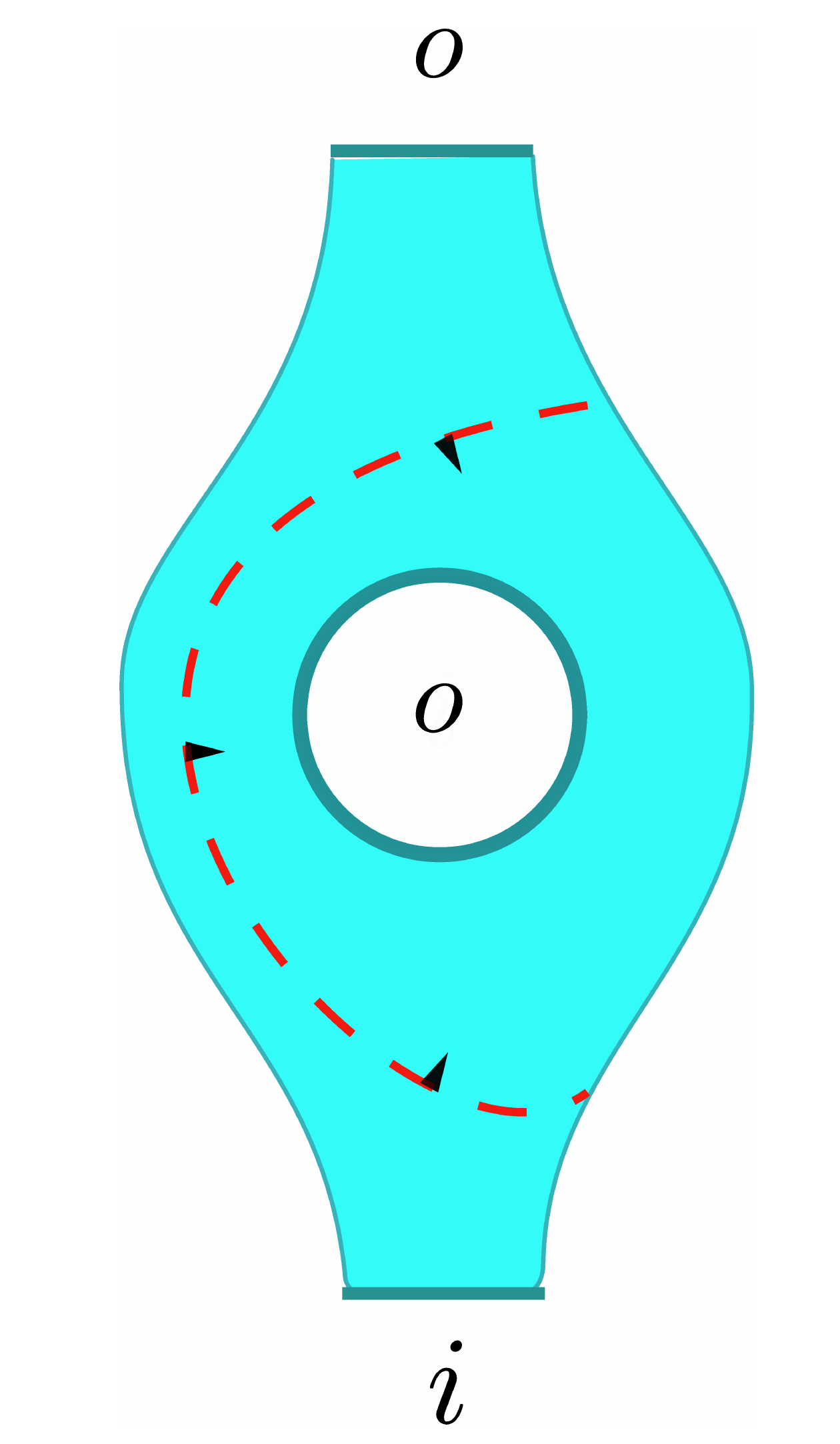}\hspace{5mm}
\figbox{0.18}{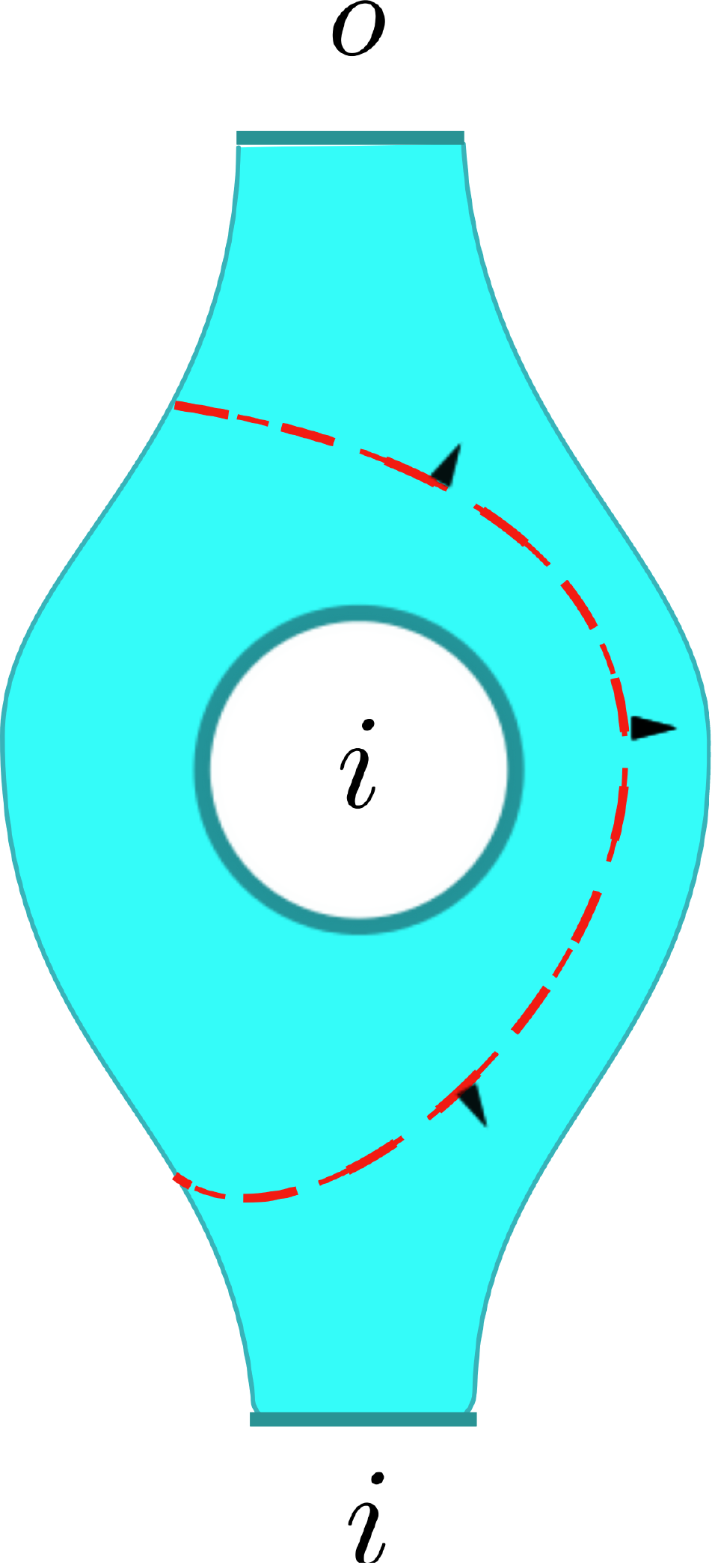} \hspace{5mm} \figbox{0.18}{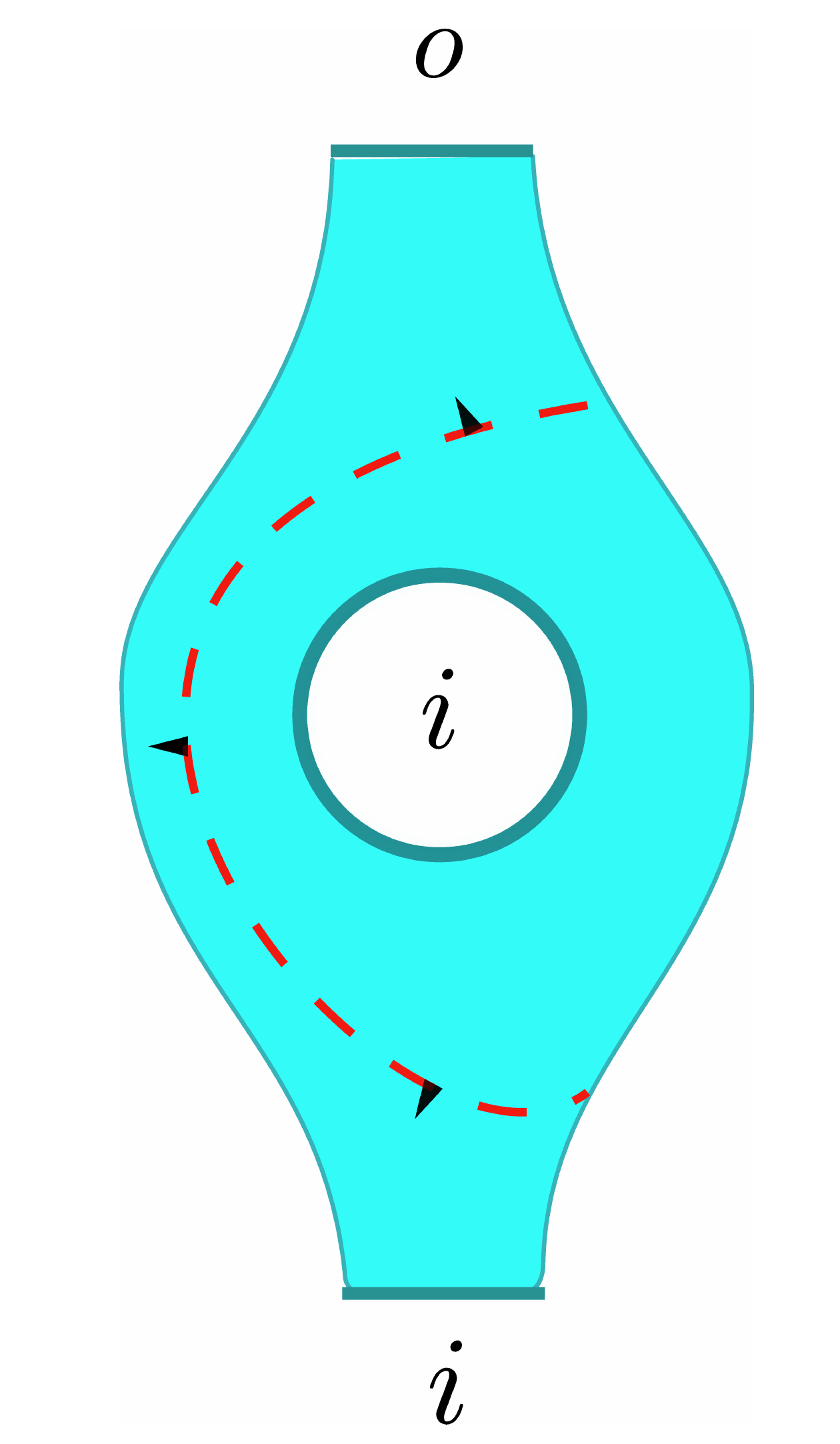}
\quad .
\end{equation}
In Definition \ref{def:wsh-subord} we will choose the first one in each pair of possible decompositions. However, making a different choice would make no difference, as each decompositions of $\Xr_\alpha$
in (\ref{pic:decom-extended-gen}) give rise to the same linear map $U(\Xr_\alpha)\rightarrow F(\Xr_\alpha)$. For
 $\alpha=oco(o,i)$ and $\alpha=oco(i,o)$, this follows from $\Rr28$ and its dual respectively. 
For $\alpha=o4$, this is not that direct. After a little thought, the reader will agree that
the following sequence of decompositions of $\Xr_{o4}$ give rise to the same
linear operator $U(\Xr_{o4})\rightarrow F(\Xr_{o4})$ using relations indicated on the arrows:
\begin{align*}
&\figbox{0.3}{Pic88}\overset{\Rr
10}{\longrightarrow}\figbox{0.3}{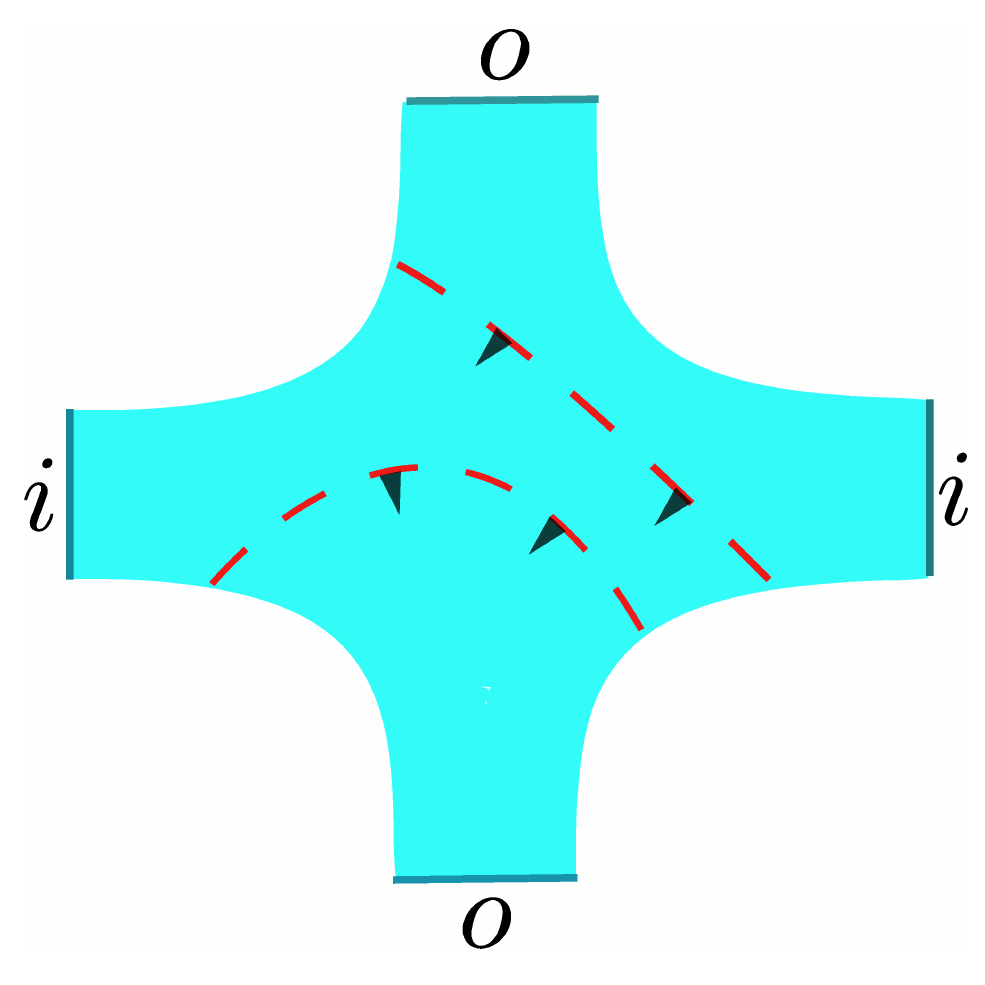}\overset{\Rr4}{\longrightarrow}\figbox{0.3}{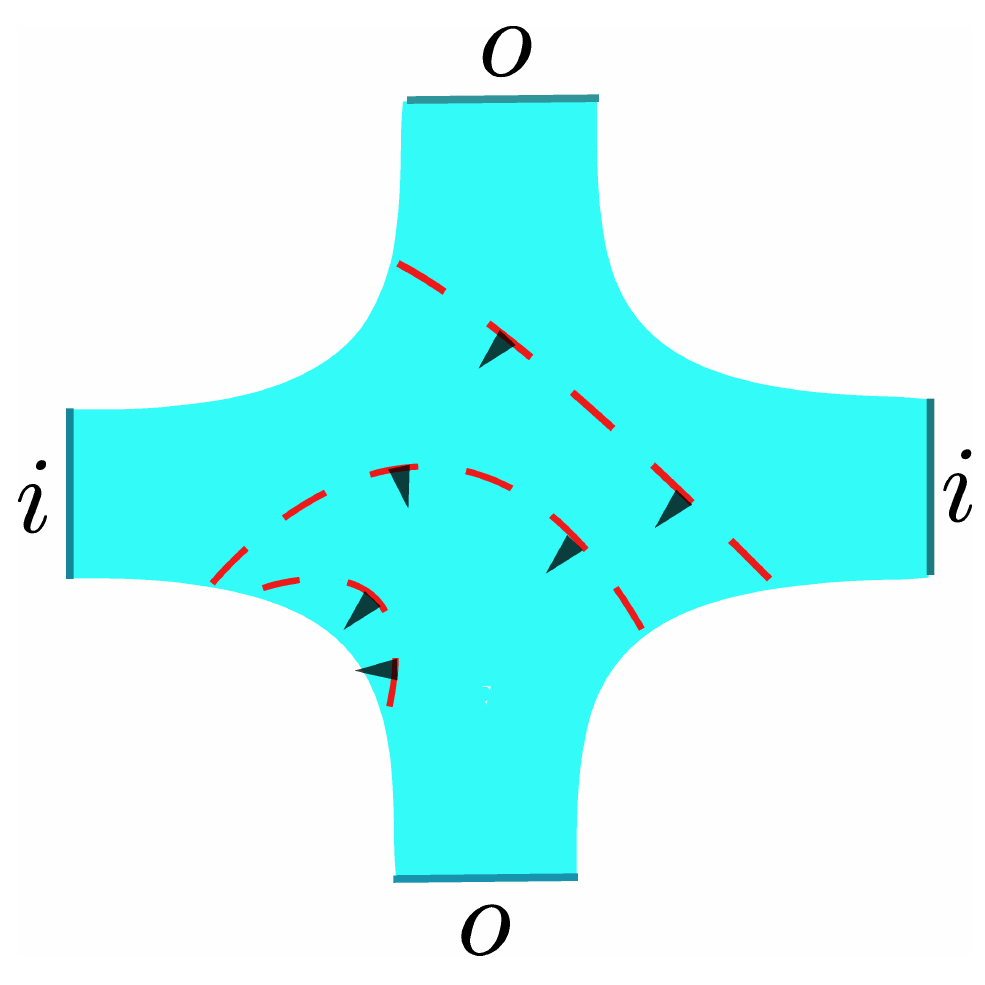}\overset{\Rr
9}{\longrightarrow}
\figbox{0.3}{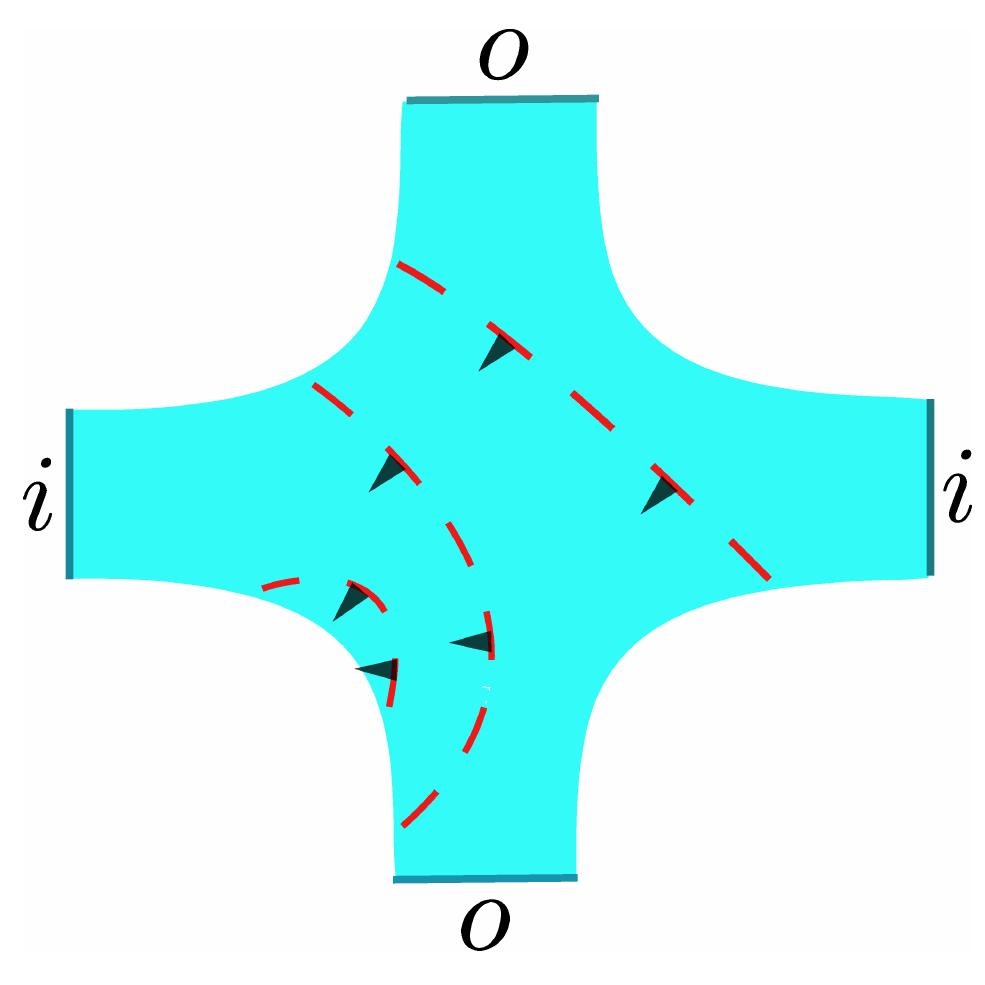}\\
\overset{\Rr 12}{\longrightarrow}&\figbox{0.3}{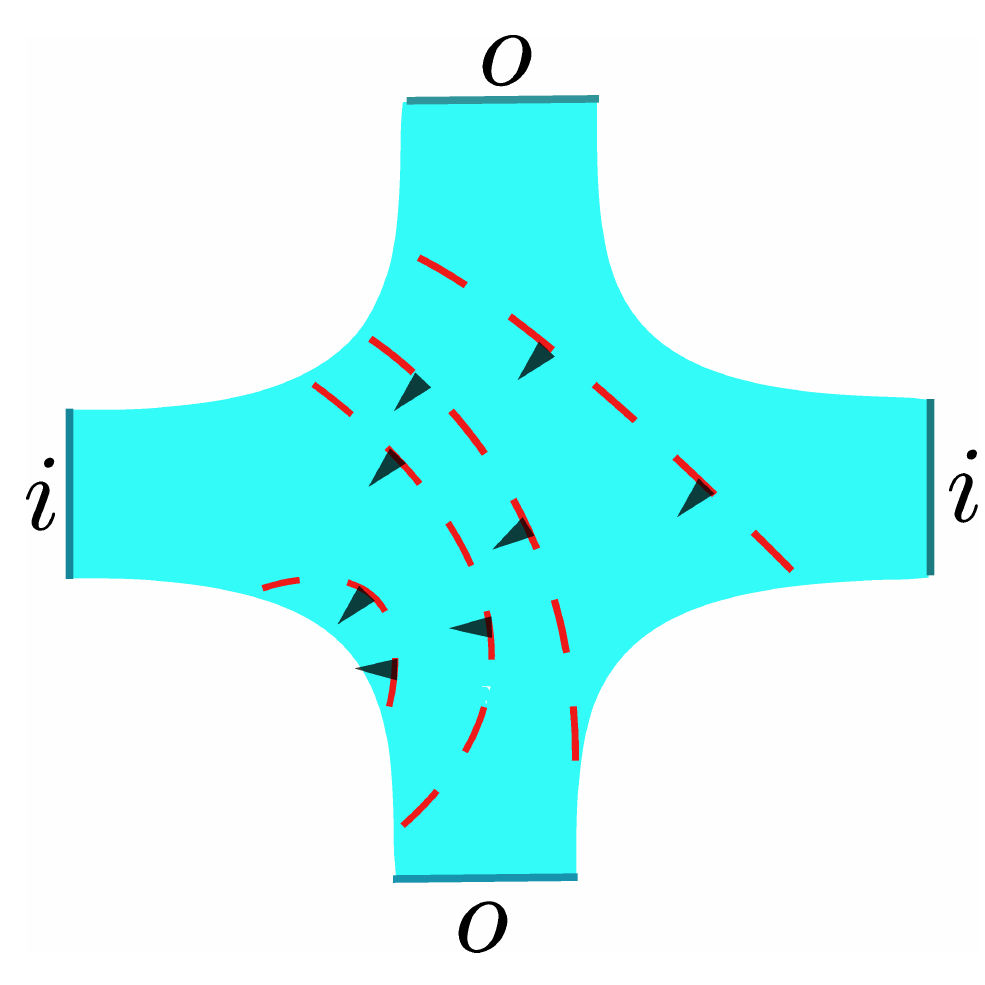}\overset{\Rr2}{\longrightarrow}\figbox{0.3}{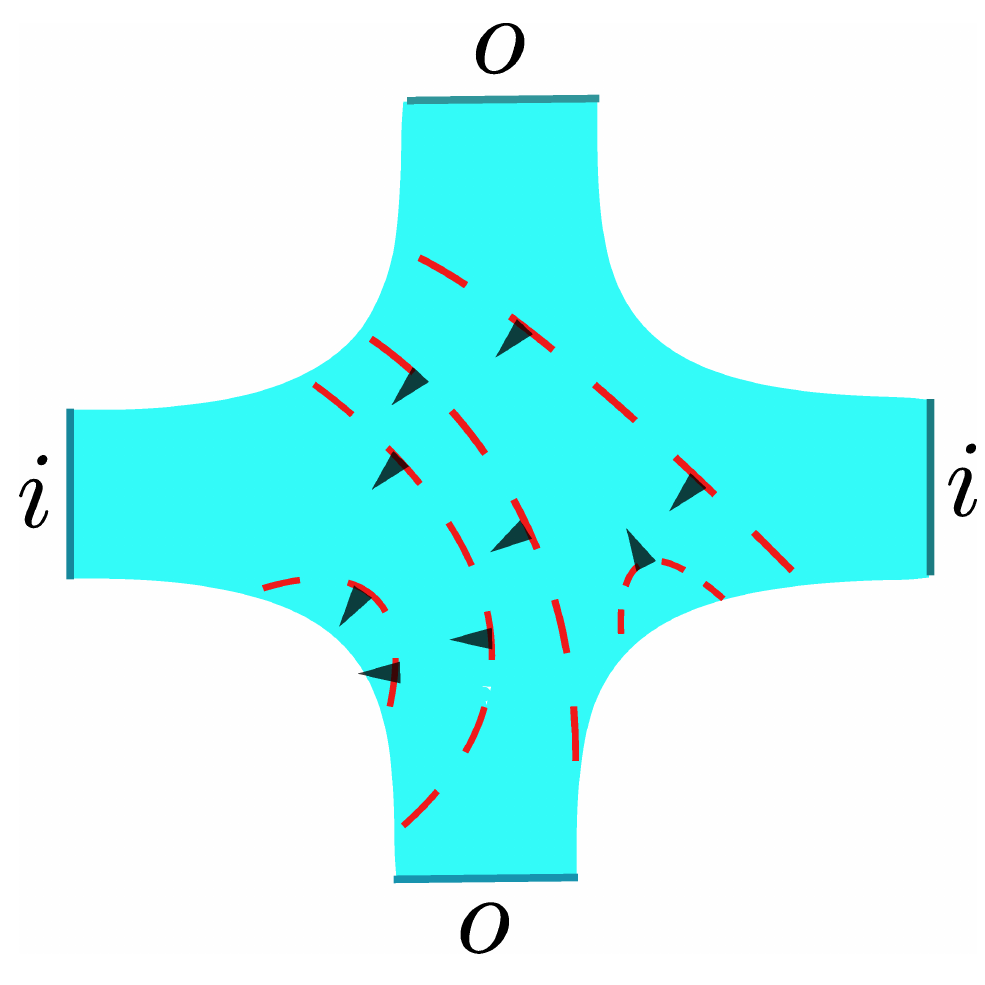}\overset{\Rr
9}{\longrightarrow}\figbox{0.3}{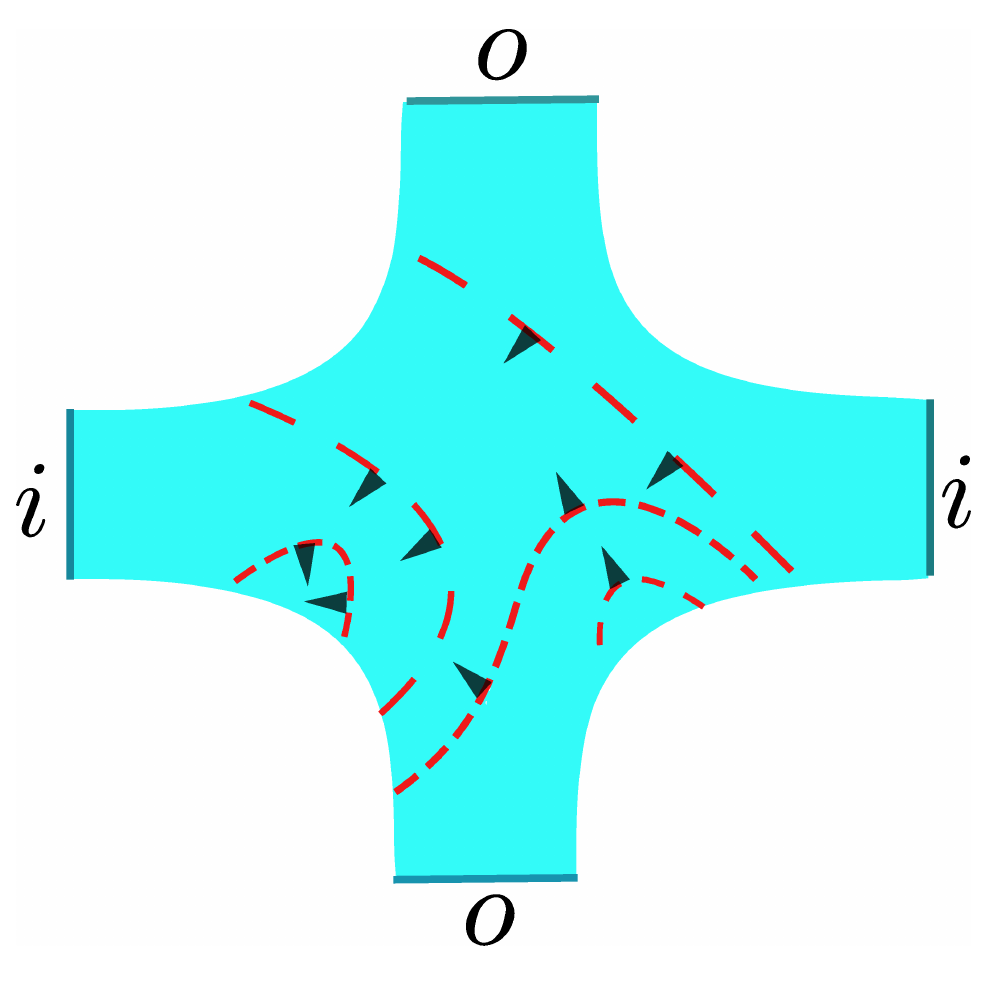}\overset{\Rr 8}{\longrightarrow}\figbox{0.3}{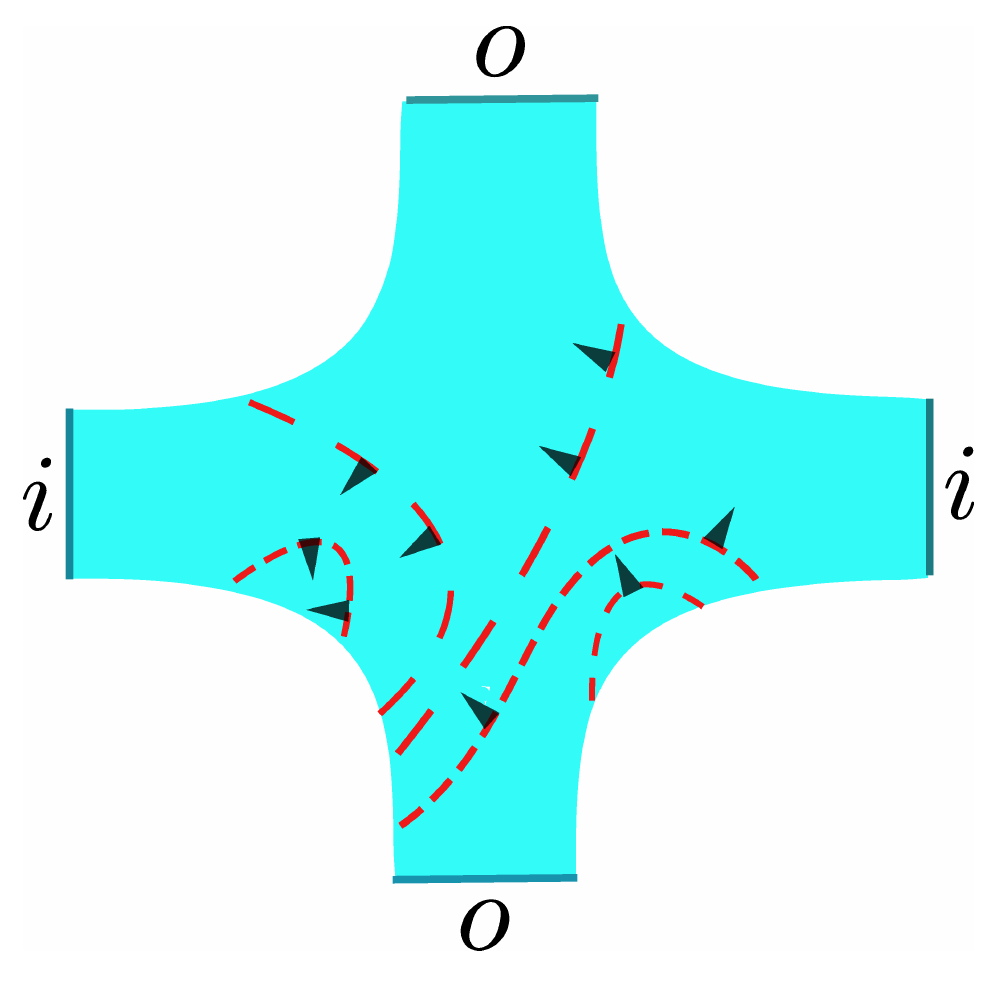}\\
\overset{\Rr
5}{\longrightarrow}&\figbox{0.3}{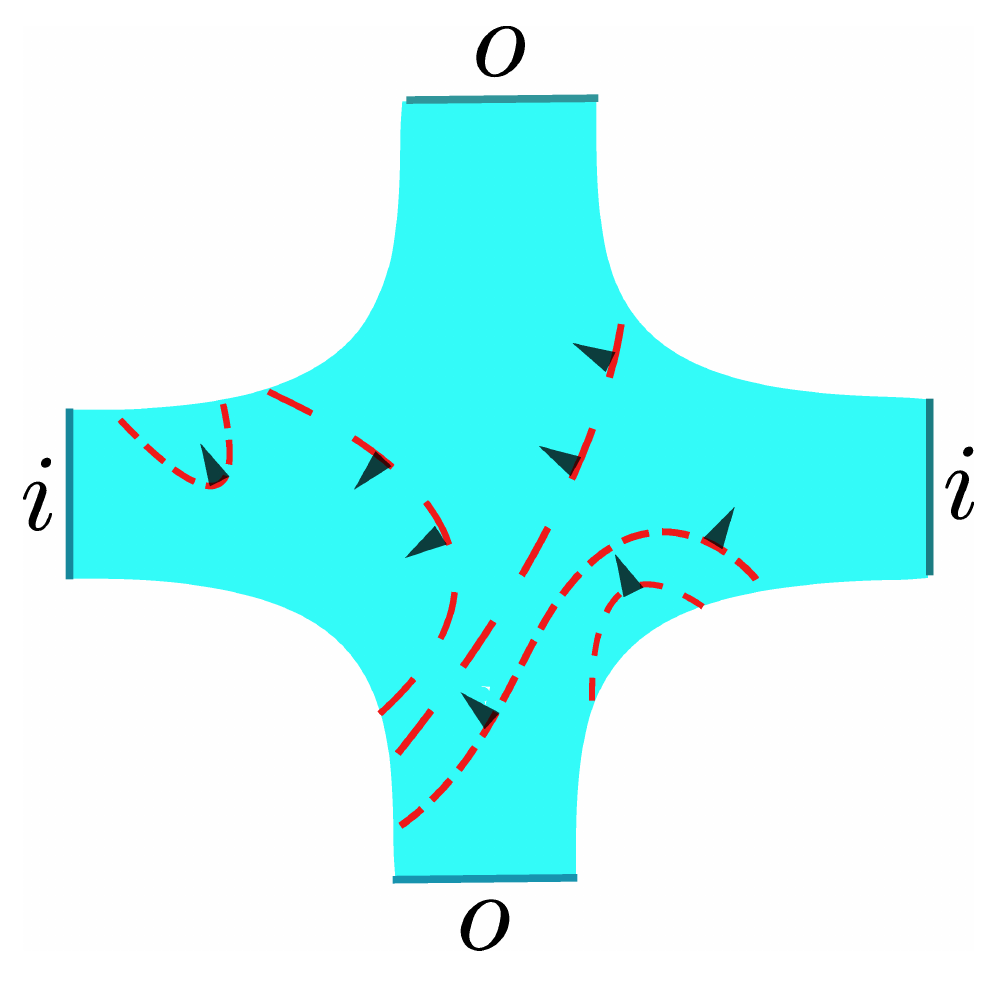}\overset{\Rr8}{\longrightarrow}\figbox{0.3}{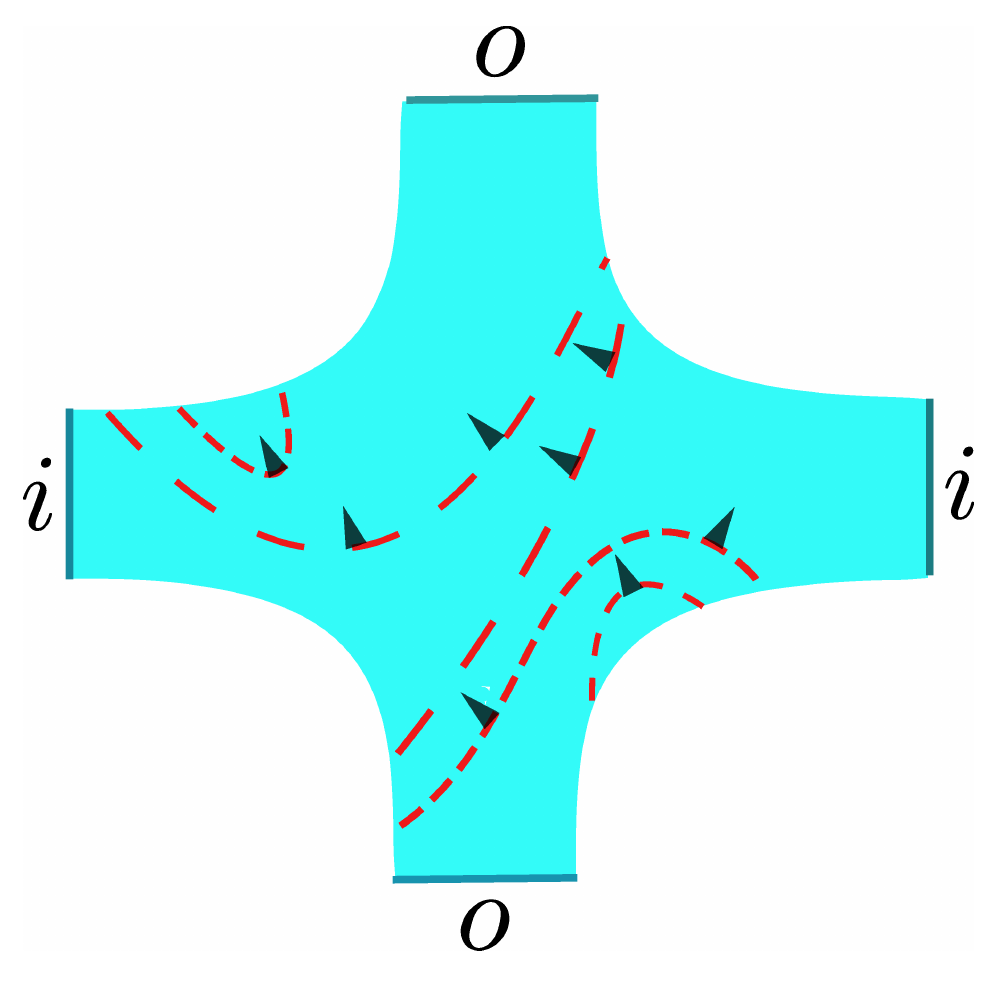}\overset{\Rr3}{\longrightarrow}
\figbox{0.3}{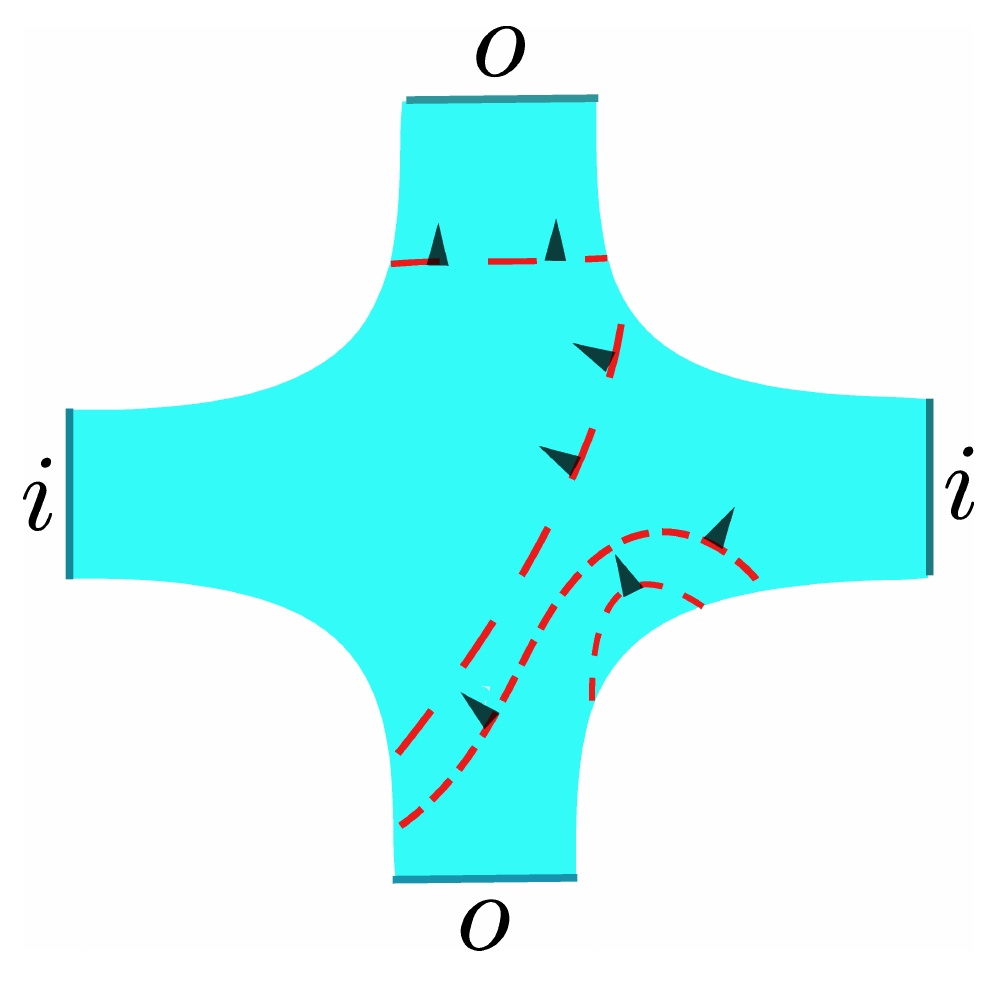}\overset{\Rr
2}{\longrightarrow}
\figbox{0.3}{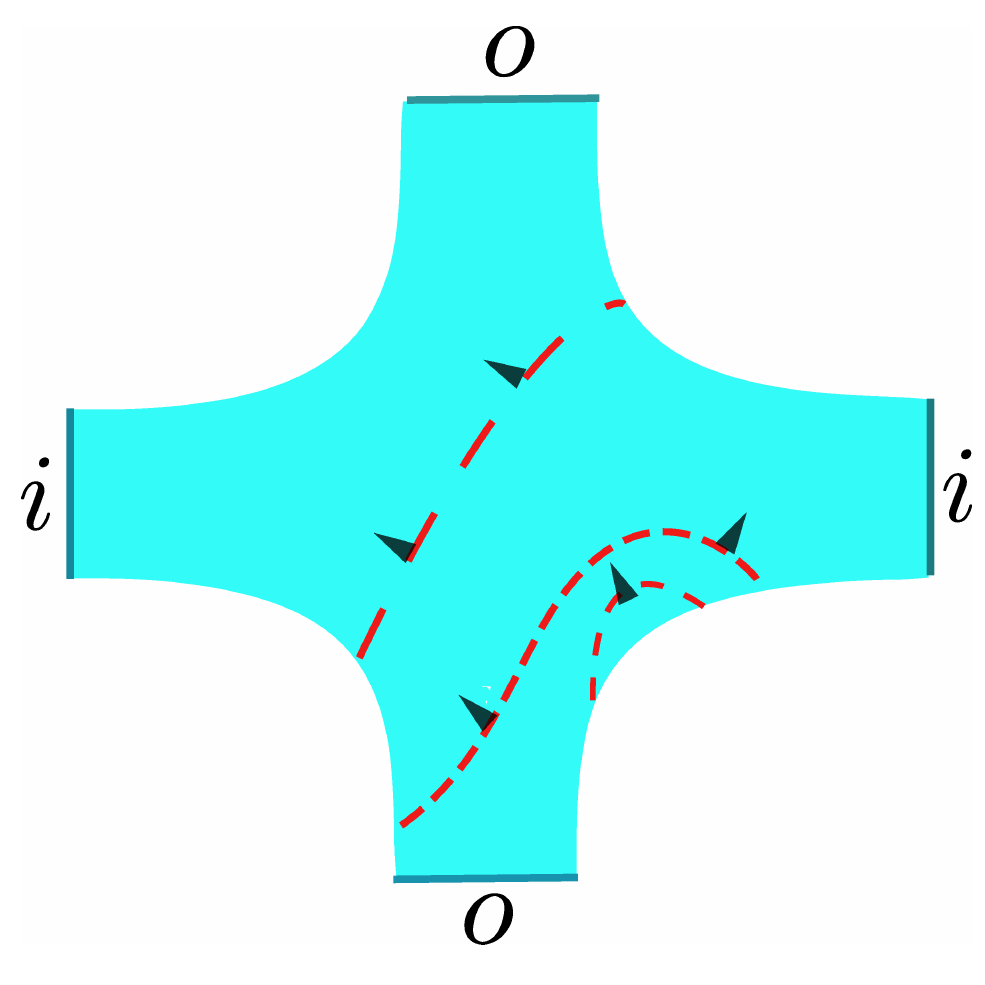}\\
\overset{\Rr
9}{\longrightarrow}&\figbox{0.3}{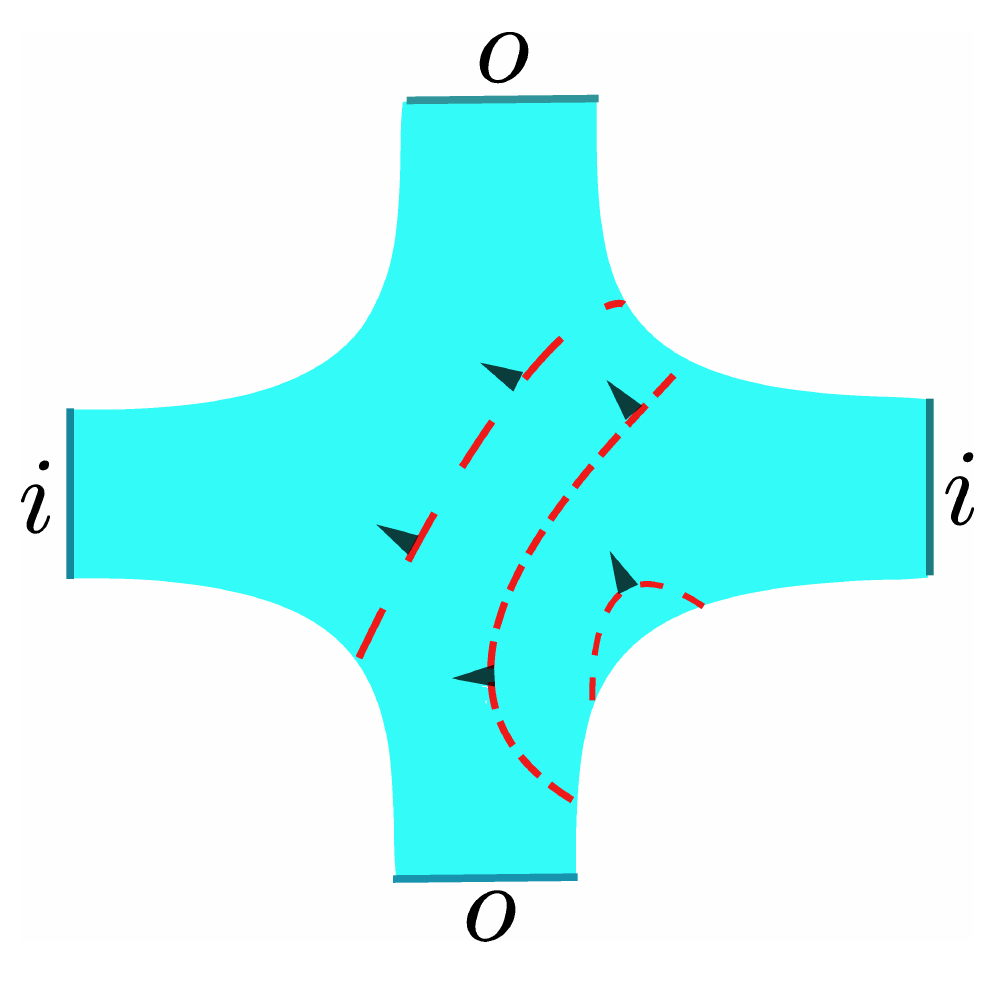}\overset{\Rr2}{\longrightarrow}\figbox{0.3}{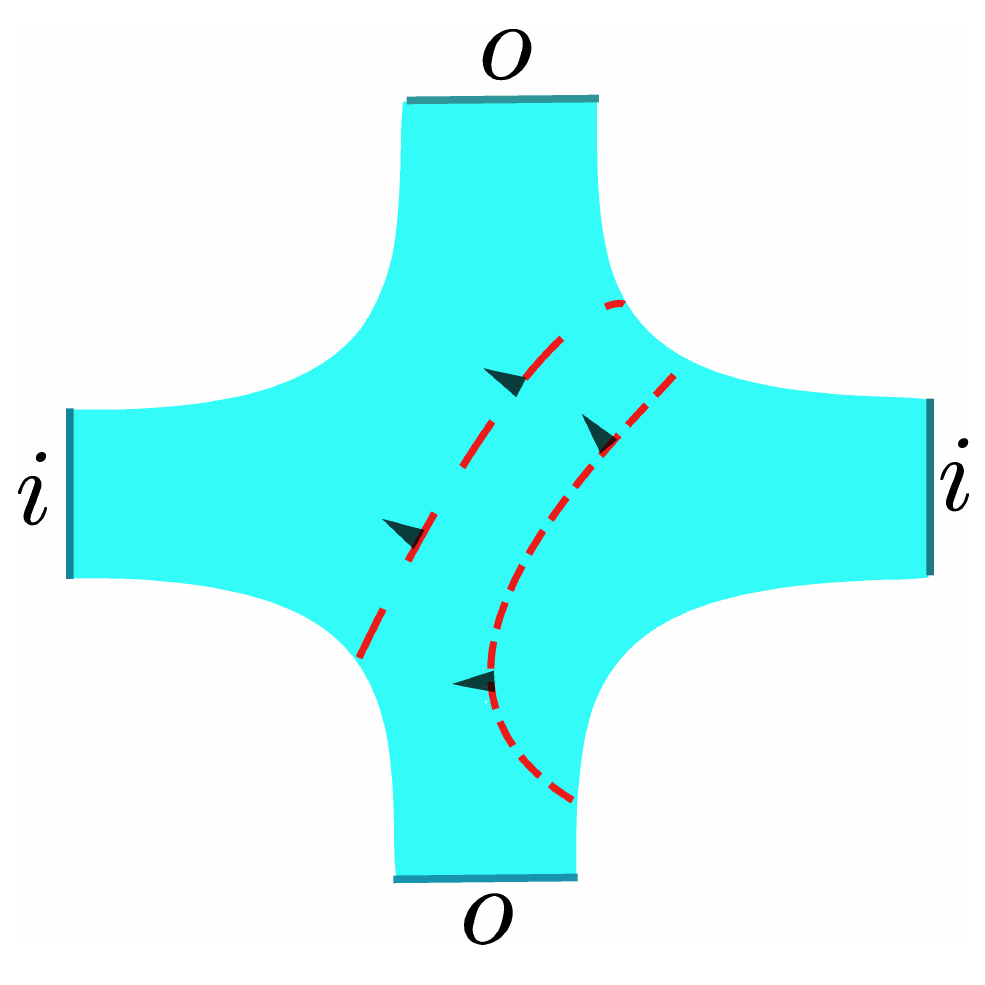}\overset{\Rr12}{\longrightarrow}
\figbox{0.3}{Pic107}
\end{align*}

\begin{defn}  \label{def:wsh-subord} {\rm
Let $\Xr$ be a world sheet and let $f$ be a generic Morse function on $\Xr$ with critical points $\{c_1,\cdots,c_{n+1}\}$. A
decomposed world sheet $\Xr^d = (\Xr,(\alpha_1,\dots,\alpha_m), \varpi)$  is called a {\em decomposition
of} $\Xr$ {\em subordinate to} $f$ if there exists a choice of values $0=t_0< f(c_1) < t_1 < \cdots <t_{n} < f(c_{n+1})
< t_{n+1}=1$ such that
 every component of  $f^{-1}([t_k,t_{k+1}])$ is either the image of one single generator $\Xr_{\alpha_i}$ under $\varpi$ for
some $\alpha_i$, or the image of two generators $\Xr_{\alpha_i}$ and $\Xr_{\alpha_j}$ under $\varpi$ as indicated in the
following pictures, where the red dotted lines denote how the two generators are sewn:
\begin{equation}\label{eq:decom-extended-gen}
 \figbox{0.3}{Pic88}\hspace{7mm} \figbox{0.18}{Pic87} \hspace{7mm}
\figbox{0.18}{Pic89}
\end{equation}

}
\end{defn}

\begin{rema}{\rm
With the above definition, every generic Morse function on a world sheet $\Xr$ has a decomposition subordinate to it.}
\end{rema}

\begin{lemma}    \label{lem:Morse-exists}  {\rm
Suppose that $C$ obeys \erf{eq:sewing-cond}. For each decomposed world sheet $\Xr^d_1
=(\Xr,(\alpha_1,\dots,\alpha_n),\varpi_1)$ there exists a Morse function $f$ on $\Xr$ and a decomposition
$\Xr^d_2 = (\Xr,(\alpha_1,\dots,\alpha_n,\beta_1,\dots,\beta_m),\varpi_2)$ subordinate to $f$ such that $C_{\Xr^d_1} =
C_{\Xr^d_2}$.
}
\end{lemma}
\noindent {\it Sketch of a proof}: 

\noindent 1) Use Lemma \ref{lem:inv-omit-cyl} to add cylinders $\Xr_{po}$ and $\Xr_{pc}$ to $\Xr^d_1$ 
until we obtain a decomposed world sheet
\be
  \Xr^d_{1a} = (\Xr,(\alpha_1,\dots,\alpha_n,po,\dots,pc,\dots),\varpi_{1a})
\ee
such that a) the sewing $\sew$ in $\varpi_{1a} = (\sew,h)$ only identifies state boundaries of $\Xr_{\alpha_i}$
with those of $\Xr_{pc}$ or $\Xr_{po}$ (i.e. never identify directly those of $\Xr_{\alpha_i}$ with those of $\Xr_{\alpha_j}$),
b) all boundaries of
each $\Xr_{\alpha_i}$ are contained in $\sew$ (i.e.\ all state boundaries of $\Xr$ come from state boundaries of
the cylinders $\Xr_{po/pc}$), and c) $C_{\Xr^d_1} = C_{\Xr^d_{1a}}$.

\smallskip
\noindent 2) On each of the generator $\dot\Xr_{\alpha_i}$, if $\alpha_i\not\in\{pc,po\}$, we pick a Morse function as the
height function depicted in \eqref{eq:standard-gen-ws}
with one critical point. On $\dot\Xr_{po}$, we take the Morse function to be the height function in the following
picture:\footnote{We take the Morse function on $\Xr_{po}$ in (\ref{eqn:Morse-po}) instead of the height function in
\eqref{eq:standard-gen-ws} to guarantee that the Morse function has at least one critical point on each generator $\Xr_{\alpha_i}$.}
\begin{equation}\label{eqn:Morse-po}
\figbox{0.64}{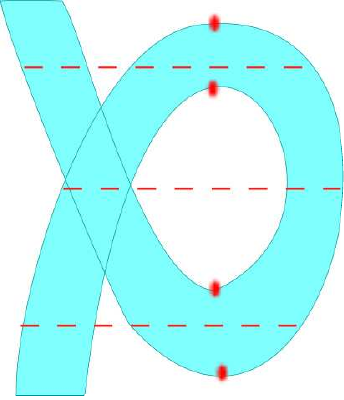}
\end{equation}

We can take a similar Morse function on $\dot\Xr_{pc}$. In all the cases, the Morse function takes values in an
interval $[x_i,x_i+\eps]$ for $0<x_i<1$ and some small $\eps>0$, such that the
Morse function takes the value $x_i$ on the in-coming state boundaries of $\dot\Xr_{\alpha_i}$
and $x_i+\eps$ on the out-going state boundaries. Choose the $x_i$ and $\eps$ such that all the intervals $[x_i,x_i+\eps]$ for
$i=1,\dots,n$ are disjoint.

\smallskip
\noindent 3) Let's extend the Morse function from step 2) to the cylinders $\dot\Xr_{po}$ and $\dot\Xr_{pc}$ added in step 1. Since
the state
boundaries of such $\Xr_{po}$ and $\Xr_{pc}$ are either identified with the state boundaries of the generators $\Xr_{\alpha_i}$ for
$1\leqslant i\leqslant n$ or part of the state boundaries of $\Xr$, the value of the Morse
function on the state boundaries of $\dot\Xr_{po}$ and $\dot\Xr_{pc}$ are uniquely determined. On
each cylinder $\dot\Xr_{po}$ or $\dot\Xr_{pc}$, if the value on the in-going boundary is less than the value on
the out-going boundary, we take the Morse function as in (\ref{eqn:Morse-po}). Otherwise we take the Morse function to be of the
form
\be \label{eq:morse-po-critical}
 \includegraphics[width=20mm]{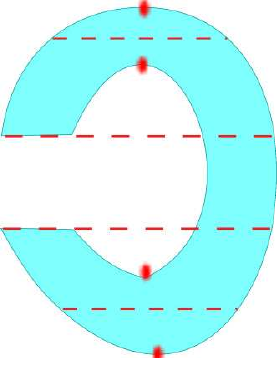}
\ee
with four critical points. Altogether we can obtain a Morse function $f$ on $\dot\Xr$ with critical points
$c_1,\dots,c_l$ for which there exists a choice of $\{ t_i \}$ with
 $0=t_0 <f(c_1) < t_1 < \cdots <t_{l-1} < f(c_l) < t_{l}=1$ such that each $\Xr_{\alpha_i}$ is isomorphic to one of the
connected components of $f^{-1}([t_k,t_{k+1}])$ for some $k$. 
\\[.3em]
4) If on a cylinder $\Xr_{po}$ the Morse function $f$ is of the form as shown in (\ref{eq:morse-po-critical}), we decompose
$\Xr_{po}$ further via a morphism $\Xr_{mo}\otimes\Xr_{\Delta o}\otimes\Xr_{\eta o}\otimes\Xr_{\eps o}\otimes\Xr_{po} \rightarrow
\Xr_{po}$ as indicated by the dotted lines in \erf{eq:morse-po-critical}, and similar for $\Xr_{pc}$. This does not affect the
value of $C$ according to Corollary \ref{cor:inv-Rn} in the case of R8 , R1 and R3 and Lemma \ref{lem:inv-omit-cyl}. We have a
similar result for the decomposition in (\ref{eqn:Morse-po}).
\\[.3em]
5) After subdividing the cylinders $\Xr_{po}$ and $\Xr_{pc}$ into more cylinders when necessary, we obtain a
decomposed world sheet
\be
  \Xr^d_{2} = (\Xr,(\alpha_1,\dots,\alpha_n,\beta_1,\dots,\beta_m),\varpi_{2})
\ee
which is subordinate to $f$. The $\beta_i \in \{  mo, \Delta o, \eta o, \eps o, mc, \Delta c, \eta c, \eps c, pc,
po \}$ resulted from steps 1,2 and 3, and from the additional subdivisions of $\Xr_{po}$ and $\Xr_{pc}$.
By construction, we have $C_{\Xr^d_1} = C_{\Xr^d_2}$.
\epf

\medskip
By Lemma \ref{lem:Morse-exists}, to prove Theorem \ref{thm:gen-rel} we only need to show that if two decompositions $\Xr^d_1$ and
$\Xr^d_2$ of a world sheet $\Xr$ are subordinate to Morse functions $f_0$ and $f_1$ respectively, then we have $C_{\Xr^d_1} =
C_{\Xr^d_2}$.
This is achieved in Proposition\,\ref{prop:Morse-equal}. We will prove this result in steps. First, we prove the cases that the
slicing sets of $\Xr^d_1$ and $\Xr^d_2$ are isotopic (Lemma \ref{lem:Morse-isotopic-level-set}); secondly, we prove the cases that
two decompositions are subordinate to the same Morse function (Lemma \ref{lem:Morse-same}); thirdly, we prove  the cases that two
decompositions are subordinate to two Morse functions that lie in the same path connecting component of generic Morse functions
(Lemma \ref{lem:path-connected generic Morse functions}); after recalling the degeneration that can occur
in Lemma \ref{lem:finite-exceptions}, we finally
prove Proposition \ref{prop:Morse-equal}.

\begin{lemma}\label{lem:Morse-isotopic-level-set} {\rm
Let $f_0$ and $f_1$ be  generic Morse functions on a world sheet $\Xr$. Let
$\Xr^d_0
=(\Xr,(\alpha_1,\cdots,\alpha_m),\varpi_0)$ be a decomposition of $\Xr$ subordinate to $f_0$ with (by definition) a
choice of values $0=s_0<s_1<\cdots<s_{k+1}=1$, where $\alpha_i\in\mathcal{S}$. Let
$\Xr^d_1=(\Xr,(\beta_1,\cdots,\beta_n),\varpi_1)$ be a decomposition of $\Xr$ subordinate to $f_1$ with a choice of
values $0=t_0<t_1<\cdots<t_{l+1}=1$, where $\beta_i\in\mathcal{S}$. Assume that there is an isotopy
$$F:(\sqcup_{0\leqslant i\leqslant k+1}f_0^{-1}(s_i))\times [0,1]\rightarrow \dot\Xr$$ such that
\begin{enumerate}
\item $F|_{\sqcup_{0\leqslant i\leqslant k+1}f_0^{-1}(s_i)\times\{0\}}$ is the embedding of $\sqcup_{0\leqslant
i\leqslant k+1}f_0^{-1}(s_i)$ in $\dot\Xr$,
\item $k=l$ and $F(f_0^{-1}(s_i),1)=f_1^{-1}(t_i)$ for each $i\in \{1, \cdots, k\}$,
\item If $x\in\sqcup_{0\leqslant i\leqslant k+1}f_0^{-1}(s_i)$ lies on the boundary of $\dot\Xr$, then $F(x,s)$ lies on
the boundary of $\dot\Xr$ for all $s\in[0,1]$. If $x\in\sqcup_{0\leqslant i\leqslant k+1}f_0^{-1}(s_i)$ lies in the
interior of $\dot\Xr$, then $F(x,s)$ lies in the interior of $\dot\Xr$ for all $s\in[0,1]$. 
\end{enumerate}
 Then we have
$C_{\Xr^d_0}=C_{\Xr^d_1}$.}
\end{lemma}

\pf
Let $M=\sqcup_{0\leqslant i\leqslant k+1} f_0^{-1}(s_i)$ and $N=\sqcup_{0\leqslant i\leqslant l+1} f_1^{-1}(t_i)$. The
third condition of $F$ guarantees that we can extend\footnote{Such an extension is possible as can be seen e.g.\ by a slight modification of \cite[Thm.\,8.1.3]{Hirsch}. That theorem does not directly apply since $F$ does not satisfy the condition that the image of $F$ is either contained in
$\dot\Xr\setminus\partial\dot\Xr$ or in $\partial\dot\Xr$. But condition 3 of $F$ in Lemma \ref{lem:Morse-isotopic-level-set} makes
the argument in the proof of \cite[Thm.\,8.1.3]{Hirsch} work also for our case.} the isotopy $F$ to obtain a family of
diffeomorphisms $H_s$ of $\dot\Xr$ for $s\in [0,1]$ such that $H_s|_M=F|_{M\times\{s\}}$, $H_0=\id_{\dot\Xr}$. In particular, we
have $H_1(f_0^{-1}(s_i))=f_1^{-1}(t_i)$ for $0\leq i \leq k+1$. 

Let $\{\Yr_{1},\cdots, \Yr_m\}$ be the set of components of $f_0^{-1}([s_0,s_1])$, and let $\{\Zr_1,\cdots,\Zr_n\}$ be the set of  components of $f_1^{-1}([t_0,t_1])$.  There always exist
$s$ and $t$ such that $s\leqslant t_0<t_1\leqslant t$ and $H_1(\Yr_{1})\subset
f_1^{-1}([s,t])$. Here we choose  $[s,t]$ to be the smallest interval such that $H_1(\Yr_{1})\subset
f_1^{-1}([s,t])$. We claim that $s=t_0$, $t=t_1$. If $t>t_{1}$, then there exists an interior point $x$ of
$H_1(\Yr_{1})$ such that $f_1(x)=t_{1}$ and $H_1^{-1}(x)\in f_0^{-1}(s_1)\subset\partial Y_1$. This  contradicts to the fact that $H_1(\partial
\Yr_{1})=\partial H_1(\Yr_{1})$ and the injectivity of $H_1$.  Similarly we can show that $s=t_{0}$. Thus there exists $j\in\{1,\cdots,n\}$ such that $H_1(\Yr_{1})\subset\Zr_{j}$. By the same argument, we can show that
$H_1^{-1}(\Zr_{j})\subset\Yr_{i}$ for some $i\in\{1,\cdots,m\}$. It is clear that $i=1$ and we must have
$H_1(\Yr_{1})=\Zr_{j}$. This argument works for all components of $f_0^{-1}([s_{0},s_{1}])$, and also works for all
$f_0^{-1}([s_{i_0},s_{i_0+1}]), 0\leqslant i_0\leqslant k$. Note that either both $\Yr_{1}$ and $\Zr_{j}$ are the images of one
single generator, or both of them are the images of two generators. Thus it is clear that $m=n$. After a permutation of indices, we can assume that
$(\alpha_1,\cdots,\alpha_m)=(\beta_1,\cdots,\beta_n)$. If the
components of $f_0^{-1}([s_{i_0},s_{i_0+1}]), 0\leqslant i_0\leqslant k$ are all images of single generators, then it is clear
that the diffeomorphism $H_1$ of $\dot\Xr$ induces morphisms of world
sheets $h_i:\Xr_{\alpha_i}\rightarrow\Xr_{\beta_i}$ such that the following diagram commutes. 
\begin{equation}\label{diagram}
\xymatrixcolsep{5pc}\xymatrix{
\Xr_{\alpha_1}\otimes\cdots\otimes\Xr_{\alpha_m} \ar[d]_{\varpi_0} \ar[r]^{h_1\otimes\cdots\otimes h_m}
&\Xr_{\beta_1}\otimes\cdots\otimes\Xr_{\beta_m} \ar[d]_{\varpi_1}\\
\Xr \ar[r]^{H_1} &\Xr }
\end{equation}
If some component of $f_0^{-1}([s_{i_0},s_{i_0+1}])$ is a world sheet in (\ref{eq:decom-extended-gen}), then so is its image
under
$H_1$. We can deform $H_1$ in the space of diffeomorphisms of $\Xr$ such that it preserves the red dotted line in
(\ref{eq:decom-extended-gen}). In particular, this deformation of $H_1$ does not change the property
that it is homotopic to the identity map. For convenience, we still denote the deformed diffeomorphism by $H_1$, and we still have
the
commutative diagram (\ref{diagram}).
Let ${\Xr^d_0}'=(\Xr,(\alpha_1,\cdots,\alpha_m),H_1\circ\varpi_0)$
and ${\Xr^d_1}'=(\Xr,(\alpha_1,\cdots,\alpha_m),\varpi_1\circ(h_1\otimes\cdots\otimes h_k))$. Then we obtain
\begin{equation*}
     C_{\Xr^d_0}\overset{(1)}{=}C_{{\Xr^d_0}'}\overset{(2)}{=}C_{{\Xr^d_1}'}\overset{(3)}{=}  C_{\Xr^d_1}~,
\end{equation*}
where step 1 uses the fact that  $H_1$ is homotopic to the identity morphism; step 2 is the commutativity of (\ref{diagram});
step 3 uses Lemma \ref{lem:inv-mapgen} and Lemma \ref{lem:change-part}.
\epf

\smallskip
\begin{lemma}\label{lem:Morse-same}{\rm
Let $f$ be a generic Morse function on $\Xr$ with critical points $\{c_1,\cdots,c_{k+1}\}$ such that $0<f(c_1)<\cdots<f(c_{k+1})<1$, and let $\Xr^d_0=(\Xr,(\alpha_1,\cdots,\alpha_m),\varpi_0)$ be a decomposition of $\Xr$ subordinate to $f$ with the choice of values $0=s_0<f(c_1)<s_1<\cdots<s_k<f(c_{k+1})<s_{k+1}=1$. Let $\Xr^d_1=(\Xr,(\beta_1,\cdots,\beta_n),\varpi_1)$ be another decomposition of $\Xr$ subordinate to $f$ with 
the choice of values $0=t_0<f(c_1)<t_1<\cdots<t_k<f(c_{k+1})<t_{k+1}=1$. Then we have $C_{\Xr^d_0}=C_{\Xr^d_1}$.  }
\end{lemma}
\pf
Let $C_i:=\left( f^{-1}([t_i,t_{i+1}])\setminus f^{-1}([s_i,s_{i+1}]) \right) \sqcup \left(
f^{-1}([s_i,s_{i+1}])\setminus f^{-1}([t_i,t_{i+1}]) \right)$. The subset $f^{-1}([s_i,s_{i+1}])$ of $\dot\Xr$ contains a single
critical point  $c_{i+1}$ of $f$, and so does  $f^{-1}([t_i,t_{i+1}])$. Hence  $C_i$ does not contain any critical point of $f$.
This implies that $C_i$ is homeomorphic to the union of several cylinders and strips, as illustrated in the following picture:
\begin{equation*}
\figbox{0.3}{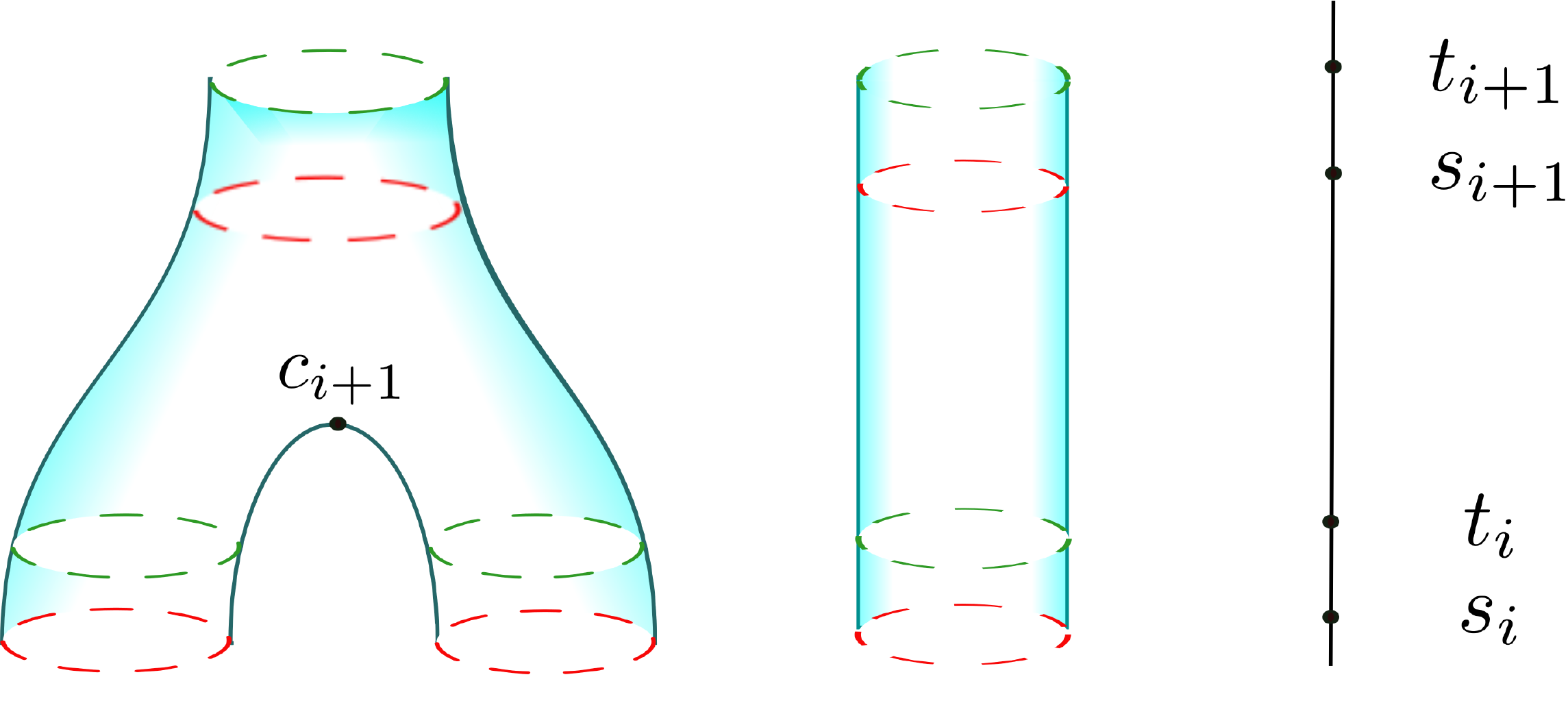}
\end{equation*}  
It is clear  that there is an isotopy from  $\sqcup_{1\leqslant i\leqslant k}f^{-1}(s_i)$ to $\sqcup_{1\leqslant
i\leqslant k}f^{-1}(t_i)$ satisfying the conditions in Lemma \ref{lem:Morse-isotopic-level-set}. Thus we have 
$C_{\Xr^d_0}=C_{\Xr^d_1}$.
\epf

\begin{lemma}  \label{lem:path-connected generic Morse functions}  {\rm
Let $\Xr^d_0=(\Xr,(\alpha_1,\cdots,\alpha_n),\varpi_0)$ and $\Xr^d_1=(\Xr,(\beta_1,\cdots,\beta_m),\varpi_1)$ be two
decompositions subordinate to $f_0$ and $f_1$ respectively. Suppose that $f_0$ and $f_1$ can be connected by a path of
generic Morse functions $f_s:\dot\Xr\rightarrow [0,1], s\in[0,1]$. Then we have $C_{\Xr^d_0}=C_{\Xr^d_1}$.
}
\end{lemma}

\pf Let $\{c_1,\cdots,c_{k+1}\}$ be the set of critical points of $f_0$ with
$0<f_0(c_1)<\cdots<f_0(c_{k+1})<1$. By
Lemma \ref{lem:Morse-family-isotopy}, we can find isotopies
$$ \varphi_s:[0,1]\rightarrow [0,1] \quad \mbox{and} \quad h_s:\dot\Xr\rightarrow \dot\Xr$$ with $\varphi_0$ and $h_0$
identity homeomorphisms and $f_0=\varphi_s\circ f_s\circ h_s$, where $\varphi_s$ is strictly increasing
since $\varphi_s$ is an orientation-preserving diffeomorphism of $[0,1]$. We have the following simple observations from
the identity $f_0=\varphi_s\circ f_s\circ h_s$:
\begin{enumerate}
\item  The diffeomorphism $h_s$ maps the critical points of $f_0$ to those of $f_s$ since the tangent maps of $h_s$ and
$\varphi_s$ are non-degenerate. It follows that the critical points $h_s(c_{i})$ form a path in $\dot\Xr$ as $s$ varies
in $[0,1]$, for each $1\leqslant i\leqslant k+1$. It is also clear that $h_s(c_i)\not=h_s(c_j)$ for $i\not=j$.

\item If $x\in f_0^{-1}(t)$, then $h_s(x)\in f_s^{-1}(\varphi_s^{-1}(t))$.
\end{enumerate}
It follows easily from the first observation that the number of critical points of $f_s$ remains constant as $s$ varies
in $[0,1]$. 

Let $\Xr^d_0$ be a decomposition of $\Xr$ subordinate to $f_0$ with the choice of values
$$0=t_0<f_0(c_1)<t_1<\cdots<t_{k}<f_0(c_{k+1})<t_{k+1}=1$$ From the above observations, it is obvious that $h_1$ maps
$f^{-1}_0([t_{i},t_{i+1}])$ homeomorphically to $f^{-1}_1([\varphi_1^{-1}(t_{i}),\varphi_1^{-1}(t_{i+1})])$ and
that
$\{h_s\}$ gives rise to an isotopy from $\sqcup_{0\leqslant i\leqslant{k+1}}f^{-1}_0(t_{i})$ to $\sqcup_{0\leqslant
i\leqslant{k+1}}f^{-1}_1(\varphi_1^{-1}(t_{i}))$. By the inequality $t_{i}<f_0(c_{i+1})<t_{i+1}$ and the fact that
$\varphi_1^{-1}$ is increasing  we obtain: $$\varphi_1^{-1}(t_{i})<f_1(h_1(c_{i+1}))<\varphi_1^{-1}(t_{i+1})$$
It follows that $f_1^{-1}([\varphi_1^{-1}(t_{i}),\varphi_1^{-1}(t_{i+1})])$ contains a single critical point
$h_1(c_{i+1})$ of $f_1$. By Lemma \ref{lem:Morse-same}, we can assume without loss of generality that the
decomposition $\Xr^d_1=(\Xr,(\beta_1,\cdots,\beta_{m}),\varpi_1)$ of $\Xr$ is subordinate to $f_1$ with the choice
of values $$0< f_1(h_1(c_1)) <\varphi^{-1}(t_1)<\cdots<\varphi^{-1}(t_{k})<f_1(h_1(c_{k+1}))<1 \ .$$ Since $\{h_s\}$ gives rise
to an isotopy from $\sqcup_{0\leqslant i\leqslant{k+1}}f^{-1}_0(t_i)$ to
$\sqcup_{0\leqslant i\leqslant{k+1}}f^{-1}_1(\varphi_1^{-1}(t_i))$ satisfying the conditions in
Lemma \ref{lem:Morse-isotopic-level-set}, we have $C_{\Xr^d_0}=C_{\Xr^d_1}$.
\epf

\medskip
The space of generic Morse functions is not path connected. We need the following result (see for example \cite{MS} for details).
\begin{lemma}  \label{lem:finite-exceptions}  {\rm
Let $f_0$ and $f_1$ be two generic Morse functions on a world sheet $\Xr$.
There exists a path of smooth functions $f_s,s\in[0,1]$ such that
 $f_s$ are all generic Morse functions except for finitely many points. There are three possibilities at
an exceptional point $s_0$:
\begin{enumerate}
 \item $f_{s_0}$ is a Morse function. However, there are two critical points $c_1\not=c_2$ of $f_{s_0}$, such that
$f_{s_0}(c_1)=f_{s_0}(c_2)$. For $i,j\geqslant 3$, we have $f_{s_0}(c_i)\not=f_{s_0}(c_1)$ and
$f_{s_0}(c_i)\not=f_{s_0}(c_j)$ if $i\not=j$.
 \item $f_{s_0}$ has only one degenerate critical point $p\in\dot\Xr$. If $p$ is in the interior of $\dot\Xr$, then
       there exists $\epsilon>0$, such that for
       $s\in[s_0-\epsilon,s_0+\epsilon]$, there exists local coordinates $(x_1,\cdots,x_n)$ centred at $p$, such that
       $f_s$ can be expressed as
       $$ f_s(x_1,x_2)=f_{s_0}(p)+x_1^3+\epsilon_1(s-s_0)x_1+\epsilon_2 x_2^2$$ where
       $\epsilon_2=\pm 1$. Similarly, if $p$ is on the boundary of $\dot\Xr$, then on the boundary, the Morse function
       can be locally express as  $$ f_s(x_1)=f_{s_0}(p)+x_1^3+\epsilon_1(s-s_0)x_1$$ The following picture shows the picture
       of such Morse functions as $s-s_0$ goes from negative to positive.
      \begin{center}
\includegraphics[width=130mm]{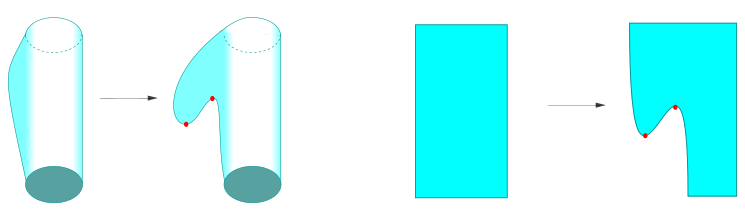}
\end{center}

\item There exists a boundary critical point of $f_{s_0}$, such that the normal derivative of $f_{s_0}$ vanishes at
this boundary critical point. A full discussion of this situation is given in the appendix of \cite{MS}. There are eight
cases, one of which is illustrated in the following picture: 
\begin{equation}        
\figbox{0.5}{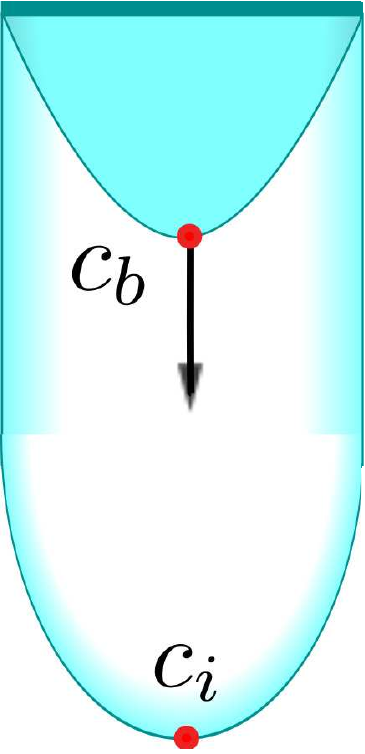}\hspace{20mm}\figbox{0.5}{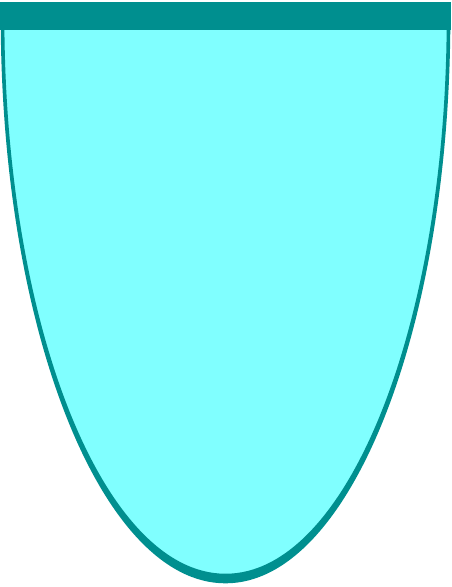}~.
\end{equation}
The Morse functions for above two surfaces are the height function. In the left picture, $c_b$ and $c_i$ are the
boundary and interior critical points of the Morse function, respectively, and the arrow in the left picture denotes the
normal direction at $c_b$. We deform the Morse function from the left picture by pulling the boundary critical point
$c_b$ downward and pushing the interior critical point upward until the height function becomes the one in the right
picture. At one point in this process, the interior critical point $c_i$ merge into the boundary critical point $c_b$.
At this point, the normal direction at the boundary critical point is horizontal with respect to the height function.
That is when the normal derivative vanishes.
\end{enumerate}
}
\end{lemma}

\begin{prop}   \label{prop:Morse-equal}  {\rm
Suppose that $C$ obeys \erf{eq:sewing-cond}. Let $f_0$ and $f_1$ be generic Morse functions on $\Xr$, and let
$\Xr^d_0$, $\Xr^d_1$ be decompositions of $\Xr$ subordinate to $f_0$ and $f_1$ respectively. Then we have $C_{\Xr^d_0} =
C_{\Xr^d_1}$.
}
\end{prop}
\noindent \pf
By Lemma \ref{lem:finite-exceptions}, there exists a family $\{ g_s \}_{0\le s \le 1}$ of smooth functions, such that 
$g_0=f_0$, $g_1=f_1$ and $g_s$ is a generic Morse function on $\Xr$ except at finitely many $s=s_0, \cdots, s_k$. By
Lemma \ref{lem:path-connected generic Morse functions}, it is enough to prove the Lemma for the case in which there is
only one point $s=s_0$ such that $g_{s_0}$ is not a generic Morse function.

Let $\Xr^{d}(s) = (\Xr,(\alpha_1,\dots,\alpha_m)|_s,\varpi(s))$ be a decomposition of $\Xr$ subordinate to $g_s$. By Lemma \ref{lem:path-connected generic Morse functions}, $(\alpha_1,\dots,\alpha_m)|_s$ and $C_{\Xr^{d}(s)}$ can change only when $s$ passes $s_0$. At the point $s=s_0$, there are three possibilities according to the three cases described in Lemma \ref{lem:finite-exceptions}. For $s\in [s_0-\epsilon, s_0+\epsilon]$ and $\epsilon$ being sufficiently small, when $s$ pass $s_0$ from one side to the other side, two critical points $c_k$ and $c_{k+1}$ of $g_s$ merge into one in case 3 of Lemma \ref{lem:finite-exceptions};
merge into one then disappear in case 2; two critical values $g_s(c_k(s))$ and $g_s(c_{k+1}(s))$ merge into one at $s=s_0$ and separate again in case 1.

Let $t_c$ denote the degenerate critical value in either  case. Let $\epsilon,\delta>0$  be small enough, such that
$g_s^{-1}([t_c-\delta,t_c+\delta])$ contains no critical points other than $c_k(s)$ and $c_{k+1}(s)$ for
$s\in[s_0-\epsilon,s_0+\epsilon]$. Let $\Yr=g^{-1}_{s_0}([0,t_c-\delta])\cup g^{-1}_{s_0}([t_c+\delta,1])$ and
$\Zr=g_{s_0}^{-1}([t_c-\delta,t_c+\delta])$. 
We can assume\footnote{If otherwise, we can always deform  the function $g_s$ only locally within a small neighbourhood
of $g_s^{-1}(t_c - \delta)$ and $g_s^{-1}(t_c - \delta)$ in its domain $\Xr$, for all $s\in [s_0-\epsilon, s_0) \cup (s_0,
s_0+\epsilon]$, so that the assumption is achieved. Moreover, this deformation of $g_s$ can be carried out within the
same path connecting component of generic Morse functions for $s\in [s_0-\epsilon, s_0+\epsilon]$.} that the restriction
of $g_s$ on $\Yr$ are generic Morse functions for all $s\in [s_0-\epsilon, s_0+\epsilon]$. Let
$\Yr_{-}^d=(\Yr,(\alpha_1,\cdots,\alpha_n),\eta_-)$ and $\Yr_{+}^d=(\Yr,(\alpha_1,\cdots,\alpha_n),\eta_+)$ be two decompositions of $\Yr$ subordinated to $g_{s_0-\epsilon}|_{\Yr}$ and $g_{s_0+\epsilon}|_{\Yr}$, respectively.
By Lemma \ref{lem:path-connected generic Morse functions}, we have $C_{\Yr_{-}^d}=C_{\Yr_{+}^d}$. Let
$\Zr_{\pm}^d =(\Zr,(\beta_1,\cdots,\beta_k),\xi_\pm)$ be
two decompositions of $\Zr$ subordinated to $g_{s_0\pm\epsilon}|_{\Zr}$, respectively.
Note that there is a morphism of world sheets $\xi:\Yr\otimes
\Zr\rightarrow \Xr$, such that
\begin{align*}
\Xr^d_- &= (\Xr,(\alpha_1,\cdots,\alpha_n,\beta_1,\cdots,\beta_k),\xi\circ(\eta_-\otimes\xi_-)) \ , \\
\Xr^d_+ &= (\Xr,(\alpha_1, \cdots,\alpha_n,\gamma_1,\cdots,\gamma_l),\xi\circ(\eta_+\otimes\xi_+))
\end{align*}
are
two decompositions of $\Xr$ subordinate to $g_{s_0-\epsilon}$ and $g_{s_0+\epsilon}$ respectively.  By Lemma
\ref{lem:change-part}, to prove that $C_{\Xr^d_-}=C_{\Xr^d_+}$, it remains to show that
$C_{\Zr_{-}^d}=C_{\Zr_{+}^d}$.

If $g_{s_0}$ is in case 2 in Lemma \ref{lem:finite-exceptions}, then $C_{\Zr_{-}^d}=C_{\Zr_{+}^d}$
follows from the unit axioms in R1-R32. Therefore, it remains to consider the case 1 and case 3 in
Lemma \ref{lem:finite-exceptions}. The proof for case 3 is similar to the argument in Appendix A.2 of \cite{MS}, which we omit.
For case 1, we consider those world sheets
which allow Morse functions with two critical points on the same level set.
All such open/closed world sheets with exactly two critical points and corresponding decompositions are studied case by case in \cite[App.\,A1\,\&\,A2]{MS} for 2-d topological field theories. The case-by-case study can be
immediately adopted here to give a proof of $C_{\Zr_{-}^d}=C_{\Zr_{+}^d}$, except for one important change, which is when $\Zr$ is given by a torus with two holes. This case is trivial in 2-d
topological field theory but nontrivial in conformal field theory. The proof of the two-holed torus case is given in Lemmas \ref{lem:intersection-number=1} and \ref{lemma:X_12} below. 

Given these results,  altogether we have
$$
C_{\Xr_0^d}
\overset{\text{Lem.\,\ref{lem:path-connected generic Morse functions}}}=
C_{\Xr^d(s_0-\eps)}
\overset{\text{Lem.\,\ref{lem:Morse-same}}}=
C_{\Xr_{-}^d}
=
C_{\Xr_{+}^d}
\overset{\text{Lem.\,\ref{lem:Morse-same}}}=
C_{\Xr^d(s_0+\eps)}
\overset{\text{Lem.\,\ref{lem:path-connected generic Morse functions}}}=
C_{\Xr_1^d} \ .
$$
\epf

\medskip
We will give, in the following two lemmas, a detailed analysis of generic Morse functions on the two-holed genus 1
world sheet, which is denoted by $\Tr$ and depicted pictorially below:
\begin{equation}
{\figbox{0.8}{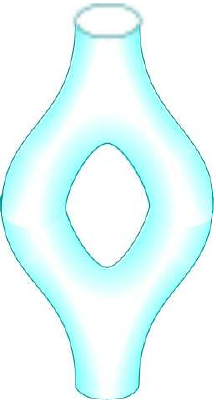}}
\end{equation}

\begin{lemma}  \label{lem:intersection-number=1}  {\rm
Let $f_0$ and $f_1$ be two generic Morse functions on $\Tr$, and let $\Tr^d_i$ be  decompositions of $\Tr$
subordinate to $f_i$ with a choice of value $t_i$ for $i=0,1$. Suppose that one component of $f_0^{-1}(t_0)$ has intersection
number $1$ with a component of
$f_1^{-1}(t_1)$. Then we have $C_{\Tr^d_0}=C_{\Tr^d_1}$.}
\end{lemma}
\pf
Since one component of $f_0^{-1}(t_0)$ has intersection number $1$ with a component of $f_1^{-1}(t_1)$, we can deform $f_1$ in the
space of generic Morse functions so that they intersect each other at exactly one point. By Lemma \ref{lem:path-connected generic Morse
functions}, we only need to prove the lemma for this situation. We can express $f_0$ and $f_1$ pictorially as
\begin{equation}  \label{eq:level-set-f-01}
  {\figbox{0.5}{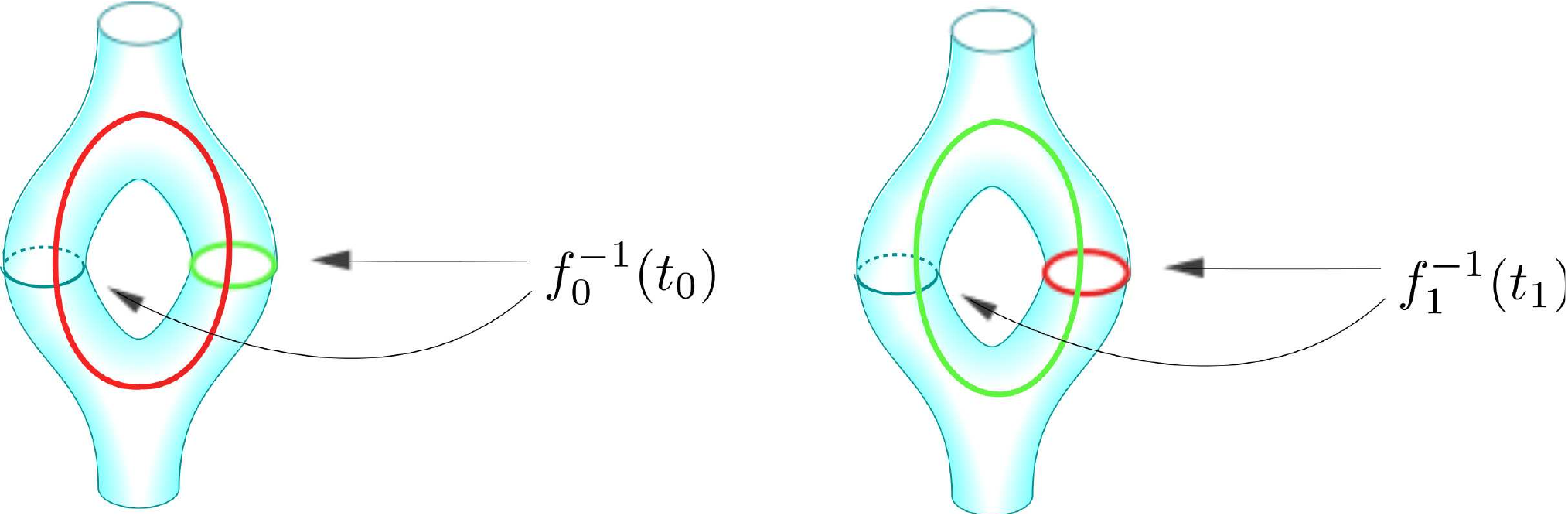}} \quad .
\end{equation}
Here, 
$f_0$ and $f_1$ are the height functions in the left and right figure, respectively. The horizontal green
circle in the left picture is a component of $f_0^{-1}(t_0)$, and the horizontal red circle in the right picture is a component of
$f_1^{-1}(t_1)$. There is a diffeomorphism from the torus in the left picture to the one in the right picture that maps the green
circle to the green circle and the red circle to the red circle. 

The proof of the identity $C_{\Tr^d_0}=C_{\Tr^d_1}$ is shown in the following pictures. The 
coloured dashed circles denote the boundaries that are sewn
in a decomposition of $\Tr$. Circles and lines in the same colour represent the same curve in
$\Tr$. 
$$
\Tr^d_0=\figbox{0.75}{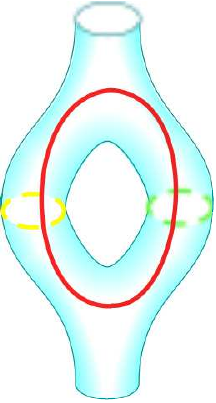}=\figbox{0.75}{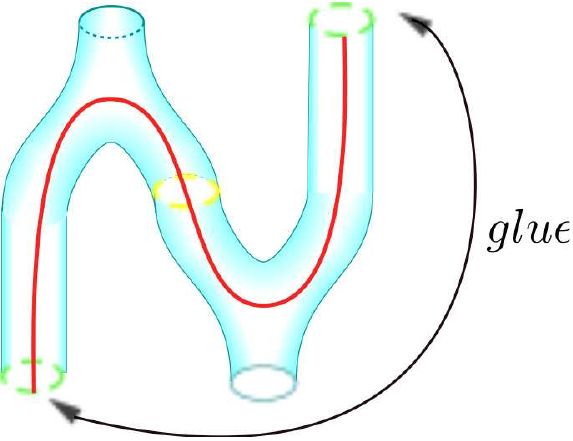}\hspace{2mm}\overset{\Rr18}{\longrightarrow}\hspace{2mm
} \figbox { 0.75 } {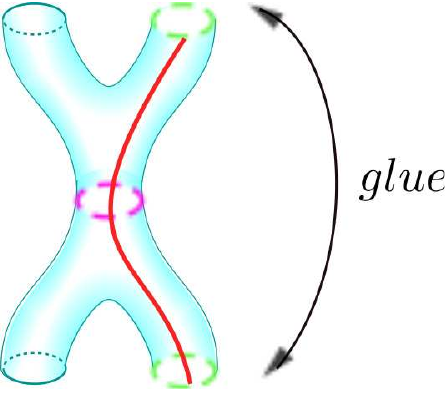}
$$
$$
\hspace{2mm}\overset{\Rr18}{\longrightarrow}\hspace{2mm}\figbox{0.75}{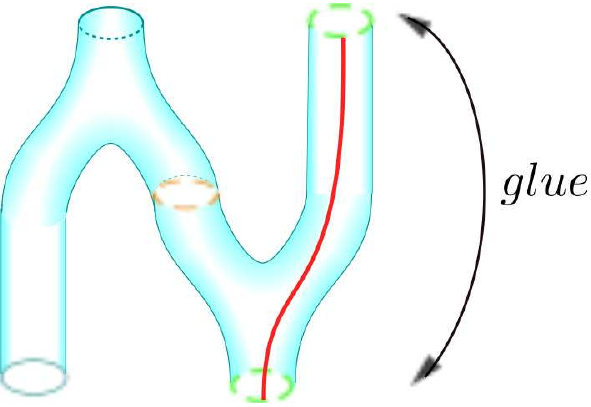}\overset{\Rr32}{\longrightarrow}\hspace
{2mm}\figbox{0.75}{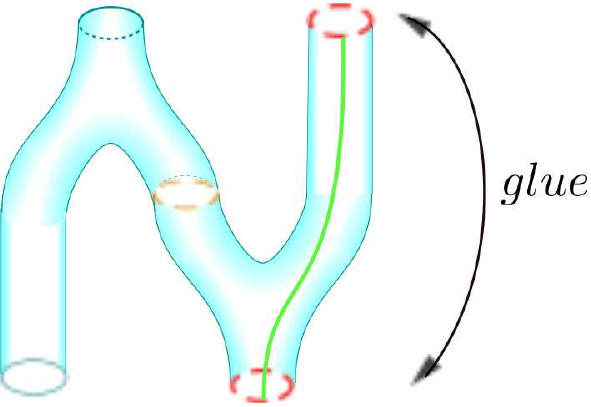}\hspace{2mm}
$$
$$
\overset{\Rr18}{\longrightarrow}
\hspace{2mm}\figbox{0.75}{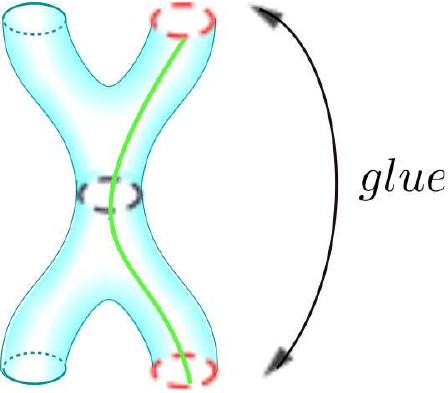}\overset{\Rr18}{\longrightarrow}\hspace{2mm}\Tr^d_1=\figbox{0.75}{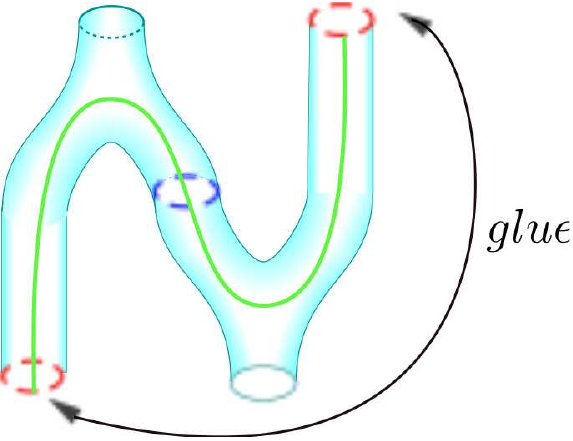}
=\figbox{0.75}{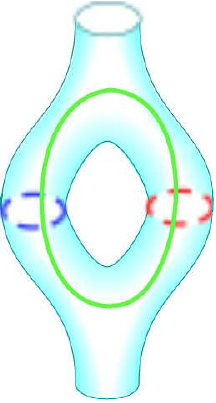}
$$
\epf

\begin{lemma} \label{lemma:X_12}  {\rm
Let $g_0$ and $g_1$ be generic Morse functions on $\Tr$, and let $\Tr^d_1, \Tr^d_2$ be decompositions of $\Tr$ subordinate to $g_0$ and $g_1$ respectively. Then we have $C_{\Tr^d_0}=C_{\Tr^d_1}$. 
}
\end{lemma}
\pf
It is clear that a generic Morse function on $\Tr$ has at least two critical points. Let $g$ be an arbitrary generic Morse
function on $\Tr$, we claim that there exists a path of smooth functions connecting $g$ with another generic Morse function with
only two critical points such that there are only finitely many non-generic Morse functions in this path of smooth functions, all
belonging to the following cases:
\begin{enumerate}
\item Case 2 in Lemma \ref{lem:finite-exceptions}
\item Two critical points of indices 0 and 0 have the same critical value
\item Two critical points of indices 0 and 1 have the same critical value
\item Two critical points of indices 0 and 2 have the same critical value
\item Two critical points of indices 1 and 2 have the same critical value
\item Two critical points of indices 2 and 2 have the same critical value
\end{enumerate}

We prove the claim by induction on the number of critical points of $g$. If $g$ has two critical points, then there is nothing to
prove.  Let $\{c_1,\cdots, c_{n+1}\}$ be the set of critical points of $g$ with $n>1$, such that $g_(c_1)<g(c_2)<\cdots<g(c_{n+1})$,
and let the decomposition be subordinate to $0=t_0<f(c_1)<t_1<f(c_2)<\cdots<f(c_{n+1})<t_{n+1}=1$.  Pick any critical point $c_k$ of
index 0 or 2. (Such a critical point exists since we assume that $g$ has at least three critical points.) We consider the adjacent
critical points of $c_k$, namely, $c_{k-1}$ and $c_{k+1}$. There are the following situations:
\begin{itemize} \setlength{\leftskip}{-1em}
\item
$c_k$
and $c_{k+1}$ lie on the same
connected component of $f^{-1}([t_{k-1},t_{k+1}])$, or $c_k$ and $c_{k-1}$ lie on the same connected component of
$f^{-1}([t_{k-2},t_{k}])$: In this situation we can deform $g$ to a non-generic Morse function as in case 2 of Lemma
\ref{lem:finite-exceptions}. After this deformation, the number of critical points decreases by 2, and we are done
by induction hypothesis. 
\item
$c_{k}$ and $c_{k+1}$ lie on different components of
$f^{-1}([t_{k-1},t_{k+1}])$, and  $c_{k}$ and $c_{k-1}$ lie on different components of $f^{-1}([t_{k-2},t_{k}])$: In this
situation, we can deform the Morse function $g$ to move the critical point $c_k$ up or down so that the value of $c_k$ is greater than
the value of $c_{k+1}$ or less than the value of $c_{k-1}$. Since $c_k$ is of index $0$ or $2$, only cases 2-6 could happen. 
\end{itemize}
After
repeating this operation finitely many times, the first situation will appear, and the claim is proved. 

Let $f$ be a generic Morse function that can be connected to $g$ by a path of smooth functions in which there is only one function
that is not generic Morse. If this degenerate function is of case 1, then the decompositions subordinate to $f$
and $g$ give rise to the same $C_{\Tr}$ by the unit axioms. If this degenerate function is of case 2, then decompositions
subordinate to $f$ and $g$ only differ by some cylinders. Thus decompositions subordinate to $f$ and $g$ give rise to the same
$C_{\Tr}$ by Lemma \ref{lem:inv-omit-cyl}. Cases 3-6 are similar to case 2. Combining this discussion with the above claim, we
conclude that for any generic Morse function $g$ on $\Tr$, there is a generic Morse function with two critical points, such that
decompositions subordinate to them give rise to the same $C_{\Tr}$. Hence it is enough to prove the statement of this Lemma for
generic Morse functions $g_i$, $i=0,1$ each with two critical points and a choice of value $t_i$ between two critical values for
$i=0,1$.

Without loss of generality, we can set $g_0=f_0$ where $f_0$ is the height function of $\Tr$ with a choice of value $t_0$ shown in the left picture in (\ref{eq:level-set-f-01}). Let us denote the homology class of the boundary component of $\Tr$ at the bottom of the left picture of (\ref{eq:level-set-f-01}) by $c$, one component of $f_0^{-1}(t_0)$ (the green curve) by $a$ (hence, the other component must be $c-a$), and one component of $f_1^{-1}(t_1)$ (the red curve) by $b$.
In general, the homology classes of the two components of $g_1^{-1}(t_1)$ must be $\alpha a+\beta b$ and $c-\alpha a-\beta b$ for some $\alpha,\beta\in\mathbb{Z}$ being coprime. We prove the Lemma by induction on $|\alpha|$: 

\medskip
\noindent 1) If $|\alpha|=0$, then one of the components of $g_1^{-1}(t_1)$ must be homologous to       
$b$ or $-b$. In either case, $g_1$ looks exactly like the height function $f_1$ shown in the right picture in (\ref{eq:level-set-f-01}). Therefore, these cases are covered by Lemma \ref{lem:intersection-number=1}. 

\medskip
\noindent 2) Assume that we have $C_{\Tr^d_0} = C_{\Tr^d_1}$
for all generic Morse functions $g_1$ such that one component of $g_1^{-1}(t_1)$ homologous to $\alpha a+\beta b$ with $|\alpha|\leqslant k$. Now we assume that $g_1$ is a generic Morse function such that the homology class of one component of $g_1^{-1}(t_1)$ is $\alpha a+\beta b$ with $|\alpha|=k+1$. Because $\alpha$ and $\beta$ are coprime, we can find $\alpha_1,\beta_1\in\mathbb{Z}$ such that $\alpha_1\beta-\beta_1\alpha=1$, which also implies that the intersection number of $\alpha_1a+\beta_1b$ with $\alpha a+\beta b$ is $1$.  It is clear that we can always find an integer $m$ such that $|\alpha_1+m\alpha|\leqslant k$ (for example, if $k=0$, take $m=-\alpha_1$). Let $h$ be a generic Morse function on $\Tr$ with two critical points (with a choice of $s$ between two critical values) such that the homology class of one component of $h^{-1}(s)$ is $\alpha_1a+\beta_1 b+m(\alpha a+\beta b)$. Let $\Tr^d_h$ be a decomposition of $\Tr$ subordinate to $h$. 
Since $|\alpha_1+m\alpha|\leqslant k$, we have $C_{\Tr^d_0}=C_{\Tr^d_h}$ by our induction hypothesis. On the other hand, it is easy to check that the intersection number of $\alpha_1a+\beta_1 b+m(\alpha a+\beta b)$ with $\alpha a+\beta b$ is $1$. By Lemma \ref{lem:intersection-number=1}, we have $C_{\Tr^d_h}=C_{\Tr^d_1}$.
\epf

\bigskip
\noindent
{\em Proof of Theorem \ref{thm:gen-rel}, part 2:}
\\[.3em]
Suppose that $C$ obeys \erf{eq:sewing-cond}. 
Let $\Xr$ be a world sheet and let $\Xr_1^d$ and $\Xr_2^d$ be decomposed world sheets which map to $\Xr$ under the
forgetful functor. By Lemma \ref{lem:Morse-exists} there exist Morse functions $f_i$, $i=1,2$ on $\dot\Xr$ and world
sheets $\Yr_i^d$, $i=1,2$ such that $\Yr_i$ is subordinate to $f_i$ and $C_{\Xr^d_i} = C_{\Yr^d_i}$. By Proposition
\ref{prop:Morse-equal}, we obtain $C_{\Yr^d_1} = C_{\Yr^d_2}$ which implies $C_{\Xr^d_1} = C_{\Xr^d_2}$. It is therefore
consistent to define
\be
  \mu_\Xr = C_{\Xr^d_1} ~.
\ee
It remains to show that $\mu$ is a monoidal natural transformation. Naturality amounts to the commutativity of the  diagram
\be
  \xymatrix{
  U(\Xr) \ar[d]^{\mu_{\Xr}} \ar[r]^{U(\xi)} & U(\Yr) \ar[d]^{\mu_{\Yr}} \\
  F(\Xr) \ar[r]^{F(\xi)} & F(\Yr)}
\labl{eq:inv-sew-aux1}
for all $\xi : \Xr \rightarrow \Yr$. Let $\Xr^d = (\Xr,(\alpha_1,\dots,\alpha_m),\varpi)$ be a decomposed world sheet.
The morphism $\xi$ gives the decomposed world sheet $\Yr^d = (\Yr,(\alpha_1,\dots,\alpha_m),\xi \circ \varpi)$, and
 $\xi : \Xr^d \rightarrow \Yr^d$ is a morphism of decomposed world sheets.
By definition, we have $\mu_\Xr = C_{\Xr^d}$ and $\mu_\Yr = C_{\Yr^d}$. Then the commutativity of \erf{eq:inv-sew-aux1} follows because $C$ itself is a natural transformation according to Lemma \ref{lem:all-nat-C}, i.e.\ the diagram
\be
  \xymatrix{
  U(\Xr^d) \ar[d]^{C_{\Xr^d}} \ar[r]^{U^\delta(\xi)} & U(\Yr^d) \ar[d]^{C_{\Yr^d}} \\
  F(\Xr^d) \ar[r]^{F^\delta(\xi)} & F(\Yr^d)}
\ee
commutes.  That $\mu$ is monoidal equally follows directly from the fact that $C$ is monoidal.
\epf

\begin{rema}  {\rm
An alternative strategy to prove the above result is to extend the work of \cite{Moore:1988qv,BK-Lego,FuGe} to surfaces with orientation-reversing involution. These works use cut-systems, markings and rigid structures of surfaces, and elementary moves between them to present modular functors in terms of generators and relations. A proof using this language would be desirable as it is closer in spirit to the definition of $\WSh$ and $\WSh^d$. It would also be interesting to see if the proof of Theorem 3.1 in \cite{KP} can be used to simplify the presentation in this appendix.
}
\end{rema}

\end{document}